%% file: main.tex
\pdfoutput=1 
\documentclass{article}[12pt]

\input{head}

\title{A $y$-ification of Khovanov homology}
\author{Taketo Sano}

\begin{document}

\maketitle

\begin{abstract}%
    \input{abst}
\end{abstract}


\input{1}
\input{acknowledgements}

\input{2}
\input{3}
\input{4}
\input{5}
\input{6}

\pagebreak
\printbibliography

\end{document}

%% file: head.tex
\usepackage{amssymb,amsmath,amsthm,amsfonts,mathrsfs}
\usepackage[hmargin=3cm,vmargin=3.5cm]{geometry}
\usepackage{graphicx} 
\usepackage{hyperref}
\usepackage{amsfonts}
\usepackage{amssymb}
\usepackage{amsmath}
\usepackage{amsthm}
\usepackage{mathtools}
\usepackage{cleveref}
\usepackage{tikz}
\usepackage{tikz-cd}
\usepackage{subcaption}
\usepackage{todonotes} 
\usepackage{enumitem}
\usepackage{fancyvrb}

\usepackage[
    backend=biber,
    style=alphabetic,
    sorting=nyt,
    doi=false,
    url=false,
    isbn=false,
    maxbibnames=10,
]{biblatex}

\hypersetup{
  colorlinks=true,
  linkcolor=blue,
  citecolor=blue,
  urlcolor=blue
}

\usetikzlibrary{tikzmark,arrows.meta,knots}

\theoremstyle{plain}
    \newtheorem{mainthm}{Theorem}
    \newtheorem{thm}{Theorem}[section]
    \newtheorem{prop}[thm]{Proposition}
    \newtheorem{cor}[thm]{Corollary}
    \newtheorem{lem}[thm]{Lemma}

\theoremstyle{definition}
    \newtheorem{defn}[thm]{Definition}
    \newtheorem{ex}[thm]{Example}

\theoremstyle{remark}
	\newtheorem{rem}[thm]{Remark}%

\crefname{mainthm}{Theorem}{Theorems}
\crefname{thm}{Theorem}{Theorems}
\crefname{prop}{Proposition}{Propositions}
\crefname{lem}{Lemma}{Lemmas}
\crefname{cor}{Corollary}{Corollaries}
\crefname{ex}{Example}{Examples}

\numberwithin{equation}{section}

\newcommand{\QQ}{\mathbb Q}
\newcommand{\ZZ}{\mathbb Z}
 
\newcommand{\CC}{\mathbb C}
\newcommand{\FF}{\mathbb F}

\newcommand{\msf}{\mathsf} 

\newcommand{\mcA}{{\mathcal A}}

\newcommand{\mcC}{{\mathcal C}}
\newcommand{\mcCA}{{\mcC\mcA}}

\newcommand{\mcF}{{\mathcal F}}
\newcommand{\mcH}{{\mathcal H}}

\newcommand{\mfS}{\mathfrak{S}}

\newcommand{\sfB}{\msf{B}}
\newcommand{\sfE}{\msf{E}}
\newcommand{\sfU}{\msf{U}}

\newcommand{\del}{\partial}
\newcommand{\ddel}[2]{\frac{\del #1}{\del #2}}

\newcommand{\CKh}{\mathit{CKh}}
\newcommand{\Kh}{\mathit{Kh}}
\newcommand{\rCKh}{\widetilde{\CKh}{}}
\newcommand{\rKh}{\widetilde{\Kh}{}}
\newcommand{\rKhscript}{\smash{\rKh}}

\newcommand{\hchi}{\hat{\chi}}

\renewcommand{\sl}{\mathfrak{sl}}
\newcommand{\sle}{\msf{e}}
\newcommand{\slf}{\msf{f}}
\newcommand{\slh}{\msf{h}}
\newcommand{\slE}{E}

\DeclareMathOperator{\Hom}{Hom}

\DeclareMathOperator{\Ima}{Im}

\DeclareMathOperator{\Kom}{Kom}

\DeclareMathOperator{\gr}{gr}
\DeclareMathOperator{\cl}{cl}

\DeclareMathOperator{\sgn}{sgn}

\newcommand{\id}{\mathsf{id}}
\newcommand{\Mod}[1]{#1\operatorname{-Mod}}

\newcommand{\isom}{\cong}
\newcommand{\htpy}{\simeq}

\DeclareMathOperator{\Cob}{Cob}

\DeclareMathOperator{\Mat}{Mat}
\DeclareMathOperator{\Kob}{Kob}
\DeclareMathOperator{\Diag}{Diag}
\DeclareMathOperator{\BDiag}{BDiag}

\renewcommand{\epsilon}{\varepsilon}
\renewcommand{\emptyset}{\varnothing}

\newcommand{\hdot}{\circ}

\newcommand{\dual}{\mathsf{D}}
\newcommand{\bfy}{\mathbf{y}}

\newcommand{\bbeta}{\bar{\beta}}

\newcommand{\up}{\uparrow}
\newcommand{\dn}{\downarrow}

\AddToHook{env/prop/begin}{\crefalias{thm}{prop}} 
\AddToHook{env/lem/begin}{\crefalias{thm}{lem}}
\AddToHook{env/cor/begin}{\crefalias{thm}{cor}}
\AddToHook{env/defn/begin}{\crefalias{thm}{defn}}
\AddToHook{env/rem/begin}{\crefalias{thm}{rem}}
\AddToHook{env/ex/begin}{\crefalias{thm}{ex}}

%% file: abst.tex
Motivated by the $y$-ification of HOMFLY--PT homology by Gorsky and Hogancamp \cite{GH:y-ification2022} and the $\mathfrak{sl}_2$-action of Gorsky, Hogancamp, and Mellit \cite{GHM:symmety-kr2024}, we construct $y$-ifications of Khovanov homology and its equivariant versions within Bar-Natan's framework for tangles, and define an action of the element $e$ in $\mathfrak{sl}_2$ on these $y$-ifications. We then prove that our construction is compatible with the previous ones under Rasmussen's spectral sequence from HOMFLY--PT homology to Khovanov homology. Our construction is elementary and well suited to diagrammatic manipulations and algorithmic implementations. As a result, we verify directly that these additional structures distinguish pairs of knots with identical Khovanov homology and HOMFLY--PT homology, in particular the Conway knot and the Kinoshita--Terasaka knot.

%% file: 1.tex
\section{Introduction}
\label{sec:intro}

Khovanov's categorification of the Jones polynomial~\cite{Khovanov:2000}, known as \textit{Khovanov homology}, assigns to an oriented link $L\subset S^3$ a bigraded homology $\Kh^{*, *}(L; \QQ)$ whose graded Euler characteristic recovers the Jones polynomial. 
A persistent theme in the subject is that adding \textit{extra structures} on Khovanov homology strengthens the invariant, and often detects phenomena invisible at the level of polynomials or bigraded vector spaces. A foundational instance is the work of Kronheimer and Mrowka showing that \textit{Khovanov homology detects the unknot}, a property still unknown for the Jones polynomial~\cite{KM:unknot}. The result is obtained by relating Khovanov homology with \textit{singular instanton Floer homology} via a spectral sequence. Subsequent works also show that Khovanov homology \textit{detects the unlink}, a property that does not hold for the Jones polynomial~\cite{Hedden-Ni:2013,Batson-Seed:2015}.\footnote{%
    See \cite{BSX:2018, BS:2022, Ma:2022, BDLLS:2021, LS:2022, BS:2024} for more detection results in this direction. 
}\footnote{%
    See \cite{OS05, B11, Da15, Sca15, Szabo:2015, lin2019bar, alishahi2023khovanov, Do24, Li24, nahm2025spectral} for more spectral sequences.
} 

Another influential direction was led by Lee's deformation of Khovanov homology~\cite{Lee:2005}, which gave Rasmussen's concordance invariant $s$ and its applications to slice genus bounds~\cite{Rasmussen:2010}. Subsequently, Bar-Natan reformulated Khovanov homology via tangles and cobordisms~\cite{BarNatan:2005}, and Khovanov unified rank-two deformations of Khovanov homology via Frobenius extensions, yielding the \textit{equivariant versions} of Khovanov homology~\cite{Khovanov:2004}. Their algebraic structures have been studied in \cite{Naot:2006,Turner:2006,MTV:2007,Shumakovitch:2014,Wigderson:2016,Turner:2020,Khovanov:2022,Khovanov-Sano:2025,Chen-Yang:2025}. In recent works, the torsion part of equivariant Khovanov homology has found topological applications, as in \cite{Alishahi:2017,Alishahi:2018,Sarkar:2020,Onkar:2020,Caprau-etal:2021,Zhuang:2022,Hayden:2023,Lewark-Marino-Zibrowius:2024,Iltgen-Lewark-Marino:2025}. Up-to-sign functoriality with respect to link cobordisms was established in \cite{Jacobsson:2002,BarNatan:2005},\footnote{%
    Strict functoriality is addressed in \cite{Cap:2007,CMW:2009,Blanc:2010,Beliakova:2019,Vog:2020,Sano:2020-b}.
}
and was applied to the detection of \textit{exotic slice disks}~\cite{Hayden-Sundberg:2021,Hayden:2023}.
Lipshitz and Sarkar gave a space-level refinement of Khovanov homology, which equips it with Steenrod squares and strictly refines the $s$-invariant~\cite{LS:2014a,LS:2014b,LS:2014c}. There are also variants adapted to additional geometry or symmetry of the link, as in~\cite{APS:2004, GLW:2018, Pol:2019, Lobb-Watson:2021, LS:2024, Sano:2025, BDMS:2025}. These developments underscore the potential of the algebraic structures that Khovanov homology admits, and the geometric properties they reflect.

In this paper, we focus on the works of Gorsky and Hogancamp \cite{GH:y-ification2022} and of Gorsky, Hogancamp and Mellit \cite{GHM:symmety-kr2024}, which are extensions of \textit{HOMFLY--PT homology}---a triply graded link homology theory introduced by Khovanov and Rozansky as a categorification of the HOMFLY--PT polynomial~\cite{KR:2008b}. In \cite{GH:y-ification2022}, Gorsky and Hogancamp gave a \textit{$y$-ification} of HOMFLY--PT homology by adjoining formal parameters $y_c$ for each component $c$ of a link $L$, giving a triply graded module over $\QQ[x_c, y_c]_{c\in\pi_0(L)}$. This extension is designed to give a \textit{link-splitting spectral sequence} analogous to the one given by Batson and Seed on Khovanov homology~\cite{Batson-Seed:2015}, while also fitting into conjectural geometric structures coming from Hilbert schemes of points in the plane. Building on this construction, in \cite{GHM:symmety-kr2024}, Gorsky, Hogancamp and Mellit further construct a family of commuting endomorphisms $F_k$ on the $y$-ified HOMFLY--PT homology, analogous to the action of tautological classes in the cohomology of character varieties. They prove a hard Lefschetz-type statement for the second operator $F_2$, which in turn extends to an $\sl_2$-action on the homology group, further proving the $Q \leftrightarrow TQ^{-1}$ symmetry conjectured by Dunfield--Gukov--Rasmussen~\cite{DGR:2006}. 

The purpose of this paper is to bring the idea of $y$-ifications to Khovanov homology and to study the resulting structure both theoretically and computationally. We construct a framework independent of the previous work, and then relate the two formulations via Rasmussen's spectral sequence from HOMFLY--PT homology to Khovanov homology~\cite{Rasmussen:2015}.

\subsection*{$y$-ification of Khovanov homology}

\begin{figure}
    \centering
    \begin{subfigure}[t]{0.3\textwidth}
        \centering
        \input{tikzpictures/B_n}
        \caption{Boundary points $\sfB_n$}
        \label{fig:Bn}
    \end{subfigure}
    \begin{subfigure}[t]{0.3\textwidth}
        \centering
        \input{tikzpictures/xi_i}
        \caption{Homotopy $\xi_i$}
        \label{fig:xi_i}
    \end{subfigure}
    \begin{subfigure}[t]{0.3\textwidth}
        \centering
        \input{tikzpictures/braid}
        \caption{Closure $\bbeta$}
        \label{fig:braid-closure}
    \end{subfigure}
    \caption{}
\end{figure}

We start by explaining the construction of the \textit{$y$-ified Khovanov complex} $y[\beta]$ for a braid diagram $\beta$ on $n$ strands. Here, for simplicity, we work in the non-equivariant setting over $\QQ$. Let $\sfB_n \subset \del I^2$ denote the set of boundary points of $\beta$, with $n$ incoming points $\{P_i\}$ on $I \times \{0\}$ and $n$ outgoing points $\{P'_i\}$ on $I \times \{1\}$ (see \Cref{fig:Bn}), and $w \in \mfS_n$ the underlying permutation of $\beta$. Regard $\beta$ as a tangle diagram with boundary $\sfB_n$ and consider the \textit{formal Khovanov bracket} $[\beta]$ of $\beta$---a complex in the additive category of \textit{dotted cobordisms} $\Cob(\sfB_n)$ defined in \cite{BarNatan:2005}. We follow the grading convention that a dot has \textit{positive} quantum degree $2$, which is opposite from \cite{Khovanov:2000,BarNatan:2005} and the same as \cite{Khovanov:2004,KR:2008a}, since it better suits with the grading convention for HOMFLY--PT homology. 

For any regular point $p$ on $\beta$ (i.e.\ a point which is not a crossing), a bidegree $(0, 2)$ endomorphism $X_p$ of $[\beta]$ is given by an array of identity cobordisms each with a dot placed on the component that contains $p$. Let $\beta_0$ be the orientation preserving resolution of $\beta$, which is given by $n$ vertical arcs. We define a bidegree $(0, 2)$ endomorphism $x_p$ of $[\beta]$ by 
\[
    x_p := \begin{cases}
        X_p & \text{if $p$ lies on the odd indexed strand of $\beta_0$}, \\ 
        -X_p & \text{otherwise}. 
    \end{cases}
\]
This sign choice is made so that the following proposition holds.
\begin{prop}
    If two points $p, q$ lie on the same strand and are separated by a single crossing $c$, then there is a bidegree $(-1, 2)$ homotopy $\chi_c$ defined by the reversal of the local differential at $c$, giving 
    \[
        [d, \chi_c] = \pm(x_p - x_q).
    \]
\end{prop}
$\chi_c$ is called a \textit{dot-sliding homotopy} for the crossing $c$. For each $i = 1, \ldots, n$, define endomorphisms corresponding to the endpoints of $\beta$,
\[
    x_i := x_{P_i}, \quad x'_i := x_{P'_i}.
\]
This gives an action of the polynomial ring $R^e_n := \QQ[x_1, \ldots, x_n; x'_1, \ldots, x'_n]$ on $[\beta]$. For each $i$-th strand of $\beta$, there is a bidegree $(-1, 2)$ homotopy $\xi_i$ giving 
\[
    [d, \xi_i] = x_i - x'_{w(i)},
\]
defined by a signed sum of dot-sliding homotopies for the crossings that lie between $P_i$ and $P'_{w(i)}$ (see \Cref{fig:xi_i}). Furthermore, there is a bidegree $(-2, 4)$ homotopy $u$ giving
\[
    [d, u] = \sum_{i = 1}^n (x_i + x'_{w(i)}) \xi_i
\]
defined by a quadratic form of dot-sliding homotopies. We claim that the endomorphisms $x_i$, $x'_i$, $\xi_i$ and $u$ super-commute (see \Cref{prop:homotopy-u}). These are the additional data that are required to give an \textit{$\mcA^w_n$-module structure} on $[\beta]$, which in turn yields a \textit{$y$-ification} of $[\beta]$ and the associated \textit{$\sle$-operator}. 

\begin{defn}
\label{def:into-y-ify}
    The \textit{$y$-ified Khovanov complex} $(y[\beta], D)$ of the braid $\beta$ is given by
    \[
        y[\beta] := [\beta] \otimes_\QQ \QQ[y_1, \ldots, y_n]
    \]
    with (curved) differential 
    \[
        D = d + \sum_i \xi_i y_i.
    \]
    Here, each $y_i$ is assigned bidegree $(2, -2)$, so that $D$ is homogeneous of bidegree $(1, 0)$. Furthermore, $y[\beta]$ is equipped with an endomorphism of bidegree $(-2, 4)$
    \[
        \sle := u + \sum_i (x_i + x'_{w(i)}) \ddel{}{y_i}
    \]
    called the \textit{$\sle$-operator} for $\beta$.\footnote{%
        We write $e, f, h$ for the standard generators of the Lie algebra $\sl_2$, and the corresponding actions by $\sle, \slf, \slh$. 
    } 
\end{defn}

Here, the differential $D$ squares to $D^2 = \sum_i (x_i - x'_{w(i)}) y_i$, which is non-zero, but chain maps and chain homotopies make sense in the usual way (see \Cref{def:y-ification}). One easily verifies from the above equations that $\sle$ is a chain endomorphism of $y[\beta]$. 
We shall prove that, 

\begin{mainthm}
\label{mainthm:y-ify-braid}
    The $\sle$-equivariant chain homotopy type of $y[\beta]$ is an invariant of the braid. 
\end{mainthm}

There is also a $y$-ification for closed braids. Observe that the closure $\bbeta$ of $\beta$ can be realized inserting $\beta$ into the $1$-input \textit{planar arc diagram} depicted in \Cref{fig:braid-closure}. From the functoriality of planar arc diagrams, the complex $[\beta]$ turns into a complex $[\bbeta]$ in the category $\Cob(\emptyset)$. The actions $x_i$ and $x'_i$ become equal on $[\bbeta]$, making $[\bbeta]$ a \textit{central $\mcA^w_n$-module}. Let $l$ be the number of components of $\bbeta$. The strand-indexed homotopies $\xi_i$ $(i = 1, \ldots, n)$ 
give rise to component-indexed homotopies $\bar{\xi}_k$ $(k = 1, \ldots, l)$, and the degree $-2$ homotopy $u$ gives rise to a homotopy $\bar{u}$, such that 
\[
    [d, \bar{\xi}_k] = 0,\quad 
    [d, \bar{u}] = 2 \sum_{k = 1}^l \bar{x}_k \bar{\xi}_k
\]
(see \Cref{prop:reduce-central-dga-action}). Here, each $\bar{x}_k$ is defined by the action $x_p$ for a chosen point $p$ on the $k$-th component of $\bbeta$. By formulas similar to the ones given in \Cref{def:into-y-ify}, we obtain the definition of the \textit{$y$-ified Khovanov complex} $(y[\bbeta], D)$ for the closed braid $\bbeta$ over the polynomial ring $\QQ[y_1, \ldots, y_l]$. For this $y$-ified complex, we have $D^2 = 0$, and similar to \Cref{mainthm:y-ify-braid}, we prove

\begin{mainthm}
\label{mainthm:y-ify-link}
    The $\sle$-equivariant chain homotopy type of $y[\bbeta]$ is an invariant of the link. 
\end{mainthm}

Finally, suppose that $\beta$ closes to a knot. In this case $l = 1$, and it follows from the definition that $\bar{\xi}_1 = 0$ (see \Cref{prop:y-complex-knot}). Therefore, we have $(y[\bbeta], D) = ([\bbeta] \otimes \QQ[y], d)$ and the underlying complex $[\bbeta]$ is a subcomplex of $y[\bbeta]$. Furthermore, since $\ddel{}{y}$ acts as zero on $[\bbeta]$, the chain endomorphism $\sle$ restricts to $\bar{u}$ on $[\bbeta]$. We also write $\sle$ for the restriction to $[\bbeta]$. We prove, 

\begin{mainthm}
\label{mainthm:y-ify-knot}
    If $\beta$ closes to a knot, the $\sle$-equivariant chain homotopy type of $[\bbeta]$ is an invariant of the knot. 
\end{mainthm}

Since $[\bbeta]$ is a complex in $\Cob(\emptyset)$, we may apply the \textit{tautological functor} $\mcF_0$ and recover the ordinary Khovanov complex $\CKh(\bbeta)$. By extending $\mcF_0$ over the $y$-ifications, we define $y\CKh(\bbeta) := \mcF_0(y[\bbeta])$ the \textit{$y$-ified Khovanov complex} of $\bbeta$, and its homology $y\Kh(\bbeta)$ the \textit{$y$-ified Khovanov homology} of $\bbeta$. Furthermore, $\mcF_0(\sle)$ gives a chain endomorphism of $y\CKh(\bbeta)$, and also induces an endomorphism of $y\Kh(\bbeta)$, both of which are called the \textit{$\sle$-operator}. 
\Cref{mainthm:y-ify-link,mainthm:y-ify-knot} immediately imply,

\begin{cor}
    Let $L$ be a link represented by a braid $\beta$. The $y$-ified Khovanov homology $y\Kh(\bbeta)$ and the $\sle$-operator are invariants of $L$. In particular, for a knot $K$, the $\sle$-operator on Khovanov homology $\Kh(K)$ is an invariant of $K$. 
\end{cor}

In \Cref{sec:y-ified-kh}, we shall formulate the above arguments in the setting of \textit{$U(2)$-equivariant Khovanov homology}, from which all other rank-two deformations can be derived. Then \Cref{mainthm:y-ify-braid} is proved by showing the invariance under the three braid relations, and \Cref{mainthm:y-ify-link} by showing the invariance under the two Markov moves. The proofs heavily rely on a general framework of \textit{dg-modules over the dga $\mcA$}, developed in \Cref{sec:formal}. $y$-ifications for the \textit{reduced versions} are defined similarly.  

We remark that the above arguments are formally identical to those given in \cite{GH:y-ification2022,GHM:symmety-kr2024}. The differences are as follows. (i) We defined the endomorphism $x_p$ depending on which strand of $\beta_0$ the point $p$ belongs to; no such adjustment is needed on the HOMFLY--PT side. (ii) We build on Bar-Natan's framework of Khovanov homology for tangles~\cite{BarNatan:2005}, and passing from $\beta$ to $\bbeta$ is simply realized by the closure functor; on the HOMFLY--PT side, they build on \textit{Rouquier complexes} for braids~\cite{Rouquier:2004} and the closure operation is realized by taking the \textit{Hochschild cohomology functor} as in the formulation of \cite{Khovanov:2007}. (iii) Our construction works over arbitrary commutative ring $R$, both unreduced and reduced; whereas their construction is for the reduced theory over $\QQ$ (and $\CC$). (iv) The invariance of $y\Kh$ and the $\sle$-operator is given by explicit chain homotopy equivalences; whereas theirs involves more sophisticated quasi-isomorphisms of dgas. (v) We need not restrict to braids, in fact, and may consider general tangle diagrams. Therefore, apart from (i), our construction is more elementary and is well suited to diagrammatic manipulations and algorithmic implementations. 

\subsection*{Connection with $y$-ified HOMFLY--PT homology}

Having observed that our constructions are formally identical to those given in \cite{GH:y-ification2022,GHM:symmety-kr2024}, it is natural to expect that they are related in a manageable way. In \cite{Rasmussen:2015}, Rasmussen constructs, for any link $L$ and for each $N \geq 1$, a spectral sequence starting from the reduced HOMFLY--PT homology $\mcH(L)$ of $L$ and converging to the reduced $\sl_N$-homology $\overline{H}_N(L)$ of $L$. In particular, for $N = 2$, it is known that $\overline{H}_2(L)$ is isomorphic to the reduced Khovanov homology $\rKhscript(L^*)$ of (the mirror of) $L$ (see \cite{KR:2008a}). In \cite{CG:structure-in-homfly2024}, Chandler and Gorsky extended Rasmussen's spectral sequences to the $y$-ified setting, giving for each $N$ a spectral sequence starting from the reduced $y$-ified HOMFLY--PT homology $y\mcH(L)$ and converging to the reduced $y$-ified $\sl_N$-homology $y\overline{H}_N(L)$. Furthermore, it is proved that the action of $e$ in $\sl_2$ on $y\mcH(L)$ (given by the $F_2$-operator) extends over the spectral sequence and induces an action on $y\overline{H}_N(L)$. Building on these works, we prove

\begin{mainthm}
\label{mainthm:homfly}
    Our constructions are compatible with the constructions of \cite{GH:y-ification2022,GHM:symmety-kr2024} under the $y$-ified spectral sequence of \cite{CG:structure-in-homfly2024}. Namely, our reduced $y$-ified Khovanov homology $y\rKh(L^*)$ is isomorphic to the reduced $y$-ified $\sl_2$-homology $y\overline{H}_2(L)$, and the associated $\sle$-operator coincides with the induced action of $e$ on $y\overline{H}_2(L)$. In particular, for a knot $K$, our $\sle$-operator on $\rKh(K^*)$ coincides with the induced action of $e$ on $\overline{H}_2(K)$.
\end{mainthm}

As stated earlier, the $y$-ified HOMFLY--PT homology gives a \textit{link-splitting spectral sequence} analogous to the one given by Batson and Seed on Khovanov homology~\cite{Batson-Seed:2015}. To state precisely, for a link $L$, let $\operatorname{split}(L)$ denote the \textit{split union} of the components of $L$. Batson and Seed give a deformation of the Khovanov complex that gives rise to a spectral sequence starting from $\Kh(L)$ and converging to $\Kh(\operatorname{split}(L))$ (\cite[Theorem 1.2]{Batson-Seed:2015}). The existence of the link-splitting spectral sequence implies that (i) it gives a lower bound for the \textit{splitting number} of $L$ (\cite[Theorem 1.2]{Batson-Seed:2015}), and (ii) the $\FF_2$-Khovanov homology detects the unlink (\cite[Theorem 1.3]{Batson-Seed:2015}). 
Analogously, a specialization of the $y$-ified HOMFLY--PT complex gives rise to a spectral sequence starting from $\mcH(L)$ and converging to $\mcH(\operatorname{split}(L))$ (\cite[Theorem 1.12]{GH:y-ification2022}). With \Cref{mainthm:homfly}, we may set this similarity as a statement. First,

\begin{prop}
\label{prop:intro-batson-seed}
    The deformed Khovanov complex given in \cite{Batson-Seed:2015}, obtained from a link diagram $D$ together with a set of weights $\nu = \{\nu_c\}$ assigned to each component $c$ of $D$, is isomorphic (as a filtered complex) to the $y$-ified Khovanov complex specialized as
    \[
        y \CKh(D; \QQ_\nu) := y \CKh(D) \otimes_{\QQ[y_1, \ldots, y_l]} \QQ_\nu
    \]
    where $\QQ_\nu$ denotes $\QQ$ regarded as a $\QQ[y_1, \ldots, y_l]$-module with $y_k$ acted on as $\nu_k$. In particular, if $\nu_c - \nu_{c'}$ is invertible for each pair of components $(c, c')$ of $D$, then the homology $y \Kh(D; \QQ_\nu)$ of $y \CKh(D; \QQ_\nu)$ is isomorphic to $\Kh(\operatorname{split}(L))$. 
\end{prop}

Now, consider the following four spectral sequences: (i) Rasmussen's spectral sequence \cite{Rasmussen:2015} for $N = 2$, (ii) its $y$-ified extension \cite{CG:structure-in-homfly2024} specialized by $y_k = \nu_k$, (iii) Batson--Seed's link-splitting spectral sequence~\cite{Batson-Seed:2015} for reduced Khovanov homology, and (iv) its HOMFLY--PT analogue~\cite{GH:y-ification2022},
\[
    \begin{tikzcd}[row sep=3em, column sep=4em]
    \mcH(L) 
        \arrow[r, "\text{(iv)}", Rightarrow] 
        \arrow[d, "\text{(i)}"', Rightarrow] 
    & \mcH(\operatorname{split}(L)) 
        \arrow[d, "\text{(ii)}", Rightarrow] \\
    \rKh(L^*) 
        \arrow[r, "\text{(iii)}", Rightarrow] 
    & \rKh(\operatorname{split}(L^*)).
    \end{tikzcd}
\]
By combining \Cref{mainthm:homfly} and \Cref{prop:intro-batson-seed}, we shall see that all of these spectral sequences arise from a double filtration on the HOMFLY--PT complex.  A precise statement is given in \Cref{prop:4-ss}. 

\subsection*{Computations}
\label{subsec:intro-computation}

In this section, we only consider reduced versions of HOMFLY--PT homology and Khovanov homology over $\QQ$, and hence omit the adjective `reduced'. 

Having established a bridge between the $\sle$-operators on HOMFLY--PT homology $\mcH$ and on Khovanov homology $\rKhscript$, we compute $\sle$ on $\rKhscript$ for concrete examples and observe how they relate to the $\sl_2$-action on $\mcH$. In \cite{Nakagane-Sano:2025}, Nakagane and the author gave an algorithm to compute the HOMFLY--PT homology, and performed computations for all prime knots with up to $11$ crossings. These computational results will be used together with the program developed by the author~\cite{Sano:YUI}, which now has the feature of computing the $\sle$-operator on Khovanov homology.\footnote{
    As of February 2026, the program is under preparation for public release. 
} 

First we fix conventions. We use Rasmussen's triple grading of $\mcH$~\cite{Rasmussen:2015}, where we denote the grading functions by $(q, a, v)$. It corresponds to the usual bigrading $(t, {q_2})$ of $\Kh$ under the spectral sequence as 
\begin{equation}
\label{eqn:homfly-kh-grading}
    t = \frac{v - a}{2}, \quad 
    {q_2} = q + 2 a.
\end{equation}
The \textit{$\Delta$-grading} on $\mcH$ is defined by 
\[
    \Delta = q + a + v
\]
which coincides with the \textit{$\delta$-grading} on $\rKh$ defined by 
\[
    \delta = 2t + {q_2}.
\]
A \textit{KR-thin} knot $K$ is such that $\mcH(K)$ is supported on a single $\Delta$-grading. For such a knot, the spectral sequence collapses at the $E_1$ page and gives $\mcH(K) \isom \rKhscript(K^*)$. In particular, a KR-thin knot is $\Kh$-thin, i.e.\ $\rKhscript(K^*)$ is supported on $\delta = \Delta$.

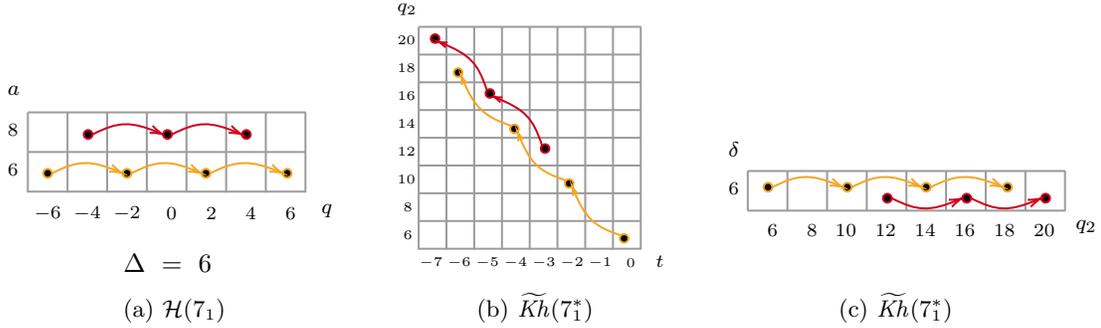
\begin{figure}[t]
    \centering
    \begin{subfigure}[b]{0.3\textwidth}
        \centering
        \input{tikzpictures/table-H-7_1}
        \caption{$\mcH(7_1)$}
        \label{fig:H-7_1}
    \end{subfigure}
    \begin{subfigure}[b]{0.3\textwidth}
        \centering
        \input{tikzpictures/table-Kh-7_1}
        \caption{$\rKh(7_1^*)$}
        \label{fig:Kh-7_1}
    \end{subfigure}
    \begin{subfigure}[b]{0.3\textwidth}
        \centering
        \input{tikzpictures/table-Kh-7_1-slice}
        \caption{$\rKh(7_1^*)$}
        \label{fig:Kh-7_1-slice}
    \end{subfigure}
    \caption{Computation of $\mcH(7_1)$ and $\rKh(7_1^*)$}
    \label{fig:7_1}
\end{figure}

Now, consider $K = 7_1$, the $(2, 7)$-torus knot, which is known to be KR-thin. \Cref{fig:H-7_1} depicts $\mcH(7_1)$ in the format proposed in \cite{CG:structure-in-homfly2024}, showing the single slice $\Delta = 6$ as a bigraded $\QQ$-module with bigrading $(q, a)$. A dot $\bullet$ in a cell at position $(q, a) = (i, j)$ represents a generator of $\mcH^{i, j, \Delta - (i + j)}(7_1)$. For visibility, dots are colored according to the $a$-grading. The arrows indicate the action of $e$, and we see two irreducible $\sl_2$-representations with $h$ acting as $\frac12 q$. The $\sl_2$-action exhibits an apparent symmetry $(q, a) \leftrightarrow (q^{-1}, a)$ on the table (\cite[Theorem 1.2]{GHM:symmety-kr2024}). 

To the right, \Cref{fig:Kh-7_1} depicts $\rKhscript(7_1^*)$ as a bigraded $\QQ$-module with bigrading $(t, {q_2})$.  \Cref{fig:Kh-7_1-slice} depicts the same homology but with an alternative bigrading $(\delta, q_2)$ with $q_2$ on the horizontal axis. Obviously, the two have the same amount of information, but we prefer to use the latter since it is more compact and exhibits the correspondence between $\mcH$ more clearly. Computation of the $\sle$-operator on $\rKhscript$ shows that $e$ acts as indicated by the arrows, and we may confirm that the bigraded $\QQ[e]$-module structure is isomorphic to that on $\mcH$ under the correspondence of \eqref{eqn:homfly-kh-grading}. (In \Cref{prop:torus-knot-str}, we determine the $\QQ[e]$-module structure for all $(2, 2k + 1)$-torus knot and see that this pattern appears in general.) However, the symmetry of $\mcH$ has been lost in $\rKhscript$, and the two $e$-strings that lie separately on $a = 6, 8$ in $\mcH$ collapse onto the same line in $\rKhscript$. This is the case for KR-thin knots, where passing from $\mcH$ to $\rKhscript$ simply collapses the $a$-grading and loses the visible $\sl_2$-symmetry. 

In \cite{Nakagane-Sano:2025}, for all prime knots with up to $11$-crossings, we identified the 96 KR-thick prime knots, and further grouped them according to having identical HOMFLY--PT homology, resulting in 44 groups containing two or more distinct knots. These are the targets of interest from the viewpoint of refining the knot invariants. Our computation showed,

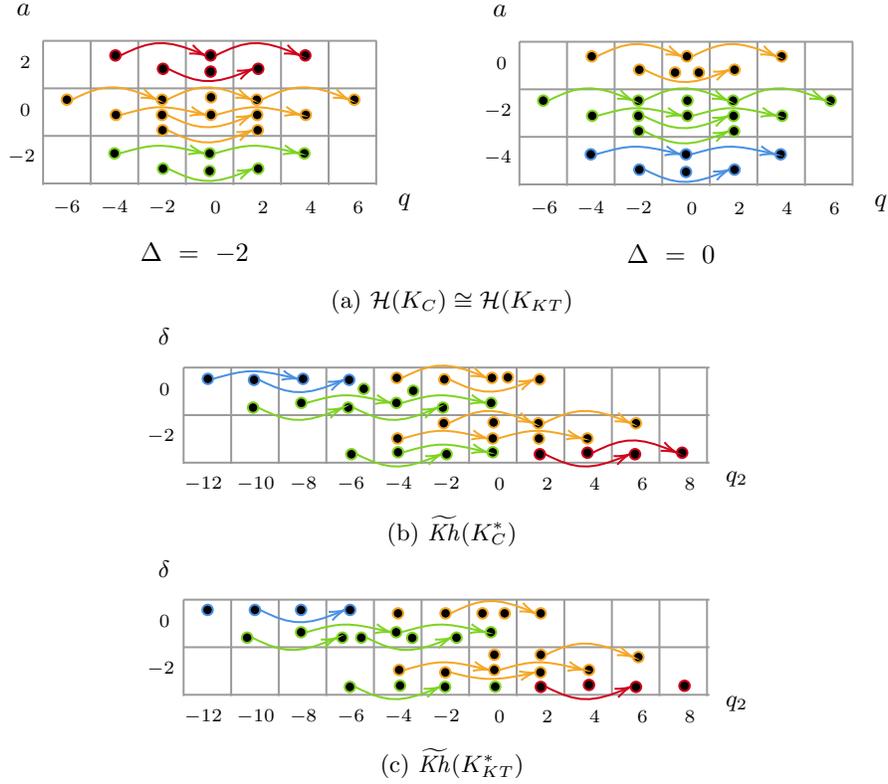
\begin{figure}[t]
    \centering
    \begin{subfigure}{\textwidth}
        \centering
        \input{tikzpictures/table-H-conway}
        \caption{$\mcH(K_C) \isom \mcH(K_{KT})$}
        \label{fig:H-conway}
    \end{subfigure}
    \begin{subfigure}[t]{\textwidth}
        \centering
        \input{tikzpictures/table-kh-conway}
        \caption{$\rKh(K_C^*)$}
        \label{fig:Kh-conway}
    \end{subfigure}
    \begin{subfigure}[t]{\textwidth}
        \centering
        \input{tikzpictures/table-kh-kt}
        \caption{$\rKh(K_{KT}^*)$}
        \label{fig:Kh-KT}
    \end{subfigure}
    \caption{Distinguishing $K_C$ and $K_{KT}$ by the action of $e$}
\end{figure}

\begin{mainthm}
\label{mainthm:conway-and-kt}
    The Conway knot $K_C = 11n_{34}$ and the Kinoshita--Terasaka knot $K_{KT} = 11n_{42}$ have identical HOMFLY--PT homology and reduced Khovanov homology, but have distinct $\sle$-operators on reduced Khovanov homology (i.e. they are not isomorphic as bigraded $\QQ[e]$-modules).
\end{mainthm}

Recall that $K_C$ and $K_{KT}$ are Conway mutant pairs, and are hard to distinguish by classical knot invariants. Both have trivial Alexander and Conway polynomials; identical Jones, HOMFLY--PT and Kauffman polynomials (since these are mutation invariant); and are algebraically slice (hence have trivial signature). Modern \textit{slice-torus invariants}, such as Rasmussen's $s$-invariant \cite{Rasmussen:2010} and Ozsv\'{a}th--Szab\'{o}'s $\tau$-invariant \cite{OS:2003} also fail to distinguish them. Computations of Khovanov homology for knots with up to $11$ crossings have been performed in \cite{Bar-Natan:2002}, and that the two knots have identical Khovanov homology was explicitly stated in \cite{Wehrli:2003}. That the two knots have identical HOMFLY--PT homology was proved in \cite{MV:08}. Invariants that distinguish the two knots include: the knot group \cite{Riley:1971}, the genus $g(K_C) = 3,\ g(K_{KT}) = 2$ \cite{Gabai:1986}, and the slice genus $g_4(K_C) \neq g_4(K_{KT}) = 0$ \cite{Piccirillo:2020}; all of which are highly non-trivial results. Knot Floer homology also distinguish these knots, since it detects the genus~\cite{OS:2003}; explicit computations are given in \cite{OS:2004,BW:2012}. As for quantum invariants, certain \textit{colored HOMFLY--PT polynomials} can distinguish them \cite{Morton-Cromwell:1996,NRV:2017}. See \cite{Hom:2022} for a detailed survey. 
\Cref{mainthm:conway-and-kt} adds Khovanov homology with $\sle$ to the list of successful methods, and in particular, it shows that the $\sle$ is \textit{not} mutation invariant. 

Let us observe how the bigraded $\QQ[e]$-module structures of the two homology groups differ. First, \Cref{fig:H-conway} depicts the $\sl_2$-module $\mcH(K_C) \isom \mcH(K_{KT})$ supported on $\Delta = -2, 0$. Since they have KR-thickness $2$, each spectral sequence collapses at the $E_2$ page.  \Cref{fig:Kh-conway,fig:Kh-KT} depict the computed $\QQ[e]$-module structures of $\rKhscript(K_C^*)$ and $\rKhscript(K_{KT}^*)$ (without the colors). Observe that, as bigraded $\QQ$-modules, the two  are isomorphic, but as bigraded $\QQ[e]$-modules, they are not. For instance, compare the left-end or the right-end of the two tables, we see that there are isolated summands in $\rKhscript(K_{KT}^*)$, which are not present in $\rKhscript(K_{C}^*)$. By further analyzing the spectral sequence over $\QQ[e]$, we may determine the triply graded $\QQ[e]$-module structures of the $E_2$ pages, hence giving the colors depicted in \Cref{fig:Kh-conway,fig:Kh-KT}. A detailed argument is given in \Cref{subsec:conway-and-kt}. 

In \Cref{prop:stronger-than-homfly}, we give three more pairs of knots that have both identical HOMFLY--PT and identical Khovanov homology but distinct $\sle$-operators on Khovanov homology, which are 
\[
    (11n_{39}, 11n_{45}), \quad
    (11n_{73}, 11n_{74}), \quad
    (11n_{151}, 11n_{152}).
\]
Still, within prime knots with up to $11$ crossings, there are 40 groups of knots that cannot be distinguished by HOMFLY--PT, Khovanov homology and $\sle$ combined; for instance the pair $(8_8, m(10_{129}))$. If we only consider Khovanov homology, then of the 90 groups that have identical Khovanov homology, 57 groups are ruled out for having distinct $\sle$-operator (\Cref{prop:stronger-than-kh}). 

\subsection*{Future directions}

\begin{enumerate}
    \item As a continuation of this work, we shall study how $y$-ified Khovanov homology and the $\sle$-operator behave under link cobordisms. It is possible that we obtain a refinement of the \textit{$s$-invariant} by following the approach taken by Lipshitz and Sarkar in \cite{LS:2014c}, where they give a refinement of $s$ via the Steenrod squares on Khovanov homology. 
    
    \item In \cite{ISST:2025a,ISST:2025b}, Imori, Sato, Taniguchi and the author studied Kronheimer and Mrowka's spectral sequence from Khovanov homology to singular instanton homology~\cite{KM:unknot}, and proved that (i) the spectral sequence is compatible with the cobordism maps, and (ii) the compatibility further extends to \textit{immersed cobordisms}. One may ask whether a similar $\sle$-operator can be defined on the singular instanton side, which is compatible with our $\sle$ and is natural with respect to cobordisms. 
    
    \item We have seen above that the $\sl_2$-symmetry on $\mcH$ collapses as we pass from $\mcH$ to $\rKhscript$. Still, we may ask whether Khovanov homology admits a natural $\sl_2$-action that is compatible (in a certain sense) with that on $\mcH$ under Rasmussen's spectral sequence. If there is such action, then it would be interesting to ask for its geometric or physical interpretations. 
    
    \item Recently, there has been interest in endowing HOMFLY--PT and (equivariant) $\mathfrak{gl}_N$ homologies with $\sl_2$-actions \cite{KR:2016, QRSW:2023, QRSW:2024}. In contrast to the $\sl_2$-action of Gorsky--Hogancamp--Mellit, these actions are homogeneous of homological degree $0$. It is interesting to ask whether such actions can also be specialized to Khovanov homology, and study how they interact with the $\sl_2$-action coming from $y$-ification.
\end{enumerate}

\subsection*{Organization}

This paper is organized as follows. In \Cref{sec:formal}, we develop a general framework of $\mcA$-modules and their $y$-ifications, adapting the approach of \cite{GHM:symmety-kr2024} to our setting. In \Cref{sec:y-ified-kh}, we construct the $y$-ified Khovanov complex for braids and links, and prove \Cref{mainthm:y-ify-braid,mainthm:y-ify-link,mainthm:y-ify-knot} with generalizations to the equivariant setting (\Cref{thm:inv-yified-cpx,thm:closed-y-invariance,thm:e-inv-knot}). In \Cref{sec:homfly}, we relate our construction to the $y$-ified HOMFLY--PT homology of \cite{GH:y-ification2022,GHM:symmety-kr2024} via Rasmussen's spectral sequence and prove \Cref{mainthm:homfly} with generalization to the deformed setting (\Cref{thm:y-compatibility-deformed}). We further prove its relation with Batson--Seed's link splitting spectral sequence~\cite{Batson-Seed:2015}. In \Cref{sec:diagrammatic-computation}, we develop diagrammatic reduction techniques for $y$-ified Khovanov homology in the spirit of \cite{BarNatan:2007}, and completely describe the $\QQ[e]$-module structure on Khovanov homology for $(2,2k+1)$-torus knots and twist knots. Finally, in \Cref{sec:direct-computation}, we describe the implementation and present computational results, including the result showing \Cref{mainthm:conway-and-kt}, with a detailed analysis in \Cref{subsec:conway-and-kt}.

%% file: tikzpictures/B_n.tex
\tikzset{every picture/.style={line width=0.75pt}} 

\begin{tikzpicture}[x=0.75pt,y=0.75pt,yscale=-.75,xscale=.75]

\draw   (20,40) -- (149,40) -- (149,169) -- (20,169) -- cycle ;
\draw    (30,176) -- (30,158) ;
\draw [shift={(30,156)}, rotate = 90] [color={rgb, 255:red, 0; green, 0; blue, 0 }  ][line width=0.75]    (10.93,-3.29) .. controls (6.95,-1.4) and (3.31,-0.3) .. (0,0) .. controls (3.31,0.3) and (6.95,1.4) .. (10.93,3.29)   ;
\draw    (59,176) -- (59,158) ;
\draw [shift={(59,156)}, rotate = 90] [color={rgb, 255:red, 0; green, 0; blue, 0 }  ][line width=0.75]    (10.93,-3.29) .. controls (6.95,-1.4) and (3.31,-0.3) .. (0,0) .. controls (3.31,0.3) and (6.95,1.4) .. (10.93,3.29)   ;
\draw    (140,176) -- (140,158) ;
\draw [shift={(140,156)}, rotate = 90] [color={rgb, 255:red, 0; green, 0; blue, 0 }  ][line width=0.75]    (10.93,-3.29) .. controls (6.95,-1.4) and (3.31,-0.3) .. (0,0) .. controls (3.31,0.3) and (6.95,1.4) .. (10.93,3.29)   ;
\draw    (30,48) -- (30,30) ;
\draw [shift={(30,28)}, rotate = 90] [color={rgb, 255:red, 0; green, 0; blue, 0 }  ][line width=0.75]    (10.93,-3.29) .. controls (6.95,-1.4) and (3.31,-0.3) .. (0,0) .. controls (3.31,0.3) and (6.95,1.4) .. (10.93,3.29)   ;
\draw    (59,48) -- (59,30) ;
\draw [shift={(59,28)}, rotate = 90] [color={rgb, 255:red, 0; green, 0; blue, 0 }  ][line width=0.75]    (10.93,-3.29) .. controls (6.95,-1.4) and (3.31,-0.3) .. (0,0) .. controls (3.31,0.3) and (6.95,1.4) .. (10.93,3.29)   ;
\draw    (140,48) -- (140,30) ;
\draw [shift={(140,28)}, rotate = 90] [color={rgb, 255:red, 0; green, 0; blue, 0 }  ][line width=0.75]    (10.93,-3.29) .. controls (6.95,-1.4) and (3.31,-0.3) .. (0,0) .. controls (3.31,0.3) and (6.95,1.4) .. (10.93,3.29)   ;

\draw (92,179.4) node [anchor=north west][inner sep=0.75pt]    {$\cdots $};
\draw (91,10.4) node [anchor=north west][inner sep=0.75pt]    {$\cdots $};
\draw (30,179.4) node [anchor=north] [inner sep=0.75pt]    {$P_{1}$};
\draw (59,179.4) node [anchor=north] [inner sep=0.75pt]    {$P_{2}$};
\draw (140,179.4) node [anchor=north] [inner sep=0.75pt]    {$P_{n}$};
\draw (30,24.6) node [anchor=south] [inner sep=0.75pt]    {$P'_{1}$};
\draw (59,24.6) node [anchor=south] [inner sep=0.75pt]    {$P'_{2}$};
\draw (140,24.6) node [anchor=south] [inner sep=0.75pt]    {$P'_{n}$};

\end{tikzpicture}

%% file: tikzpictures/xi_i.tex
\tikzset{every picture/.style={line width=0.75pt}} 

\begin{tikzpicture}[x=0.75pt,y=0.75pt,yscale=-.75,xscale=.75]

\draw   (20,40) -- (149,40) -- (149,169) -- (20,169) -- cycle ;
\draw    (49.5,177) -- (49.5,159) ;
\draw [shift={(49.5,157)}, rotate = 90] [color={rgb, 255:red, 0; green, 0; blue, 0 }  ][line width=0.75]    (10.93,-3.29) .. controls (6.95,-1.4) and (3.31,-0.3) .. (0,0) .. controls (3.31,0.3) and (6.95,1.4) .. (10.93,3.29)   ;
\draw    (120.5,48) -- (120.5,30) ;
\draw [shift={(120.5,28)}, rotate = 90] [color={rgb, 255:red, 0; green, 0; blue, 0 }  ][line width=0.75]    (10.93,-3.29) .. controls (6.95,-1.4) and (3.31,-0.3) .. (0,0) .. controls (3.31,0.3) and (6.95,1.4) .. (10.93,3.29)   ;
\draw  [dash pattern={on 0.84pt off 2.51pt}]  (49.5,157) .. controls (50.2,122) and (119.8,90.4) .. (120.5,48) ;
\draw    (50,135) -- (60,145) ;
\draw    (60,135) -- (50,145) ;

\draw    (66,115) -- (76,125) ;
\draw    (76,115) -- (66,125) ;

\draw (49.5,181.4) node [anchor=north] [inner sep=0.75pt]    {$P_{i}$};
\draw (120.5,24.6) node [anchor=south] [inner sep=0.75pt]    {$P'_{w( i)}$};
\draw (44,138) node [anchor=east] [inner sep=0.75pt]  [font=\footnotesize]  {$c_{1}$};
\draw (61,117) node [anchor=east] [inner sep=0.75pt]  [font=\footnotesize]  {$c_{2}$};
\draw (89,97.6) node [anchor=south east] [inner sep=0.75pt]    {$\xi _{i}$};

\end{tikzpicture}

%% file: tikzpictures/braid.tex
\tikzset{every picture/.style={line width=0.75pt}} 

\begin{tikzpicture}[x=0.75pt,y=0.75pt,yscale=-.75,xscale=.75]

\draw   (30,60) -- (110,60) -- (110,140) -- (30,140) -- cycle ;
\draw    (40.15,42) -- (40,60) ;
\draw  [draw opacity=0] (40.15,42) .. controls (40.05,41.34) and (40,40.67) .. (40,40) .. controls (40,23.43) and (71.34,10) .. (110,10) .. controls (148.66,10) and (180,23.43) .. (180,40) -- (110,40) -- cycle ; \draw   (40.15,42) .. controls (40.05,41.34) and (40,40.67) .. (40,40) .. controls (40,23.43) and (71.34,10) .. (110,10) .. controls (148.66,10) and (180,23.43) .. (180,40) ;  
\draw    (180,40) -- (180,160) ;
\draw [shift={(180,106)}, rotate = 270] [color={rgb, 255:red, 0; green, 0; blue, 0 }  ][line width=0.75]    (10.93,-3.29) .. controls (6.95,-1.4) and (3.31,-0.3) .. (0,0) .. controls (3.31,0.3) and (6.95,1.4) .. (10.93,3.29)   ;
\draw  [draw opacity=0] (40.15,158) .. controls (40.05,158.66) and (40,159.33) .. (40,160) .. controls (40,176.57) and (71.34,190) .. (110,190) .. controls (148.66,190) and (180,176.57) .. (180,160) -- (110,160) -- cycle ; \draw   (40.15,158) .. controls (40.05,158.66) and (40,159.33) .. (40,160) .. controls (40,176.57) and (71.34,190) .. (110,190) .. controls (148.66,190) and (180,176.57) .. (180,160) ;  
\draw    (40,140) -- (40.15,158) ;
\draw    (100,40) -- (100,60) ;
\draw    (99.85,142) -- (100,160) ;
\draw  [draw opacity=0] (100,40) .. controls (99.99,39.86) and (99.98,39.71) .. (99.98,39.57) .. controls (99.98,34.6) and (106.73,30.57) .. (115.06,30.57) .. controls (123.39,30.57) and (130.14,34.6) .. (130.14,39.57) .. controls (130.14,39.7) and (130.13,39.83) .. (130.12,39.95) -- (115.06,39.57) -- cycle ; \draw   (100,40) .. controls (99.99,39.86) and (99.98,39.71) .. (99.98,39.57) .. controls (99.98,34.6) and (106.73,30.57) .. (115.06,30.57) .. controls (123.39,30.57) and (130.14,34.6) .. (130.14,39.57) .. controls (130.14,39.7) and (130.13,39.83) .. (130.12,39.95) ;  
\draw    (130,39) -- (130,159) ;
\draw [shift={(130,105)}, rotate = 270] [color={rgb, 255:red, 0; green, 0; blue, 0 }  ][line width=0.75]    (10.93,-3.29) .. controls (6.95,-1.4) and (3.31,-0.3) .. (0,0) .. controls (3.31,0.3) and (6.95,1.4) .. (10.93,3.29)   ;
\draw  [draw opacity=0] (100,160) .. controls (99.99,160.14) and (99.98,160.29) .. (99.98,160.43) .. controls (99.98,165.4) and (106.73,169.43) .. (115.06,169.43) .. controls (123.39,169.43) and (130.14,165.4) .. (130.14,160.43) .. controls (130.14,159.78) and (130.02,159.15) .. (129.8,158.54) -- (115.06,160.43) -- cycle ; \draw   (100,160) .. controls (99.99,160.14) and (99.98,160.29) .. (99.98,160.43) .. controls (99.98,165.4) and (106.73,169.43) .. (115.06,169.43) .. controls (123.39,169.43) and (130.14,165.4) .. (130.14,160.43) .. controls (130.14,159.78) and (130.02,159.15) .. (129.8,158.54) ;  

\draw (70,100) node    {$\beta $};
\draw (60,32.4) node [anchor=north west][inner sep=0.75pt]    {$\cdots $};
\draw (140,92.4) node [anchor=north west][inner sep=0.75pt]    {$\cdots $};

\end{tikzpicture}

%% file: tikzpictures/table-H-7_1.tex
\tikzset{every picture/.style={line width=0.75pt}} 

\begin{tikzpicture}[x=0.75pt,y=0.75pt,yscale=-1,xscale=1]

\draw  [draw opacity=0] (50,19.55) -- (190.2,19.55) -- (190.2,60) -- (50,60) -- cycle ; \draw  [color={rgb, 255:red, 155; green, 155; blue, 155 }  ,draw opacity=1 ] (50,19.55) -- (50,60)(70,19.55) -- (70,60)(90,19.55) -- (90,60)(110,19.55) -- (110,60)(130,19.55) -- (130,60)(150,19.55) -- (150,60)(170,19.55) -- (170,60)(190,19.55) -- (190,60) ; \draw  [color={rgb, 255:red, 155; green, 155; blue, 155 }  ,draw opacity=1 ] (50,19.55) -- (190.2,19.55)(50,39.55) -- (190.2,39.55)(50,59.55) -- (190.2,59.55) ; \draw  [color={rgb, 255:red, 155; green, 155; blue, 155 }  ,draw opacity=1 ]  ;
\draw  [color={rgb, 255:red, 245; green, 166; blue, 35 }  ,draw opacity=1 ][fill={rgb, 255:red, 0; green, 0; blue, 0 }  ,fill opacity=1 ] (57.8,50.3) .. controls (57.8,49.08) and (58.78,48.1) .. (60,48.1) .. controls (61.22,48.1) and (62.2,49.08) .. (62.2,50.3) .. controls (62.2,51.52) and (61.22,52.5) .. (60,52.5) .. controls (58.78,52.5) and (57.8,51.52) .. (57.8,50.3) -- cycle ;
\draw  [color={rgb, 255:red, 245; green, 166; blue, 35 }  ,draw opacity=1 ][fill={rgb, 255:red, 0; green, 0; blue, 0 }  ,fill opacity=1 ] (97.8,50.3) .. controls (97.8,49.08) and (98.78,48.1) .. (100,48.1) .. controls (101.22,48.1) and (102.2,49.08) .. (102.2,50.3) .. controls (102.2,51.52) and (101.22,52.5) .. (100,52.5) .. controls (98.78,52.5) and (97.8,51.52) .. (97.8,50.3) -- cycle ;
\draw  [color={rgb, 255:red, 245; green, 166; blue, 35 }  ,draw opacity=1 ][fill={rgb, 255:red, 0; green, 0; blue, 0 }  ,fill opacity=1 ] (137.6,50.3) .. controls (137.6,49.08) and (138.58,48.1) .. (139.8,48.1) .. controls (141.02,48.1) and (142,49.08) .. (142,50.3) .. controls (142,51.52) and (141.02,52.5) .. (139.8,52.5) .. controls (138.58,52.5) and (137.6,51.52) .. (137.6,50.3) -- cycle ;
\draw  [color={rgb, 255:red, 245; green, 166; blue, 35 }  ,draw opacity=1 ][fill={rgb, 255:red, 0; green, 0; blue, 0 }  ,fill opacity=1 ] (178.6,50.3) .. controls (178.6,49.08) and (179.58,48.1) .. (180.8,48.1) .. controls (182.02,48.1) and (183,49.08) .. (183,50.3) .. controls (183,51.52) and (182.02,52.5) .. (180.8,52.5) .. controls (179.58,52.5) and (178.6,51.52) .. (178.6,50.3) -- cycle ;
\draw [color={rgb, 255:red, 245; green, 166; blue, 35 }  ,draw opacity=1 ]   (62.2,50.3) .. controls (75.32,43.64) and (81.37,43.41) .. (96.16,49.6) ;
\draw [shift={(97.8,50.3)}, rotate = 203.33] [color={rgb, 255:red, 245; green, 166; blue, 35 }  ,draw opacity=1 ][line width=0.75]    (6.56,-1.97) .. controls (4.17,-0.84) and (1.99,-0.18) .. (0,0) .. controls (1.99,0.18) and (4.17,0.84) .. (6.56,1.97)   ;
\draw [color={rgb, 255:red, 245; green, 166; blue, 35 }  ,draw opacity=1 ]   (102.2,50.3) .. controls (115.32,43.64) and (121.18,43.41) .. (135.96,49.6) ;
\draw [shift={(137.6,50.3)}, rotate = 203.33] [color={rgb, 255:red, 245; green, 166; blue, 35 }  ,draw opacity=1 ][line width=0.75]    (6.56,-1.97) .. controls (4.17,-0.84) and (1.99,-0.18) .. (0,0) .. controls (1.99,0.18) and (4.17,0.84) .. (6.56,1.97)   ;
\draw [color={rgb, 255:red, 245; green, 166; blue, 35 }  ,draw opacity=1 ]   (142,50.3) .. controls (155.12,43.64) and (162.1,43.41) .. (176.95,49.6) ;
\draw [shift={(178.6,50.3)}, rotate = 203.33] [color={rgb, 255:red, 245; green, 166; blue, 35 }  ,draw opacity=1 ][line width=0.75]    (6.56,-1.97) .. controls (4.17,-0.84) and (1.99,-0.18) .. (0,0) .. controls (1.99,0.18) and (4.17,0.84) .. (6.56,1.97)   ;

\draw [color={rgb, 255:red, 208; green, 2; blue, 27 }  ,draw opacity=1 ]   (82.6,30.7) .. controls (95.72,24.04) and (101.77,23.81) .. (116.56,30) ;
\draw [shift={(118.2,30.7)}, rotate = 203.33] [color={rgb, 255:red, 208; green, 2; blue, 27 }  ,draw opacity=1 ][line width=0.75]    (6.56,-1.97) .. controls (4.17,-0.84) and (1.99,-0.18) .. (0,0) .. controls (1.99,0.18) and (4.17,0.84) .. (6.56,1.97)   ;
\draw  [color={rgb, 255:red, 208; green, 2; blue, 27 }  ,draw opacity=1 ][fill={rgb, 255:red, 0; green, 0; blue, 0 }  ,fill opacity=1 ] (78.2,30.7) .. controls (78.2,29.48) and (79.18,28.5) .. (80.4,28.5) .. controls (81.62,28.5) and (82.6,29.48) .. (82.6,30.7) .. controls (82.6,31.92) and (81.62,32.9) .. (80.4,32.9) .. controls (79.18,32.9) and (78.2,31.92) .. (78.2,30.7) -- cycle ;
\draw  [color={rgb, 255:red, 208; green, 2; blue, 27 }  ,draw opacity=1 ][fill={rgb, 255:red, 0; green, 0; blue, 0 }  ,fill opacity=1 ] (118.2,30.7) .. controls (118.2,29.48) and (119.18,28.5) .. (120.4,28.5) .. controls (121.62,28.5) and (122.6,29.48) .. (122.6,30.7) .. controls (122.6,31.92) and (121.62,32.9) .. (120.4,32.9) .. controls (119.18,32.9) and (118.2,31.92) .. (118.2,30.7) -- cycle ;
\draw  [color={rgb, 255:red, 208; green, 2; blue, 27 }  ,draw opacity=1 ][fill={rgb, 255:red, 0; green, 0; blue, 0 }  ,fill opacity=1 ] (158,30.7) .. controls (158,29.48) and (158.98,28.5) .. (160.2,28.5) .. controls (161.42,28.5) and (162.4,29.48) .. (162.4,30.7) .. controls (162.4,31.92) and (161.42,32.9) .. (160.2,32.9) .. controls (158.98,32.9) and (158,31.92) .. (158,30.7) -- cycle ;
\draw [color={rgb, 255:red, 208; green, 2; blue, 27 }  ,draw opacity=1 ]   (122.6,30.7) .. controls (135.72,24.04) and (141.58,23.81) .. (156.36,30) ;
\draw [shift={(158,30.7)}, rotate = 203.33] [color={rgb, 255:red, 208; green, 2; blue, 27 }  ,draw opacity=1 ][line width=0.75]    (6.56,-1.97) .. controls (4.17,-0.84) and (1.99,-0.18) .. (0,0) .. controls (1.99,0.18) and (4.17,0.84) .. (6.56,1.97)   ;

\draw (42.5,28.5) node  [font=\scriptsize]  {$8$};
\draw (42.5,48.5) node  [font=\scriptsize]  {$6$};
\draw (120.5,95.5) node    {$\Delta \ =\ 6$};
\draw (42.82,8.4) node  [font=\footnotesize]  {$a$};
\draw (200.82,69.2) node  [font=\footnotesize]  {$q$};
\draw (122.5,70.5) node  [font=\scriptsize]  {$0$};
\draw (100,70.5) node  [font=\scriptsize]  {$-2$};
\draw (162.5,70.5) node  [font=\scriptsize]  {$4$};
\draw (142.5,70.5) node  [font=\scriptsize]  {$2$};
\draw (182.5,70.5) node  [font=\scriptsize]  {$6$};
\draw (80,70.5) node  [font=\scriptsize]  {$-4$};
\draw (60,70.5) node  [font=\scriptsize]  {$-6$};

\end{tikzpicture}

%% file: tikzpictures/table-Kh-7_1.tex
\tikzset{every picture/.style={line width=0.75pt}} 

\begin{tikzpicture}[x=0.75pt,y=0.75pt,yscale=-.7,xscale=.7]

\draw  [draw opacity=0] (50,20) -- (210.2,20) -- (210.2,180.6) -- (50,180.6) -- cycle ; \draw  [color={rgb, 255:red, 155; green, 155; blue, 155 }  ,draw opacity=1 ] (50,20) -- (50,180.6)(70,20) -- (70,180.6)(90,20) -- (90,180.6)(110,20) -- (110,180.6)(130,20) -- (130,180.6)(150,20) -- (150,180.6)(170,20) -- (170,180.6)(190,20) -- (190,180.6)(210,20) -- (210,180.6) ; \draw  [color={rgb, 255:red, 155; green, 155; blue, 155 }  ,draw opacity=1 ] (50,20) -- (210.2,20)(50,40) -- (210.2,40)(50,60) -- (210.2,60)(50,80) -- (210.2,80)(50,100) -- (210.2,100)(50,120) -- (210.2,120)(50,140) -- (210.2,140)(50,160) -- (210.2,160)(50,180) -- (210.2,180) ; \draw  [color={rgb, 255:red, 155; green, 155; blue, 155 }  ,draw opacity=1 ]  ;
\draw  [color={rgb, 255:red, 245; green, 166; blue, 35 }  ,draw opacity=1 ][fill={rgb, 255:red, 0; green, 0; blue, 0 }  ,fill opacity=1 ] (200.19,174.72) .. controls (198.99,175.93) and (197.04,175.93) .. (195.84,174.72) .. controls (194.63,173.52) and (194.63,171.57) .. (195.84,170.37) .. controls (197.04,169.16) and (198.99,169.16) .. (200.19,170.37) .. controls (201.4,171.57) and (201.4,173.52) .. (200.19,174.72) -- cycle ;
\draw  [color={rgb, 255:red, 245; green, 166; blue, 35 }  ,draw opacity=1 ][fill={rgb, 255:red, 0; green, 0; blue, 0 }  ,fill opacity=1 ] (160.59,135.12) .. controls (159.39,136.32) and (157.44,136.32) .. (156.23,135.12) .. controls (155.03,133.92) and (155.03,131.97) .. (156.23,130.76) .. controls (157.44,129.56) and (159.39,129.56) .. (160.59,130.76) .. controls (161.79,131.97) and (161.79,133.92) .. (160.59,135.12) -- cycle ;
\draw  [color={rgb, 255:red, 245; green, 166; blue, 35 }  ,draw opacity=1 ][fill={rgb, 255:red, 0; green, 0; blue, 0 }  ,fill opacity=1 ] (121.19,95.71) .. controls (119.98,96.92) and (118.03,96.92) .. (116.83,95.71) .. controls (115.63,94.51) and (115.63,92.56) .. (116.83,91.36) .. controls (118.03,90.15) and (119.98,90.15) .. (121.19,91.36) .. controls (122.39,92.56) and (122.39,94.51) .. (121.19,95.71) -- cycle ;
\draw  [color={rgb, 255:red, 245; green, 166; blue, 35 }  ,draw opacity=1 ][fill={rgb, 255:red, 0; green, 0; blue, 0 }  ,fill opacity=1 ] (80.59,55.12) .. controls (79.39,56.32) and (77.44,56.32) .. (76.24,55.12) .. controls (75.03,53.92) and (75.03,51.97) .. (76.24,50.76) .. controls (77.44,49.56) and (79.39,49.56) .. (80.59,50.76) .. controls (81.79,51.97) and (81.79,53.92) .. (80.59,55.12) -- cycle ;
\draw [color={rgb, 255:red, 245; green, 166; blue, 35 }  ,draw opacity=1 ][line width=0.75]    (195.84,170.37) .. controls (176.05,163.9) and (169.91,158.09) .. (161.26,136.79) ;
\draw [shift={(160.59,135.12)}, rotate = 68.33] [color={rgb, 255:red, 245; green, 166; blue, 35 }  ,draw opacity=1 ][line width=0.75]    (6.56,-1.97) .. controls (4.17,-0.84) and (1.99,-0.18) .. (0,0) .. controls (1.99,0.18) and (4.17,0.84) .. (6.56,1.97)   ;
\draw [color={rgb, 255:red, 245; green, 166; blue, 35 }  ,draw opacity=1 ][line width=0.75]    (156.23,130.76) .. controls (136.45,124.3) and (130.49,118.67) .. (121.86,97.38) ;
\draw [shift={(121.19,95.71)}, rotate = 68.33] [color={rgb, 255:red, 245; green, 166; blue, 35 }  ,draw opacity=1 ][line width=0.75]    (6.56,-1.97) .. controls (4.17,-0.84) and (1.99,-0.18) .. (0,0) .. controls (1.99,0.18) and (4.17,0.84) .. (6.56,1.97)   ;
\draw [color={rgb, 255:red, 245; green, 166; blue, 35 }  ,draw opacity=1 ][line width=0.75]    (116.83,91.36) .. controls (97.04,84.89) and (89.96,78.14) .. (81.26,56.79) ;
\draw [shift={(80.59,55.12)}, rotate = 68.33] [color={rgb, 255:red, 245; green, 166; blue, 35 }  ,draw opacity=1 ][line width=0.75]    (6.56,-1.97) .. controls (4.17,-0.84) and (1.99,-0.18) .. (0,0) .. controls (1.99,0.18) and (4.17,0.84) .. (6.56,1.97)   ;

\draw [color={rgb, 255:red, 208; green, 2; blue, 27 }  ,draw opacity=1 ][line width=0.75]    (138.89,105.57) .. controls (132.4,85.71) and (126.57,79.55) .. (105.2,70.88) ;
\draw [shift={(103.53,70.21)}, rotate = 21.67] [color={rgb, 255:red, 208; green, 2; blue, 27 }  ,draw opacity=1 ][line width=0.75]    (6.56,-1.97) .. controls (4.17,-0.84) and (1.99,-0.18) .. (0,0) .. controls (1.99,0.18) and (4.17,0.84) .. (6.56,1.97)   ;
\draw  [color={rgb, 255:red, 208; green, 2; blue, 27 }  ,draw opacity=1 ][fill={rgb, 255:red, 0; green, 0; blue, 0 }  ,fill opacity=1 ] (143.26,109.94) .. controls (144.47,108.73) and (144.47,106.78) .. (143.26,105.57) .. controls (142.05,104.36) and (140.1,104.36) .. (138.89,105.57) .. controls (137.68,106.78) and (137.68,108.73) .. (138.89,109.94) .. controls (140.1,111.15) and (142.05,111.15) .. (143.26,109.94) -- cycle ;
\draw  [color={rgb, 255:red, 208; green, 2; blue, 27 }  ,draw opacity=1 ][fill={rgb, 255:red, 0; green, 0; blue, 0 }  ,fill opacity=1 ] (103.53,70.21) .. controls (104.73,69) and (104.73,67.04) .. (103.53,65.83) .. controls (102.32,64.63) and (100.36,64.63) .. (99.16,65.83) .. controls (97.95,67.04) and (97.95,69) .. (99.16,70.21) .. controls (100.36,71.41) and (102.32,71.41) .. (103.53,70.21) -- cycle ;
\draw  [color={rgb, 255:red, 208; green, 2; blue, 27 }  ,draw opacity=1 ][fill={rgb, 255:red, 0; green, 0; blue, 0 }  ,fill opacity=1 ] (63.99,30.67) .. controls (65.2,29.46) and (65.2,27.51) .. (63.99,26.3) .. controls (62.78,25.09) and (60.83,25.09) .. (59.62,26.3) .. controls (58.41,27.51) and (58.41,29.46) .. (59.62,30.67) .. controls (60.83,31.88) and (62.78,31.88) .. (63.99,30.67) -- cycle ;
\draw [color={rgb, 255:red, 208; green, 2; blue, 27 }  ,draw opacity=1 ][line width=0.75]    (99.16,65.83) .. controls (92.67,45.98) and (87.03,40.01) .. (65.67,31.34) ;
\draw [shift={(63.99,30.67)}, rotate = 21.67] [color={rgb, 255:red, 208; green, 2; blue, 27 }  ,draw opacity=1 ][line width=0.75]    (6.56,-1.97) .. controls (4.17,-0.84) and (1.99,-0.18) .. (0,0) .. controls (1.99,0.18) and (4.17,0.84) .. (6.56,1.97)   ;

\draw (41.2,109.84) node  [font=\tiny]  {$12$};
\draw (41.2,129.95) node  [font=\tiny]  {$10$};
\draw (41.2,69.62) node  [font=\tiny]  {$16$};
\draw (41.2,89.73) node  [font=\tiny]  {$14$};
\draw (41.2,49.51) node  [font=\tiny]  {$18$};
\draw (42.2,150.06) node  [font=\tiny]  {$8$};
\draw (42.2,170.2) node  [font=\tiny]  {$6$};
\draw (41.2,29.4) node  [font=\tiny]  {$20$};
\draw (119.87,190) node  [font=\tiny]  {$-4$};
\draw (99.58,190) node  [font=\tiny]  {$-5$};
\draw (160.45,190) node  [font=\tiny]  {$-2$};
\draw (140.16,190) node  [font=\tiny]  {$-3$};
\draw (180.74,190) node  [font=\tiny]  {$-1$};
\draw (79.29,190) node  [font=\tiny]  {$-6$};
\draw (59,190) node  [font=\tiny]  {$-7$};
\draw (202.8,190) node  [font=\tiny]  {$0$};
\draw (223.9,188.8) node  [font=\scriptsize]  {$t$};
\draw (41.7,7.6) node  [font=\scriptsize]  {$q_{2}$};

\end{tikzpicture}

%% file: tikzpictures/table-Kh-7_1-slice.tex
\tikzset{every picture/.style={line width=0.75pt}} 

\begin{tikzpicture}[x=0.75pt,y=0.75pt,yscale=-1,xscale=1]

\draw  [draw opacity=0] (50,30) -- (210.2,30) -- (210.2,50.1) -- (50,50.1) -- cycle ; \draw  [color={rgb, 255:red, 155; green, 155; blue, 155 }  ,draw opacity=1 ] (50,30) -- (50,50.1)(70,30) -- (70,50.1)(90,30) -- (90,50.1)(110,30) -- (110,50.1)(130,30) -- (130,50.1)(150,30) -- (150,50.1)(170,30) -- (170,50.1)(190,30) -- (190,50.1)(210,30) -- (210,50.1) ; \draw  [color={rgb, 255:red, 155; green, 155; blue, 155 }  ,draw opacity=1 ] (50,30) -- (210.2,30)(50,50) -- (210.2,50) ; \draw  [color={rgb, 255:red, 155; green, 155; blue, 155 }  ,draw opacity=1 ]  ;
\draw  [color={rgb, 255:red, 245; green, 166; blue, 35 }  ,draw opacity=1 ][fill={rgb, 255:red, 0; green, 0; blue, 0 }  ,fill opacity=1 ] (58,38.17) .. controls (58,36.96) and (58.98,35.97) .. (60.2,35.97) .. controls (61.42,35.97) and (62.4,36.96) .. (62.4,38.17) .. controls (62.4,39.39) and (61.42,40.38) .. (60.2,40.38) .. controls (58.98,40.38) and (58,39.39) .. (58,38.17) -- cycle ;
\draw  [color={rgb, 255:red, 245; green, 166; blue, 35 }  ,draw opacity=1 ][fill={rgb, 255:red, 0; green, 0; blue, 0 }  ,fill opacity=1 ] (98,38.17) .. controls (98,36.96) and (98.98,35.97) .. (100.2,35.97) .. controls (101.42,35.97) and (102.4,36.96) .. (102.4,38.17) .. controls (102.4,39.39) and (101.42,40.38) .. (100.2,40.38) .. controls (98.98,40.38) and (98,39.39) .. (98,38.17) -- cycle ;
\draw  [color={rgb, 255:red, 245; green, 166; blue, 35 }  ,draw opacity=1 ][fill={rgb, 255:red, 0; green, 0; blue, 0 }  ,fill opacity=1 ] (137.8,38.17) .. controls (137.8,36.96) and (138.78,35.97) .. (140,35.97) .. controls (141.22,35.97) and (142.2,36.96) .. (142.2,38.17) .. controls (142.2,39.39) and (141.22,40.38) .. (140,40.38) .. controls (138.78,40.38) and (137.8,39.39) .. (137.8,38.17) -- cycle ;
\draw  [color={rgb, 255:red, 245; green, 166; blue, 35 }  ,draw opacity=1 ][fill={rgb, 255:red, 0; green, 0; blue, 0 }  ,fill opacity=1 ] (178.8,38.17) .. controls (178.8,36.96) and (179.78,35.97) .. (181,35.97) .. controls (182.22,35.97) and (183.2,36.96) .. (183.2,38.17) .. controls (183.2,39.39) and (182.22,40.38) .. (181,40.38) .. controls (179.78,40.38) and (178.8,39.39) .. (178.8,38.17) -- cycle ;
\draw [color={rgb, 255:red, 245; green, 166; blue, 35 }  ,draw opacity=1 ]   (62.4,38.17) .. controls (75.52,31.52) and (81.57,31.28) .. (96.36,37.48) ;
\draw [shift={(98,38.17)}, rotate = 203.33] [color={rgb, 255:red, 245; green, 166; blue, 35 }  ,draw opacity=1 ][line width=0.75]    (6.56,-1.97) .. controls (4.17,-0.84) and (1.99,-0.18) .. (0,0) .. controls (1.99,0.18) and (4.17,0.84) .. (6.56,1.97)   ;
\draw [color={rgb, 255:red, 245; green, 166; blue, 35 }  ,draw opacity=1 ]   (102.4,38.17) .. controls (115.52,31.52) and (121.38,31.28) .. (136.16,37.48) ;
\draw [shift={(137.8,38.17)}, rotate = 203.33] [color={rgb, 255:red, 245; green, 166; blue, 35 }  ,draw opacity=1 ][line width=0.75]    (6.56,-1.97) .. controls (4.17,-0.84) and (1.99,-0.18) .. (0,0) .. controls (1.99,0.18) and (4.17,0.84) .. (6.56,1.97)   ;
\draw [color={rgb, 255:red, 245; green, 166; blue, 35 }  ,draw opacity=1 ]   (142.2,38.17) .. controls (155.32,31.52) and (162.3,31.28) .. (177.15,37.48) ;
\draw [shift={(178.8,38.17)}, rotate = 203.33] [color={rgb, 255:red, 245; green, 166; blue, 35 }  ,draw opacity=1 ][line width=0.75]    (6.56,-1.97) .. controls (4.17,-0.84) and (1.99,-0.18) .. (0,0) .. controls (1.99,0.18) and (4.17,0.84) .. (6.56,1.97)   ;

\draw [color={rgb, 255:red, 208; green, 2; blue, 27 }  ,draw opacity=1 ]   (122.6,43.72) .. controls (135.72,50.38) and (141.77,50.62) .. (156.56,44.42) ;
\draw [shift={(158.2,43.72)}, rotate = 156.67] [color={rgb, 255:red, 208; green, 2; blue, 27 }  ,draw opacity=1 ][line width=0.75]    (6.56,-1.97) .. controls (4.17,-0.84) and (1.99,-0.18) .. (0,0) .. controls (1.99,0.18) and (4.17,0.84) .. (6.56,1.97)   ;
\draw  [color={rgb, 255:red, 208; green, 2; blue, 27 }  ,draw opacity=1 ][fill={rgb, 255:red, 0; green, 0; blue, 0 }  ,fill opacity=1 ] (118.2,43.73) .. controls (118.2,44.94) and (119.18,45.93) .. (120.4,45.93) .. controls (121.62,45.93) and (122.6,44.94) .. (122.6,43.73) .. controls (122.6,42.51) and (121.62,41.53) .. (120.4,41.53) .. controls (119.18,41.53) and (118.2,42.51) .. (118.2,43.73) -- cycle ;
\draw  [color={rgb, 255:red, 208; green, 2; blue, 27 }  ,draw opacity=1 ][fill={rgb, 255:red, 0; green, 0; blue, 0 }  ,fill opacity=1 ] (158.2,43.72) .. controls (158.2,44.94) and (159.18,45.93) .. (160.4,45.93) .. controls (161.62,45.93) and (162.6,44.94) .. (162.6,43.72) .. controls (162.6,42.51) and (161.62,41.52) .. (160.4,41.52) .. controls (159.18,41.52) and (158.2,42.51) .. (158.2,43.72) -- cycle ;
\draw  [color={rgb, 255:red, 208; green, 2; blue, 27 }  ,draw opacity=1 ][fill={rgb, 255:red, 0; green, 0; blue, 0 }  ,fill opacity=1 ] (198,43.72) .. controls (198,44.94) and (198.98,45.93) .. (200.2,45.93) .. controls (201.42,45.93) and (202.4,44.94) .. (202.4,43.72) .. controls (202.4,42.51) and (201.42,41.52) .. (200.2,41.52) .. controls (198.98,41.52) and (198,42.51) .. (198,43.72) -- cycle ;
\draw [color={rgb, 255:red, 208; green, 2; blue, 27 }  ,draw opacity=1 ]   (162.6,43.72) .. controls (175.72,50.38) and (181.58,50.62) .. (196.36,44.42) ;
\draw [shift={(198,43.72)}, rotate = 156.67] [color={rgb, 255:red, 208; green, 2; blue, 27 }  ,draw opacity=1 ][line width=0.75]    (6.56,-1.97) .. controls (4.17,-0.84) and (1.99,-0.18) .. (0,0) .. controls (1.99,0.18) and (4.17,0.84) .. (6.56,1.97)   ;

\draw (119.5,60) node  [font=\scriptsize]  {$12$};
\draw (99.5,60) node  [font=\scriptsize]  {$10$};
\draw (159.5,60) node  [font=\scriptsize]  {$16$};
\draw (139.5,60) node  [font=\scriptsize]  {$14$};
\draw (179.5,60) node  [font=\scriptsize]  {$18$};
\draw (82.5,60) node  [font=\scriptsize]  {$8$};
\draw (62.5,60) node  [font=\scriptsize]  {$6$};
\draw (199.5,60) node  [font=\scriptsize]  {$20$};
\draw (43,19) node  [font=\footnotesize]  {$\delta $};
\draw (220.8,58.2) node  [font=\footnotesize]  {$q_{2}$};
\draw (43,39) node  [font=\scriptsize]  {$6$};

\end{tikzpicture}

%% file: tikzpictures/table-H-conway.tex
\tikzset{every picture/.style={line width=0.75pt}} 

\begin{tikzpicture}[x=0.75pt,y=0.75pt,yscale=-1.2,xscale=1.2]

\draw  [draw opacity=0] (50,19.55) -- (190.2,19.55) -- (190.2,80.2) -- (50,80.2) -- cycle ; \draw  [color={rgb, 255:red, 155; green, 155; blue, 155 }  ,draw opacity=1 ] (50,19.55) -- (50,80.2)(70,19.55) -- (70,80.2)(90,19.55) -- (90,80.2)(110,19.55) -- (110,80.2)(130,19.55) -- (130,80.2)(150,19.55) -- (150,80.2)(170,19.55) -- (170,80.2)(190,19.55) -- (190,80.2) ; \draw  [color={rgb, 255:red, 155; green, 155; blue, 155 }  ,draw opacity=1 ] (50,19.55) -- (190.2,19.55)(50,39.55) -- (190.2,39.55)(50,59.55) -- (190.2,59.55)(50,79.55) -- (190.2,79.55) ; \draw  [color={rgb, 255:red, 155; green, 155; blue, 155 }  ,draw opacity=1 ]  ;
\draw  [color={rgb, 255:red, 245; green, 166; blue, 35 }  ,draw opacity=1 ][fill={rgb, 255:red, 0; green, 0; blue, 0 }  ,fill opacity=1 ] (57.8,44.3) .. controls (57.8,43.08) and (58.78,42.1) .. (60,42.1) .. controls (61.22,42.1) and (62.2,43.08) .. (62.2,44.3) .. controls (62.2,45.52) and (61.22,46.5) .. (60,46.5) .. controls (58.78,46.5) and (57.8,45.52) .. (57.8,44.3) -- cycle ;
\draw  [color={rgb, 255:red, 245; green, 166; blue, 35 }  ,draw opacity=1 ][fill={rgb, 255:red, 0; green, 0; blue, 0 }  ,fill opacity=1 ] (97.8,44.3) .. controls (97.8,43.08) and (98.78,42.1) .. (100,42.1) .. controls (101.22,42.1) and (102.2,43.08) .. (102.2,44.3) .. controls (102.2,45.52) and (101.22,46.5) .. (100,46.5) .. controls (98.78,46.5) and (97.8,45.52) .. (97.8,44.3) -- cycle ;
\draw  [color={rgb, 255:red, 245; green, 166; blue, 35 }  ,draw opacity=1 ][fill={rgb, 255:red, 0; green, 0; blue, 0 }  ,fill opacity=1 ] (137.6,44.3) .. controls (137.6,43.08) and (138.58,42.1) .. (139.8,42.1) .. controls (141.02,42.1) and (142,43.08) .. (142,44.3) .. controls (142,45.52) and (141.02,46.5) .. (139.8,46.5) .. controls (138.58,46.5) and (137.6,45.52) .. (137.6,44.3) -- cycle ;
\draw  [color={rgb, 255:red, 245; green, 166; blue, 35 }  ,draw opacity=1 ][fill={rgb, 255:red, 0; green, 0; blue, 0 }  ,fill opacity=1 ] (178.6,44.3) .. controls (178.6,43.08) and (179.58,42.1) .. (180.8,42.1) .. controls (182.02,42.1) and (183,43.08) .. (183,44.3) .. controls (183,45.52) and (182.02,46.5) .. (180.8,46.5) .. controls (179.58,46.5) and (178.6,45.52) .. (178.6,44.3) -- cycle ;
\draw [color={rgb, 255:red, 245; green, 166; blue, 35 }  ,draw opacity=1 ]   (62.2,44.3) .. controls (75.32,37.64) and (81.37,37.41) .. (96.16,43.6) ;
\draw [shift={(97.8,44.3)}, rotate = 203.33] [color={rgb, 255:red, 245; green, 166; blue, 35 }  ,draw opacity=1 ][line width=0.75]    (6.56,-1.97) .. controls (4.17,-0.84) and (1.99,-0.18) .. (0,0) .. controls (1.99,0.18) and (4.17,0.84) .. (6.56,1.97)   ;
\draw [color={rgb, 255:red, 245; green, 166; blue, 35 }  ,draw opacity=1 ]   (102.2,44.3) .. controls (115.32,37.64) and (121.18,37.41) .. (135.96,43.6) ;
\draw [shift={(137.6,44.3)}, rotate = 203.33] [color={rgb, 255:red, 245; green, 166; blue, 35 }  ,draw opacity=1 ][line width=0.75]    (6.56,-1.97) .. controls (4.17,-0.84) and (1.99,-0.18) .. (0,0) .. controls (1.99,0.18) and (4.17,0.84) .. (6.56,1.97)   ;
\draw [color={rgb, 255:red, 245; green, 166; blue, 35 }  ,draw opacity=1 ]   (142,44.3) .. controls (155.12,37.64) and (162.1,37.41) .. (176.95,43.6) ;
\draw [shift={(178.6,44.3)}, rotate = 203.33] [color={rgb, 255:red, 245; green, 166; blue, 35 }  ,draw opacity=1 ][line width=0.75]    (6.56,-1.97) .. controls (4.17,-0.84) and (1.99,-0.18) .. (0,0) .. controls (1.99,0.18) and (4.17,0.84) .. (6.56,1.97)   ;

\draw [color={rgb, 255:red, 208; green, 2; blue, 27 }  ,draw opacity=1 ]   (82.6,25.7) .. controls (95.72,19.04) and (101.77,18.81) .. (116.56,25) ;
\draw [shift={(118.2,25.7)}, rotate = 203.33] [color={rgb, 255:red, 208; green, 2; blue, 27 }  ,draw opacity=1 ][line width=0.75]    (6.56,-1.97) .. controls (4.17,-0.84) and (1.99,-0.18) .. (0,0) .. controls (1.99,0.18) and (4.17,0.84) .. (6.56,1.97)   ;
\draw  [color={rgb, 255:red, 208; green, 2; blue, 27 }  ,draw opacity=1 ][fill={rgb, 255:red, 0; green, 0; blue, 0 }  ,fill opacity=1 ] (78.2,25.7) .. controls (78.2,24.48) and (79.18,23.5) .. (80.4,23.5) .. controls (81.62,23.5) and (82.6,24.48) .. (82.6,25.7) .. controls (82.6,26.92) and (81.62,27.9) .. (80.4,27.9) .. controls (79.18,27.9) and (78.2,26.92) .. (78.2,25.7) -- cycle ;
\draw  [color={rgb, 255:red, 208; green, 2; blue, 27 }  ,draw opacity=1 ][fill={rgb, 255:red, 0; green, 0; blue, 0 }  ,fill opacity=1 ] (118.2,25.7) .. controls (118.2,24.48) and (119.18,23.5) .. (120.4,23.5) .. controls (121.62,23.5) and (122.6,24.48) .. (122.6,25.7) .. controls (122.6,26.92) and (121.62,27.9) .. (120.4,27.9) .. controls (119.18,27.9) and (118.2,26.92) .. (118.2,25.7) -- cycle ;
\draw  [color={rgb, 255:red, 208; green, 2; blue, 27 }  ,draw opacity=1 ][fill={rgb, 255:red, 0; green, 0; blue, 0 }  ,fill opacity=1 ] (158,25.7) .. controls (158,24.48) and (158.98,23.5) .. (160.2,23.5) .. controls (161.42,23.5) and (162.4,24.48) .. (162.4,25.7) .. controls (162.4,26.92) and (161.42,27.9) .. (160.2,27.9) .. controls (158.98,27.9) and (158,26.92) .. (158,25.7) -- cycle ;
\draw [color={rgb, 255:red, 208; green, 2; blue, 27 }  ,draw opacity=1 ]   (122.6,25.7) .. controls (135.72,19.04) and (141.58,18.81) .. (156.36,25) ;
\draw [shift={(158,25.7)}, rotate = 203.33] [color={rgb, 255:red, 208; green, 2; blue, 27 }  ,draw opacity=1 ][line width=0.75]    (6.56,-1.97) .. controls (4.17,-0.84) and (1.99,-0.18) .. (0,0) .. controls (1.99,0.18) and (4.17,0.84) .. (6.56,1.97)   ;

\draw  [color={rgb, 255:red, 208; green, 2; blue, 27 }  ,draw opacity=1 ][fill={rgb, 255:red, 0; green, 0; blue, 0 }  ,fill opacity=1 ] (118.2,32.5) .. controls (118.2,31.28) and (119.18,30.3) .. (120.4,30.3) .. controls (121.62,30.3) and (122.6,31.28) .. (122.6,32.5) .. controls (122.6,33.72) and (121.62,34.7) .. (120.4,34.7) .. controls (119.18,34.7) and (118.2,33.72) .. (118.2,32.5) -- cycle ;
\draw [color={rgb, 255:red, 208; green, 2; blue, 27 }  ,draw opacity=1 ]   (102.6,31.33) .. controls (115.72,37.98) and (121.77,38.22) .. (136.56,32.02) ;
\draw [shift={(138.2,31.33)}, rotate = 156.67] [color={rgb, 255:red, 208; green, 2; blue, 27 }  ,draw opacity=1 ][line width=0.75]    (6.56,-1.97) .. controls (4.17,-0.84) and (1.99,-0.18) .. (0,0) .. controls (1.99,0.18) and (4.17,0.84) .. (6.56,1.97)   ;
\draw  [color={rgb, 255:red, 208; green, 2; blue, 27 }  ,draw opacity=1 ][fill={rgb, 255:red, 0; green, 0; blue, 0 }  ,fill opacity=1 ] (98.2,31.33) .. controls (98.2,32.54) and (99.18,33.53) .. (100.4,33.53) .. controls (101.62,33.53) and (102.6,32.54) .. (102.6,31.33) .. controls (102.6,30.11) and (101.62,29.13) .. (100.4,29.13) .. controls (99.18,29.13) and (98.2,30.11) .. (98.2,31.33) -- cycle ;
\draw  [color={rgb, 255:red, 208; green, 2; blue, 27 }  ,draw opacity=1 ][fill={rgb, 255:red, 0; green, 0; blue, 0 }  ,fill opacity=1 ] (138.2,31.33) .. controls (138.2,32.54) and (139.18,33.53) .. (140.4,33.53) .. controls (141.62,33.53) and (142.6,32.54) .. (142.6,31.33) .. controls (142.6,30.11) and (141.62,29.13) .. (140.4,29.13) .. controls (139.18,29.13) and (138.2,30.11) .. (138.2,31.33) -- cycle ;

\draw [color={rgb, 255:red, 245; green, 166; blue, 35 }  ,draw opacity=1 ]   (102.4,57.13) .. controls (115.52,63.78) and (121.57,64.02) .. (136.36,57.82) ;
\draw [shift={(138,57.13)}, rotate = 156.67] [color={rgb, 255:red, 245; green, 166; blue, 35 }  ,draw opacity=1 ][line width=0.75]    (6.56,-1.97) .. controls (4.17,-0.84) and (1.99,-0.18) .. (0,0) .. controls (1.99,0.18) and (4.17,0.84) .. (6.56,1.97)   ;
\draw  [color={rgb, 255:red, 245; green, 166; blue, 35 }  ,draw opacity=1 ][fill={rgb, 255:red, 0; green, 0; blue, 0 }  ,fill opacity=1 ] (98,57.13) .. controls (98,58.34) and (98.98,59.33) .. (100.2,59.33) .. controls (101.42,59.33) and (102.4,58.34) .. (102.4,57.13) .. controls (102.4,55.91) and (101.42,54.93) .. (100.2,54.93) .. controls (98.98,54.93) and (98,55.91) .. (98,57.13) -- cycle ;
\draw  [color={rgb, 255:red, 245; green, 166; blue, 35 }  ,draw opacity=1 ][fill={rgb, 255:red, 0; green, 0; blue, 0 }  ,fill opacity=1 ] (138,57.13) .. controls (138,58.34) and (138.98,59.33) .. (140.2,59.33) .. controls (141.42,59.33) and (142.4,58.34) .. (142.4,57.13) .. controls (142.4,55.91) and (141.42,54.93) .. (140.2,54.93) .. controls (138.98,54.93) and (138,55.91) .. (138,57.13) -- cycle ;

\draw [color={rgb, 255:red, 245; green, 166; blue, 35 }  ,draw opacity=1 ]   (82.8,50.7) .. controls (97.28,46.55) and (100.02,46.41) .. (116.55,50.26) ;
\draw [shift={(118.4,50.7)}, rotate = 193.29] [color={rgb, 255:red, 245; green, 166; blue, 35 }  ,draw opacity=1 ][line width=0.75]    (6.56,-1.97) .. controls (4.17,-0.84) and (1.99,-0.18) .. (0,0) .. controls (1.99,0.18) and (4.17,0.84) .. (6.56,1.97)   ;
\draw  [color={rgb, 255:red, 245; green, 166; blue, 35 }  ,draw opacity=1 ][fill={rgb, 255:red, 0; green, 0; blue, 0 }  ,fill opacity=1 ] (78.4,50.7) .. controls (78.4,49.48) and (79.38,48.5) .. (80.6,48.5) .. controls (81.82,48.5) and (82.8,49.48) .. (82.8,50.7) .. controls (82.8,51.92) and (81.82,52.9) .. (80.6,52.9) .. controls (79.38,52.9) and (78.4,51.92) .. (78.4,50.7) -- cycle ;
\draw  [color={rgb, 255:red, 245; green, 166; blue, 35 }  ,draw opacity=1 ][fill={rgb, 255:red, 0; green, 0; blue, 0 }  ,fill opacity=1 ] (118.4,50.7) .. controls (118.4,49.48) and (119.38,48.5) .. (120.6,48.5) .. controls (121.82,48.5) and (122.8,49.48) .. (122.8,50.7) .. controls (122.8,51.92) and (121.82,52.9) .. (120.6,52.9) .. controls (119.38,52.9) and (118.4,51.92) .. (118.4,50.7) -- cycle ;
\draw  [color={rgb, 255:red, 245; green, 166; blue, 35 }  ,draw opacity=1 ][fill={rgb, 255:red, 0; green, 0; blue, 0 }  ,fill opacity=1 ] (158.2,50.7) .. controls (158.2,49.48) and (159.18,48.5) .. (160.4,48.5) .. controls (161.62,48.5) and (162.6,49.48) .. (162.6,50.7) .. controls (162.6,51.92) and (161.62,52.9) .. (160.4,52.9) .. controls (159.18,52.9) and (158.2,51.92) .. (158.2,50.7) -- cycle ;
\draw [color={rgb, 255:red, 245; green, 166; blue, 35 }  ,draw opacity=1 ]   (122.8,50.7) .. controls (137.28,46.55) and (139.83,46.41) .. (156.35,50.26) ;
\draw [shift={(158.2,50.7)}, rotate = 193.29] [color={rgb, 255:red, 245; green, 166; blue, 35 }  ,draw opacity=1 ][line width=0.75]    (6.56,-1.97) .. controls (4.17,-0.84) and (1.99,-0.18) .. (0,0) .. controls (1.99,0.18) and (4.17,0.84) .. (6.56,1.97)   ;

\draw [color={rgb, 255:red, 245; green, 166; blue, 35 }  ,draw opacity=1 ]   (102.2,50.73) .. controls (115.32,57.38) and (121.37,57.62) .. (136.16,51.42) ;
\draw [shift={(137.8,50.73)}, rotate = 156.67] [color={rgb, 255:red, 245; green, 166; blue, 35 }  ,draw opacity=1 ][line width=0.75]    (6.56,-1.97) .. controls (4.17,-0.84) and (1.99,-0.18) .. (0,0) .. controls (1.99,0.18) and (4.17,0.84) .. (6.56,1.97)   ;
\draw  [color={rgb, 255:red, 245; green, 166; blue, 35 }  ,draw opacity=1 ][fill={rgb, 255:red, 0; green, 0; blue, 0 }  ,fill opacity=1 ] (97.8,50.73) .. controls (97.8,51.94) and (98.78,52.93) .. (100,52.93) .. controls (101.22,52.93) and (102.2,51.94) .. (102.2,50.73) .. controls (102.2,49.51) and (101.22,48.53) .. (100,48.53) .. controls (98.78,48.53) and (97.8,49.51) .. (97.8,50.73) -- cycle ;
\draw  [color={rgb, 255:red, 245; green, 166; blue, 35 }  ,draw opacity=1 ][fill={rgb, 255:red, 0; green, 0; blue, 0 }  ,fill opacity=1 ] (137.8,50.73) .. controls (137.8,51.94) and (138.78,52.93) .. (140,52.93) .. controls (141.22,52.93) and (142.2,51.94) .. (142.2,50.73) .. controls (142.2,49.51) and (141.22,48.53) .. (140,48.53) .. controls (138.78,48.53) and (137.8,49.51) .. (137.8,50.73) -- cycle ;

\draw [color={rgb, 255:red, 126; green, 211; blue, 33 }  ,draw opacity=1 ]   (82.2,66.9) .. controls (96.68,62.75) and (99.42,62.61) .. (115.95,66.46) ;
\draw [shift={(117.8,66.9)}, rotate = 193.29] [color={rgb, 255:red, 126; green, 211; blue, 33 }  ,draw opacity=1 ][line width=0.75]    (6.56,-1.97) .. controls (4.17,-0.84) and (1.99,-0.18) .. (0,0) .. controls (1.99,0.18) and (4.17,0.84) .. (6.56,1.97)   ;
\draw  [color={rgb, 255:red, 126; green, 211; blue, 33 }  ,draw opacity=1 ][fill={rgb, 255:red, 0; green, 0; blue, 0 }  ,fill opacity=1 ] (77.8,66.9) .. controls (77.8,65.68) and (78.78,64.7) .. (80,64.7) .. controls (81.22,64.7) and (82.2,65.68) .. (82.2,66.9) .. controls (82.2,68.12) and (81.22,69.1) .. (80,69.1) .. controls (78.78,69.1) and (77.8,68.12) .. (77.8,66.9) -- cycle ;
\draw  [color={rgb, 255:red, 126; green, 211; blue, 33 }  ,draw opacity=1 ][fill={rgb, 255:red, 0; green, 0; blue, 0 }  ,fill opacity=1 ] (117.8,66.9) .. controls (117.8,65.68) and (118.78,64.7) .. (120,64.7) .. controls (121.22,64.7) and (122.2,65.68) .. (122.2,66.9) .. controls (122.2,68.12) and (121.22,69.1) .. (120,69.1) .. controls (118.78,69.1) and (117.8,68.12) .. (117.8,66.9) -- cycle ;
\draw  [color={rgb, 255:red, 126; green, 211; blue, 33 }  ,draw opacity=1 ][fill={rgb, 255:red, 0; green, 0; blue, 0 }  ,fill opacity=1 ] (157.6,66.9) .. controls (157.6,65.68) and (158.58,64.7) .. (159.8,64.7) .. controls (161.02,64.7) and (162,65.68) .. (162,66.9) .. controls (162,68.12) and (161.02,69.1) .. (159.8,69.1) .. controls (158.58,69.1) and (157.6,68.12) .. (157.6,66.9) -- cycle ;
\draw [color={rgb, 255:red, 126; green, 211; blue, 33 }  ,draw opacity=1 ]   (122.2,66.9) .. controls (136.68,62.75) and (139.23,62.61) .. (155.75,66.46) ;
\draw [shift={(157.6,66.9)}, rotate = 193.29] [color={rgb, 255:red, 126; green, 211; blue, 33 }  ,draw opacity=1 ][line width=0.75]    (6.56,-1.97) .. controls (4.17,-0.84) and (1.99,-0.18) .. (0,0) .. controls (1.99,0.18) and (4.17,0.84) .. (6.56,1.97)   ;

\draw  [color={rgb, 255:red, 126; green, 211; blue, 33 }  ,draw opacity=1 ][fill={rgb, 255:red, 0; green, 0; blue, 0 }  ,fill opacity=1 ] (117.8,74.3) .. controls (117.8,73.08) and (118.78,72.1) .. (120,72.1) .. controls (121.22,72.1) and (122.2,73.08) .. (122.2,74.3) .. controls (122.2,75.52) and (121.22,76.5) .. (120,76.5) .. controls (118.78,76.5) and (117.8,75.52) .. (117.8,74.3) -- cycle ;
\draw [color={rgb, 255:red, 126; green, 211; blue, 33 }  ,draw opacity=1 ]   (102.6,73.33) .. controls (115.72,79.98) and (121.77,80.22) .. (136.56,74.02) ;
\draw [shift={(138.2,73.33)}, rotate = 156.67] [color={rgb, 255:red, 126; green, 211; blue, 33 }  ,draw opacity=1 ][line width=0.75]    (6.56,-1.97) .. controls (4.17,-0.84) and (1.99,-0.18) .. (0,0) .. controls (1.99,0.18) and (4.17,0.84) .. (6.56,1.97)   ;
\draw  [color={rgb, 255:red, 126; green, 211; blue, 33 }  ,draw opacity=1 ][fill={rgb, 255:red, 0; green, 0; blue, 0 }  ,fill opacity=1 ] (98.2,73.33) .. controls (98.2,74.54) and (99.18,75.53) .. (100.4,75.53) .. controls (101.62,75.53) and (102.6,74.54) .. (102.6,73.33) .. controls (102.6,72.11) and (101.62,71.13) .. (100.4,71.13) .. controls (99.18,71.13) and (98.2,72.11) .. (98.2,73.33) -- cycle ;
\draw  [color={rgb, 255:red, 126; green, 211; blue, 33 }  ,draw opacity=1 ][fill={rgb, 255:red, 0; green, 0; blue, 0 }  ,fill opacity=1 ] (138.2,73.33) .. controls (138.2,74.54) and (139.18,75.53) .. (140.4,75.53) .. controls (141.62,75.53) and (142.6,74.54) .. (142.6,73.33) .. controls (142.6,72.11) and (141.62,71.13) .. (140.4,71.13) .. controls (139.18,71.13) and (138.2,72.11) .. (138.2,73.33) -- cycle ;

\draw  [color={rgb, 255:red, 245; green, 166; blue, 35 }  ,draw opacity=1 ][fill={rgb, 255:red, 0; green, 0; blue, 0 }  ,fill opacity=1 ] (118.4,43.3) .. controls (118.4,42.08) and (119.38,41.1) .. (120.6,41.1) .. controls (121.82,41.1) and (122.8,42.08) .. (122.8,43.3) .. controls (122.8,44.52) and (121.82,45.5) .. (120.6,45.5) .. controls (119.38,45.5) and (118.4,44.52) .. (118.4,43.3) -- cycle ;
\draw  [draw opacity=0] (250.2,19.95) -- (390.4,19.95) -- (390.4,80.6) -- (250.2,80.6) -- cycle ; \draw  [color={rgb, 255:red, 155; green, 155; blue, 155 }  ,draw opacity=1 ] (250.2,19.95) -- (250.2,80.6)(270.2,19.95) -- (270.2,80.6)(290.2,19.95) -- (290.2,80.6)(310.2,19.95) -- (310.2,80.6)(330.2,19.95) -- (330.2,80.6)(350.2,19.95) -- (350.2,80.6)(370.2,19.95) -- (370.2,80.6)(390.2,19.95) -- (390.2,80.6) ; \draw  [color={rgb, 255:red, 155; green, 155; blue, 155 }  ,draw opacity=1 ] (250.2,19.95) -- (390.4,19.95)(250.2,39.95) -- (390.4,39.95)(250.2,59.95) -- (390.4,59.95)(250.2,79.95) -- (390.4,79.95) ; \draw  [color={rgb, 255:red, 155; green, 155; blue, 155 }  ,draw opacity=1 ]  ;
\draw  [color={rgb, 255:red, 126; green, 211; blue, 33 }  ,draw opacity=1 ][fill={rgb, 255:red, 0; green, 0; blue, 0 }  ,fill opacity=1 ] (258,44.7) .. controls (258,43.48) and (258.98,42.5) .. (260.2,42.5) .. controls (261.42,42.5) and (262.4,43.48) .. (262.4,44.7) .. controls (262.4,45.92) and (261.42,46.9) .. (260.2,46.9) .. controls (258.98,46.9) and (258,45.92) .. (258,44.7) -- cycle ;
\draw  [color={rgb, 255:red, 126; green, 211; blue, 33 }  ,draw opacity=1 ][fill={rgb, 255:red, 0; green, 0; blue, 0 }  ,fill opacity=1 ] (298,44.7) .. controls (298,43.48) and (298.98,42.5) .. (300.2,42.5) .. controls (301.42,42.5) and (302.4,43.48) .. (302.4,44.7) .. controls (302.4,45.92) and (301.42,46.9) .. (300.2,46.9) .. controls (298.98,46.9) and (298,45.92) .. (298,44.7) -- cycle ;
\draw  [color={rgb, 255:red, 126; green, 211; blue, 33 }  ,draw opacity=1 ][fill={rgb, 255:red, 0; green, 0; blue, 0 }  ,fill opacity=1 ] (337.8,44.7) .. controls (337.8,43.48) and (338.78,42.5) .. (340,42.5) .. controls (341.22,42.5) and (342.2,43.48) .. (342.2,44.7) .. controls (342.2,45.92) and (341.22,46.9) .. (340,46.9) .. controls (338.78,46.9) and (337.8,45.92) .. (337.8,44.7) -- cycle ;
\draw  [color={rgb, 255:red, 126; green, 211; blue, 33 }  ,draw opacity=1 ][fill={rgb, 255:red, 0; green, 0; blue, 0 }  ,fill opacity=1 ] (378.8,44.7) .. controls (378.8,43.48) and (379.78,42.5) .. (381,42.5) .. controls (382.22,42.5) and (383.2,43.48) .. (383.2,44.7) .. controls (383.2,45.92) and (382.22,46.9) .. (381,46.9) .. controls (379.78,46.9) and (378.8,45.92) .. (378.8,44.7) -- cycle ;
\draw [color={rgb, 255:red, 126; green, 211; blue, 33 }  ,draw opacity=1 ]   (262.4,44.7) .. controls (275.52,38.04) and (281.57,37.81) .. (296.36,44) ;
\draw [shift={(298,44.7)}, rotate = 203.33] [color={rgb, 255:red, 126; green, 211; blue, 33 }  ,draw opacity=1 ][line width=0.75]    (6.56,-1.97) .. controls (4.17,-0.84) and (1.99,-0.18) .. (0,0) .. controls (1.99,0.18) and (4.17,0.84) .. (6.56,1.97)   ;
\draw [color={rgb, 255:red, 126; green, 211; blue, 33 }  ,draw opacity=1 ]   (302.4,44.7) .. controls (315.52,38.04) and (321.38,37.81) .. (336.16,44) ;
\draw [shift={(337.8,44.7)}, rotate = 203.33] [color={rgb, 255:red, 126; green, 211; blue, 33 }  ,draw opacity=1 ][line width=0.75]    (6.56,-1.97) .. controls (4.17,-0.84) and (1.99,-0.18) .. (0,0) .. controls (1.99,0.18) and (4.17,0.84) .. (6.56,1.97)   ;
\draw [color={rgb, 255:red, 126; green, 211; blue, 33 }  ,draw opacity=1 ]   (342.2,44.7) .. controls (355.32,38.04) and (362.3,37.81) .. (377.15,44) ;
\draw [shift={(378.8,44.7)}, rotate = 203.33] [color={rgb, 255:red, 126; green, 211; blue, 33 }  ,draw opacity=1 ][line width=0.75]    (6.56,-1.97) .. controls (4.17,-0.84) and (1.99,-0.18) .. (0,0) .. controls (1.99,0.18) and (4.17,0.84) .. (6.56,1.97)   ;

\draw [color={rgb, 255:red, 245; green, 166; blue, 35 }  ,draw opacity=1 ]   (282.8,26.1) .. controls (295.92,19.44) and (301.97,19.21) .. (316.76,25.4) ;
\draw [shift={(318.4,26.1)}, rotate = 203.33] [color={rgb, 255:red, 245; green, 166; blue, 35 }  ,draw opacity=1 ][line width=0.75]    (6.56,-1.97) .. controls (4.17,-0.84) and (1.99,-0.18) .. (0,0) .. controls (1.99,0.18) and (4.17,0.84) .. (6.56,1.97)   ;
\draw  [color={rgb, 255:red, 245; green, 166; blue, 35 }  ,draw opacity=1 ][fill={rgb, 255:red, 0; green, 0; blue, 0 }  ,fill opacity=1 ] (278.4,26.1) .. controls (278.4,24.88) and (279.38,23.9) .. (280.6,23.9) .. controls (281.82,23.9) and (282.8,24.88) .. (282.8,26.1) .. controls (282.8,27.32) and (281.82,28.3) .. (280.6,28.3) .. controls (279.38,28.3) and (278.4,27.32) .. (278.4,26.1) -- cycle ;
\draw  [color={rgb, 255:red, 245; green, 166; blue, 35 }  ,draw opacity=1 ][fill={rgb, 255:red, 0; green, 0; blue, 0 }  ,fill opacity=1 ] (318.4,26.1) .. controls (318.4,24.88) and (319.38,23.9) .. (320.6,23.9) .. controls (321.82,23.9) and (322.8,24.88) .. (322.8,26.1) .. controls (322.8,27.32) and (321.82,28.3) .. (320.6,28.3) .. controls (319.38,28.3) and (318.4,27.32) .. (318.4,26.1) -- cycle ;
\draw  [color={rgb, 255:red, 245; green, 166; blue, 35 }  ,draw opacity=1 ][fill={rgb, 255:red, 0; green, 0; blue, 0 }  ,fill opacity=1 ] (358.2,26.1) .. controls (358.2,24.88) and (359.18,23.9) .. (360.4,23.9) .. controls (361.62,23.9) and (362.6,24.88) .. (362.6,26.1) .. controls (362.6,27.32) and (361.62,28.3) .. (360.4,28.3) .. controls (359.18,28.3) and (358.2,27.32) .. (358.2,26.1) -- cycle ;
\draw [color={rgb, 255:red, 245; green, 166; blue, 35 }  ,draw opacity=1 ]   (322.8,26.1) .. controls (335.92,19.44) and (341.78,19.21) .. (356.56,25.4) ;
\draw [shift={(358.2,26.1)}, rotate = 203.33] [color={rgb, 255:red, 245; green, 166; blue, 35 }  ,draw opacity=1 ][line width=0.75]    (6.56,-1.97) .. controls (4.17,-0.84) and (1.99,-0.18) .. (0,0) .. controls (1.99,0.18) and (4.17,0.84) .. (6.56,1.97)   ;

\draw  [color={rgb, 255:red, 245; green, 166; blue, 35 }  ,draw opacity=1 ][fill={rgb, 255:red, 0; green, 0; blue, 0 }  ,fill opacity=1 ] (313.4,32.9) .. controls (313.4,31.68) and (314.38,30.7) .. (315.6,30.7) .. controls (316.82,30.7) and (317.8,31.68) .. (317.8,32.9) .. controls (317.8,34.12) and (316.82,35.1) .. (315.6,35.1) .. controls (314.38,35.1) and (313.4,34.12) .. (313.4,32.9) -- cycle ;
\draw [color={rgb, 255:red, 245; green, 166; blue, 35 }  ,draw opacity=1 ]   (302.8,31.72) .. controls (315.92,38.38) and (321.97,38.62) .. (336.76,32.42) ;
\draw [shift={(338.4,31.72)}, rotate = 156.67] [color={rgb, 255:red, 245; green, 166; blue, 35 }  ,draw opacity=1 ][line width=0.75]    (6.56,-1.97) .. controls (4.17,-0.84) and (1.99,-0.18) .. (0,0) .. controls (1.99,0.18) and (4.17,0.84) .. (6.56,1.97)   ;
\draw  [color={rgb, 255:red, 245; green, 166; blue, 35 }  ,draw opacity=1 ][fill={rgb, 255:red, 0; green, 0; blue, 0 }  ,fill opacity=1 ] (298.4,31.72) .. controls (298.4,32.94) and (299.38,33.92) .. (300.6,33.92) .. controls (301.82,33.92) and (302.8,32.94) .. (302.8,31.72) .. controls (302.8,30.51) and (301.82,29.52) .. (300.6,29.52) .. controls (299.38,29.52) and (298.4,30.51) .. (298.4,31.72) -- cycle ;
\draw  [color={rgb, 255:red, 245; green, 166; blue, 35 }  ,draw opacity=1 ][fill={rgb, 255:red, 0; green, 0; blue, 0 }  ,fill opacity=1 ] (338.4,31.72) .. controls (338.4,32.94) and (339.38,33.92) .. (340.6,33.92) .. controls (341.82,33.92) and (342.8,32.94) .. (342.8,31.72) .. controls (342.8,30.51) and (341.82,29.52) .. (340.6,29.52) .. controls (339.38,29.52) and (338.4,30.51) .. (338.4,31.72) -- cycle ;

\draw [color={rgb, 255:red, 126; green, 211; blue, 33 }  ,draw opacity=1 ]   (302.6,57.53) .. controls (315.72,64.18) and (321.77,64.42) .. (336.56,58.22) ;
\draw [shift={(338.2,57.53)}, rotate = 156.67] [color={rgb, 255:red, 126; green, 211; blue, 33 }  ,draw opacity=1 ][line width=0.75]    (6.56,-1.97) .. controls (4.17,-0.84) and (1.99,-0.18) .. (0,0) .. controls (1.99,0.18) and (4.17,0.84) .. (6.56,1.97)   ;
\draw  [color={rgb, 255:red, 126; green, 211; blue, 33 }  ,draw opacity=1 ][fill={rgb, 255:red, 0; green, 0; blue, 0 }  ,fill opacity=1 ] (298.2,57.53) .. controls (298.2,58.74) and (299.18,59.73) .. (300.4,59.73) .. controls (301.62,59.73) and (302.6,58.74) .. (302.6,57.53) .. controls (302.6,56.31) and (301.62,55.33) .. (300.4,55.33) .. controls (299.18,55.33) and (298.2,56.31) .. (298.2,57.53) -- cycle ;
\draw  [color={rgb, 255:red, 126; green, 211; blue, 33 }  ,draw opacity=1 ][fill={rgb, 255:red, 0; green, 0; blue, 0 }  ,fill opacity=1 ] (338.2,57.53) .. controls (338.2,58.74) and (339.18,59.73) .. (340.4,59.73) .. controls (341.62,59.73) and (342.6,58.74) .. (342.6,57.53) .. controls (342.6,56.31) and (341.62,55.33) .. (340.4,55.33) .. controls (339.18,55.33) and (338.2,56.31) .. (338.2,57.53) -- cycle ;
\draw [color={rgb, 255:red, 126; green, 211; blue, 33 }  ,draw opacity=1 ]   (283,51.1) .. controls (297.48,46.95) and (300.22,46.81) .. (316.75,50.66) ;
\draw [shift={(318.6,51.1)}, rotate = 193.29] [color={rgb, 255:red, 126; green, 211; blue, 33 }  ,draw opacity=1 ][line width=0.75]    (6.56,-1.97) .. controls (4.17,-0.84) and (1.99,-0.18) .. (0,0) .. controls (1.99,0.18) and (4.17,0.84) .. (6.56,1.97)   ;
\draw  [color={rgb, 255:red, 126; green, 211; blue, 33 }  ,draw opacity=1 ][fill={rgb, 255:red, 0; green, 0; blue, 0 }  ,fill opacity=1 ] (278.6,51.1) .. controls (278.6,49.88) and (279.58,48.9) .. (280.8,48.9) .. controls (282.02,48.9) and (283,49.88) .. (283,51.1) .. controls (283,52.32) and (282.02,53.3) .. (280.8,53.3) .. controls (279.58,53.3) and (278.6,52.32) .. (278.6,51.1) -- cycle ;
\draw  [color={rgb, 255:red, 126; green, 211; blue, 33 }  ,draw opacity=1 ][fill={rgb, 255:red, 0; green, 0; blue, 0 }  ,fill opacity=1 ] (318.6,51.1) .. controls (318.6,49.88) and (319.58,48.9) .. (320.8,48.9) .. controls (322.02,48.9) and (323,49.88) .. (323,51.1) .. controls (323,52.32) and (322.02,53.3) .. (320.8,53.3) .. controls (319.58,53.3) and (318.6,52.32) .. (318.6,51.1) -- cycle ;
\draw  [color={rgb, 255:red, 126; green, 211; blue, 33 }  ,draw opacity=1 ][fill={rgb, 255:red, 0; green, 0; blue, 0 }  ,fill opacity=1 ] (358.4,51.1) .. controls (358.4,49.88) and (359.38,48.9) .. (360.6,48.9) .. controls (361.82,48.9) and (362.8,49.88) .. (362.8,51.1) .. controls (362.8,52.32) and (361.82,53.3) .. (360.6,53.3) .. controls (359.38,53.3) and (358.4,52.32) .. (358.4,51.1) -- cycle ;
\draw [color={rgb, 255:red, 126; green, 211; blue, 33 }  ,draw opacity=1 ]   (323,51.1) .. controls (337.48,46.95) and (340.03,46.81) .. (356.55,50.66) ;
\draw [shift={(358.4,51.1)}, rotate = 193.29] [color={rgb, 255:red, 126; green, 211; blue, 33 }  ,draw opacity=1 ][line width=0.75]    (6.56,-1.97) .. controls (4.17,-0.84) and (1.99,-0.18) .. (0,0) .. controls (1.99,0.18) and (4.17,0.84) .. (6.56,1.97)   ;

\draw [color={rgb, 255:red, 126; green, 211; blue, 33 }  ,draw opacity=1 ]   (302.4,51.13) .. controls (315.52,57.78) and (321.57,58.02) .. (336.36,51.82) ;
\draw [shift={(338,51.13)}, rotate = 156.67] [color={rgb, 255:red, 126; green, 211; blue, 33 }  ,draw opacity=1 ][line width=0.75]    (6.56,-1.97) .. controls (4.17,-0.84) and (1.99,-0.18) .. (0,0) .. controls (1.99,0.18) and (4.17,0.84) .. (6.56,1.97)   ;
\draw  [color={rgb, 255:red, 126; green, 211; blue, 33 }  ,draw opacity=1 ][fill={rgb, 255:red, 0; green, 0; blue, 0 }  ,fill opacity=1 ] (298,51.13) .. controls (298,52.34) and (298.98,53.33) .. (300.2,53.33) .. controls (301.42,53.33) and (302.4,52.34) .. (302.4,51.13) .. controls (302.4,49.91) and (301.42,48.93) .. (300.2,48.93) .. controls (298.98,48.93) and (298,49.91) .. (298,51.13) -- cycle ;
\draw  [color={rgb, 255:red, 126; green, 211; blue, 33 }  ,draw opacity=1 ][fill={rgb, 255:red, 0; green, 0; blue, 0 }  ,fill opacity=1 ] (338,51.13) .. controls (338,52.34) and (338.98,53.33) .. (340.2,53.33) .. controls (341.42,53.33) and (342.4,52.34) .. (342.4,51.13) .. controls (342.4,49.91) and (341.42,48.93) .. (340.2,48.93) .. controls (338.98,48.93) and (338,49.91) .. (338,51.13) -- cycle ;
\draw [color={rgb, 255:red, 74; green, 144; blue, 226 }  ,draw opacity=1 ]   (282.4,67.3) .. controls (296.88,63.15) and (299.62,63.01) .. (316.15,66.86) ;
\draw [shift={(318,67.3)}, rotate = 193.29] [color={rgb, 255:red, 74; green, 144; blue, 226 }  ,draw opacity=1 ][line width=0.75]    (6.56,-1.97) .. controls (4.17,-0.84) and (1.99,-0.18) .. (0,0) .. controls (1.99,0.18) and (4.17,0.84) .. (6.56,1.97)   ;
\draw  [color={rgb, 255:red, 74; green, 144; blue, 226 }  ,draw opacity=1 ][fill={rgb, 255:red, 0; green, 0; blue, 0 }  ,fill opacity=1 ] (278,67.3) .. controls (278,66.08) and (278.98,65.1) .. (280.2,65.1) .. controls (281.42,65.1) and (282.4,66.08) .. (282.4,67.3) .. controls (282.4,68.52) and (281.42,69.5) .. (280.2,69.5) .. controls (278.98,69.5) and (278,68.52) .. (278,67.3) -- cycle ;
\draw  [color={rgb, 255:red, 74; green, 144; blue, 226 }  ,draw opacity=1 ][fill={rgb, 255:red, 0; green, 0; blue, 0 }  ,fill opacity=1 ] (318,67.3) .. controls (318,66.08) and (318.98,65.1) .. (320.2,65.1) .. controls (321.42,65.1) and (322.4,66.08) .. (322.4,67.3) .. controls (322.4,68.52) and (321.42,69.5) .. (320.2,69.5) .. controls (318.98,69.5) and (318,68.52) .. (318,67.3) -- cycle ;
\draw  [color={rgb, 255:red, 74; green, 144; blue, 226 }  ,draw opacity=1 ][fill={rgb, 255:red, 0; green, 0; blue, 0 }  ,fill opacity=1 ] (357.8,67.3) .. controls (357.8,66.08) and (358.78,65.1) .. (360,65.1) .. controls (361.22,65.1) and (362.2,66.08) .. (362.2,67.3) .. controls (362.2,68.52) and (361.22,69.5) .. (360,69.5) .. controls (358.78,69.5) and (357.8,68.52) .. (357.8,67.3) -- cycle ;
\draw [color={rgb, 255:red, 74; green, 144; blue, 226 }  ,draw opacity=1 ]   (322.4,67.3) .. controls (336.88,63.15) and (339.43,63.01) .. (355.95,66.86) ;
\draw [shift={(357.8,67.3)}, rotate = 193.29] [color={rgb, 255:red, 74; green, 144; blue, 226 }  ,draw opacity=1 ][line width=0.75]    (6.56,-1.97) .. controls (4.17,-0.84) and (1.99,-0.18) .. (0,0) .. controls (1.99,0.18) and (4.17,0.84) .. (6.56,1.97)   ;

\draw  [color={rgb, 255:red, 74; green, 144; blue, 226 }  ,draw opacity=1 ][fill={rgb, 255:red, 0; green, 0; blue, 0 }  ,fill opacity=1 ] (318,74.7) .. controls (318,73.48) and (318.98,72.5) .. (320.2,72.5) .. controls (321.42,72.5) and (322.4,73.48) .. (322.4,74.7) .. controls (322.4,75.92) and (321.42,76.9) .. (320.2,76.9) .. controls (318.98,76.9) and (318,75.92) .. (318,74.7) -- cycle ;
\draw [color={rgb, 255:red, 74; green, 144; blue, 226 }  ,draw opacity=1 ]   (302.8,73.73) .. controls (315.92,80.38) and (321.97,80.62) .. (336.76,74.42) ;
\draw [shift={(338.4,73.73)}, rotate = 156.67] [color={rgb, 255:red, 74; green, 144; blue, 226 }  ,draw opacity=1 ][line width=0.75]    (6.56,-1.97) .. controls (4.17,-0.84) and (1.99,-0.18) .. (0,0) .. controls (1.99,0.18) and (4.17,0.84) .. (6.56,1.97)   ;
\draw  [color={rgb, 255:red, 74; green, 144; blue, 226 }  ,draw opacity=1 ][fill={rgb, 255:red, 0; green, 0; blue, 0 }  ,fill opacity=1 ] (298.4,73.73) .. controls (298.4,74.94) and (299.38,75.93) .. (300.6,75.93) .. controls (301.82,75.93) and (302.8,74.94) .. (302.8,73.73) .. controls (302.8,72.51) and (301.82,71.53) .. (300.6,71.53) .. controls (299.38,71.53) and (298.4,72.51) .. (298.4,73.73) -- cycle ;
\draw  [color={rgb, 255:red, 74; green, 144; blue, 226 }  ,draw opacity=1 ][fill={rgb, 255:red, 0; green, 0; blue, 0 }  ,fill opacity=1 ] (338.4,73.73) .. controls (338.4,74.94) and (339.38,75.93) .. (340.6,75.93) .. controls (341.82,75.93) and (342.8,74.94) .. (342.8,73.73) .. controls (342.8,72.51) and (341.82,71.53) .. (340.6,71.53) .. controls (339.38,71.53) and (338.4,72.51) .. (338.4,73.73) -- cycle ;

\draw  [color={rgb, 255:red, 126; green, 211; blue, 33 }  ,draw opacity=1 ][fill={rgb, 255:red, 0; green, 0; blue, 0 }  ,fill opacity=1 ] (318.6,44.7) .. controls (318.6,43.48) and (319.58,42.5) .. (320.8,42.5) .. controls (322.02,42.5) and (323,43.48) .. (323,44.7) .. controls (323,45.92) and (322.02,46.9) .. (320.8,46.9) .. controls (319.58,46.9) and (318.6,45.92) .. (318.6,44.7) -- cycle ;
\draw  [color={rgb, 255:red, 245; green, 166; blue, 35 }  ,draw opacity=1 ][fill={rgb, 255:red, 0; green, 0; blue, 0 }  ,fill opacity=1 ] (323.4,32.9) .. controls (323.4,31.68) and (324.38,30.7) .. (325.6,30.7) .. controls (326.82,30.7) and (327.8,31.68) .. (327.8,32.9) .. controls (327.8,34.12) and (326.82,35.1) .. (325.6,35.1) .. controls (324.38,35.1) and (323.4,34.12) .. (323.4,32.9) -- cycle ;

\draw (42.5,28.5) node  [font=\scriptsize]  {$2$};
\draw (42.5,48.5) node  [font=\scriptsize]  {$0$};
\draw (113.7,108.7) node    {$\Delta \ =\ -2$};
\draw (41.82,6.4) node    {$a$};
\draw (201.82,87.4) node    {$q$};
\draw (122.5,89.7) node  [font=\scriptsize]  {$0$};
\draw (100,89.7) node  [font=\scriptsize]  {$-2$};
\draw (162.5,89.7) node  [font=\scriptsize]  {$4$};
\draw (142.5,89.7) node  [font=\scriptsize]  {$2$};
\draw (182.5,89.7) node  [font=\scriptsize]  {$6$};
\draw (80,89.7) node  [font=\scriptsize]  {$-4$};
\draw (60,89.7) node  [font=\scriptsize]  {$-6$};
\draw (242.7,28.9) node  [font=\scriptsize]  {$0$};
\draw (241.2,48.9) node  [font=\scriptsize]  {$-2$};
\draw (313.9,109.1) node    {$\Delta \ =\ 0$};
\draw (242.02,6.8) node    {$a$};
\draw (402.02,87.8) node    {$q$};
\draw (322.7,90.1) node  [font=\scriptsize]  {$0$};
\draw (300.2,90.1) node  [font=\scriptsize]  {$-2$};
\draw (362.7,90.1) node  [font=\scriptsize]  {$4$};
\draw (342.7,90.1) node  [font=\scriptsize]  {$2$};
\draw (382.7,90.1) node  [font=\scriptsize]  {$6$};
\draw (280.2,90.1) node  [font=\scriptsize]  {$-4$};
\draw (260.6,90.1) node  [font=\scriptsize]  {$-6$};
\draw (41,68.1) node  [font=\scriptsize]  {$-2$};
\draw (241.2,68.9) node  [font=\scriptsize]  {$-4$};

\end{tikzpicture}

%% file: tikzpictures/table-kh-conway.tex
\tikzset{every picture/.style={line width=0.75pt}} 

\begin{tikzpicture}[x=0.75pt,y=0.75pt,yscale=-1.2,xscale=1.2]

\draw [color={rgb, 255:red, 245; green, 166; blue, 35 }  ,draw opacity=1 ]   (142.13,49.37) .. controls (156.61,45.22) and (159.35,45.07) .. (175.88,48.93) ;
\draw [shift={(177.73,49.37)}, rotate = 193.29] [color={rgb, 255:red, 245; green, 166; blue, 35 }  ,draw opacity=1 ][line width=0.75]    (6.56,-1.97) .. controls (4.17,-0.84) and (1.99,-0.18) .. (0,0) .. controls (1.99,0.18) and (4.17,0.84) .. (6.56,1.97)   ;
\draw  [color={rgb, 255:red, 245; green, 166; blue, 35 }  ,draw opacity=1 ][fill={rgb, 255:red, 0; green, 0; blue, 0 }  ,fill opacity=1 ] (137.73,49.37) .. controls (137.73,48.15) and (138.72,47.17) .. (139.93,47.17) .. controls (141.15,47.17) and (142.13,48.15) .. (142.13,49.37) .. controls (142.13,50.58) and (141.15,51.57) .. (139.93,51.57) .. controls (138.72,51.57) and (137.73,50.58) .. (137.73,49.37) -- cycle ;
\draw  [color={rgb, 255:red, 245; green, 166; blue, 35 }  ,draw opacity=1 ][fill={rgb, 255:red, 0; green, 0; blue, 0 }  ,fill opacity=1 ] (177.73,49.37) .. controls (177.73,48.15) and (178.72,47.17) .. (179.93,47.17) .. controls (181.15,47.17) and (182.13,48.15) .. (182.13,49.37) .. controls (182.13,50.58) and (181.15,51.57) .. (179.93,51.57) .. controls (178.72,51.57) and (177.73,50.58) .. (177.73,49.37) -- cycle ;
\draw  [color={rgb, 255:red, 245; green, 166; blue, 35 }  ,draw opacity=1 ][fill={rgb, 255:red, 0; green, 0; blue, 0 }  ,fill opacity=1 ] (217.53,49.37) .. controls (217.53,48.15) and (218.52,47.17) .. (219.73,47.17) .. controls (220.95,47.17) and (221.93,48.15) .. (221.93,49.37) .. controls (221.93,50.58) and (220.95,51.57) .. (219.73,51.57) .. controls (218.52,51.57) and (217.53,50.58) .. (217.53,49.37) -- cycle ;
\draw [color={rgb, 255:red, 245; green, 166; blue, 35 }  ,draw opacity=1 ]   (182.13,49.37) .. controls (196.61,45.22) and (199.16,45.07) .. (215.68,48.93) ;
\draw [shift={(217.53,49.37)}, rotate = 193.29] [color={rgb, 255:red, 245; green, 166; blue, 35 }  ,draw opacity=1 ][line width=0.75]    (6.56,-1.97) .. controls (4.17,-0.84) and (1.99,-0.18) .. (0,0) .. controls (1.99,0.18) and (4.17,0.84) .. (6.56,1.97)   ;

\draw  [draw opacity=0] (50,19.55) -- (271,19.55) -- (271,60.07) -- (50,60.07) -- cycle ; \draw  [color={rgb, 255:red, 155; green, 155; blue, 155 }  ,draw opacity=1 ] (50,19.55) -- (50,60.07)(70,19.55) -- (70,60.07)(90,19.55) -- (90,60.07)(110,19.55) -- (110,60.07)(130,19.55) -- (130,60.07)(150,19.55) -- (150,60.07)(170,19.55) -- (170,60.07)(190,19.55) -- (190,60.07)(210,19.55) -- (210,60.07)(230,19.55) -- (230,60.07)(250,19.55) -- (250,60.07)(270,19.55) -- (270,60.07) ; \draw  [color={rgb, 255:red, 155; green, 155; blue, 155 }  ,draw opacity=1 ] (50,19.55) -- (271,19.55)(50,39.55) -- (271,39.55)(50,59.55) -- (271,59.55) ; \draw  [color={rgb, 255:red, 155; green, 155; blue, 155 }  ,draw opacity=1 ]  ;
\draw  [color={rgb, 255:red, 245; green, 166; blue, 35 }  ,draw opacity=1 ][fill={rgb, 255:red, 0; green, 0; blue, 0 }  ,fill opacity=1 ] (157.13,42.97) .. controls (157.13,41.75) and (158.12,40.77) .. (159.33,40.77) .. controls (160.55,40.77) and (161.53,41.75) .. (161.53,42.97) .. controls (161.53,44.18) and (160.55,45.17) .. (159.33,45.17) .. controls (158.12,45.17) and (157.13,44.18) .. (157.13,42.97) -- cycle ;
\draw  [color={rgb, 255:red, 245; green, 166; blue, 35 }  ,draw opacity=1 ][fill={rgb, 255:red, 0; green, 0; blue, 0 }  ,fill opacity=1 ] (196.93,42.97) .. controls (196.93,41.75) and (197.92,40.77) .. (199.13,40.77) .. controls (200.35,40.77) and (201.33,41.75) .. (201.33,42.97) .. controls (201.33,44.18) and (200.35,45.17) .. (199.13,45.17) .. controls (197.92,45.17) and (196.93,44.18) .. (196.93,42.97) -- cycle ;
\draw  [color={rgb, 255:red, 245; green, 166; blue, 35 }  ,draw opacity=1 ][fill={rgb, 255:red, 0; green, 0; blue, 0 }  ,fill opacity=1 ] (237.93,42.97) .. controls (237.93,41.75) and (238.92,40.77) .. (240.13,40.77) .. controls (241.35,40.77) and (242.33,41.75) .. (242.33,42.97) .. controls (242.33,44.18) and (241.35,45.17) .. (240.13,45.17) .. controls (238.92,45.17) and (237.93,44.18) .. (237.93,42.97) -- cycle ;
\draw [color={rgb, 255:red, 245; green, 166; blue, 35 }  ,draw opacity=1 ]   (161.53,42.97) .. controls (174.66,36.31) and (180.52,36.08) .. (195.29,42.27) ;
\draw [shift={(196.93,42.97)}, rotate = 203.33] [color={rgb, 255:red, 245; green, 166; blue, 35 }  ,draw opacity=1 ][line width=0.75]    (6.56,-1.97) .. controls (4.17,-0.84) and (1.99,-0.18) .. (0,0) .. controls (1.99,0.18) and (4.17,0.84) .. (6.56,1.97)   ;
\draw [color={rgb, 255:red, 245; green, 166; blue, 35 }  ,draw opacity=1 ]   (201.33,42.97) .. controls (214.46,36.31) and (221.44,36.08) .. (236.29,42.27) ;
\draw [shift={(237.93,42.97)}, rotate = 203.33] [color={rgb, 255:red, 245; green, 166; blue, 35 }  ,draw opacity=1 ][line width=0.75]    (6.56,-1.97) .. controls (4.17,-0.84) and (1.99,-0.18) .. (0,0) .. controls (1.99,0.18) and (4.17,0.84) .. (6.56,1.97)   ;

\draw  [color={rgb, 255:red, 208; green, 2; blue, 27 }  ,draw opacity=1 ][fill={rgb, 255:red, 0; green, 0; blue, 0 }  ,fill opacity=1 ] (217.53,55.37) .. controls (217.53,54.15) and (218.52,53.17) .. (219.73,53.17) .. controls (220.95,53.17) and (221.93,54.15) .. (221.93,55.37) .. controls (221.93,56.58) and (220.95,57.57) .. (219.73,57.57) .. controls (218.52,57.57) and (217.53,56.58) .. (217.53,55.37) -- cycle ;
\draw  [color={rgb, 255:red, 208; green, 2; blue, 27 }  ,draw opacity=1 ][fill={rgb, 255:red, 0; green, 0; blue, 0 }  ,fill opacity=1 ] (257.33,55.37) .. controls (257.33,54.15) and (258.32,53.17) .. (259.53,53.17) .. controls (260.75,53.17) and (261.73,54.15) .. (261.73,55.37) .. controls (261.73,56.58) and (260.75,57.57) .. (259.53,57.57) .. controls (258.32,57.57) and (257.33,56.58) .. (257.33,55.37) -- cycle ;
\draw [color={rgb, 255:red, 208; green, 2; blue, 27 }  ,draw opacity=1 ]   (221.93,55.37) .. controls (235.06,48.71) and (240.92,48.48) .. (255.69,54.67) ;
\draw [shift={(257.33,55.37)}, rotate = 203.33] [color={rgb, 255:red, 208; green, 2; blue, 27 }  ,draw opacity=1 ][line width=0.75]    (6.56,-1.97) .. controls (4.17,-0.84) and (1.99,-0.18) .. (0,0) .. controls (1.99,0.18) and (4.17,0.84) .. (6.56,1.97)   ;

\draw [color={rgb, 255:red, 208; green, 2; blue, 27 }  ,draw opacity=1 ]   (201.93,55.99) .. controls (215.06,62.65) and (221.1,62.88) .. (235.89,56.69) ;
\draw [shift={(237.53,55.99)}, rotate = 156.67] [color={rgb, 255:red, 208; green, 2; blue, 27 }  ,draw opacity=1 ][line width=0.75]    (6.56,-1.97) .. controls (4.17,-0.84) and (1.99,-0.18) .. (0,0) .. controls (1.99,0.18) and (4.17,0.84) .. (6.56,1.97)   ;
\draw  [color={rgb, 255:red, 208; green, 2; blue, 27 }  ,draw opacity=1 ][fill={rgb, 255:red, 0; green, 0; blue, 0 }  ,fill opacity=1 ] (197.53,55.99) .. controls (197.53,57.21) and (198.52,58.19) .. (199.73,58.19) .. controls (200.95,58.19) and (201.93,57.21) .. (201.93,55.99) .. controls (201.93,54.78) and (200.95,53.79) .. (199.73,53.79) .. controls (198.52,53.79) and (197.53,54.78) .. (197.53,55.99) -- cycle ;
\draw  [color={rgb, 255:red, 208; green, 2; blue, 27 }  ,draw opacity=1 ][fill={rgb, 255:red, 0; green, 0; blue, 0 }  ,fill opacity=1 ] (237.53,55.99) .. controls (237.53,57.21) and (238.52,58.19) .. (239.73,58.19) .. controls (240.95,58.19) and (241.93,57.21) .. (241.93,55.99) .. controls (241.93,54.78) and (240.95,53.79) .. (239.73,53.79) .. controls (238.52,53.79) and (237.53,54.78) .. (237.53,55.99) -- cycle ;

\draw  [color={rgb, 255:red, 245; green, 166; blue, 35 }  ,draw opacity=1 ][fill={rgb, 255:red, 0; green, 0; blue, 0 }  ,fill opacity=1 ] (197.13,49.39) .. controls (197.13,50.61) and (198.12,51.59) .. (199.33,51.59) .. controls (200.55,51.59) and (201.53,50.61) .. (201.53,49.39) .. controls (201.53,48.18) and (200.55,47.19) .. (199.33,47.19) .. controls (198.12,47.19) and (197.13,48.18) .. (197.13,49.39) -- cycle ;
\draw  [color={rgb, 255:red, 126; green, 211; blue, 33 }  ,draw opacity=1 ][fill={rgb, 255:red, 0; green, 0; blue, 0 }  ,fill opacity=1 ] (137.93,55.17) .. controls (137.93,53.95) and (138.92,52.97) .. (140.13,52.97) .. controls (141.35,52.97) and (142.33,53.95) .. (142.33,55.17) .. controls (142.33,56.38) and (141.35,57.37) .. (140.13,57.37) .. controls (138.92,57.37) and (137.93,56.38) .. (137.93,55.17) -- cycle ;
\draw  [color={rgb, 255:red, 126; green, 211; blue, 33 }  ,draw opacity=1 ][fill={rgb, 255:red, 0; green, 0; blue, 0 }  ,fill opacity=1 ] (177.73,55.17) .. controls (177.73,53.95) and (178.72,52.97) .. (179.93,52.97) .. controls (181.15,52.97) and (182.13,53.95) .. (182.13,55.17) .. controls (182.13,56.38) and (181.15,57.37) .. (179.93,57.37) .. controls (178.72,57.37) and (177.73,56.38) .. (177.73,55.17) -- cycle ;
\draw [color={rgb, 255:red, 126; green, 211; blue, 33 }  ,draw opacity=1 ]   (142.33,55.17) .. controls (156.81,51.02) and (159.36,50.87) .. (175.88,54.73) ;
\draw [shift={(177.73,55.17)}, rotate = 193.29] [color={rgb, 255:red, 126; green, 211; blue, 33 }  ,draw opacity=1 ][line width=0.75]    (6.56,-1.97) .. controls (4.17,-0.84) and (1.99,-0.18) .. (0,0) .. controls (1.99,0.18) and (4.17,0.84) .. (6.56,1.97)   ;

\draw [color={rgb, 255:red, 126; green, 211; blue, 33 }  ,draw opacity=1 ]   (122.73,55.99) .. controls (135.86,62.65) and (141.9,62.88) .. (156.69,56.69) ;
\draw [shift={(158.33,55.99)}, rotate = 156.67] [color={rgb, 255:red, 126; green, 211; blue, 33 }  ,draw opacity=1 ][line width=0.75]    (6.56,-1.97) .. controls (4.17,-0.84) and (1.99,-0.18) .. (0,0) .. controls (1.99,0.18) and (4.17,0.84) .. (6.56,1.97)   ;
\draw  [color={rgb, 255:red, 126; green, 211; blue, 33 }  ,draw opacity=1 ][fill={rgb, 255:red, 0; green, 0; blue, 0 }  ,fill opacity=1 ] (118.33,55.99) .. controls (118.33,57.21) and (119.32,58.19) .. (120.53,58.19) .. controls (121.75,58.19) and (122.73,57.21) .. (122.73,55.99) .. controls (122.73,54.78) and (121.75,53.79) .. (120.53,53.79) .. controls (119.32,53.79) and (118.33,54.78) .. (118.33,55.99) -- cycle ;
\draw  [color={rgb, 255:red, 126; green, 211; blue, 33 }  ,draw opacity=1 ][fill={rgb, 255:red, 0; green, 0; blue, 0 }  ,fill opacity=1 ] (158.33,55.99) .. controls (158.33,57.21) and (159.32,58.19) .. (160.53,58.19) .. controls (161.75,58.19) and (162.73,57.21) .. (162.73,55.99) .. controls (162.73,54.78) and (161.75,53.79) .. (160.53,53.79) .. controls (159.32,53.79) and (158.33,54.78) .. (158.33,55.99) -- cycle ;

\draw  [color={rgb, 255:red, 245; green, 166; blue, 35 }  ,draw opacity=1 ][fill={rgb, 255:red, 0; green, 0; blue, 0 }  ,fill opacity=1 ] (178.07,42.63) .. controls (178.07,41.42) and (179.05,40.43) .. (180.27,40.43) .. controls (181.48,40.43) and (182.47,41.42) .. (182.47,42.63) .. controls (182.47,43.85) and (181.48,44.83) .. (180.27,44.83) .. controls (179.05,44.83) and (178.07,43.85) .. (178.07,42.63) -- cycle ;
\draw  [color={rgb, 255:red, 126; green, 211; blue, 33 }  ,draw opacity=1 ][fill={rgb, 255:red, 0; green, 0; blue, 0 }  ,fill opacity=1 ] (77,36.39) .. controls (77,37.61) and (77.98,38.59) .. (79.2,38.59) .. controls (80.42,38.59) and (81.4,37.61) .. (81.4,36.39) .. controls (81.4,35.18) and (80.42,34.19) .. (79.2,34.19) .. controls (77.98,34.19) and (77,35.18) .. (77,36.39) -- cycle ;
\draw  [color={rgb, 255:red, 126; green, 211; blue, 33 }  ,draw opacity=1 ][fill={rgb, 255:red, 0; green, 0; blue, 0 }  ,fill opacity=1 ] (117,36.39) .. controls (117,37.61) and (117.98,38.59) .. (119.2,38.59) .. controls (120.42,38.59) and (121.4,37.61) .. (121.4,36.39) .. controls (121.4,35.18) and (120.42,34.19) .. (119.2,34.19) .. controls (117.98,34.19) and (117,35.18) .. (117,36.39) -- cycle ;
\draw  [color={rgb, 255:red, 126; green, 211; blue, 33 }  ,draw opacity=1 ][fill={rgb, 255:red, 0; green, 0; blue, 0 }  ,fill opacity=1 ] (156.8,36.39) .. controls (156.8,37.61) and (157.78,38.59) .. (159,38.59) .. controls (160.22,38.59) and (161.2,37.61) .. (161.2,36.39) .. controls (161.2,35.18) and (160.22,34.19) .. (159,34.19) .. controls (157.78,34.19) and (156.8,35.18) .. (156.8,36.39) -- cycle ;
\draw [color={rgb, 255:red, 126; green, 211; blue, 33 }  ,draw opacity=1 ]   (81.4,36.39) .. controls (94.52,43.05) and (100.57,43.28) .. (115.36,37.09) ;
\draw [shift={(117,36.39)}, rotate = 156.67] [color={rgb, 255:red, 126; green, 211; blue, 33 }  ,draw opacity=1 ][line width=0.75]    (6.56,-1.97) .. controls (4.17,-0.84) and (1.99,-0.18) .. (0,0) .. controls (1.99,0.18) and (4.17,0.84) .. (6.56,1.97)   ;
\draw [color={rgb, 255:red, 126; green, 211; blue, 33 }  ,draw opacity=1 ]   (121.4,36.39) .. controls (134.52,43.05) and (140.38,43.28) .. (155.16,37.09) ;
\draw [shift={(156.8,36.39)}, rotate = 156.67] [color={rgb, 255:red, 126; green, 211; blue, 33 }  ,draw opacity=1 ][line width=0.75]    (6.56,-1.97) .. controls (4.17,-0.84) and (1.99,-0.18) .. (0,0) .. controls (1.99,0.18) and (4.17,0.84) .. (6.56,1.97)   ;

\draw [color={rgb, 255:red, 245; green, 166; blue, 35 }  ,draw opacity=1 ]   (141.8,23.9) .. controls (154.92,17.24) and (160.97,17.01) .. (175.76,23.2) ;
\draw [shift={(177.4,23.9)}, rotate = 203.33] [color={rgb, 255:red, 245; green, 166; blue, 35 }  ,draw opacity=1 ][line width=0.75]    (6.56,-1.97) .. controls (4.17,-0.84) and (1.99,-0.18) .. (0,0) .. controls (1.99,0.18) and (4.17,0.84) .. (6.56,1.97)   ;
\draw  [color={rgb, 255:red, 245; green, 166; blue, 35 }  ,draw opacity=1 ][fill={rgb, 255:red, 0; green, 0; blue, 0 }  ,fill opacity=1 ] (137.4,23.9) .. controls (137.4,22.68) and (138.38,21.7) .. (139.6,21.7) .. controls (140.82,21.7) and (141.8,22.68) .. (141.8,23.9) .. controls (141.8,25.12) and (140.82,26.1) .. (139.6,26.1) .. controls (138.38,26.1) and (137.4,25.12) .. (137.4,23.9) -- cycle ;
\draw  [color={rgb, 255:red, 245; green, 166; blue, 35 }  ,draw opacity=1 ][fill={rgb, 255:red, 0; green, 0; blue, 0 }  ,fill opacity=1 ] (177.4,23.9) .. controls (177.4,22.68) and (178.38,21.7) .. (179.6,21.7) .. controls (180.82,21.7) and (181.8,22.68) .. (181.8,23.9) .. controls (181.8,25.12) and (180.82,26.1) .. (179.6,26.1) .. controls (178.38,26.1) and (177.4,25.12) .. (177.4,23.9) -- cycle ;

\draw [color={rgb, 255:red, 245; green, 166; blue, 35 }  ,draw opacity=1 ]   (161.8,24.53) .. controls (174.92,31.18) and (180.97,31.42) .. (195.76,25.22) ;
\draw [shift={(197.4,24.53)}, rotate = 156.67] [color={rgb, 255:red, 245; green, 166; blue, 35 }  ,draw opacity=1 ][line width=0.75]    (6.56,-1.97) .. controls (4.17,-0.84) and (1.99,-0.18) .. (0,0) .. controls (1.99,0.18) and (4.17,0.84) .. (6.56,1.97)   ;
\draw  [color={rgb, 255:red, 245; green, 166; blue, 35 }  ,draw opacity=1 ][fill={rgb, 255:red, 0; green, 0; blue, 0 }  ,fill opacity=1 ] (157.4,24.52) .. controls (157.4,25.74) and (158.38,26.72) .. (159.6,26.72) .. controls (160.82,26.72) and (161.8,25.74) .. (161.8,24.52) .. controls (161.8,23.31) and (160.82,22.32) .. (159.6,22.32) .. controls (158.38,22.32) and (157.4,23.31) .. (157.4,24.52) -- cycle ;
\draw  [color={rgb, 255:red, 245; green, 166; blue, 35 }  ,draw opacity=1 ][fill={rgb, 255:red, 0; green, 0; blue, 0 }  ,fill opacity=1 ] (197.4,24.52) .. controls (197.4,25.74) and (198.38,26.72) .. (199.6,26.72) .. controls (200.82,26.72) and (201.8,25.74) .. (201.8,24.52) .. controls (201.8,23.31) and (200.82,22.32) .. (199.6,22.32) .. controls (198.38,22.32) and (197.4,23.31) .. (197.4,24.52) -- cycle ;

\draw [color={rgb, 255:red, 126; green, 211; blue, 33 }  ,draw opacity=1 ]   (101.67,34.43) .. controls (116.14,30.28) and (118.88,30.14) .. (135.41,34) ;
\draw [shift={(137.27,34.43)}, rotate = 193.29] [color={rgb, 255:red, 126; green, 211; blue, 33 }  ,draw opacity=1 ][line width=0.75]    (6.56,-1.97) .. controls (4.17,-0.84) and (1.99,-0.18) .. (0,0) .. controls (1.99,0.18) and (4.17,0.84) .. (6.56,1.97)   ;
\draw  [color={rgb, 255:red, 126; green, 211; blue, 33 }  ,draw opacity=1 ][fill={rgb, 255:red, 0; green, 0; blue, 0 }  ,fill opacity=1 ] (97.27,34.43) .. controls (97.27,33.22) and (98.25,32.23) .. (99.47,32.23) .. controls (100.68,32.23) and (101.67,33.22) .. (101.67,34.43) .. controls (101.67,35.65) and (100.68,36.63) .. (99.47,36.63) .. controls (98.25,36.63) and (97.27,35.65) .. (97.27,34.43) -- cycle ;
\draw [color={rgb, 255:red, 126; green, 211; blue, 33 }  ,draw opacity=1 ]   (141.67,34.43) .. controls (156.14,30.28) and (158.7,30.14) .. (175.21,34) ;
\draw [shift={(177.07,34.43)}, rotate = 193.29] [color={rgb, 255:red, 126; green, 211; blue, 33 }  ,draw opacity=1 ][line width=0.75]    (6.56,-1.97) .. controls (4.17,-0.84) and (1.99,-0.18) .. (0,0) .. controls (1.99,0.18) and (4.17,0.84) .. (6.56,1.97)   ;
\draw  [color={rgb, 255:red, 126; green, 211; blue, 33 }  ,draw opacity=1 ][fill={rgb, 255:red, 0; green, 0; blue, 0 }  ,fill opacity=1 ] (137.27,34.43) .. controls (137.27,33.22) and (138.25,32.23) .. (139.47,32.23) .. controls (140.68,32.23) and (141.67,33.22) .. (141.67,34.43) .. controls (141.67,35.65) and (140.68,36.63) .. (139.47,36.63) .. controls (138.25,36.63) and (137.27,35.65) .. (137.27,34.43) -- cycle ;
\draw  [color={rgb, 255:red, 126; green, 211; blue, 33 }  ,draw opacity=1 ][fill={rgb, 255:red, 0; green, 0; blue, 0 }  ,fill opacity=1 ] (177.07,34.43) .. controls (177.07,33.22) and (178.05,32.23) .. (179.27,32.23) .. controls (180.48,32.23) and (181.47,33.22) .. (181.47,34.43) .. controls (181.47,35.65) and (180.48,36.63) .. (179.27,36.63) .. controls (178.05,36.63) and (177.07,35.65) .. (177.07,34.43) -- cycle ;

\draw [color={rgb, 255:red, 74; green, 144; blue, 226 }  ,draw opacity=1 ]   (62.4,24.3) .. controls (76.87,20.15) and (79.62,20.01) .. (96.15,23.86) ;
\draw [shift={(98,24.3)}, rotate = 193.29] [color={rgb, 255:red, 74; green, 144; blue, 226 }  ,draw opacity=1 ][line width=0.75]    (6.56,-1.97) .. controls (4.17,-0.84) and (1.99,-0.18) .. (0,0) .. controls (1.99,0.18) and (4.17,0.84) .. (6.56,1.97)   ;
\draw  [color={rgb, 255:red, 74; green, 144; blue, 226 }  ,draw opacity=1 ][fill={rgb, 255:red, 0; green, 0; blue, 0 }  ,fill opacity=1 ] (58,24.3) .. controls (58,23.08) and (58.98,22.1) .. (60.2,22.1) .. controls (61.42,22.1) and (62.4,23.08) .. (62.4,24.3) .. controls (62.4,25.52) and (61.42,26.5) .. (60.2,26.5) .. controls (58.98,26.5) and (58,25.52) .. (58,24.3) -- cycle ;
\draw  [color={rgb, 255:red, 74; green, 144; blue, 226 }  ,draw opacity=1 ][fill={rgb, 255:red, 0; green, 0; blue, 0 }  ,fill opacity=1 ] (98,24.3) .. controls (98,23.08) and (98.98,22.1) .. (100.2,22.1) .. controls (101.42,22.1) and (102.4,23.08) .. (102.4,24.3) .. controls (102.4,25.52) and (101.42,26.5) .. (100.2,26.5) .. controls (98.98,26.5) and (98,25.52) .. (98,24.3) -- cycle ;

\draw [color={rgb, 255:red, 74; green, 144; blue, 226 }  ,draw opacity=1 ]   (81.8,24.73) .. controls (94.92,31.38) and (100.97,31.62) .. (115.76,25.42) ;
\draw [shift={(117.4,24.73)}, rotate = 156.67] [color={rgb, 255:red, 74; green, 144; blue, 226 }  ,draw opacity=1 ][line width=0.75]    (6.56,-1.97) .. controls (4.17,-0.84) and (1.99,-0.18) .. (0,0) .. controls (1.99,0.18) and (4.17,0.84) .. (6.56,1.97)   ;
\draw  [color={rgb, 255:red, 74; green, 144; blue, 226 }  ,draw opacity=1 ][fill={rgb, 255:red, 0; green, 0; blue, 0 }  ,fill opacity=1 ] (77.4,24.73) .. controls (77.4,25.94) and (78.38,26.93) .. (79.6,26.93) .. controls (80.82,26.93) and (81.8,25.94) .. (81.8,24.73) .. controls (81.8,23.51) and (80.82,22.53) .. (79.6,22.53) .. controls (78.38,22.53) and (77.4,23.51) .. (77.4,24.73) -- cycle ;
\draw  [color={rgb, 255:red, 74; green, 144; blue, 226 }  ,draw opacity=1 ][fill={rgb, 255:red, 0; green, 0; blue, 0 }  ,fill opacity=1 ] (117.4,24.73) .. controls (117.4,25.94) and (118.38,26.93) .. (119.6,26.93) .. controls (120.82,26.93) and (121.8,25.94) .. (121.8,24.73) .. controls (121.8,23.51) and (120.82,22.53) .. (119.6,22.53) .. controls (118.38,22.53) and (117.4,23.51) .. (117.4,24.73) -- cycle ;

\draw  [color={rgb, 255:red, 126; green, 211; blue, 33 }  ,draw opacity=1 ][fill={rgb, 255:red, 0; green, 0; blue, 0 }  ,fill opacity=1 ] (144.33,29.3) .. controls (144.33,28.08) and (145.32,27.1) .. (146.53,27.1) .. controls (147.75,27.1) and (148.73,28.08) .. (148.73,29.3) .. controls (148.73,30.52) and (147.75,31.5) .. (146.53,31.5) .. controls (145.32,31.5) and (144.33,30.52) .. (144.33,29.3) -- cycle ;
\draw  [color={rgb, 255:red, 245; green, 166; blue, 35 }  ,draw opacity=1 ][fill={rgb, 255:red, 0; green, 0; blue, 0 }  ,fill opacity=1 ] (184.07,23.7) .. controls (184.07,22.48) and (185.05,21.5) .. (186.27,21.5) .. controls (187.48,21.5) and (188.47,22.48) .. (188.47,23.7) .. controls (188.47,24.92) and (187.48,25.9) .. (186.27,25.9) .. controls (185.05,25.9) and (184.07,24.92) .. (184.07,23.7) -- cycle ;
\draw  [color={rgb, 255:red, 126; green, 211; blue, 33 }  ,draw opacity=1 ][fill={rgb, 255:red, 0; green, 0; blue, 0 }  ,fill opacity=1 ] (123.53,28.5) .. controls (123.53,27.28) and (124.52,26.3) .. (125.73,26.3) .. controls (126.95,26.3) and (127.93,27.28) .. (127.93,28.5) .. controls (127.93,29.72) and (126.95,30.7) .. (125.73,30.7) .. controls (124.52,30.7) and (123.53,29.72) .. (123.53,28.5) -- cycle ;

\draw (42,28.5) node  [font=\scriptsize]  {$0$};
\draw (40.5,48.5) node  [font=\scriptsize]  {$-2$};
\draw (41.82,6.4) node  [font=\small]  {$\delta $};
\draw (282.49,66.47) node  [font=\small]  {$q_{2}$};
\draw (182.83,68.6) node  [font=\scriptsize]  {$0$};
\draw (160.33,68.6) node  [font=\scriptsize]  {$-2$};
\draw (222.83,68.6) node  [font=\scriptsize]  {$4$};
\draw (202.83,68.6) node  [font=\scriptsize]  {$2$};
\draw (242.83,68.6) node  [font=\scriptsize]  {$6$};
\draw (140.33,68.6) node  [font=\scriptsize]  {$-4$};
\draw (120.33,68.6) node  [font=\scriptsize]  {$-6$};
\draw (100.33,68.6) node  [font=\scriptsize]  {$-8$};
\draw (80.33,68.6) node  [font=\scriptsize]  {$-10$};
\draw (58.33,68.6) node  [font=\scriptsize]  {$-12$};
\draw (262.83,68.6) node  [font=\scriptsize]  {$8$};

\end{tikzpicture}

%% file: tikzpictures/table-kh-kt.tex
\tikzset{every picture/.style={line width=0.75pt}} 

\begin{tikzpicture}[x=0.75pt,y=0.75pt,yscale=-1.2,xscale=1.2]

\draw  [draw opacity=0] (49,31.04) -- (270,31.04) -- (270,71.57) -- (49,71.57) -- cycle ; \draw  [color={rgb, 255:red, 155; green, 155; blue, 155 }  ,draw opacity=1 ] (49,31.04) -- (49,71.57)(69,31.04) -- (69,71.57)(89,31.04) -- (89,71.57)(109,31.04) -- (109,71.57)(129,31.04) -- (129,71.57)(149,31.04) -- (149,71.57)(169,31.04) -- (169,71.57)(189,31.04) -- (189,71.57)(209,31.04) -- (209,71.57)(229,31.04) -- (229,71.57)(249,31.04) -- (249,71.57)(269,31.04) -- (269,71.57) ; \draw  [color={rgb, 255:red, 155; green, 155; blue, 155 }  ,draw opacity=1 ] (49,31.04) -- (270,31.04)(49,51.04) -- (270,51.04)(49,71.04) -- (270,71.04) ; \draw  [color={rgb, 255:red, 155; green, 155; blue, 155 }  ,draw opacity=1 ]  ;
\draw  [color={rgb, 255:red, 245; green, 166; blue, 35 }  ,draw opacity=1 ][fill={rgb, 255:red, 0; green, 0; blue, 0 }  ,fill opacity=1 ] (196.93,54.19) .. controls (196.93,52.98) and (197.92,51.99) .. (199.13,51.99) .. controls (200.35,51.99) and (201.33,52.98) .. (201.33,54.19) .. controls (201.33,55.41) and (200.35,56.39) .. (199.13,56.39) .. controls (197.92,56.39) and (196.93,55.41) .. (196.93,54.19) -- cycle ;
\draw  [color={rgb, 255:red, 245; green, 166; blue, 35 }  ,draw opacity=1 ][fill={rgb, 255:red, 0; green, 0; blue, 0 }  ,fill opacity=1 ] (237.93,55.19) .. controls (237.93,53.98) and (238.92,52.99) .. (240.13,52.99) .. controls (241.35,52.99) and (242.33,53.98) .. (242.33,55.19) .. controls (242.33,56.41) and (241.35,57.39) .. (240.13,57.39) .. controls (238.92,57.39) and (237.93,56.41) .. (237.93,55.19) -- cycle ;
\draw [color={rgb, 255:red, 245; green, 166; blue, 35 }  ,draw opacity=1 ]   (201.33,54.19) .. controls (214.46,47.54) and (221.44,48.23) .. (236.29,54.49) ;
\draw [shift={(237.93,55.19)}, rotate = 203.33] [color={rgb, 255:red, 245; green, 166; blue, 35 }  ,draw opacity=1 ][line width=0.75]    (6.56,-1.97) .. controls (4.17,-0.84) and (1.99,-0.18) .. (0,0) .. controls (1.99,0.18) and (4.17,0.84) .. (6.56,1.97)   ;

\draw  [color={rgb, 255:red, 208; green, 2; blue, 27 }  ,draw opacity=1 ][fill={rgb, 255:red, 0; green, 0; blue, 0 }  ,fill opacity=1 ] (257.33,67.26) .. controls (257.33,66.05) and (258.32,65.06) .. (259.53,65.06) .. controls (260.75,65.06) and (261.73,66.05) .. (261.73,67.26) .. controls (261.73,68.48) and (260.75,69.46) .. (259.53,69.46) .. controls (258.32,69.46) and (257.33,68.48) .. (257.33,67.26) -- cycle ;
\draw  [color={rgb, 255:red, 208; green, 2; blue, 27 }  ,draw opacity=1 ][fill={rgb, 255:red, 0; green, 0; blue, 0 }  ,fill opacity=1 ] (217.2,66.89) .. controls (217.2,65.68) and (218.18,64.69) .. (219.4,64.69) .. controls (220.62,64.69) and (221.6,65.68) .. (221.6,66.89) .. controls (221.6,68.11) and (220.62,69.09) .. (219.4,69.09) .. controls (218.18,69.09) and (217.2,68.11) .. (217.2,66.89) -- cycle ;
\draw [color={rgb, 255:red, 208; green, 2; blue, 27 }  ,draw opacity=1 ]   (201.33,67.72) .. controls (214.46,74.38) and (220.5,74.61) .. (235.29,68.42) ;
\draw [shift={(236.93,67.72)}, rotate = 156.67] [color={rgb, 255:red, 208; green, 2; blue, 27 }  ,draw opacity=1 ][line width=0.75]    (6.56,-1.97) .. controls (4.17,-0.84) and (1.99,-0.18) .. (0,0) .. controls (1.99,0.18) and (4.17,0.84) .. (6.56,1.97)   ;
\draw  [color={rgb, 255:red, 208; green, 2; blue, 27 }  ,draw opacity=1 ][fill={rgb, 255:red, 0; green, 0; blue, 0 }  ,fill opacity=1 ] (196.93,67.72) .. controls (196.93,68.93) and (197.92,69.92) .. (199.13,69.92) .. controls (200.35,69.92) and (201.33,68.93) .. (201.33,67.72) .. controls (201.33,66.5) and (200.35,65.52) .. (199.13,65.52) .. controls (197.92,65.52) and (196.93,66.5) .. (196.93,67.72) -- cycle ;
\draw  [color={rgb, 255:red, 208; green, 2; blue, 27 }  ,draw opacity=1 ][fill={rgb, 255:red, 0; green, 0; blue, 0 }  ,fill opacity=1 ] (236.93,67.72) .. controls (236.93,68.93) and (237.92,69.92) .. (239.13,69.92) .. controls (240.35,69.92) and (241.33,68.93) .. (241.33,67.72) .. controls (241.33,66.5) and (240.35,65.52) .. (239.13,65.52) .. controls (237.92,65.52) and (236.93,66.5) .. (236.93,67.72) -- cycle ;

\draw [color={rgb, 255:red, 245; green, 166; blue, 35 }  ,draw opacity=1 ]   (141.8,60.59) .. controls (156.28,56.44) and (159.02,56.3) .. (175.55,60.16) ;
\draw [shift={(177.4,60.59)}, rotate = 193.29] [color={rgb, 255:red, 245; green, 166; blue, 35 }  ,draw opacity=1 ][line width=0.75]    (6.56,-1.97) .. controls (4.17,-0.84) and (1.99,-0.18) .. (0,0) .. controls (1.99,0.18) and (4.17,0.84) .. (6.56,1.97)   ;
\draw  [color={rgb, 255:red, 245; green, 166; blue, 35 }  ,draw opacity=1 ][fill={rgb, 255:red, 0; green, 0; blue, 0 }  ,fill opacity=1 ] (137.4,60.59) .. controls (137.4,59.38) and (138.38,58.39) .. (139.6,58.39) .. controls (140.82,58.39) and (141.8,59.38) .. (141.8,60.59) .. controls (141.8,61.81) and (140.82,62.79) .. (139.6,62.79) .. controls (138.38,62.79) and (137.4,61.81) .. (137.4,60.59) -- cycle ;
\draw  [color={rgb, 255:red, 245; green, 166; blue, 35 }  ,draw opacity=1 ][fill={rgb, 255:red, 0; green, 0; blue, 0 }  ,fill opacity=1 ] (177.4,60.59) .. controls (177.4,59.38) and (178.38,58.39) .. (179.6,58.39) .. controls (180.82,58.39) and (181.8,59.38) .. (181.8,60.59) .. controls (181.8,61.81) and (180.82,62.79) .. (179.6,62.79) .. controls (178.38,62.79) and (177.4,61.81) .. (177.4,60.59) -- cycle ;
\draw  [color={rgb, 255:red, 245; green, 166; blue, 35 }  ,draw opacity=1 ][fill={rgb, 255:red, 0; green, 0; blue, 0 }  ,fill opacity=1 ] (217.2,60.59) .. controls (217.2,59.38) and (218.18,58.39) .. (219.4,58.39) .. controls (220.62,58.39) and (221.6,59.38) .. (221.6,60.59) .. controls (221.6,61.81) and (220.62,62.79) .. (219.4,62.79) .. controls (218.18,62.79) and (217.2,61.81) .. (217.2,60.59) -- cycle ;
\draw [color={rgb, 255:red, 245; green, 166; blue, 35 }  ,draw opacity=1 ]   (181.8,60.59) .. controls (196.28,56.44) and (198.83,56.3) .. (215.35,60.16) ;
\draw [shift={(217.2,60.59)}, rotate = 193.29] [color={rgb, 255:red, 245; green, 166; blue, 35 }  ,draw opacity=1 ][line width=0.75]    (6.56,-1.97) .. controls (4.17,-0.84) and (1.99,-0.18) .. (0,0) .. controls (1.99,0.18) and (4.17,0.84) .. (6.56,1.97)   ;

\draw [color={rgb, 255:red, 245; green, 166; blue, 35 }  ,draw opacity=1 ]   (161.2,61.62) .. controls (175.41,65.47) and (181.2,65.01) .. (195.03,62.01) ;
\draw [shift={(196.8,61.62)}, rotate = 167.53] [color={rgb, 255:red, 245; green, 166; blue, 35 }  ,draw opacity=1 ][line width=0.75]    (6.56,-1.97) .. controls (4.17,-0.84) and (1.99,-0.18) .. (0,0) .. controls (1.99,0.18) and (4.17,0.84) .. (6.56,1.97)   ;
\draw  [color={rgb, 255:red, 245; green, 166; blue, 35 }  ,draw opacity=1 ][fill={rgb, 255:red, 0; green, 0; blue, 0 }  ,fill opacity=1 ] (156.8,61.62) .. controls (156.8,62.83) and (157.78,63.82) .. (159,63.82) .. controls (160.22,63.82) and (161.2,62.83) .. (161.2,61.62) .. controls (161.2,60.4) and (160.22,59.42) .. (159,59.42) .. controls (157.78,59.42) and (156.8,60.4) .. (156.8,61.62) -- cycle ;
\draw  [color={rgb, 255:red, 245; green, 166; blue, 35 }  ,draw opacity=1 ][fill={rgb, 255:red, 0; green, 0; blue, 0 }  ,fill opacity=1 ] (196.8,61.62) .. controls (196.8,62.83) and (197.78,63.82) .. (199,63.82) .. controls (200.22,63.82) and (201.2,62.83) .. (201.2,61.62) .. controls (201.2,60.4) and (200.22,59.42) .. (199,59.42) .. controls (197.78,59.42) and (196.8,60.4) .. (196.8,61.62) -- cycle ;

\draw  [color={rgb, 255:red, 126; green, 211; blue, 33 }  ,draw opacity=1 ][fill={rgb, 255:red, 0; green, 0; blue, 0 }  ,fill opacity=1 ] (177.77,67.63) .. controls (177.77,66.41) and (178.75,65.43) .. (179.97,65.43) .. controls (181.18,65.43) and (182.17,66.41) .. (182.17,67.63) .. controls (182.17,68.84) and (181.18,69.83) .. (179.97,69.83) .. controls (178.75,69.83) and (177.77,68.84) .. (177.77,67.63) -- cycle ;
\draw  [color={rgb, 255:red, 126; green, 211; blue, 33 }  ,draw opacity=1 ][fill={rgb, 255:red, 0; green, 0; blue, 0 }  ,fill opacity=1 ] (137.8,67.19) .. controls (137.8,65.98) and (138.78,64.99) .. (140,64.99) .. controls (141.22,64.99) and (142.2,65.98) .. (142.2,67.19) .. controls (142.2,68.41) and (141.22,69.39) .. (140,69.39) .. controls (138.78,69.39) and (137.8,68.41) .. (137.8,67.19) -- cycle ;
\draw [color={rgb, 255:red, 126; green, 211; blue, 33 }  ,draw opacity=1 ]   (121.1,67.72) .. controls (134.22,74.38) and (140.27,74.61) .. (155.06,68.42) ;
\draw [shift={(156.7,67.72)}, rotate = 156.67] [color={rgb, 255:red, 126; green, 211; blue, 33 }  ,draw opacity=1 ][line width=0.75]    (6.56,-1.97) .. controls (4.17,-0.84) and (1.99,-0.18) .. (0,0) .. controls (1.99,0.18) and (4.17,0.84) .. (6.56,1.97)   ;
\draw  [color={rgb, 255:red, 126; green, 211; blue, 33 }  ,draw opacity=1 ][fill={rgb, 255:red, 0; green, 0; blue, 0 }  ,fill opacity=1 ] (116.7,67.72) .. controls (116.7,68.93) and (117.68,69.92) .. (118.9,69.92) .. controls (120.12,69.92) and (121.1,68.93) .. (121.1,67.72) .. controls (121.1,66.5) and (120.12,65.52) .. (118.9,65.52) .. controls (117.68,65.52) and (116.7,66.5) .. (116.7,67.72) -- cycle ;
\draw  [color={rgb, 255:red, 126; green, 211; blue, 33 }  ,draw opacity=1 ][fill={rgb, 255:red, 0; green, 0; blue, 0 }  ,fill opacity=1 ] (156.7,67.72) .. controls (156.7,68.93) and (157.68,69.92) .. (158.9,69.92) .. controls (160.12,69.92) and (161.1,68.93) .. (161.1,67.72) .. controls (161.1,66.5) and (160.12,65.52) .. (158.9,65.52) .. controls (157.68,65.52) and (156.7,66.5) .. (156.7,67.72) -- cycle ;

\draw  [color={rgb, 255:red, 245; green, 166; blue, 35 }  ,draw opacity=1 ][fill={rgb, 255:red, 0; green, 0; blue, 0 }  ,fill opacity=1 ] (177.4,54.19) .. controls (177.4,52.98) and (178.38,51.99) .. (179.6,51.99) .. controls (180.82,51.99) and (181.8,52.98) .. (181.8,54.19) .. controls (181.8,55.41) and (180.82,56.39) .. (179.6,56.39) .. controls (178.38,56.39) and (177.4,55.41) .. (177.4,54.19) -- cycle ;
\draw  [color={rgb, 255:red, 126; green, 211; blue, 33 }  ,draw opacity=1 ][fill={rgb, 255:red, 0; green, 0; blue, 0 }  ,fill opacity=1 ] (73.5,47.02) .. controls (73.5,48.23) and (74.48,49.22) .. (75.7,49.22) .. controls (76.92,49.22) and (77.9,48.23) .. (77.9,47.02) .. controls (77.9,45.8) and (76.92,44.82) .. (75.7,44.82) .. controls (74.48,44.82) and (73.5,45.8) .. (73.5,47.02) -- cycle ;
\draw  [color={rgb, 255:red, 126; green, 211; blue, 33 }  ,draw opacity=1 ][fill={rgb, 255:red, 0; green, 0; blue, 0 }  ,fill opacity=1 ] (113.5,47.02) .. controls (113.5,48.23) and (114.48,49.22) .. (115.7,49.22) .. controls (116.92,49.22) and (117.9,48.23) .. (117.9,47.02) .. controls (117.9,45.8) and (116.92,44.82) .. (115.7,44.82) .. controls (114.48,44.82) and (113.5,45.8) .. (113.5,47.02) -- cycle ;
\draw [color={rgb, 255:red, 126; green, 211; blue, 33 }  ,draw opacity=1 ]   (77.9,47.02) .. controls (91.02,53.68) and (97.07,53.91) .. (111.86,47.72) ;
\draw [shift={(113.5,47.02)}, rotate = 156.67] [color={rgb, 255:red, 126; green, 211; blue, 33 }  ,draw opacity=1 ][line width=0.75]    (6.56,-1.97) .. controls (4.17,-0.84) and (1.99,-0.18) .. (0,0) .. controls (1.99,0.18) and (4.17,0.84) .. (6.56,1.97)   ;

\draw  [color={rgb, 255:red, 245; green, 166; blue, 35 }  ,draw opacity=1 ][fill={rgb, 255:red, 0; green, 0; blue, 0 }  ,fill opacity=1 ] (136.93,36.79) .. controls (136.93,35.58) and (137.92,34.59) .. (139.13,34.59) .. controls (140.35,34.59) and (141.33,35.58) .. (141.33,36.79) .. controls (141.33,38.01) and (140.35,38.99) .. (139.13,38.99) .. controls (137.92,38.99) and (136.93,38.01) .. (136.93,36.79) -- cycle ;
\draw  [color={rgb, 255:red, 245; green, 166; blue, 35 }  ,draw opacity=1 ][fill={rgb, 255:red, 0; green, 0; blue, 0 }  ,fill opacity=1 ] (181.73,36.79) .. controls (181.73,35.58) and (182.72,34.59) .. (183.93,34.59) .. controls (185.15,34.59) and (186.13,35.58) .. (186.13,36.79) .. controls (186.13,38.01) and (185.15,38.99) .. (183.93,38.99) .. controls (182.72,38.99) and (181.73,38.01) .. (181.73,36.79) -- cycle ;
\draw [color={rgb, 255:red, 245; green, 166; blue, 35 }  ,draw opacity=1 ]   (161.3,36.79) .. controls (174.42,30.14) and (180.47,29.9) .. (195.26,36.09) ;
\draw [shift={(196.9,36.79)}, rotate = 203.33] [color={rgb, 255:red, 245; green, 166; blue, 35 }  ,draw opacity=1 ][line width=0.75]    (6.56,-1.97) .. controls (4.17,-0.84) and (1.99,-0.18) .. (0,0) .. controls (1.99,0.18) and (4.17,0.84) .. (6.56,1.97)   ;
\draw  [color={rgb, 255:red, 245; green, 166; blue, 35 }  ,draw opacity=1 ][fill={rgb, 255:red, 0; green, 0; blue, 0 }  ,fill opacity=1 ] (156.9,36.79) .. controls (156.9,35.58) and (157.88,34.59) .. (159.1,34.59) .. controls (160.32,34.59) and (161.3,35.58) .. (161.3,36.79) .. controls (161.3,38.01) and (160.32,38.99) .. (159.1,38.99) .. controls (157.88,38.99) and (156.9,38.01) .. (156.9,36.79) -- cycle ;
\draw  [color={rgb, 255:red, 245; green, 166; blue, 35 }  ,draw opacity=1 ][fill={rgb, 255:red, 0; green, 0; blue, 0 }  ,fill opacity=1 ] (196.9,36.79) .. controls (196.9,35.58) and (197.88,34.59) .. (199.1,34.59) .. controls (200.32,34.59) and (201.3,35.58) .. (201.3,36.79) .. controls (201.3,38.01) and (200.32,38.99) .. (199.1,38.99) .. controls (197.88,38.99) and (196.9,38.01) .. (196.9,36.79) -- cycle ;

\draw [color={rgb, 255:red, 126; green, 211; blue, 33 }  ,draw opacity=1 ]   (100.57,44.73) .. controls (115.04,40.58) and (117.78,40.43) .. (134.31,44.29) ;
\draw [shift={(136.17,44.73)}, rotate = 193.29] [color={rgb, 255:red, 126; green, 211; blue, 33 }  ,draw opacity=1 ][line width=0.75]    (6.56,-1.97) .. controls (4.17,-0.84) and (1.99,-0.18) .. (0,0) .. controls (1.99,0.18) and (4.17,0.84) .. (6.56,1.97)   ;
\draw  [color={rgb, 255:red, 126; green, 211; blue, 33 }  ,draw opacity=1 ][fill={rgb, 255:red, 0; green, 0; blue, 0 }  ,fill opacity=1 ] (96.17,44.73) .. controls (96.17,43.51) and (97.15,42.53) .. (98.37,42.53) .. controls (99.58,42.53) and (100.57,43.51) .. (100.57,44.73) .. controls (100.57,45.94) and (99.58,46.93) .. (98.37,46.93) .. controls (97.15,46.93) and (96.17,45.94) .. (96.17,44.73) -- cycle ;
\draw  [color={rgb, 255:red, 126; green, 211; blue, 33 }  ,draw opacity=1 ][fill={rgb, 255:red, 0; green, 0; blue, 0 }  ,fill opacity=1 ] (136.17,44.73) .. controls (136.17,43.51) and (137.15,42.53) .. (138.37,42.53) .. controls (139.58,42.53) and (140.57,43.51) .. (140.57,44.73) .. controls (140.57,45.94) and (139.58,46.93) .. (138.37,46.93) .. controls (137.15,46.93) and (136.17,45.94) .. (136.17,44.73) -- cycle ;
\draw  [color={rgb, 255:red, 126; green, 211; blue, 33 }  ,draw opacity=1 ][fill={rgb, 255:red, 0; green, 0; blue, 0 }  ,fill opacity=1 ] (175.97,44.73) .. controls (175.97,43.51) and (176.95,42.53) .. (178.17,42.53) .. controls (179.38,42.53) and (180.37,43.51) .. (180.37,44.73) .. controls (180.37,45.94) and (179.38,46.93) .. (178.17,46.93) .. controls (176.95,46.93) and (175.97,45.94) .. (175.97,44.73) -- cycle ;
\draw [color={rgb, 255:red, 126; green, 211; blue, 33 }  ,draw opacity=1 ]   (140.57,44.73) .. controls (155.04,40.58) and (157.6,40.43) .. (174.11,44.29) ;
\draw [shift={(175.97,44.73)}, rotate = 193.29] [color={rgb, 255:red, 126; green, 211; blue, 33 }  ,draw opacity=1 ][line width=0.75]    (6.56,-1.97) .. controls (4.17,-0.84) and (1.99,-0.18) .. (0,0) .. controls (1.99,0.18) and (4.17,0.84) .. (6.56,1.97)   ;

\draw [color={rgb, 255:red, 126; green, 211; blue, 33 }  ,draw opacity=1 ]   (125.9,47.02) .. controls (139.02,53.68) and (145.07,53.91) .. (159.86,47.72) ;
\draw [shift={(161.5,47.02)}, rotate = 156.67] [color={rgb, 255:red, 126; green, 211; blue, 33 }  ,draw opacity=1 ][line width=0.75]    (6.56,-1.97) .. controls (4.17,-0.84) and (1.99,-0.18) .. (0,0) .. controls (1.99,0.18) and (4.17,0.84) .. (6.56,1.97)   ;
\draw  [color={rgb, 255:red, 126; green, 211; blue, 33 }  ,draw opacity=1 ][fill={rgb, 255:red, 0; green, 0; blue, 0 }  ,fill opacity=1 ] (121.5,47.02) .. controls (121.5,48.23) and (122.48,49.22) .. (123.7,49.22) .. controls (124.92,49.22) and (125.9,48.23) .. (125.9,47.02) .. controls (125.9,45.8) and (124.92,44.82) .. (123.7,44.82) .. controls (122.48,44.82) and (121.5,45.8) .. (121.5,47.02) -- cycle ;
\draw  [color={rgb, 255:red, 126; green, 211; blue, 33 }  ,draw opacity=1 ][fill={rgb, 255:red, 0; green, 0; blue, 0 }  ,fill opacity=1 ] (161.5,47.02) .. controls (161.5,48.23) and (162.48,49.22) .. (163.7,49.22) .. controls (164.92,49.22) and (165.9,48.23) .. (165.9,47.02) .. controls (165.9,45.8) and (164.92,44.82) .. (163.7,44.82) .. controls (162.48,44.82) and (161.5,45.8) .. (161.5,47.02) -- cycle ;

\draw  [color={rgb, 255:red, 74; green, 144; blue, 226 }  ,draw opacity=1 ][fill={rgb, 255:red, 0; green, 0; blue, 0 }  ,fill opacity=1 ] (56.9,35.35) .. controls (56.9,34.14) and (57.88,33.15) .. (59.1,33.15) .. controls (60.32,33.15) and (61.3,34.14) .. (61.3,35.35) .. controls (61.3,36.57) and (60.32,37.55) .. (59.1,37.55) .. controls (57.88,37.55) and (56.9,36.57) .. (56.9,35.35) -- cycle ;
\draw  [color={rgb, 255:red, 74; green, 144; blue, 226 }  ,draw opacity=1 ][fill={rgb, 255:red, 0; green, 0; blue, 0 }  ,fill opacity=1 ] (96.17,35.35) .. controls (96.17,34.14) and (97.15,33.15) .. (98.37,33.15) .. controls (99.58,33.15) and (100.57,34.14) .. (100.57,35.35) .. controls (100.57,36.57) and (99.58,37.55) .. (98.37,37.55) .. controls (97.15,37.55) and (96.17,36.57) .. (96.17,35.35) -- cycle ;
\draw [color={rgb, 255:red, 74; green, 144; blue, 226 }  ,draw opacity=1 ]   (81.17,35.35) .. controls (94.29,42.01) and (100.34,42.24) .. (115.12,36.05) ;
\draw [shift={(116.77,35.35)}, rotate = 156.67] [color={rgb, 255:red, 74; green, 144; blue, 226 }  ,draw opacity=1 ][line width=0.75]    (6.56,-1.97) .. controls (4.17,-0.84) and (1.99,-0.18) .. (0,0) .. controls (1.99,0.18) and (4.17,0.84) .. (6.56,1.97)   ;
\draw  [color={rgb, 255:red, 74; green, 144; blue, 226 }  ,draw opacity=1 ][fill={rgb, 255:red, 0; green, 0; blue, 0 }  ,fill opacity=1 ] (76.77,35.35) .. controls (76.77,36.57) and (77.75,37.55) .. (78.97,37.55) .. controls (80.18,37.55) and (81.17,36.57) .. (81.17,35.35) .. controls (81.17,34.14) and (80.18,33.15) .. (78.97,33.15) .. controls (77.75,33.15) and (76.77,34.14) .. (76.77,35.35) -- cycle ;
\draw  [color={rgb, 255:red, 74; green, 144; blue, 226 }  ,draw opacity=1 ][fill={rgb, 255:red, 0; green, 0; blue, 0 }  ,fill opacity=1 ] (116.77,35.35) .. controls (116.77,36.57) and (117.75,37.55) .. (118.97,37.55) .. controls (120.18,37.55) and (121.17,36.57) .. (121.17,35.35) .. controls (121.17,34.14) and (120.18,33.15) .. (118.97,33.15) .. controls (117.75,33.15) and (116.77,34.14) .. (116.77,35.35) -- cycle ;

\draw  [color={rgb, 255:red, 126; green, 211; blue, 33 }  ,draw opacity=1 ][fill={rgb, 255:red, 0; green, 0; blue, 0 }  ,fill opacity=1 ] (142.9,47.13) .. controls (142.9,45.91) and (143.88,44.93) .. (145.1,44.93) .. controls (146.32,44.93) and (147.3,45.91) .. (147.3,47.13) .. controls (147.3,48.34) and (146.32,49.33) .. (145.1,49.33) .. controls (143.88,49.33) and (142.9,48.34) .. (142.9,47.13) -- cycle ;
\draw  [color={rgb, 255:red, 245; green, 166; blue, 35 }  ,draw opacity=1 ][fill={rgb, 255:red, 0; green, 0; blue, 0 }  ,fill opacity=1 ] (172.4,36.79) .. controls (172.4,35.58) and (173.38,34.59) .. (174.6,34.59) .. controls (175.82,34.59) and (176.8,35.58) .. (176.8,36.79) .. controls (176.8,38.01) and (175.82,38.99) .. (174.6,38.99) .. controls (173.38,38.99) and (172.4,38.01) .. (172.4,36.79) -- cycle ;

\draw (41,39.99) node  [font=\scriptsize]  {$0$};
\draw (39.5,59.99) node  [font=\scriptsize]  {$-2$};
\draw (40.82,17.89) node  [font=\small]  {$\delta $};
\draw (281.49,74.96) node  [font=\small]  {$q_{2}$};
\draw (181.83,80.09) node  [font=\scriptsize]  {$0$};
\draw (159.33,80.09) node  [font=\scriptsize]  {$-2$};
\draw (221.83,80.09) node  [font=\scriptsize]  {$4$};
\draw (201.83,80.09) node  [font=\scriptsize]  {$2$};
\draw (241.83,80.09) node  [font=\scriptsize]  {$6$};
\draw (139.33,80.09) node  [font=\scriptsize]  {$-4$};
\draw (119.33,80.09) node  [font=\scriptsize]  {$-6$};
\draw (99.33,80.09) node  [font=\scriptsize]  {$-8$};
\draw (79.33,80.09) node  [font=\scriptsize]  {$-10$};
\draw (57.33,80.09) node  [font=\scriptsize]  {$-12$};
\draw (261.83,80.09) node  [font=\scriptsize]  {$8$};

\end{tikzpicture}

%% file: acknowledgements.tex
\subsection*{Acknowledgement}

The author thanks Mikhail Khovanov for suggesting the study of $\sl_2$-actions on Khovanov homology. 
The author also thanks Motoo Tange and Tetsuya Abe for giving him the opportunity to present this work at the conference \textit{Differential Topology '26}. 
The author was supported by JSPS KAKENHI Grant Number 23K12982 and academist crowdfunding Project No.\ 121.

%% file: 2.tex
\section
[A-modules and y-ifications]
{$\mcA$-modules and $y$-ifications}
\label{sec:formal}

\subsection
[Differential graded algebra A]
{Differential graded algebra $\mcA$}
\label{subsec:dga-A}

Here, we introduce a simplified version of the differential graded algebra (dga) $\mcA$ defined in \cite[Definitions 3.4, 3.15]{GHM:symmety-kr2024}. Throughout this section, we let $R$ be an arbitrary commutative ring and $n$ a positive integer. Let 
\[
    R_n := R[x_1, \ldots, x_n]
\]
and 
\[
    R^e_n := R[x_1, \ldots, x_n; x'_1, \ldots, x'_n]. 
\]
Also, consider the quotient of $R^e_n$ given by 
\[
    B_n := \frac{R^e_n}{\left(\sum_{i=1}^n x_i^k = \sum_{i=1}^n (x'_i)^k\ ;\ k = 1, \ldots, n\right)}.
\]
Both $R^e_n$ and $B_n$ are naturally $R_n$-$R_n$ bimodules. In particular, when $R$ is a field of characteristic $0$, one may alternatively write
\[
    B_n = R^n \otimes_{R_n^{S}} R^n
\]
where $R_n^{S}$ is the algebra of symmetric polynomials in $n$ variables. Rings $R_n, R^e_n, B_n$ are given graded ring structures by assigning \textit{quantum degrees}
\[
    \deg_q(x_i) = \deg_q(x'_i) = 2.
\]

\begin{defn}
\label{def:dga-A}
    Let $w \in \mfS_n$ be a permutation. A dga $\mcA^w_n$ over $B_n$ is defined as follows. Let $\mcA^w_n$ be the graded-commutative $B_n$-algebra freely generated by elements $\xi_i$ $(i = 1, \dots, n)$ of degree $-1$ and $u$ of degree $-2$, modulo the ideal generated by the squares of $\xi_i$. The differential $d$ of degree $1$ is defined by
    \[
        d(\xi_i) = x_i - x'_{w(i)},\quad 
        d(u) = \sum_{i = 1}^n (x_i + x'_{w(i)}) \xi_i.
    \]
    $\mcA^w_n$ is also equipped with the \textit{quantum grading}, compatible with that of $B_n$, by additionally imposing $\deg_q(\xi_i) = 2$ and $\deg_q(u) = 4$.
\end{defn}

That $d$ is indeed a differential can be seen by   
\begin{align*}
    d^2(u) 
        &= \sum_i (x_i + x'_{w(i)}) d(\xi_i) \\
        &= \sum_i (x^2_i - x'^2_{w(i)}) \\
        &= \sum_i x^2_i - \sum_i x'^2_i \\
        &= 0.
\end{align*}

\begin{rem}
    Compared with the definition of $\mcA^w_n$ given in \cite[Definitions 3.4, 3.15]{GHM:symmety-kr2024}, we have dropped the generators $u_k\ (k = 1, \ldots, n)$ except for $u = u_2$.
\end{rem}

The primary grading of $\mcA^w_n$ is called the \textit{homological grading}, and its degree function is denoted $\deg_t$. For brevity, we write $|a| = \deg_t(a)$ for $a \in \mcA^w_n$. Throughout this paper, we make use of the \textit{super-commutator}
\[
    [a, b] = a b - (-1)^{|a||b|} b a.
\]

\begin{defn}
\label{def:mcA-module}
    Let $C$ be a bigraded chain complex over $B_n$ with differential $d$ of degree $(1, 0)$. $C$ is an \textit{$\mcA^w_n$-module} if it admits an action of $\mcA^w_n$ such that for each element $a \in \mcA^w_n$, the action of $d(a)$ coincides with the action of the super-commutator
    \[
        [d, a] = d \circ a - (-1)^{|a|} a \circ d.
    \]
    Let $\Kom_{\mcA^w_n}(\mcC)$ denote the category of \textit{bounded} $\mcA^w_n$-modules and $\mcA^w_n$-linear chain maps, i.e.\ chain maps that commute with the actions of $\mcA^w_n$. 
\end{defn}

\subsection
[y-ifications from A-modules]
{$y$-ifications from $\mcA$-modules}

Following \cite[Section 2.5]{GH:y-ification2022}, we define a \textit{$y$-ification} of a complex $C$ in an arbitrary graded $R^e_n$-linear category $\mcC$. For a non-negative multi-exponent $\alpha = (\alpha_1, \ldots, \alpha_n) \in (\ZZ_{\geq 0})^n$, we let $y^\alpha$ denote the monomial $\Pi_i y_i^{\alpha_i} \in R[y_1, \ldots, y_n]$. 

\begin{defn}
\label{def:formal-tensor-product}
    Let $\mcC$ be a graded $R$-linear category and $C = \{C_p\}_{p \in \ZZ}$ a bounded sequence of objects in $\mcC$. The \textit{formal tensor product} of $C$ and the polynomial ring $R[\bfy] = R[y_1, \ldots, y_n]$, denoted $C[\mathbf{y}] = C \otimes R[y_1, \ldots, y_n]$, is an infinite sequence of objects in $\mcC$, whose $p$-th object is given by the finite direct sum
    \[
        (C[\mathbf{y}])^p := \bigoplus_{k \geq 0} \bigoplus_{|\alpha| = k} C^{p - 2k} y^\alpha.
    \]
    Here, $k$ ranges over non-negative integers, and $\alpha$ ranges over multi-exponents of total degree
    \[
        |\alpha| = \sum_i \alpha_i = k.
    \]
    Each $y_i$ is assigned bidegree $(2, -2)$, and $C^{p - 2k} y^\alpha$ is a copy of $C^{p - 2k}$ with bidegree shifted by $(2k, -2k)$. For a sequence of maps $f = \{f_p\}_{p \in \ZZ}\colon C \to C'$ of degree $a$ and for any multi-exponent $\gamma$, a map
    \[
        f \otimes y^\gamma\colon C[\bfy] \to C'[\bfy]
    \]
    of degree $a + 2|\gamma|$ is defined by the sequence of morphisms
    \[
    \begin{tikzcd}
        C^p y^\alpha \arrow[r, "f^p"] & C'^{p+a} y^{\alpha+\gamma}.
    \end{tikzcd}
    \]
    Hereafter, $f \otimes y^\gamma$ will be denoted $f y^\gamma$ when there is no risk of confusion.
\end{defn}

\begin{defn}
\label{def:y-ification}
    Let $\mcC$ be a graded $R^e_n$-linear category, and $C$ a bounded chain complex in $\mcC$ with differential $d$. A \textit{$y$-ification of $C$} is a triple $(C[\mathbf{y}], w, D)$ satisfying the following conditions: (i) $w \in \mfS_n$ is a permutation, (ii) $D$ is an endomorphism of $C[\mathbf{y}]$ such that $D \equiv d \bmod{(\mathbf{y})}$, and (iii) $(C[\mathbf{y}], D)$ is a \textit{curved complex} with curvature
    \[
        D^2 = \sum_{i=1}^n (x_i - x'_{w(i)}) y_i. 
    \]
    Let $y\Kom(\mcC)$ denote the category of $y$-ified complexes and (curved) chain maps. 
\end{defn}

Note that the curvature $D^2$ of $D$ is central in $\mcC$, i.e.\ it commutes with all morphisms in $\mcC$. See \cite[Appendix A]{GH:y-ification2022} for a general treatment of curved complexes. We call an object of $y\Kom(\mcC)$ a \textit{$y$-complex}, and a morphism in $y\Kom(\mcC)$ a \textit{$y$-chain map}. 


\begin{prop}
\label{prop:y-ification-from-dga}
    There is a homotopy invariant faithful functor
    \[
        y\colon \Kom_{\mcA^w_n}(\mcC) \to y\Kom(\mcC)
    \]
    such that for each $\mcA^w_n$-module $C$, it assigns a $y$-ification $y(C) := (C[\mathbf{y}], w, D)$ of $C$ with differential
    \[
        D := d + \sum_i \xi_i y_i,
    \]
    and for each $\mcA^w_n$-homomorphism $f\colon C \to C'$, it assigns a $y$-chain map 
    \[
        y(f) := f \otimes 1 \colon C[\bfy] \to C'[\bfy].
    \]
\end{prop}

\begin{proof}
    We verify that $y(C)$ is a $y$-complex. Conditions (i), (ii) of \Cref{def:y-ification} are obvious. For condition (iii), we have
    \[
        D^2 = d^2 + \sum_i (d \xi_i + \xi_i d) y_i + \sum_{i, j} \xi_i \xi_j y_i y_j.
    \]
    The first term is $0$, and the second term is exactly the curvature, and the third term is 
    \[
        \sum_i \xi_i^2 y_i ^2 + \sum_{i < j} [\xi_i, \xi_j] y_i y_j = 0
    \]
    from the conditions given in \Cref{def:dga-A}. 
\end{proof}

We call $y(C)$ in \Cref{prop:y-ification-from-dga} the \textit{$y$-ification associated to the $\mcA^w_n$-module $C$}. Next, following the formula given in \cite[Theorem 1.9]{GHM:symmety-kr2024}, we define an operator $\sle$ on $y(C)$ using the action of $u$.\footnote{
    In \cite{GHM:symmety-kr2024}, the corresponding operator is denoted $F_2$. In \cite{CG:structure-in-homfly2024}, it is replaced with $E$ to emphasize that it corresponds to the raising operator for an $\sl_2$-action (see Remark 3.24 therein).
}

\begin{defn}
    We retain assumptions of \Cref{def:y-ification}. The $i$-th \textit{formal partial derivative}
    \[
        \ddel{}{y_i}\colon C[\bfy] \to C[\bfy]
    \]
    is defined by an array of maps
    \[
    \begin{tikzcd}
        C^p y^\alpha \arrow[r, "\alpha_i"] & C^p y^{\alpha - e_i}.
    \end{tikzcd}
    \]
    Here, $e_i$ denotes the multi-exponent whose $i$-th component is $1$ and all others are $0$. 
\end{defn}

\begin{lem}
    For maps $f, g\colon C \to C'$ satisfying $[f, g] = 0$, we have 
    \[
        [f \ddel{}{y_i}, g y_j] = fg \delta_{i, j}.
    \]
\end{lem}


\begin{prop}
\label{prop:action-e}
    For an $\mcA^w_n$-module $C$, 
    \begin{equation}
    \label{eqn:e-def}
        \sle := u + \sum_i (x_i + x'_{w(i)}) \ddel{}{y_i}
    \end{equation}
    gives a chain endomorphism of $y(C)$ of bidegree $(-2, 4)$.
\end{prop}

\begin{proof}
    Compute
    \[
        [D, \sle] 
            = [d, u] 
            + \sum_i [\xi_i, u] y_i
            + \sum_i [d, x_i + x'_{w(i)}] \ddel{}{y_i}
            + \sum_{i, j} \xi_i (x_j + x'_{w(j)}) [y_i, \ddel{}{y_j}].
    \]
    The second and third terms vanish, and the first and fourth terms cancel. 
\end{proof}

We call the endomorphism $\sle$ in \Cref{prop:action-e} the \textit{$\sle$-operator associated to the $\mcA^w_n$-module $C$}. To formulate \Cref{prop:action-e} categorically, we introduce the category of \textit{$\sle$-equivariant $y$-complexes}, where the equivariance is considered up to homotopy. 

\begin{defn}
\label{def:e-equiv-y-complexes}
    The \textit{category of $\sle$-equivariant $y$-complexes}, denoted $y\Kom_\sle(\mcC)$, is defined as follows. An object of $y\Kom_\sle(\mcC)$ is a pair $(C, \sle)$ of a $y$-complex $C$ in $y\Kom(\mcC)$ and an endomorphism $\sle$ of bidegree $(-2, 4)$. A morphism $(F, \tilde{F})$ between objects $(C, \sle), (C', \sle')$ is a pair of a $y$-chain map $F: C \to C'$ and a degree $-3$ homotopy $\tilde{F}$ satisfying 
    \[
        [D, \tilde{F}] + [\sle, F] = 0.
    \]
    The identity morphism of an object $C$ is given by $(I_C, 0)$, and the composition map is defined by
    \[
        (G, \tilde{G}) \circ (F, \tilde{F}) := 
        (GF,\ \tilde{G} F + G \tilde{F}). 
    \]
    A morphism $(F, \tilde{F})$ in $y\Kom_\sle(\mcC)$ is called an \textit{$\sle$-equivariant $y$-chain map}. In particular, if $\tilde{F} = 0$, it is called a \textit{strictly $\sle$-equivariant $y$-chain map}. 
\end{defn}

One easily sees that $y\Kom_\sle(\mcC)$ is indeed a category. Note that there is an obvious forgetful functor $y\Kom_{\sle}(\mcC) \to y\Kom(\mcC)$.

\begin{defn}
\label{def:e-equiv-homotopy}
    Let $(F, \tilde{F}), (G, \tilde{G}')\colon C \to C'$ be two $\sle$-equivariant $y$-chain maps. We say $(F, \tilde{F}), (G, \tilde{G}')$ are \textit{(equivariantly) homotopic} if there is a pair of homotopies $(H, \tilde{H})$, where $H$ is a degree $-1$ homotopy giving
    \[
        F - G = [D, H]
    \]
    and $\tilde{H}$ is a degree $-4$ homotopy satisfying
    \[
        \tilde{F} - \tilde{G} = [D, \tilde{H}] + [\sle, H].
    \]
\end{defn}

The degree $-4$ homotopy $\tilde{H}$ in the \Cref{def:e-equiv-homotopy} can be understood as a $3$-cell filling the following cylinder consisting of $2$-cells made of four degree $-3$ homotopies: 
\[
    \begin{tikzcd}[column sep=5em, row sep=3.5em]
    C
        \arrow[d, equal]
        \arrow[r, "\sle F"{name=ef}, bend left]
        \arrow[r, "F \sle"'{name=fe}, bend right] 
        \arrow[Rightarrow, from=ef, to=fe, "\tilde{F}", shorten <= 2pt, shorten >= 2pt]
    & C' 
        \arrow[d, equal] \\
    C
        \arrow[r, "\sle G"{name=eg}, bend left]
        \arrow[r, "G \sle"'{name=ge}, bend right] 
        \arrow[Rightarrow, from=eg, to=ge, "\tilde{G}", shorten <= 2pt, shorten >= 2pt]
    & C'
    \end{tikzcd}
\]
One verifies that homotopies give an equivalence relation.

\begin{prop}
\label{prop:y-functor-with-e}
    The functor $y$ of \Cref{prop:y-ification-from-dga} lifts to a homotopy invariant faithful functor
    \[
        y\colon \Kom_{\mcA^w_n}(\mcC) \to y\Kom_\sle(\mcC). 
    \]
\end{prop}

\subsection
[Reduction of central A-modules]
{Reduction of central $\mcA$-modules}
\label{subsec:reduce-dga}

Next, we consider an algebraic analogue of the braid closing operator $\beta \mapsto \bar{\beta}$. We follow arguments given in \cite[Section 5.1]{GHM:symmety-kr2024}. 

\begin{defn}
\label{def:central}
    An $R_n$-$R_n$ bimodule $X$ is \textit{central} if the left and right actions of $x_i \in R_n$ are equal.
\end{defn}

\begin{defn}
\label{def:central-dga}
    Let $w \in \mfS_n$ be a permutation. A dga $\mcCA^w_n$ over $R_n$ is defined as follows. Let $\mcCA^w_n$ be the graded-commutative $R_n$-algebra freely generated by elements $\xi_i$ $(i = 1, \dots, n)$ of degree $-1$ and $u$ of degree $-2$, modulo the ideal generated by the squares of $\xi_i$. The differential $d$ is defined by 
    \[
        d(\xi_i) = 0,\quad 
        d(u) = \sum_{i = 1}^n (x_i + x_{w(i)}) \xi_i.
    \]
    $\mcCA^w_n$ is also endowed with the quantum grading as in \Cref{def:dga-A}. 
\end{defn}

The dga $\mcCA^w_n$ can alternatively be defined as the quotient $\mcA^w_n / I \mcA^w_n$, where $I \subset R^e_n$ is the ideal generated by $x_i - x'_i$. A central $\mcA^w_n$-module is nothing but a $\mcCA^w_n$-module, and the category of $\mcCA^w_n$-modules $\Kom_{\mcCA^w_n}(\mcC)$ may be regarded as a full subcategory of $\Kom_{\mcA^w_n}(\mcC)$. We write $\mcCA_n = \mcCA^{\id}_n$ when $w = \id_n$. 
We shall prove that any central $\mcA^w_n$-module can be reduced to a smaller $\mcCA_l$-module. 

For a permutation $w \in \mfS_n$, consider a cycle decomposition
\begin{equation}
\label{eqn:w-cycle-decomp}
    w = (i_{1, 1} \cdots i_{1, n_1}) \cdots (i_{l, 1} \cdots i_{l, n_l}).
\end{equation}
Here, $l$ is the number of cycles and $n_k$ is the length of the $k$-th cycle, so that $\sum_k n_k = n$, $w(i_{k, j}) = i_{k, j + 1}$ $(j \in \ZZ/n_k)$ and $\{i_{k, j}\}_{k, j} = \{1, \ldots, n\}$. Having fixed such cycle decomposition of $w$, for a $\mcA^w_n$-module $C$, we write the endomorphisms $x_{k, j} := x_{i_{k, j}}$ and $\xi_{k, j} := \xi_{i_{k, j}}$.

\begin{prop}
\label{prop:reduce-central-dga-action}
    Consider the cycle decomposition \eqref{eqn:w-cycle-decomp} of $w$. With the ring inclusion
    \[
        i\colon R_l \hookrightarrow R_n,\quad 
        x_k \mapsto x_{k, 1}\ 
        (k = 1, \ldots, l),
    \]
    there is a dga-morphism over $R_l$
    \[
        \bar{i}: \mcCA_l \to i^* \mcCA^w_n
    \]
    given by mapping 
    \[
        \xi_k \mapsto \sum_{j = 1}^{n_k} \xi_{k, j}\ (1 \leq k \leq l)
    \]
    and 
    \[
        u \mapsto u + \sum_{k = 1}^{l} \tilde{u}_k,\quad 
        \tilde{u}_k := \sum_{1 \leq j < j' \leq n_k} \xi_{k, j} \xi_{k, j'}.
    \]
\end{prop}

\begin{proof}
    We write $\bar{x}_k = i(x_k)$, $\bar{\xi}_k = \bar{i}(\xi_k)$ and $\bar{u} = \bar{i}(u)$. It is easy to see that 
    \[
        \bar{\xi}_k^{\ 2} = 0,\quad 
        [\bar{\xi}_k, \bar{\xi}_{k'}] = 0,\quad 
        [\bar{\xi}_k, \bar{u}] = 0
    \]
    so $\bar{i}$ is well-defined as a graded-algebra map. It remains to prove that $\bar{x}_k, \bar{u}$ satisfy the conditions of \Cref{def:central-dga}. First, 
    \[
        d(\bar{\xi}_k) = \sum_{j = 1}^{n_k} d(\xi_{k, j}) = \sum_{j = 1}^{n_k} (x_{k, j} - x_{k, j + 1}) = 0. 
    \]
    Next, for each $k$, 
    \begin{align*}
        d(\tilde{u}_k)
            &= \sum_{j < j'} d(\xi_{k, j}) \xi_{k, j'} 
             - \sum_{j < j'} \xi_{k, j} d(\xi_{k, j'}) \\ 
            &= \sum_{j < j'} (x_{k, j} - x_{k, j + 1}) \xi_{k, j'} 
             - \sum_{j < j'} \xi_{k, j} (x_{k, j'} - x_{k, j' + 1}) \\
            &= \sum_{j = 1}^{n_k} (2 x_{k, 1} - (x_{k, j} + x_{k, j + 1})) \xi_{k, j}.
    \end{align*}
    Thus,
    \[
        d(\bar{u}) = d(u) + \sum_k d(\tilde{u}_k) = \sum_{k, j} (2 x_{k, 1}) \xi_{k, j} = 2\sum_{k} \bar{x}_k \bar{\xi}_k.
        \qedhere
    \]
\end{proof}

\begin{prop}
\label{prop:reduce-functor} 
    The dga-morphism of \Cref{prop:reduce-central-dga-action} gives rise to a functor
    \[
        \bar{i}^*\colon 
        \Kom_{\mcCA^w_n}(\mcC) \to \Kom_{\mcCA_{l}}(i^* \mcC).
    \]
    Here, $i^* \mcC$ denotes the $R_l$-linear category obtained by pulling-back $\mcC$ by the ring homomorphism $i$. 
\end{prop}

For a central $\mcA^w_n$-module $C$, the $\mcCA_l$-module $\bar{C} := \bar{i}^*C$ is called the \textit{reduction of $C$} with respect to the cycle decomposition \eqref{eqn:w-cycle-decomp} of $w$. 

\begin{prop}
\label{prop:reduced-y-ification}
    Let $C$ be a central $\mcA^w_n$-module, and $\bar{C}$ be the reduction of $C$ by a cycle decomposition of $w$. The associated $y$-ification $y(\bar{C}) = (C[y_1, \ldots, y_l], \id_l, D)$ has differential
    \[
        D = d + \sum_{k = 1}^l \bar{\xi}_k y_k
    \]
    and the associated $\sle$-operator
    \[
        \sle = \bar{u} + 2 \sum_{k = 1}^l \bar{x}_k \ddel{}{y_k}.
    \]
\end{prop}

Now, suppose we choose another cycle decomposition of $w$ as 
\[
    w = (i'_{1, 1} \cdots i'_{1, n_1}) \cdots (i'_{l, 1} \cdots i'_{l, n_l}),
\]
such that the cycles are indexed the same as \Cref{eqn:w-cycle-decomp}, and the indices in each cycle differ by a cyclic permutation. From \Cref{prop:reduce-central-dga-action}, this gives another ring inclusion
\[
    i'\colon R_l \hookrightarrow R_n
\]
and a dga-morphism
\[
    \bar{i}'\colon \mcCA_l \to i'^* \mcCA^w_n.
\]
We write $\bar{x}'_k = i'(x_k)$, $\bar{\xi}'_k = \bar{i}'(\xi_k)$ and $\bar{u}' = \bar{i}'(u)$. For each $k$, let $t_k$ be the unique index such that
\[
    i'_{k, 1} = i_{k, t_k}.
\]

\begin{lem}
\label{prop:cycle-decomp-difference}
    \[
        \bar{x}_k - \bar{x}'_k = d(\xi_{k, 1} + \cdots + \xi_{k, t_k - 1}), \quad
        \bar{\xi}_k = \bar{\xi}'_k,\quad
        \bar{u} - \bar{u}' = \sum_k (\tilde{u}_k - \tilde{u}'_k)
    \]
    where 
    \[
        \tilde{u}_k - \tilde{u}'_k = 2(\xi_{k, 1} + \cdots + \xi_{k, t_k - 1})(\xi_{k, t_k} + \cdots + \xi_{k, n_k}).
    \]
\end{lem}

\begin{proof}
    The first two equations are easy. For the third, 
    \begin{align*}
        \tilde{u}_k - \tilde{u}'_k 
            &= \sum_{1 \leq j < j' \leq n_k} \xi_{k, j} \xi_{k, j'} 
             - \sum_{1 \leq j < j' \leq n_k} \xi_{k, t_k + j} \xi_{k, t_k + j'} \\ 
            &= \sum_{\substack{1 \leq j < t_k \\ t_k \leq j' \leq n_k}} \xi_{k, j} \xi_{k, j'} 
             - \sum_{\substack{t_k \leq j \leq n_k \\ 1 \leq j' < t_k}} \xi_{k, j} \xi_{k, j'} \\ 
            &= 2\sum_{\substack{1 \leq j < t_k \\ t_k \leq j' \leq n_k}} \xi_{k, j} \xi_{k, j'}.
        \qedhere
    \end{align*}
\end{proof}

\begin{prop}
\label{prop:reduced-y-ification-difference} 
    Let $C$ be a central $\mcA^w_n$-module, and $\bar{C}$, $\bar{C'}$ be reductions of $C$ obtained from two cycle decompositions of $w$. Then the two associated $y$-ifications $y(\bar{C})$ and $y(\bar{C}')$ have identical complexes, homotopic $R_l$-actions, and homotopic $\sle$-operators, regarded over $R$. 
\end{prop}

\begin{proof}
    From \Cref{prop:cycle-decomp-difference}, it follows that $y(\bar{C})$ and $y(\bar{C}')$ have identical complexes and homotopic $R_l$-actions. The difference of the two $\sle$-operators are given by 
    \[
        \sle - \sle' = \sum_{k = 1}^l \left ( (\tilde{u}_k - \tilde{u}'_k) + 2(\bar{x}_k - \bar{x}'_k) \ddel{}{y_k} \right).
    \]
    Consider a degree $-3$ homotopy 
    \[
        h_k := 2(\xi_{k, 1} + \cdots + \xi_{k, t_k - 1}) \ddel{}{y_k}.
    \]
    We have 
    \[
        [d, h_k] = 2(\bar{x}_k - \bar{x}'_k) \ddel{}{y_k}
    \]
    and 
    \begin{align*}
        [\bar{\xi}_k y_k, h_k] 
            &= 2(\xi_{k, 1} + \cdots + \xi_{k, t_k - 1}) \bar{\xi}_k \\ 
            &= 2(\xi_{k, 1} + \cdots + \xi_{k, t_k - 1}) (\xi_{k, t_k} + \cdots + \xi_{k, n_k}) \\
            &= \tilde{u}_k - \tilde{u}'_k.
    \end{align*}
    Thus, $h = \sum_k h_k$ gives $\sle \htpy \sle'$. 
\end{proof}

\Cref{prop:reduced-y-ification-difference} can be restated as saying that the two $y$-complexes $(y(\bar{C}), \sle)$ and $(y(\bar{C}'), \sle')$ are $\sle$-equivariantly isomorphic by $(F, \tilde{F}) = (\id, h)$. We summarize,

\begin{prop}
    The composition of the functors of \Cref{prop:reduce-functor,prop:y-functor-with-e}
    \[
    \begin{tikzcd}
        {\Kom_{\mcCA^w_n}(\mcC)} 
            \arrow[r, "\text{red.}"] 
        & {\Kom_{\mcCA_{l}}(i^* \mcC)}
            \arrow[r, "y"] 
        & {y\Kom_\sle(\bar{\mcC}),}
    \end{tikzcd}
    \]
    is independent, up to natural isomorphism, of the choice of the cycle decomposition of $w$. Here, $\bar{\mcC}$ denotes the $R$-linear category obtained from $\mcC$ by forgetting the polynomial action. 
\end{prop}

\subsection
[Homotopy coherent A-morphisms]
{Homotopy coherent $\mcA$-morphisms}
\label{subsec:hA-morphisms}

Suppose we are given $\mcA^w_n$-modules $C, C'$ and an $R^e_n$-linear chain map $f\colon C \to C'$, which is not necessarily $\mcA^w_n$-linear. We consider the problem of extending $f$ to a $y$-chain map between the two associated $y$-complexes. Here, we follow the idea given in \cite[Section 2.10]{GH:y-ification2022}. 

Suppose $f$ extends to a $y$-chain map as a polynomial 
\[
    F = \sum_\alpha f_\alpha y^\alpha \colon y(C) \to y(C') 
\]
where each coefficient $f_\alpha$ of $y^\alpha$ is an array of bidegree $(-2|\alpha|, 2|\alpha|)$ morphisms
\[
    f_\alpha^i\colon C^i \to C'^{i - 2|\alpha|}
\]
in $\mcC$. For $F$ to be an extension of $f$, we must have $f_0 = f$. For $F$ to be a $y$-chain map, we must have 
\begin{align*}
    [D, F] 
        &= \sum_\alpha \left( [d, f_\alpha] y^\alpha + \sum_{i = 1}^n [\xi_i, f_\alpha] y^{\alpha + e_i} \right) \\
        &= \sum_\alpha \left( [d, f_\alpha] + \sum_i [\xi_i, f_{\alpha - e_i}] \right) y^\alpha \\
        & = 0.
\end{align*}
Here, by convention, we set $f_\alpha = 0$ if any component of $\alpha$ is negative. The above equation requires that, for each $\alpha$, 
\begin{equation}
\label{eqn:f-extension}
    [d, f_\alpha] + \sum_i [\xi_i, f_{\alpha - e_i}] = 0.
\end{equation}
In particular, for $\alpha = 0$, the equation \eqref{eqn:f-extension} specializes to $[d, f] = 0$. For $\alpha = e_i$, we have 
\begin{equation}
\label{eqn:f-extension-lowest1}
    [d, f_i] + [\xi_i, f] = 0
\end{equation}
where $f_i := f_{e_i}$. Thus, $f$ is required to commute with $\xi_i$ up to homotopy $f_i$. 

If we further require $F$ to be $\sle$-equivariant by a homotopy
\[
    \tilde{F} = \sum_\alpha \tilde{f}_\alpha y^\alpha
\]
giving
\[
    [D, \tilde{F}] + [\sle, F] = 0,
\]
then from 
\begin{align*}
    [\sle, F] 
        &= [u + \sum_i (x_i + x'_{w(i)}) \ddel{}{y_i}, \sum_\alpha f_\alpha y^\alpha] \\ 
        &= \sum_\alpha \left( [u, f_\alpha] y^\alpha + \sum_i \alpha_i (x_i + x'_{w(i)}) f_\alpha y^{\alpha - e_i} \right) \\ 
        &= \sum_\alpha \left( [u, f_\alpha] + \sum_i (\alpha_i + 1) (x_i + x'_{w(i)}) f_{\alpha + e_i} \right) y^\alpha
\end{align*}
we must have, for each $\alpha$,
\begin{equation}
\label{eqn:f-extension2}
    [d, \tilde{f}_\alpha] + \sum_i [\xi_i, \tilde{f}_{\alpha - e_i}] + [u, f_\alpha] + \sum_i (\alpha_i + 1) (x_i + x'_{w(i)}) f_{\alpha + e_i} = 0.
\end{equation}
In particular, for $\alpha = 0$, \eqref{eqn:f-extension2} specializes to
\begin{equation}
\label{eqn:f-extension-lowest2}
    [d, f_u] + [u, f] + \sum_i (x_i + x'_{w(i)}) f_i = 0
\end{equation}
where $f_u := \tilde{f}_0$. The homotopy $f_u$ exhibits the up-to-homotopy commutativity of $u$ and $f$ together with terms coming from $f_i$. This observation leads us to the following definition. 

\begin{defn}
\label{def:homotopy-A-linear}
    Let $C, C'$ be $\mcA^w_n$-modules. A \textit{homotopy coherent $\mcA^w_n$-morphism} (or \textit{$h\mcA^w_n$-morphism}) $F = (\{f_\alpha\}, \{\tilde{f}_\alpha\})$ is a pair of collections of maps indexed by non-negative multi-exponents
    \[
        f_\alpha, \tilde{f}_\alpha \colon C \to C'
    \]
    having bidegrees
    \[
        \deg(f_\alpha) = (-2|\alpha|, 2|\alpha|),\quad 
        \deg(\tilde{f}_\alpha) = (-2|\alpha| - 3, 2|\alpha| + 4),
    \]
    and satisfying the equations \eqref{eqn:f-extension}, \eqref{eqn:f-extension2}. Hereafter, we write $F = (f_\alpha, \tilde{f}_\alpha)$ to simplify the notation. 
\end{defn}

The composition of two $h\mcA^w_n$-morphisms is defined so that it is compatible with the composition of the corresponding $\sle$-equivariant $y$-chain maps. To be precise, for two $\mcA^w_n$-modules $C, C'$, let $\Hom_{h\mcA^w_n}(C, C')$ denote the set of $h\mcA^w_n$-morphisms from $C$ to $C'$. For $F = (f_\alpha, \tilde{f}_\alpha) \in \Hom_{h\mcA^w_n}(C, C')$, an $\sle$-equivariant $y$-chain map $y(F) = (y(f_\alpha), y(\tilde{f}_\alpha))$ is defined by
\[
    y(f_\alpha) = \sum_\alpha f_\alpha y^\alpha \colon y(C) \to y(C')
\]
and
\[
    y(\tilde{f}_\alpha) = \sum_\alpha \tilde{f}_\alpha y^\alpha \colon y(C) \to y(C').
\]
This gives an injection 
\[
    \Hom_{h\mcA^w_n}(C, C') \hookrightarrow \Hom_{y\Kom_\sle(\mcC)}(y(C), y(C')).
\]
Now, for $h\mcA^w_n$-morphisms $F = (f_\alpha, \tilde{f}_\alpha)\colon C \to C'$ and $G = (g_\alpha, \tilde{g}_\alpha)\colon C' \to C''$, the composition $GF = (k_\alpha, \tilde{k}_\alpha)$ is defined by
\[
    k_\alpha := \sum_{\beta + \gamma = \alpha} g_\beta f_\gamma,\quad 
    \tilde{k}_\alpha = \sum_{\beta + \gamma = \alpha} \tilde{g}_\beta f_\gamma + g_\beta \tilde{f}_\gamma.
\]
Compared with \Cref{def:e-equiv-y-complexes}, we see that $y(GF) = y(G)y(F)$. Thus, the above embedding shows that the composition satisfies the unit and the associativity laws, and we obtain the category $\Kom_{h\mcA^w_n}(\mcC)$ with objects $\mcA^w_n$-modules and morphisms $h\mcA^w_n$-morphisms. 

Similarly, a \textit{homotopy} $H = (h_\alpha, \tilde{h}_\alpha)$ between two $h\mcA^w_n$-morphisms $F = (f_\alpha, \tilde{f}_\alpha),\ G = (g_\alpha, \tilde{g}_\alpha)$ is defined so that $y(H) = (y(h_\alpha), y(\tilde{h}_\alpha))$ gives a homotopy between $y(F)$ and $y(G)$, in the sense of \Cref{def:e-equiv-homotopy}. In particular, it must satisfy 
\begin{equation}
 \label{eqn:homotopy-hi}
    f - g = [d, h], \quad
    f_i - g_i = [d, h_i] + [\xi_i, h]
\end{equation}
and
\begin{equation}
\label{eqn:homotopy-hu}
    f_u - g_u = [d, h_u] + [u, h] + \sum_i (x_i + x'_{w(i)}) h_i
\end{equation}
where $h := h_0$, $h_i := h_{e_i}$ and $h_u = \tilde{h}_0$. 

\begin{prop}
\label{prop:homotopy-y-ify-functor}
    The functor $y$ of \Cref{prop:y-functor-with-e} extends to a homotopy invariant faithful functor
    \[
        y\colon \Kom_{h\mcA^w_n}(\mcC) \to y\Kom_\sle(\mcC).
    \]
\end{prop}

\begin{defn}
\label{def:truncated-weak-A-linear}
    A $h\mcA^w_n$-morphism $f = (f_\alpha, \tilde{f}_\alpha)$ is \textit{$k$-truncated} if $f_\alpha = 0$ for $|\alpha| > k$ and $\tilde{f}_\alpha = 0$ for $|\alpha| \geq k$. It is \textit{strongly $k$-truncated} if in addition $f_\alpha$ is $\mcA^w_n$-linear for $|\alpha| = k$ and $\tilde{f}_\alpha$ is $\mcA^w_n$-linear for $|\alpha| = k - 1$. In particular, we say $F$ is \textit{strict} if $F$ is (strongly) $0$-truncated. 
\end{defn}

The composition of a $k$-truncated morphism and an $l$-truncated morphism is a $(k + l)$-truncated morphism. Note that strongness is not preserved by composition unless one of the maps is strict. When a $h\mcA^w_n$-morphism $F$ is strict, we identify it with the chain map $F = f_0$. When $F$ is $1$-truncated, we write it as a triple $F = (f, f_i, f_u)$ satisfying \eqref{eqn:f-extension-lowest1}, \eqref{eqn:f-extension-lowest2}. 

Similarly, a \textit{strict homotopy} between strict $h\mcA^w_n$-morphisms is nothing but a homotopy between the corresponding $\mcA^w_n$-morphisms, and \textit{a strongly $1$-truncated homotopy} $H = (h, h_i, h_u)$ between strongly $1$-truncated $\mcA^w_n$-morphisms $f = (f, f_i, f_u)$, $g = (g, g_i, g_u)$ is such that satisfies \eqref{eqn:homotopy-hi}, \eqref{eqn:homotopy-hu} and $h_i, h_u$ are $\mcA^w_n$-linear. 

\subsection
[Transferring A-actions]
{Transferring $\mcA$-actions}
\label{subsec:pushforwad-by-strong-def-retract}

Suppose we are given an $\mcA^w_n$-module $C$ and a homotopy equivalence $C \to C'$ over $R^e_n$. We consider the problem of \textit{transferring} the $\mcA^w_n$-module structure to $C'$, so that the equivalence extends over $\mcA^w_n$ or $h\mcA^w_n$. Here, we focus on the particular case where $C'$ is a \textit{strong deformation retract} of $C$. 

\begin{defn}
\label{def:strong-def-retr}
    Let $C, C'$ be chain complexes in any additive category. We say $C'$ is a \textit{strong deformation retract} of $C$ if there is a triple $(g, f, h)$ of a chain map $g\colon C \to C'$ (the \textit{strong deformation retraction}), a chain map $f\colon C' \rightarrow C$ (the \textit{inclusion}) and a homotopy $h$ on $C$ satisfying the following conditions: (i) $fg - 1 = [d, h]$, (ii) $gf = 1$, (iii) $gh = 0$, (iv) $hf = 0$, (v) $h^2 = 0$.
\end{defn}

We shall depict the maps of \Cref{def:strong-def-retr} as
\[
    \begin{tikzcd}
    C 
        \arrow[r, "g", shift left] 
        \arrow["h"', loop, distance=2em, in=305, out=235] & 
    C'.
        \arrow[l, "f", shift left]
    \end{tikzcd}        
\]
Symbols $f, g, h$ are chosen so that they match with \cite[Definition 4.3]{BarNatan:2005}. The following argument will be used repeatedly in \Cref{subsec:y-ified-khovanov-braids,subsec:y-ified-khovanov-links} to prove the invariance of $y$-ified Khovanov homology.

\begin{prop}
\label{prop:strong-defr-induced-dga}
    Let $C$ be an $\mcA^w_n$-module and $C'$ a strong deformation retract of $C$ over $R^e_n$. Define endomorphisms of $C'$ by 
    \[
        \xi'_i := g \xi_i f,\quad 
        u' := g u f.
    \]
    Then, 
    \begin{align*}
        (\xi'_i)^2 &= -[d, g (\xi_i h \xi_i) f], \\
        [\xi'_i, \xi'_j] &= -[d, g (\xi_i h \xi_j + \xi_j h \xi_i ) f], \\
        [\xi'_i, u'] &= -[d, g (\xi_i h u + u h \xi_i) f] \\
        &\quad - \sum_j (x_j + x'_{w(j)}) g (\xi_i h \xi_j + \xi_j h \xi_i) f
    \end{align*}
    and 
    \[
        [d, \xi'_i] = 0,\quad [d, u'] = \sum_i (x_i + x'_{w(i)}) \xi'_i.
    \]
\end{prop}

\begin{proof}
    We have
    \begin{align*}
        (\xi'_i)^2
            &= (g \xi_i f) (g \xi_i f) \\
            &= g \xi_i (1 + [d, h]) \xi_i f \\ 
            &= g \xi_i [d, h] \xi_i f \\ 
            &= g [d, -\xi_i h \xi_i ] f \\ 
            &= -[d, g \xi_i h \xi_i f ].
    \end{align*}
    In the fourth equality, we used 
    \[
        [d, \xi_i h \xi_i ] = [d, \xi_i] h \xi_i - \xi_i [d, h] \xi_i + \xi_i h [d, \xi_i]
    \]
    together with $[d, \xi_i] \in R^e_n$ and $gh = hf = 0$. The other equations are proved similarly.
\end{proof}

\Cref{prop:strong-defr-induced-dga} implies that the endomorphisms $\xi'_i, u'$ do not necessarily commute strictly. In case they do give a strict $\mcA^w_n$-action on $C'$, we call it the \textit{$\mcA^w_n$-action transferred by} the strong deformation retraction $(g, f, h)$. 

\begin{cor}
    In addition to the assumption of \Cref{prop:strong-defr-induced-dga}, if $[\xi_i, h] = 0$, then 
    \[
        (\xi'_i)^2 = [\xi'_i, \xi'_j] = [\xi'_i, u'] = 0
    \]
    and the $\mcA^w_n$-action on $C$ can be transferred to $C'$. 
\end{cor}

\begin{prop}
\label{prop:strong-defr-induced-dga2}
    Under the assumption of \Cref{prop:strong-defr-induced-dga}, define maps
    \begin{gather*}
        f_i := h \xi_i f,\quad 
        g_i := g \xi_i h,\quad 
        h_i := h \xi_i h,\\
        f_u := h u f,\quad 
        g_u := -g u h,\quad 
        h_u := h u h.
    \end{gather*}
    Then, with the endomorphisms $\xi'_i, u'$ defined in \Cref{prop:strong-defr-induced-dga}, we have 
    \begin{align*}
        [d, f_i] + [\xi_i, f] &= 0,\\
        [d, g_i] + [\xi_i, g] &= 0,\\
        [d, h_i] + [\xi_i, h] &= f_i g + f g_i
    \end{align*}
    and
    \begin{align*}
        [d, f_u] + [u, f] + \sum_i(x_i + x'_{w(i)})f_i &= 0,\\
        [d, g_u] + [u, g] + \sum_i(x_i + x'_{w(i)})g_i &= 0,\\
        [d, h_u] + [u, h] + \sum_i(x_i + x'_{w(i)})h_i &= f_u g + f g_u.
    \end{align*}
\end{prop}

\begin{proof}
    Similar to the proof of \Cref{prop:strong-defr-induced-dga}. 
\end{proof}

\begin{cor}
\label{cor:strong-defr-linear}
    In addition to the assumption of \Cref{prop:strong-defr-induced-dga2}, if $[\xi_i, h] = 0$, then
    \[
        f_i = g_i = h_i = 0.
    \]
    Furthermore, if $[u, h] = 0$, then
    \[
        f_u = g_u = h_u = 0.
    \]
    If both $[\xi_i, h] = [u, h] = 0$, then the $\mcA^w_n$-action on $C$ can be transferred to $C'$, and the maps $f, g, h$ become $\mcA^w_n$-linear. Therefore, $C'$ is a strong deformation retract of $C$ in the category $\Kom_{\mcA^w_n}(\mcC)$. 
\end{cor}

\Cref{cor:strong-defr-linear} shows that the extension problem is primarily governed by $[\xi_i, h]$ and $[u, h]$. Even when these commutators are non-zero, there are cases where the $\mcA^w_n$-action can be transferred and the maps extend as $h\mcA^w_n$-morphisms. First, note that the equations of $f_i, f_u$ (and of $g_i, g_u$) in \Cref{prop:strong-defr-induced-dga2} correspond to the equations \eqref{eqn:f-extension-lowest1}, \eqref{eqn:f-extension-lowest2}, and those of $h_i, h_u$ correspond to \eqref{eqn:homotopy-hi}, \eqref{eqn:homotopy-hu}. In addition to the assumption of \Cref{prop:strong-defr-induced-dga2}, suppose that the $\mcA^w_n$-action on $C$ can be transferred to $C'$ and $g_i, g_u$ are $\mcA^w_n$-linear. Then $G = (g, g_i, g_u)$ gives a strongly $1$-truncated $h\mcA^w_n$-morphism $C \to C'$. Furthermore, if $f$, $h_i$ and $h_u$ are $\mcA^w_n$-linear, then the equations of  \Cref{prop:strong-defr-induced-dga2} show that $H = (h, h_i, h_u)$ gives a strongly $1$-truncated homotopy $fG \htpy_H I$. For the other direction, we have $Gf = I$, and thus $C'$ is a strong deformation retract of $C$ in $\Kom_{h\mcA^w_n}(\mcC)$. Passing to $y$-ifications, $y(C')$ is an $\sle$-equivariant strong deformation retract of $y(C)$. We shall encounter such a situation in the proof of \Cref{prop:y-cpx-braid-rel3}. 

%% file: 3.tex
\section
[y-ification of Khovanov homology]
{$y$-ification of Khovanov homology}
\label{sec:y-ified-kh}

\subsection{Preliminaries}
\label{subsec:preliminary}

We assume that the reader is familiar with the construction of Khovanov homology and its equivariant versions for links and tangles~\cite{Khovanov:2000,Lee:2005,BarNatan:2005,Khovanov:2004,Khovanov:2022}. Here, we briefly review the basic concepts, and fix the notation that will be used throughout this paper. 

\subsubsection{The category of dotted cobordisms}
\label{subsubsec:cob-dot}

Let $B$ be a finite set of signed points on $\partial D^2$ with equal numbers of positive and negative points. Let $\Cob(B)$ denote the \textit{category of dotted cobordisms} in $D^2 \times I$ with boundary points $B$. Namely, objects of $\Cob(B)$ are (unoriented) planar tangle diagrams $T$ in $D^2$ with $\del T = B$ (ignoring the signs on $B$), and morphisms between objects $T, T'$ are formal $\ZZ$-linear combinations of (unoriented) dotted cobordisms $S$ in $D^2 \times I$ with $\del S = T \times \{0\} \cup B \times I \cup T' \times \{1\}$, regarded up to boundary fixing isotopy, modulo the three \textit{local relations}:
\begin{center}
    \input{tikzpictures/loc-relation-dot}
\end{center}
$\Cob(B)$ is endowed with the \textit{quantum grading}\footnote{
    The quantum grading is negated compared to \cite[Definition 6.2]{BarNatan:2005}, due to the change of convention. 
}, where a cobordism $C \in \Hom_{\Cob(B)}(X, Y)$ with $k$ is declared to have 
\[
    \deg_q(C) := \frac{|B|}{2} - \chi(C) + 2k.
\]
Here, $\chi(C)$ denotes the Euler number of $C$. Note that the identity cobordism $C = X \times I$ has $\deg_q(C) = 0$, and that $\deg_q$ is additive under vertical and horizontal compositions (see \cite[Exercise 6.3]{BarNatan:2005}). For an object $X$, let $q^mX$ denote the copy of $X$ with its quantum grading shifted by $m$. 

The \textit{tautological functor}
\[
    \mcF = \Hom_{\Cob(\emptyset)}(\emptyset, -)\colon
    \Cob(\emptyset) \longrightarrow \Mod{R_{U(2)}}
\]
recovers the Frobenius extension $(R_{U(2)}, A_{U(2)})$ for the \textit{$U(2)$-equivariant Khovanov homology}. Here, $R_{U(2)}$ is the graded ring $\ZZ[h, t]$ with degrees $\deg h = 2$, $\deg t = 4$, and $A_{U(2)}$ is the Frobenius algebra $R_{U(2)}[X]/(X^2 - hX - t)$ determined by the counit $\epsilon(1) = 0,\ \epsilon(X) = 1$ and graded with $\deg X = 2$. Indeed, there is a graded ring isomorphism
\[
    \mcF(\emptyset) = \Hom_{\Cob(\emptyset)}(\emptyset, \emptyset) \isom R_{U(2)}
\]
under which $h, t \in R_{U(2)}$ correspond to closed cobordisms
\begin{center}
    \input{tikzpictures/cob-ht}
\end{center}
and a graded $R_{U(2)}$-module isomorphism
\[
    \mcF(\bigcirc) = \Hom_{\Cob(\emptyset)}(\emptyset, \bigcirc) \isom q^{-1} A_{U(2)}
\]
under which the two generators $1, X$ of $A_{U(2)}$ as an $R_{U(2)}$-module correspond to cobordisms
\begin{center}
    \input{tikzpictures/1X-cob}
\end{center}
Furthermore, the Frobenius algebra operations of $A_{U(2)}$, multiplication $m$, unit $\iota$, comultiplication $\Delta$, and counit $\epsilon$, correspond to the following cobordisms:
\begin{center}
    \input{tikzpictures/frobenius-ops-cob}
\end{center}
In this sense, we shall refer to $\mcF$ as a \textit{TQFT} for the $U(2)$-equivariant theory. As proved in \cite[Proposition 5]{Khovanov:2004}, the $U(2)$-equivariant theory is universal among all Khovanov-type link homology theories. For instance, the following equivariant theories are by obtained specializing the variables $(h, t)$ to
\begin{align*}
    t = 0 &\rightsquigarrow 
        & (R_{U(1)}, A_{U(1)}) &= (\ZZ[h],\ \ZZ[h, X]/(X^2 - hX)),\\
    h = 0 &\rightsquigarrow 
        & (R_{SU(2)}, A_{SU(2)}) &= (\ZZ[t],\ \ZZ[t, X]/(X^2 - t)),\\
    h = t = 0 &\rightsquigarrow 
        & (R_{\{e\}}, A_{\{e\}}) &= (\ZZ,\ \ZZ[X]/(X^2)).    
\end{align*}
Alternatively, deformed versions are obtained by additionally imposing relations in $\Cob(B)$, and then applying the tautological functor $\mcF$. The non-equivariant Frobenius extension $(R_{\{e\}}, A_{\{e\}})$ corresponds to Khovanov's original theory given in \cite{Khovanov:2000}, and the corresponding local relations are given in \cite[Section 11.2]{BarNatan:2005}. For an arbitrary commutative ring $R$ and elements $h, t \in R$, we write $A_{h, t} := R[X]/(X^2 - hX - t)$, and the corresponding TQFT by $\mcF_{h, t}$. In particular, we write $\mcF_0 := \mcF_{0, 0}$.

\subsubsection{Khovanov homology for tangles}
\label{subsubsec:preliminary-khovanov}

Let $\Diag(B)$ denote the \textit{category of oriented tangle diagrams} with boundary $B$. Namely, objects of $\Diag(B)$ are oriented tangle diagrams $T$ with $\del T = B$ (respecting the signs of $B$), and morphisms between diagrams $T, T'$ are oriented movies between $T, T'$ fixing the boundary $B$. In particular, when $B = \emptyset$, the objects of $\Diag(\emptyset)$ are oriented link diagrams and morphisms are movies between link diagrams. 

Given an oriented tangle diagram $T \in \Diag(B)$, the \textit{formal Khovanov bracket} $[T]$ is obtained by the standard construction of the Khovanov complex, but in the category $\Mat(\Cob(B))$, i.e.\ the additive closure of the preadditive category $\Cob(B)$. Following \cite{BarNatan:2005}, the category of complexes in $\Mat(\Cob(B))$ will be denoted $\Kob(B) := \Kom(\Mat(\Cob(B)))$. \cite[Theorem 1]{Bar-Natan:2002} states that the chain homotopy type of $[T]$ is an invariant of the tangle. In particular, when $B = \emptyset$, applying $\mcF$ to the formal Khovanov bracket $[T]$ recovers the \textit{$U(2)$-equivariant Khovanov complex} $\CKh_{U(2)}(T)$, and its homology give the \textit{$U(2)$-equivariant Khovanov homology} $\Kh_{U(2)}(T)$. Therefore, the Khovanov homology functor $\Kh_{U(2)}$ decomposes into a sequence of functors, 
\[
\begin{tikzcd}
    \Diag(B) \arrow[r, "{[\cdot]}"] & \Kob(B) \arrow[r, "\mathcal{F}", "\text{when $B = \emptyset$}"', dotted] &[4em] \Kom(\Mod{R_{U(2)}}) \arrow[r, "H"] & {\Mod{R_{U(2)}}}.
\end{tikzcd}
\]
Other variants can be obtained by specializing $(h, t)$ before taking the homology. Therefore, working in  $\Kob(B)$ is the most general approach to Khovanov homology. 

\subsubsection{Signing convention}

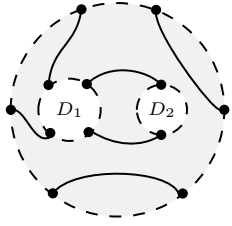
\begin{figure}
    \centering
    \input{tikzpictures/planar-diagram}
    \caption{A planar diagram}
    \label{fig:planar-diagram}
\end{figure}

One of the advantages of Bar-Natan's framework is that these categories are functorial under \textit{horizontal compositions} by \textit{planar arc diagrams}. This allows local treatments of knots and links via tangle pieces (\Cref{fig:planar-diagram}). Here, we fix the signing convention for the horizontal composition in the category $\Kob(B)$. Suppose $\Omega_i$ are objects of $\Kob(B_i)$ $(i = 1, 2)$ and are composable by a planar arc diagram $D$. Given homogeneous maps (not necessarily chain maps) $f\colon \Omega_1 \to \Omega_1$ and $g\colon \Omega_2 \to \Omega_2$ of degrees $|f|$ and $|g|$ respectively, we define maps on $\Omega = D(\Omega_1, \Omega_2)$ by 
\[
    D(f, \id) := \bigoplus_{i, j} \left( D(f^i, \id)\colon D(\Omega^i_1, \Omega^j_2) \to D(\Omega^{i + a}_1, \Omega^j_2) \right)
\]
and 
\[
    D(\id, g) := \bigoplus_{i, j} (-1)^{i|g|} \left( D(\id, g^j)\colon 
    D(\Omega^i_1, \Omega^j_2) \to D(\Omega^a_1, \Omega^{j + b}_2) \right).
\]
Note the sign $(-1)^{i|g|}$ on the second equation. With this signing convention, we always have 
\[
    [D(f, \id), D(\id, g)] = 0
\]
where $[\cdot, \cdot]$ is the super-commutator. The differential $d$ on $\Omega$ can be written as 
\[
    d = D(d_1, \id) + D(\id, d_2)
\]
as defined in \cite[Section 5]{BarNatan:2005}. Furthermore, if we have
\[
    f - f' = [d_1, h], \quad 
    g - g' = [d_2, k]
\]
then
\[
    D(f, \id) - D(f', \id) = [d, D(h, \id)],\quad
    D(\id, g) - D(\id, g') = [d, D(\id, k)]. 
\]
Therefore, with this signing convention, the super-commutator is compatible with horizontal compositions. Hereafter, we simply write $f, g$ for $D(f, \id), D(\id, g)$ when there is no risk of confusion.

\subsubsection{Braid diagrams}
\label{subsubsec:braid}

Here, we fix the notation regarding braids and their closures. Let $\sfB_n$ denote the set of boundary points on $\del I^2$ with $n$ negative points on $I \times \{0\}$ and $n$ positive points on $I \times \{1\}$. We label the bottom points of $\sfB_n$ by $P_i \in I \times \{0\}$ and the top points by $P'_i  \in I \times \{1\}$ (see \Cref{fig:Bn} in \Cref{sec:intro}).
We let $\BDiag(n)$ denote the subcategory of $\Diag(\sfB_n)$, whose objects are braid diagrams on $n$ strands represented as braid words, and morphisms between braid diagrams are \textit{braid movies}, as introduced in \cite{CS:1996,CS:1998,Khoavnov-Thomas:2007}\footnote{Movies of braid diagrams and tangle diagrams are not be treated in this paper.}. In our convention, a braid diagram $\beta \in \BDiag(n)$ is oriented from the bottom to the top. For two points $p, q$ on a braid diagram $\beta$, we write $p \sim q$ if they belong to the same strand. A braid diagram $\beta \in \BDiag(n)$ has the \textit{underlying permutation} $w \in \mfS_n$ such that $P_i \sim P'_{w(i)}$ for all $i$. The product $\beta_2 \beta_1$ of two braid diagrams $\beta_1, \beta_2 \in \BDiag(n)$ is defined by the braid diagram obtained by stacking $\beta_2$ on top of $\beta_1$. 
We let $\id_n$ denote the object represented by the empty word, $\sigma_i \in \BDiag(n)$ the object with a single positive crossing between the $i$-th and $(i+1)$-st strands, and $\sigma_i^{-1}$ its mirror image. Note that $\sigma_i \sigma_i^{-1}$ and $\sigma_i^{-1} \sigma_i$ are isomorphic to $\id_n$, but are not regarded as equal. 

A braid diagram $\beta \in \BDiag(n)$ can be closed to form a link diagram $\bbeta \in \Diag(\emptyset)$. This operation can be regarded as a functor 
\[
    \cl_n\colon \BDiag(n) \longrightarrow \Diag(\emptyset)
\]
given by the $1$-input planar arc diagram of \Cref{fig:braid-closure}. From the functoriality of planar arc diagrams, we have the following commutative diagram of functors:
\[
\begin{tikzcd}[column sep=4em]
    \BDiag(n) \arrow[r, "{[\cdot]}"] \arrow[d, "\cl_n"] & \Kob(\sfB_n) \arrow[d, "\cl_n"] \\
    \Diag(\emptyset) \arrow[r, "{[\cdot]}"] & \Kob(\emptyset).
\end{tikzcd}
\]

\subsection{Dot-sliding homotopies}
\label{subsec:dot-cross-homotopy}

Let $C$ be an object of $\Cob(B)$. For any point $p$ on $C$, an endomorphism $X_p$ of $C$ is defined by the cobordism made of the identity cobordism $C \times I$ with a dot placed on the component that contains $p$. For a tangle diagram $T$ and a regular point $p$ on $T$ (i.e.\ a point on $T$ which is not a double point), the endomorphism $X_p$ is extended as a chain endomorphism of $[T]$ in the obvious way. 

\begin{prop}
\label{prop:homotopy-crossing}
    Let $T$ be a tangle diagram with two strands and a single crossing $c$. If $p, q$ are two points that lie on the same strand separated by the crossing, then $X_p$ and $h - X_q$ are homotopic in $[T]$ via the map $\chi$ defined by the reversal of the differential $d$. 
\end{prop}

\begin{proof}
    Let $T_0, T_1$ denote the crossing-less diagrams obtained from the $0$- and $1$-resolution of $T$. The maps are given by 
    \[
\begin{tikzcd}
T_0 \arrow[r, "d", shift left] \arrow["{X_p, X_q}"', loop, distance=2em, in=305, out=235] & T_1. \arrow[l, "\chi", dashed, shift left] \arrow["{X_p, X_q}"', loop, distance=2em, in=305, out=235]
\end{tikzcd}
    \]
    It suffices to show that 
    \[
        X_p + X_q - h = \chi d
    \]
    on $T_0$, and 
    \[
        X_p + X_q - h = d \chi
    \]
    on $T_1$. Both $\chi d$ and $d \chi$ are the obvious genus $1$ endo-cobordisms on $T_0$ and $T_1$. Thus, the equations are immediate from the neck cutting relation. 
\end{proof}

The homotopy $\chi_c$ of \Cref{prop:homotopy-crossing} is called a \textit{dot-sliding homotopy} for $c$. It has appeared in \cite[Lemma 2.3]{Hedden-Ni:2013}, \cite[Proposition 6.4]{Batson-Seed:2015}, \cite[Lemma 2.1]{Alishahi:2017}, \cite[Lemma 3.3]{Alishahi:2018} in the context of Khovanov homology. 

\begin{figure}[t]
    \centering
    \begin{subfigure}[t]{.3\textwidth}
        \centering
        \input{tikzpictures/crossing-endpoints}
        \caption{}
        \label{fig:crossing-endpoints}
    \end{subfigure}
    \begin{subfigure}[t]{.3\textwidth}
        \centering
        \input{tikzpictures/crossing-admissible}
        \caption{}
        \label{fig:crossing-admissible}
    \end{subfigure}
    \caption{Endpoints of a crossing $c$ and an admissible coloring.}
\end{figure}
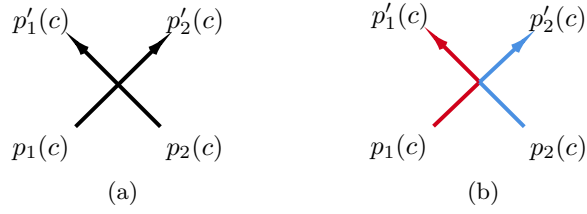

Now, let $T$ be a tangle diagram. For each crossing $c$, let $p_1(c), p_2(c), p'_1(c), p'_2(c)$ denote the four points around $c$ as in \Cref{fig:crossing-endpoints}. An \textit{$ab$-coloring} $\theta$ on $T$ is an assignment of either $a$ or $b$ to each edge of the underlying 4-valent graph $G(T)$ of $T$. For a regular point $p$ on $T$, we write $\theta(p)$ to denote the color of the edge containing $p$. We say $\theta$ is \textit{admissible} if, for every crossing $c$ of $T$,
\[
    \theta(p_1(c)) = \theta(p'_1(c)) \neq \theta(p_2(c)) = \theta(p'_2(c))
\]
(see \Cref{fig:crossing-admissible}). Observe that, if $G(T)$ is a connected graph, then $T$ admits exactly two admissible $ab$-colorings. This follows from the fact that an admissible coloring corresponds bijectively to a checkerboard coloring of the complement of the orientation preserving resolution $T_0$ of $T$ in $D^2$. Indeed, each component of $T_0$ inherits an orientation from $T$, and we may color it $a$ if it sees a black region to its left, or $b$ otherwise. More generally, if $G(T)$ has $c$ connected components, then $T$ admits exactly $2^c$ admissible colorings. 

For two tangle diagrams $T, T'$ that are composable by a planar arc diagram $D$, we say two admissible colorings of $T, T'$ are \textit{compatible} with respect to $D$ if they extend to an admissible coloring of $D(T, T')$. 
For a braid diagram $\beta$, the \textit{standard $ab$-coloring} of $\beta$ is the admissible $ab$-coloring obtained by the alternating coloring of the orientation preserving resolution $\beta_0$, which consists of $n$ vertical arcs, with the leftmost arc colored $a$. Obviously, the standard $ab$-coloring is compatible with respect to multiplication of braids and the closure operation. 

Hereafter, fix an admissible $ab$-coloring $\theta$ for a tangle diagram $T$. Then, for each regular point $p$ on $T$, we define an endomorphism $x_p$ on $[T]$ by
\[
    x_p = \begin{cases}
        X_p & \text{ if $\theta(p) = a$, } \\
        h - X_p & \text{ if $\theta(p) = b$.}
    \end{cases}
\]
Note that we have 
\[
    (x_p)^2 = h x_p + t
\]
regardless of $\theta(p)$. If two regular points $p, q$ are separated by a crossing $c$, then we have $\theta(p) \neq \theta(q)$ from the admissibility of $\theta$. From \Cref{prop:homotopy-crossing}, the corresponding endomorphisms $x_p$ and $x_q$ are homotopic by the homotopy $\chi_c$ given by the reversal of the local differential for $c$, giving
\[
    x_p - x_q = \begin{cases}
        [d, \chi_c] & \text{ if $\theta(p) = a$, } \\
        -[d, \chi_c] & \text{ if $\theta(p) = b$.}
    \end{cases}
\]
In order to unite the two cases, we define a signing function $\epsilon$ by 
\[
    \epsilon_p := \begin{cases}
        +1 & \text{ if $\theta(p) = a$, } \\
        -1 & \text{ if $\theta(p) = b$}
    \end{cases}
\]
and the signed version of $\chi_c$ by 
\[
    \hchi_c := \epsilon_{p_1(c)} \chi_c.
\]
Then we have 
\[
    x_{p_1(c)} - x_{p'_2(c)} 
    = -(x_{p_2(c)} - x_{p'_1(c)})
    = [d, \hchi_c]. 
\]
The following proposition is immediate from the definition. 

\begin{prop}
\label{prop:x-and-chi-commute}
    We have 
    \[
        [x_p, \chi_c] = 0,\quad 
        \chi_c^2 = 0, \quad 
        [\chi_c, \chi_{c'}] = 0.
    \]
    where $p$ is a regular point of $T$, and $c, c'$ are crossings of $T$.
\end{prop}

Let $p, q$ be two regular points on $T$. A \textit{path} $\gamma$ from $p$ to $q$ in $T$, denoted
\begin{equation}
\label{eqn:path}
    \gamma\colon 
    p = p_1 \xrightarrow{c_1} p_2 \xrightarrow{c_2} \cdots \xrightarrow{c_r} p_{r + 1} = q    
\end{equation}
is a sequence of regular points on $\beta$, such that each consecutive points $p_j$, $p_{j + 1}$ are separated by a single crossing $c_j$. For such path $\gamma$, write $\epsilon_{p_j} = \epsilon_j$, $\chi_j = \chi_{c_j}$  and define a degree $-1$ homotopy $\xi_\gamma$ by
\begin{equation}
\label{def:xi-gamma}
    \xi_\gamma := \sum_{j = 1}^r \epsilon_j \chi_j.
\end{equation}

\begin{prop}
\label{prop:xi-commute}
    $\xi_\gamma$ satisfies
    \[
        [d, \xi_\gamma] = x_p - x_q.
    \]
    Furthermore, 
    \[
        [x_p, \xi_\gamma] = 0,\quad 
        \xi_\gamma^2 = 0,\quad 
        [\xi_\gamma, \xi_{\gamma'}] = 0.
    \]
    where $p$ is a regular point on $T$, and  $\gamma, \gamma'$ are two arbitrary paths on $T$. 
\end{prop}

\begin{proof}
    The first equation is immediate from the definition. The second is immediate from \Cref{prop:x-and-chi-commute}. For the third, 
    \begin{align*}
        \xi_\gamma^2 
            &= (\sum_j \epsilon_j \chi_j)^2 \\
            &= \sum_j \chi_j^2 + \sum_{{j} \neq {j'}} \epsilon_{j}\epsilon_{j'} \chi_{j} \chi_{j'} \\
            &= \sum_j \chi_j^2 + \sum_{{j} < {j'}} \epsilon_{j} \epsilon_{j'} [\chi_{j},  \chi_{j'}] \\ 
            &= 0.
    \end{align*}
    The last equation is proved similarly.
\end{proof}

For a path $\gamma$ in the form of \Cref{eqn:path}, we also define a degree $-2$ homotopy $u_\gamma$ by 
\begin{equation}
\label{def:u-gamma}
    u_\gamma := \sum_{1 \leq j < j' \leq r} (\epsilon_j \chi_j) (\epsilon_{j'} \chi_{j'}).
\end{equation}

\begin{lem}
\label{lem:u-gamma}
    $u_\gamma$ satisfies
    \[
        [d, u_\gamma] = (x_p + x_q) \xi_\gamma - \sum_j \epsilon_{j} (x_{p_j} + x_{p_{j + 1}}) \chi_j
    \]
    and commutes with $x_p$ and $\xi_{\gamma'}$. 
\end{lem}

\begin{proof}
    Compute
    \begin{align*}
        [d, u_\gamma] 
        &= \sum_{j < j'} ([d, \epsilon_j \chi_j] \epsilon_{j'} \chi_{j'} - \epsilon_j \chi_j [d, \epsilon_{j'} \chi_{j'}]) \\ 
        &= \sum_{j < j'} ((x_{p_j} - x_{p_{j + 1}}) \epsilon_{j'} \chi_{j'} - \epsilon_j \chi_j (x_{p_{j'}} - x_{p_{j' + 1}})) \\ 
        &= \sum_j ((x_{p_1} - x_{p_j}) \epsilon_{j} \chi_j - \epsilon_j \chi_j (x_{p_{j + 1}} - x_{p_{r + 1}})) \\ 
        &= (x_p + x_q) \xi_\gamma - \sum_j \epsilon_{j} (x_{p_j} + x_{p_{j + 1}}) \chi_j.
        \qedhere
    \end{align*}
\end{proof}

\subsection
[A-action on Khovanov complex]
{$\mcA$-action on Khovanov complex}
\label{subsec:dg-module-on-kh-cpx}

Here, we focus on braid diagrams, and prove that the complex $[\beta]$ admits an $\mcA$-action defined in \Cref{subsec:dga-A}. Take any braid diagram $\beta \in \BDiag(n)$, and let $w$ denote the underlying permutation of $\beta$. Take the standard coloring $\theta$ of $\beta$. Define chain endomorphisms $x_i$ and $x'_i$ ($i = 1, \ldots, n$) of $[\beta]$ by 
\[
    x_i := x_{P_i},\quad 
    x'_i := x_{P'_i}.
\]
These endomorphisms obviously commute with each other, and hence give an $R^e_n$-module structure on $[\beta]$. Let $\gamma_i$ denote the path from $P_i$ to $P'_{w(i)}$, given by 
\[
    \gamma_i\colon 
    P_i = p_{i, 1} \xrightarrow{c_{i, 1}} p_{i, 2} \xrightarrow{c_{i, 2}} \cdots \xrightarrow{c_{i, r}} p_{i, {r + 1}} = P'_{w(i)}.
\]
Define endomorphisms
\[
    \xi_i := \xi_{\gamma_i},\quad 
    u := \sum_i u_{\gamma_i}. 
\]

\begin{prop}
\label{prop:homotopy-u}
    We have 
    \[
        [d, \xi_i] = x_i - x'_{w(i)},\quad 
        [d, u] = \sum_i \left( x_i + x'_{w(i)} \right) \xi_i.
    \]
    Furthermore, both $\xi_i, u$ commute with $x_i, x'_i$, and 
    \[
        \xi_i^2 = 0,\quad 
        [\xi_i, \xi_j] = 0,\quad 
        [\xi_i, u] = 0.
    \]
\end{prop}

\begin{proof}
    The equation for $\xi_i$ follows immediately from \Cref{prop:xi-commute}. For $u$, from \Cref{lem:u-gamma}, the difference of the two sides of the desired equation is given by the sum
    \[
        -\sum_i \sum_{j} \epsilon_{p_{i, j}} (x_{p_{i, j}} + x_{p_{i, {j + 1}}}) \chi_{c_{i, j}}.
    \]
    For each crossing $c$, a term with $\chi_c$ will appear exactly twice, one from a path passing $c$ from the left and another from a path passing from the right. From $\theta(p_1(c)) = \theta(p'_1(c))$ and $\theta(p_2(c)) = \theta(p'_2(c))$, the two terms add to
    \[
        \pm (x_{p_1(c)} + x_{p'_2(c)} - x_{p_2(c)} - x_{p'_1(c)}) \chi_c.
    \]
    Since $\chi_c$ acts as a saddle cobordism, it follows that
    \[
        (x_{p_1(c)} - x_{p'_1(c)}) \chi_c = 
        (x_{p_2(c)} - x_{p'_2(c)}) \chi_c = 0
    \]
    and thus the sum is zero.
\end{proof}

\begin{rem}
\label{rem:u-description}
    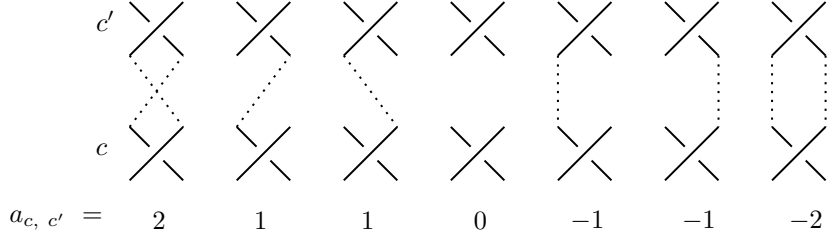
\begin{figure}
        \centering
        \input{tikzpictures/a-coeffs}
        \caption{Coefficients $a_{c, c'}$.}
        \label{fig:a-coeffs}
    \end{figure}
    An analogous description of $u = u_2$ for the Rouquier complex is given in \cite[Theorem 1.6]{GHM:symmety-kr2024}. Also, following \cite[Lemma 3.21]{CG:structure-in-homfly2024}, if $\beta = \beta_2 \beta_1$ and each $\beta_a$ is given the $\mcA^{w_a}_n$-action by endomorphisms $\xi^a_{i}, u^a$, then we may write
    \begin{equation}
    \label{eqn:u-by-braid}
        u = u^1 + u^2 + \sum_i \xi^1_i \xi^2_{w_1(i)}. 
    \end{equation}
    Furthermore, following \cite[Remark 3.22]{CG:structure-in-homfly2024}, we may also write
    \begin{equation}
    \label{eqn:u-by-chi-crossing}
        u = \sum_{c \prec c'} a_{cc'} \hchi_c \hchi_{c'},
    \end{equation}
    where $c \prec c'$ denotes the relation that $c$ precedes $c'$ in the vertical direction, and the coefficients $a_{c c'} \in \{0, \pm 1, \pm 2\}$ are given by
    \[
        a_{c c'} = \begin{cases}
            2 & \text{$p'_1(c) \sim p_2(c')$ and $p'_2(c) \sim p_1(c')$,}  \\
            1 & \text{either $p'_1(c) \sim p_2(c')$ or $p'_2(c) \sim p_1(c')$,}  \\
            0 & \text{none of the other,} \\
            -1 & \text{either $p'_1(c) \sim p_1(c')$ or $p'_2(c) \sim p_2(c')$,}  \\
            -2 & \text{$p'_1(c) \sim p_1(c')$ and $p'_2(c) \sim p_2(c')$.}
        \end{cases}
    \]
    See \Cref{fig:a-coeffs} for a pictorial description of the coefficients $a_{c, c'}$.
\end{rem}

\begin{rem}
    At this point, the reader may refer \Cref{ex:T1,ex:T2-para}, where the operators $\xi_i$ and $u$ are described explicitly for the braids $\sigma_1$ and $\sigma_1^2$. 
\end{rem}

\begin{prop}
\label{prop:Qn-module-str}
    On $[\beta]$, we have 
    \[
        \sum_{i=0}^n (x_i)^k = \sum_{i=0}^n (x'_i)^k
    \]
    for any $k \geq 0$. 
\end{prop}

\begin{proof}
    From \Cref{prop:homotopy-u}, 
    \begin{align*}
        0 &= [d, [d, u]] \\ 
        &= \sum_i (x_i + x'_{w(i)}) (x_i - x'_{w(i)}) \\
        &= \sum_i (x_i)^2 - \sum_i (x'_i)^2
    \end{align*}
    Since $x_p^2 = h x_p + t$ and $h, t$ are non-zero divisors, the above equation implies
    \[
        \sum_i x_i = \sum_i x'_i.
    \]
    Furthermore, since each monomial $x_i^k$ can be reduced to a polynomial of degree at most $1$, the claim follows. 
\end{proof}

Combining \Cref{prop:xi-commute,prop:homotopy-u,prop:Qn-module-str}, we conclude

\begin{prop}
\label{prop:beta-dg-mod-structure}
    For any braid diagram $\beta \in \BDiag(n)$ with underlying permutation $w$, the complex $[\beta]$ admits an $\mcA^w_n$-module structure. 
\end{prop}

\begin{prop}
\label{prop:ab-coloring-dependence}
    If we choose the opposite $ab$-coloring of $\theta$, then the $\mcA^w_n$-action is twisted by
    \begin{align*}
        x_i &\mapsto h - x_i, & 
        x'_i &\mapsto h - x'_i, \\ 
        \xi_i &\mapsto -\xi_i,&
        u &\mapsto u.
    \end{align*}
\end{prop}

We let $\Kob_{\mcA^w_n}(\sfB_n) := \Kom_{\mcA^w_n}(\Mat(\Cob(\sfB_n)))$ denote the category whose objects are bounded complexes in $\Kob(\sfB_n)$ equipped with $\mcA^w_n$-actions, and morphisms are $\mcA^w_n$-linear chain maps. The category $\Kob_{h\mcA^w_n}(\sfB_n)$ is defined similarly. 

\subsection
[y-ified Khovanov complex for braids]
{$y$-ified Khovanov complex for braids}
\label{subsec:y-ified-khovanov-braids}

Combining \Cref{prop:y-functor-with-e} and \Cref{prop:beta-dg-mod-structure}, we obtain the following definition. 

\begin{defn}
    Let $\beta \in \BDiag(n)$ be a braid diagram  with underlying permutation $w$. The \textit{$y$-ified Khovanov complex} of $\beta$, denoted $y[\beta]$, is the $y$-complex associated to the $\mcA^w_n$-module $[\beta]$. 
\end{defn}

Recall from \Cref{prop:y-ification-from-dga} that the differential $D$ of $y[\beta]$ is given by 
\[
    D = d + \sum_i \xi_i y_i,
\]
which has curvature
\[
    D^2 = \sum_i [d, \xi_i] y_i = \sum_i (x_i - x'_{w(i)}) y_i.
\]
The differential $D$ can alternatively be written as 
\begin{equation}
\label{eqn:y-ify-Kh-alt}
    D = d + \sum_c \hchi(c) (y_{i_1(c)} - y_{i_2(c)}),
\end{equation}
where the sum ranges over all crossings $c$ of $\beta$, and indices $i_1(c), i_2(c)$ are such that $P_{i_1(c)} \sim p_1(c)$ and $P_{i_2(c)} \sim p_2(c)$. Furthermore, from \Cref{prop:action-e}, the associated $\sle$-operator on $y[\beta]$ is given by 
\[
    \sle = u + \sum_i(x_i + x'_{w(i)}) \ddel{}{y_i}.
\]

\begin{thm}
\label{thm:inv-yified-cpx}
    The $\sle$-equivariant chain homotopy type of the $y$-ified complex $y[\beta]$ is invariant under the three braid relations. 
\end{thm}

The following three subsections are dedicated to prove \Cref{thm:inv-yified-cpx}. For each of the braid relations, we explicitly construct a (strict or homotopy coherent) $\mcA^w_n$-linear homotopy equivalence between the corresponding complexes. We then apply \Cref{prop:y-functor-with-e,prop:homotopy-y-ify-functor} to deduce \Cref{thm:inv-yified-cpx}. 

\subsubsection{Invariance under braid relation I}

\begin{prop}
\label{prop:y-cpx-braid-rel1}
    Let $\beta, \beta'$ be braid diagrams represented as 
    \[
        \beta = \beta_2 \beta_1, \quad \beta' = \beta_2 (\sigma^{-1}_k \sigma_k) \beta_1
    \]
    with $1 \leq k \leq n - 1$. Then $[\beta]$ is a strong deformation retract of $[\beta']$ as $\mcA^w_n$-modules.
\end{prop}

\begin{proof}
    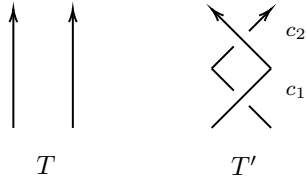
\begin{figure}[t]
        \centering
        \input{tikzpictures/R2}
        \caption{Braid relation I}
        \label{fig:R2}
    \end{figure}
    
    Let $T, T'$ denote the local tangle diagrams inside $\beta, \beta'$, such that $\beta \setminus T$ and $\beta' \setminus T'$ are identical, and $T$ and $T'$ are related by an R2-move (see \Cref{fig:R2}). As described in \cite[Figure 6]{BarNatan:2005}, there is a strong deformation retraction from $[T']$ to $[T]$, which extends as a strong deformation retraction
    \[
        \begin{tikzcd}
        {[\beta]} \arrow[r, "f", shift left] & {[\beta']}. \arrow[l, "g", shift left] \arrow["h"', loop, distance=2em, in=305, out=235]
        \end{tikzcd}
    \]
    We claim these maps are also strictly $\mcA^w_n$-linear. 

    Let $c_1, c_2$ denote the two additional crossings of $\beta'$ corresponding to the generators $\sigma_k, \sigma_k^{-1}$ respectively, and put $\hchi_i := \hchi_{c_i}$. With $\epsilon := \epsilon_{p_1(c_1)}$, we may rewrite $\hchi_i$ using the unsigned homotopies $\hchi_1 = \epsilon \chi_1$ and $\hchi_2 = \epsilon \chi_2$. 
    Let $\xi_i, u$ denote the endomorphisms of the $\mcA^w_n$-action on $[\beta]$, and $\xi'_i, u'$ denote those of $[\beta']$. 
    From the diagram, we can see that 
    \[
        \xi'_i = \begin{cases}
            \xi_i + \epsilon(\chi_1 - \chi_2) & \text{if $i = w_1^{-1}(k)$}, \\ 
            \xi_i - \epsilon(\chi_1 - \chi_2) & \text{if $i = w_1^{-1}(k + 1)$}, \\ 
            \xi_i & \text{otherwise}. 
        \end{cases}
    \]
    Furthermore, from \eqref{eqn:u-by-chi-crossing}, we have
    \begin{align*}
        u' = u &+ \sum_{c \prec c_1} a_{c, 1} \hchi_{c}\hchi_1
             + \sum_{c_1 \prec c} a_{1, c} \hchi_1\hchi_{c} \\
            &+ \sum_{c \prec c_2} a_{c, 2} \hchi_{c}\hchi_2
             + \sum_{c_2 \prec c} a_{2, c} \hchi_2\hchi_{c} \\ 
            &+ a_{1, 2} \hchi_1\hchi_2. 
    \end{align*}
    Here, we have abbreviated the coefficients as $a_{c, i} := a_{c, c_i}$ and $a_{1, 2} := a_{c_1, c_2}$. One can see from the diagram that 
    \[
        a_{c, 1} = -a_{c, 2},\quad 
        a_{1, c} = -a_{2, c}, \quad 
        a_{1, 2} = -2
    \]
    so we have
    \begin{align*}
        u' &= u + \epsilon \sum_{c \prec c_1} a_{c, 1} \hchi_{c}(\chi_1 - \chi_2) \\
            &- \epsilon \sum_{c_2 \prec c} a_{2, c} (\chi_1 - \chi_2)\hchi_{c} \\
            &- 2\chi_1\chi_2.
    \end{align*}

    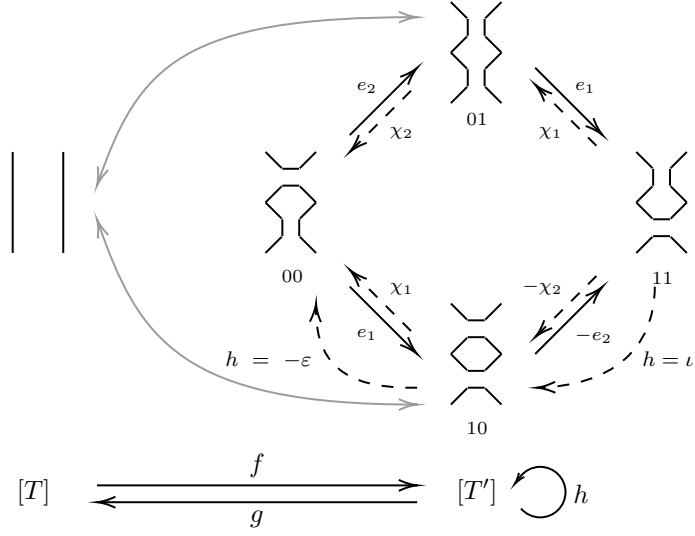
\begin{figure}[t]
        \centering
        \input{tikzpictures/R2-maps}
        \caption{Maps $f, g, h$ corresponding to the R2 move}
        \label{fig:R2-maps}
    \end{figure}

    \Cref{fig:R2-maps} depicts the maps $f, g, h$ and the homotopies $\chi_i$. Observe that the complexes $[T], [T']$ have homological width $0$ and $2$ respectively, and that the maps $h, \chi_2, \chi_3$ have homological degree $-1$. Therefore, if more than two of these maps appear in a composite map, we may immediately conclude that it is zero. Similarly, maps of the form $g \eta f$, $\eta \eta' f$ and $g \eta \eta'$ with $\deg_t(\eta), \deg_t(\eta') < 0$ are also zero. Hereafter, we refer to such reasoning as \textit{degree reasoning}. 
    
    Now, from degree reasons, we immediately obtain $g \chi_i f = 0$ and $g \chi_1 \chi_2 f = 0$, and from the explicit descriptions of $\xi'_i$ and $u'$, we have 
    \[
        g \xi'_i f = \xi_i,\quad 
        g u' f = u.
    \]
    This shows that the $\mcA^w_n$-action on $[T]$ coincides with the one transferred from $[T']$ by the strong deformation retraction. Furthermore, one sees from \Cref{fig:R2-maps} that $T_{11}$ and $T_{00}$ are identical, and
    \begin{align*}
        \chi_1 h |_{11} &= I, & \chi_2 h |_{11} &= 0, \\ 
        h \chi_1 |_{11} &= 0, & h \chi_2 |_{11} &= I.
    \end{align*}
    This gives
    \[
        [\chi_1 - \chi_2, h] = 0.
    \]
    We also have 
    \[
        \chi_1 \chi_2 h = h \chi_1 \chi_2 = 0
    \]
    from degree reasons. Therefore, 
    \[
        [\xi_i, h] = [u, h] = 0
    \]
    and the claim follows from \Cref{cor:strong-defr-linear}. 
\end{proof}

\subsubsection{Invariance under braid relation II}

\begin{prop}
\label{prop:y-cpx-braid-rel2}
    Let $\beta, \beta'$ be braid diagrams represented as 
    \[
        \beta = \beta_2 (\sigma_i \sigma_j) \beta_1, \quad \beta' = \beta_2 (\sigma_j \sigma_i) \beta_1
    \]
    where $|i - j| > 2$. Then $[\beta]$ and $[\beta']$ are isomorphic as $\mcA^w_n$-modules. 
\end{prop}

\begin{proof}
    Since the two braid diagrams are planar isotopic, and the corresponding complexes and maps are identical. 
\end{proof}

\subsubsection{Invariance under braid relation III}

\begin{prop}
\label{prop:y-cpx-braid-rel3}
    Let $\beta, \beta'$ be braids that are represented as 
    \[
        \beta = \beta_2 (\sigma_i \sigma_{i+1} \sigma_i) \beta_1, \quad 
        \beta' = \beta_2 (\sigma_{i+1} \sigma_i \sigma_{i+1}) \beta_1
    \]
    where $1 \leq i \leq n - 2$. Then $[\beta]$ and $[\beta']$ are homotopy equivalent in the category $\Kob_{h\mcA^w_n}(\sfB_n)$.
    
    More strongly, the homotopy equivalences
    \[
        \begin{tikzcd}
        {[\beta]} 
            \arrow[r, "\Psi", shift left] 
        & {[\beta']} 
            \arrow[l, "\Psi'", shift left] 
        \end{tikzcd}
    \]
    can be chosen so that (i) $\Psi$, $\Psi'$  are strongly $1$-truncated, (ii) the compositions $\Psi' \Psi$, $\Psi \Psi'$ are are strongly $1$-truncated, and (iii) $\Psi' \Psi \htpy I, \Psi \Psi' \htpy I$ are given by strongly $1$-truncated homotopies.  
\end{prop}

The additional properties stated in \Cref{prop:y-cpx-braid-rel3} will be used in the proof of \Cref{prop:y-cpx-cl-braid-rel3}. 

\begin{proof}
    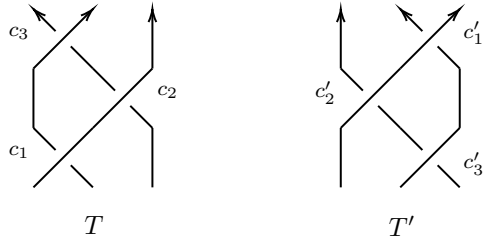
\begin{figure}[t]
        \centering
        \input{tikzpictures/R3}
        \caption{Braid relation III}
        \label{fig:R3}
    \end{figure}
    
    Let $T, T'$ denote the local tangle diagrams inside $\beta, \beta'$, such that $\beta \setminus T$ and $\beta' \setminus T'$ are identical, and $T$ and $T'$ are related by an R3-move (see \Cref{fig:R3}). Recall that, in the proof of \cite[Theorem 1]{BarNatan:2005} for the R3 move, they construct strong deformation retracts $E, E'$ of $[T], [T']$ using the maps for the R2 move, and then prove that $E \isom E'$. The chain homotopy equivalence between $[T]$ and $[T']$ is then given by the composition
    \[
    \begin{tikzcd}
        {[T]} 
            \arrow[r, "\htpy"] 
        & E 
            \arrow[r, leftrightarrow, "\isom"] 
        & E' 
        & {[T']} 
            \arrow[l, "\htpy"'].
    \end{tikzcd}
    \]
    Here, we extend these maps to give the desired equivalence. The proof is rather technical and proceeds in seven steps. Throughout the proof, we rely on the general framework developed in \Cref{subsec:pushforwad-by-strong-def-retract}, and repeatedly use degree reasoning, as we did in the proof of \Cref{prop:y-cpx-braid-rel1}.

    \bigskip 
    \noindent \textbf{Step 1.} Describe the $\mcA^w_n$-actions on $[T]$ and $[T']$.  
    \bigskip 
    
    Observe that $T$ and $T'$ differ by a 180 degree rotation and reversal of orientations. Correspondingly, let $c_1 \prec c_2 \prec c_3$ denote the crossings of $T$, and $c'_1 \succ c'_2 \succ c'_3$ those of $T'$. Here, the subscripts of $c_i$ and $c'_i$ indicate the order of the crossings, whereas $\prec$ indicates the vertical order. Put $\hchi_i = \hchi(c_i)$ and $\hchi'_i = \hchi(c'_i)$. With 
    \[
        \epsilon := \epsilon_{p_1(c_1)},
    \]
    we see from the diagram that
    \[
        \hchi_1 = \epsilon \chi_1,\quad \hchi_2 = -\epsilon \chi_2,\quad \hchi_3 = \epsilon \chi_3
    \]
    and
    \[
        \hchi'_1 = -\epsilon \chi'_1,\quad \hchi_2 = \epsilon \chi'_2,\quad \hchi'_3 = -\epsilon \chi'_3.
    \]
    
    Let $\xi_i, u$ denote the actions of $\mcA^w_n$ on $[\beta]$, and $\xi'_i, u'$ those on $[\beta']$. We may reorder the indices of the strands so that $\xi_1, \xi_2, \xi_3$ corresponds to the three strands appearing in $T$, indexed from left to right. Apply the same reordering for $\beta'$. Then we see that 
    \begin{alignat}{3}
        \xi_1 &= \hchi_1 + \hchi_2 + (\cdots) 
            &&= \epsilon(\chi_1 - \chi_2) + (\cdots), 
            \notag \\
        \xi_2 &= -\hchi_1 + \hchi_3 + (\cdots)
            &&= -\epsilon(\chi_1 - \chi_3) + (\cdots), 
            \label{eq:R3-xi} \\
        \xi_3 &= -\hchi_2 - \hchi_3 + (\cdots)
            &&= \epsilon(\chi_2 - \chi_3) + (\cdots)
            \notag
    \end{alignat}
    and 
    \begin{alignat}{3}
        \xi'_1 &= \hchi'_2 + \hchi'_1 + (\cdots)
            &&= -\epsilon(\chi'_1 - \chi'_2) + (\cdots), 
            \notag \\
        \xi'_2 &= \hchi'_3 - \hchi'_1 + (\cdots)
            &&= \epsilon(\chi'_1 - \chi'_3) + (\cdots), 
            \label{eq:R3-xi'} \\
        \xi'_3 &= -\hchi'_3 - \hchi'_2 + (\cdots)
            &&= -\epsilon(\chi'_2 - \chi'_3) + (\cdots).
            \notag
    \end{alignat}
    Here, each $(\cdots)$ part denotes a signed sum of homotopies $\hchi_c$ corresponding to crossings outside $T$ and $T'$, so that for each $i = 1, 2, 3$, these parts are identical for $\xi_i$ and $\xi'_i$. For the other strands, we have $\xi_i = \xi'_i$ $(i > 3)$. 

    Next, from \eqref{eqn:u-by-chi-crossing} we have
    \begin{align*}
        u
            &= \sum_{\substack{c \prec c' \prec c_1 \\ c_3 \prec c \prec c'}} a_{c, c'} \hchi_{c} \hchi_{c'} \\ 
            &+ \sum_{c \prec c_1} \hchi_{c} \left\{ a_{c, 1} \hchi_1 + a_{c, 2} \hchi_2 + a_{c, 3} \hchi_3 \right\} \\
            &+ \sum_{c_3 \prec c} \left\{ a_{1, c} \hchi_1 + a_{2, c} \hchi_2 + a_{3, c} \hchi_3 \right\} \hchi_{c} \\
            &+ a_{1, 2} \hchi_1 \hchi_2 + a_{1, 3} \hchi_1 \hchi_3 + a_{2, 3} \hchi_2 \hchi_3
    \end{align*}
    and 
    \begin{align*}
        u'
            &= \sum_{\substack{c \prec c' \prec c'_3 \\ c'_1 \prec c \prec c'}} a_{c, c'} \hchi_{c} \hchi_{c'} \\ 
            &+ \sum_{c \prec c'_3} \hchi_{c} \left\{ a'_{c, 3} \hchi'_3 + a'_{c, 2} \hchi'_2 + a'_{c, 1} \hchi'_1  \right\} \\
            &+ \sum_{c'_1 \prec c} \left\{ a'_{3, c} \hchi'_3 + a'_{2, c} \hchi'_2 + a'_{1, c} \hchi'_1 \right\} \hchi_{c} \\
            &+ a'_{3, 2} \hchi'_3 \hchi'_2 + a'_{3, 1} \hchi'_3 \hchi'_1 + a'_{2, 1} \hchi'_2 \hchi'_1.
    \end{align*}
    Here, we used the same abbreviation as we did in the proof of \Cref{prop:y-cpx-braid-rel1}. 
    We see that 
    \[
        a_{c, i} = a'_{c, i}, \quad
        a_{i, c} = a'_{i, c}
    \]
    and 
    \[
        a'_{1, 2} = a'_{2, 1} = 1, \quad 
        a'_{1, 3} = a'_{3, 1} = -1, \quad 
        a'_{2, 3} = a'_{3, 2} = 1
    \]
    Furthermore, by counting the connections, we see that for each crossing $c \prec c_1$, 
    \[
        a_{c, 1} - a_{c, 2} + a_{c, 3} = 0  
    \]
    and similarly, for each crossing $c \succ c_3$,
    \[
        a_{1, c} - a_{2, c} + a_{3, c} = 0.
    \]
    Finally, put 
    \begin{gather*}
        A_c = a_{c, 2},\ B_c = a_{c, 3},\ 
        C_c = a_{2, c},\ D_c = a_{3, c}.
    \end{gather*}
    Combined together, we may write 
    \begin{align}
        u
            &= (\cdots) \notag \\
            &+ \epsilon \sum_{c \prec c_1} \hchi_{c} \left\{ A_c (\chi_1 - \chi_2) - B_c (\chi_1 - \chi_3) \right\} 
            \notag \\
            &+ \epsilon \sum_{c_3 \prec c} \left\{ C_c (\chi_1 - \chi_2) - D_c (\chi_1 - \chi_3) \right\} \hchi_{c} 
            \notag \\
            &- \chi_1 \chi_2 - \chi_1 \chi_3 - \chi_2 \chi_3
            \label{eqn:R3-u}
    \end{align}
    and 
    \begin{align}
        u'
            &= (\cdots) \notag \\ 
            &- \epsilon \sum_{c \prec c'_3} \hchi_{c} \left\{ A_c (\chi'_1 - \chi'_2) - B_c(\chi'_1 - \chi'_3) \right\} 
            \notag \\
            &- \epsilon \sum_{c'_1 \prec c} \left\{ C_c(\chi'_1 - \chi'_2) - D_c(\chi'_1 - \chi'_3) \right\} \hchi_{c} 
            \notag \\
            &+ \chi'_1 \chi'_2 + \chi'_1 \chi'_3 + \chi'_2 \chi'_3.
            \label{eqn:R3-u'}
    \end{align}
    
    \bigskip 
    \noindent \textbf{Step 2.} Construct $E, E'$ and describe the relevant maps. 
    \bigskip 

    Let $T_0, T_1$ denote the tangle obtained from $T$ by resolving the crossing $c_1$ accordingly. Then the complex $[T]$ can be described as the mapping cone of the chain map $e\colon [T_0] \to [T_1]$ corresponding to the resolution change at $c_1$. Namely, we have 
    \[
        [T] = [T_0] \oplus t^1 [T_1]
    \]
    with the differential given by 
    \[
        \begin{pmatrix}
            d_0 & \\ 
            e & -d_1
        \end{pmatrix}.
    \]
    Observe that $T_1$ transforms into a planar tangle diagram $T'_1$ by an R2 move. Thus, there is a strong deformation retraction
    \[
    \begin{tikzcd}
        {[T_1]} 
            \arrow[r, "g", shift left] 
            \arrow["h"', loop, distance=2em, in=305, out=235]
        & {E_1}
            \arrow[l, "f", shift left]
    \end{tikzcd}
    \]
    as described in \Cref{fig:R2-maps}. From \cite[Lemma 4.5]{BarNatan:2005}, this yields a strong deformation retract $E$ of $[T]$, where
    \[
        E = [T_0] \oplus t^1 E_1
    \]
    with differential 
    \[
        \begin{pmatrix}
            d_0 & \\ 
            ge & -d_1
        \end{pmatrix}.
    \]

    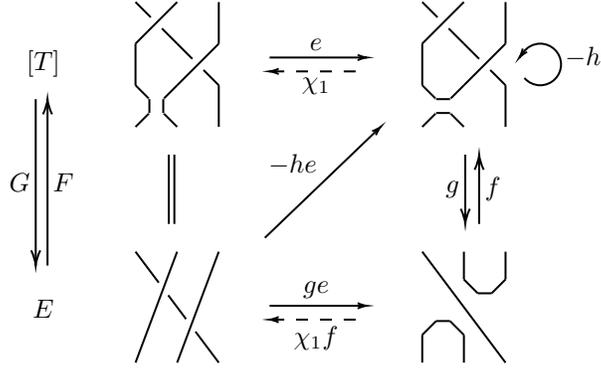
\begin{figure}[t]
        \centering
        \input{tikzpictures/R3-maps}
        \vspace{1em}
        \caption{Strong deformation retract of $[T]$.}
        \label{fig:R3-maps}
    \end{figure}

    \noindent
    The explicit maps
    \[
    \begin{tikzcd}
        {[T]} 
            \arrow[r, "G", shift left] 
            \arrow["H"', loop, distance=2em, in=305, out=235]
        & {E}
            \arrow[l, "F", shift left]
    \end{tikzcd}
    \]
    are given by 
    \begin{equation}
    \label{eqn:R3-FGH}        
        F = \begin{pmatrix}
            I & \\ 
            -he & f
        \end{pmatrix},\quad 
        G = \begin{pmatrix}
            I & \\ 
            & g
        \end{pmatrix},\quad 
        H = \begin{pmatrix}
            0 \\ 
            & -h
        \end{pmatrix}.
    \end{equation}
    See \Cref{fig:R3-maps}. 
    Similarly, there is a strong deformation retract $E'$ of $[T']$ with maps $F', G', H'$ having identical descriptions. 
    
    \bigskip 
    \noindent \textbf{Step 3.} Transfer the $\mcA^w_n$-actions to $E, E'$.
    \bigskip 

    Recall \Cref{prop:strong-defr-induced-dga2}, and consider maps
    \begin{gather*}
        F_i := H \xi_i F,\quad 
        F_u := H u F.
    \end{gather*}
    We claim that these maps are both zero. From \Cref{prop:strong-defr-induced-dga}, it follows that $E$ admits the transferred $\mcA^w_n$-action
    \[
        \bar{\xi}_i := G \xi_i F,\quad 
        \bar{u} := G u F,
    \]
    and that $F$ is strictly $\mcA^w_n$-linear. First, let us describe these maps in the form of matrices. 

    The homotopies $\chi_i$ on $[T]$ are described as 
    \begin{equation}
    \label{eqn:R3-chi}
        \chi_1 = \begin{pmatrix}
            0 & \chi_1 \\ 
            & 0
        \end{pmatrix},
        \quad 
        \chi_i = \begin{pmatrix}
            \chi_i & \\ 
            & -\chi_i
        \end{pmatrix} \ (i = 2, 3).
    \end{equation}
    With \eqref{eqn:R3-FGH}, \eqref{eqn:R3-chi}, we have 
    \[
        H \chi_1 F = \begin{pmatrix}
            0 & 0 \\ 
            0 & 0
        \end{pmatrix},\quad
        H \chi_i F = \begin{pmatrix}
            0 & 0 \\ 
            h \chi_i h e & h \chi_i f
        \end{pmatrix}
    \]
    Here, one can see that the right-hand side of the second equation is also trivial from degree reason. Furthermore, we have 
    \[
        H \chi_1 \chi_i F = \begin{pmatrix}
            0 & 0 \\ 
            0 & 0
        \end{pmatrix},\quad
        H \chi_2 \chi_3 F = 
        \begin{pmatrix}
            0 & 0 \\ 
            -h \chi_2 \chi_3 h e & -h \chi_2 \chi_3 f
        \end{pmatrix}
    \]
    Again, we see that these maps are zero from degree reasons. From the explicit descriptions of $\xi_i$ and $u$ given in \eqref{eq:R3-xi}, \eqref{eqn:R3-u}, the above equations give $F_i = F_u = 0$. Note that we ignored the $(\cdots)$ parts of $\xi_i$ and $u$, since they consist of dot-sliding homotopies $\hchi_c$ outside of $T$, thus commuting with $H, F$, and we have $HF = 0$. The same argument applies to $E'$. 

    \bigskip 
    \noindent \textbf{Step 4.} Extend $G, G'$ as strongly $1$-truncated $h\mcA^w_n$-morphisms.
    \bigskip 

    Again, from \Cref{prop:strong-defr-induced-dga2}, consider maps
    \begin{gather*}
        G_i := G \xi_i H,\quad 
        G_u := -G u H.
    \end{gather*}
    We claim that these maps are $\mcA^w_n$-linear, which implies that $\tilde{G} = (G, G_i, G_u)$ gives a strongly $1$-truncated $h\mcA^w_n$-morphism from $E$ to $[T]$. From \eqref{eqn:R3-FGH}, \eqref{eqn:R3-chi}, we have 
    \[
        G \chi_1 H = \begin{pmatrix}
            0 & -\chi_1 h\\ 
            & 0
        \end{pmatrix}, 
        \quad
        G \chi_i H = \begin{pmatrix}
            0 & \\ 
            & g \chi_i h
        \end{pmatrix}
    \]
    and 
    \[
        G \chi_1 \chi_i H = \begin{pmatrix}
            0 & \chi_1 \chi_i h\\ 
            & 0
        \end{pmatrix} 
        \quad
        G \chi_2 \chi_3 H = \begin{pmatrix}
            0 & \\ 
            & -g \chi_2 \chi_3 h
        \end{pmatrix}
    \]
    where $i = 2, 3$. Again, $G \chi_i H = G \chi_2 \chi_3 H = 0$ from degree reasons. Thus, it suffices to show that $G \chi_1 H$ and $G \chi_1 \chi_i H$ are $\mcA^w_n$-linear. 

    Recall that the transferred $\mcA^w_n$-action on $E$ is given by 
    \[
        \bar{\xi}_i := G \xi_i F,\quad 
        \bar{u} := G u F.        
    \]
    To avoid notational mess, put
    \[
        \bar{\chi}_i := G \chi_i F,\quad 
        \overline{\chi_i \chi_j} = G \chi_i \chi_j F.
    \]
    Explicitly,  
    \[
        \bar{\chi}_1 = \begin{pmatrix}
            -\chi_1 h e & \chi_1 f \\ 
            & 0
        \end{pmatrix},\quad 
        \bar{\chi}_i = \begin{pmatrix}
            \chi_i & \\ 
            & 0
        \end{pmatrix}
    \]
    and 
    \[
        \overline{\chi_1 \chi_i} = \begin{pmatrix}
            \chi_1 \chi_i h e & -\chi_1 \chi_i f \\ 
            & 0
        \end{pmatrix},\quad 
        \overline{\chi_2 \chi_3} = \begin{pmatrix}
            \chi_2 \chi_3 & \\ 
            & 0
        \end{pmatrix}.
    \]
    Note that $\overline{\chi_i \chi_j} = \bar{\chi}_i \bar{\chi}_j$ does not hold. For either $\eta = G \chi_1 H$ or $G \chi_1 \chi_i H$, we agree to write
    \[
        [\chi_j, \eta] := \bar{\chi}_j \eta - (-1)^{|\eta|} \eta \chi_j
    \]
    and 
    \[
        [\chi_j \chi_k, \eta] := \overline{\chi_j \chi_k} \eta - \eta \chi_j \chi_k.
    \]
    
    With this notation, we have
    \[
        [\chi_j, G \chi_1 H]
        = \begin{cases}
            \begin{pmatrix}
            0 & \chi_1 h e \chi_1 h \\ 
            & 0
            \end{pmatrix} 
            & (j = 1) \\ 
            \\
            \begin{pmatrix}
            0 & -\chi_1 [\chi_j, h] \\ 
            & 0
            \end{pmatrix}
            & (j = 2, 3).
        \end{cases}
    \]
    From the explicit descriptions of $e, h$ and $\chi_i$ depicted in \Cref{fig:R2-maps,fig:R3}, one can easily show that
    \[
        (h e) (\chi_1 h) |_{11} = (-I)(-I) = I.
    \]
    Furthermore, as in the proof of \Cref{prop:y-cpx-braid-rel1}, we have
    \begin{align*}
        \chi_2 h |_{11} &= -I, & \chi_3 h |_{11} &= 0, \\ 
        h \chi_2 |_{11} &= 0, & h \chi_3 |_{11} &= -I
    \end{align*}
    (The difference of signs are due to the ordering of the two crossings.) Therefore, we have
    \[
        h e \chi_1 h = -[\chi_2, h] = -[\chi_3, h]
    \]
    which gives 
    \[
        [\chi_1 - \chi_2, G \chi_1 H] = [\chi_1 - \chi_3, G \chi_1 H] = [\chi_2 - \chi_3, G \chi_1 H] = 0
    \]
    and from \eqref{eq:R3-xi} we have 
    \[
        [\xi_i, G \chi_1 H] = 0 \quad (i = 1, 2, 3).
    \]
    Furthermore, from degree reasons, we have
    \[
        [\chi_j \chi_k, G \chi_1 H] = [\chi_j, G \chi_1 \chi_i H] = [\chi_j \chi_k, G \chi_1 \chi_i H] = 0
    \]
    for $i = 2, 3$ and $1 \leq j < k \leq 3$. Therefore, from \eqref{eq:R3-xi}, \eqref{eqn:R3-u}, we obtain
    \[
        [u, G \chi_1 H] 
        = [\xi_j, G \chi_1 \chi_i H] 
        = [u, G \chi_1 \chi_i H] = 0
    \]
    proving that $G \chi_1 H$ and $G \chi_1 \chi_i H$ are $\mcA^w_n$-linear. 
    The exact same argument applies to $G'$. 

    \bigskip 
    \noindent \textbf{Step 5.} Verify that $\tilde{G}, \tilde{G}'$ are strong deformation retractions.
    \bigskip 

    Again from \Cref{prop:strong-defr-induced-dga2}, consider the following maps
    \begin{gather*}
        H_i := H \xi_i H,\quad 
        H_u := H u H.
    \end{gather*}
    One easily sees that 
    \[
        H \chi_i H = 0 \quad (i = 1, 2, 3)
    \]
    and 
    \[
        H \chi_i \chi_j F = 0 \quad (1 \leq i < j \leq 3)
    \]
    so $H_i = H_u = 0$ from \eqref{eq:R3-xi}, \eqref{eqn:R3-u}. 
    From \Cref{prop:strong-defr-induced-dga2}, it follows that 
    \[
        FG - I = [d, H],\quad 
        FG_i = [\xi_i, H],\quad 
        FG_u = [u, H]
    \]
    and we see from \eqref{eqn:homotopy-hi}, \eqref{eqn:homotopy-hu} that $H = (H, 0, 0)$ gives a strongly $1$-truncated homotopy between strongly $1$-truncated $h\mcA^w_n$-morphisms
    \[
        F \tilde{G} = (FG, FG_i, FG_u) 
        \ \htpy \ 
        I = (I, 0, 0).
    \]
    For the other direction, we have $\tilde{G} F = (GF, G_i F, G_u F) = (I, 0, 0)$, so we obtain a strong deformation retraction in the category of $h\mcA^w_n$-modules,
    \[
    \begin{tikzcd}
        {[T]} 
            \arrow[r, "\tilde{G}", shift left] 
            \arrow["H"', loop, distance=2em, in=305, out=235]
        & {E}
            \arrow[l, "F", shift left].
    \end{tikzcd}
    \]
    The exact same argument applies to $H'$.     

    \bigskip 
    \noindent \textbf{Step 6.} Prove $E \isom E'$ as $\mcA^w_n$-modules. 
    \medskip 

    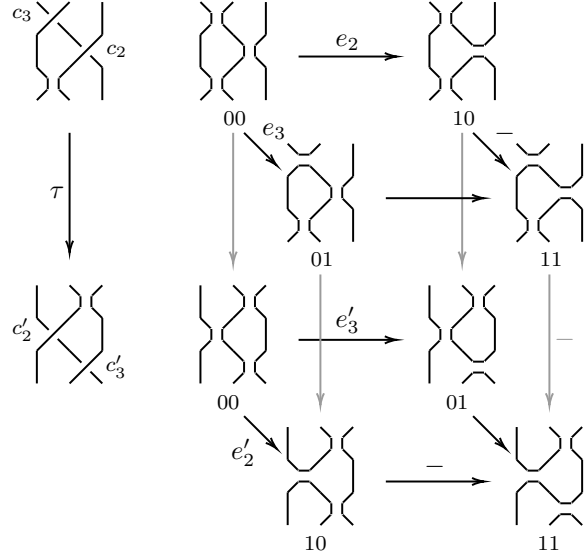
\begin{figure}[t]
        \centering
        \input{tikzpictures/R3-tau}
        \caption{Descriptions of $[T_0], [T'_0]$ and $\tau$.}
        \label{fig:R3-tau}
    \end{figure}

    We have $E = [T_0] \oplus E_1$ and $E' = [T'_0] \oplus E'_1$ from the above constructions. Note that the two diagrams $T_0$ and $T'_0$ are isotopic, but the ordering of the two crossings are opposite. Let 
    \[
        \tau\colon [T_0] \to [T'_0]
    \]
    be the chain isomorphism corresponding to the reordering, which acts as $I$ on each vertex, except at the vertex $(11)$ where it acts by $-I$. See \Cref{fig:R3-tau}. Objects $E_1$ and $E'_1$ are isotopic, so we identify the two. Then, the isomorphism $\Phi\colon E \to E'$ is given by
    \[
        \Phi = \begin{pmatrix}
            \tau & \\ 
            & -I 
        \end{pmatrix}. 
    \]
    We claim that $\Phi$ is an $\mcA^w_n$-linear isomorphism. 

    First, consider the following diagram 
    \[
\begin{tikzcd}[row sep=5em, column sep=5em]
{[T_0]} \arrow[r, equal] \arrow[d, "e"', shift right]                                                            & {[T_0]} \arrow[r, "\tau"] \arrow[d, "ge"', shift right] \arrow["{\chi_2, \chi_3}"', loop, distance=2em, in=125, out=55] & {[T_0']} \arrow[r, equal] \arrow[d, "g'e'"', shift right] \arrow["{\chi'_2, \chi'_3}"', loop, distance=2em, in=125, out=55] & {[T_0']} \arrow[d, "e'"', shift right]                                                                                       \\
{[T_1]} \arrow[u, "\chi_1"', dashed, shift right] \arrow[r, "g", shift left] \arrow["h"', loop, distance=2em, in=305, out=235] & E_1 \arrow[l, "f", shift left] \arrow[r, "-I"] \arrow[u, "\chi_1 f"', dashed, shift right]                              & E'_1 \arrow[r, "f'"', shift right] \arrow[u, "\chi_1' f'"', dashed, shift right]                                                          & {[T_1']} \arrow[u, "\chi'_1"', shift right] \arrow[l, "g'"', shift right] \arrow["h'"', loop, distance=2em, in=305, out=235]
\end{tikzcd}
    \]    
    The following equations are verified by unraveling the explicit descriptions of \Cref{fig:R3-maps,fig:R3-tau}:
    \begin{align*}
        \tau \chi_i = \chi'_j \tau 
            &\colon [T_0] \longrightarrow [T'_0], \\ 
        \tau \chi_1 f = \chi'_1 f'
            &\colon E_1 \longrightarrow [T'_0], \\ 
        ge = - g'e' \tau 
            &\colon [T_0] \longrightarrow E'_1
    \end{align*}
    where $(i, j) \in \{(2, 3), (3, 2)\}$ in the first equation. Furthermore, on $[T_0]$, we have 
    \begin{align*}
        \chi_1 he |_{00} &= 0, & \chi_1 he |_{10} &= 0,\\ 
        \chi_1 he |_{01} &= -\chi_3 |_{01}, & \chi_1 he |_{11} &= -\chi_2 |_{11}
    \end{align*}
    and on $[T'_0]$,
    \begin{align*}
        \chi'_1 h'e' |_{00} &= 0, & \chi'_1 h'e' |_{01} &= -\chi'_3 |_{01},\\ 
        \chi'_1 h'e' |_{10} &= 0, & \chi'_1 h'e' |_{11} &= -\chi'_2 |_{11}.
    \end{align*}
    
    In Step 4, we defined homotopies $\bar{\chi}_i = G \chi_i F$ and $\overline{\chi_i \chi_j} = G \chi_i \chi_j F$ on $E$. Similarly, define $\bar{\chi}'_i := G' \chi'_i F'$ and $\overline{\chi_i \chi_j}' := G' \chi'_i \chi'_j F'$ on $E'$. Using the above equations, one verifies that
    \begin{align*}
        \Phi (\bar{\chi}_1 - \bar{\chi}_i) 
        &= \begin{pmatrix}
            -\tau(\chi_1 h e + \chi_i) & \tau \chi_1 f \\ 
            & 0
        \end{pmatrix} \\ 
        &= \begin{pmatrix}
            (\chi'_1 h' e' + \chi'_i)\tau & \chi'_1 f' \\
            & 0
        \end{pmatrix} \\ 
        &= -(\bar{\chi}_1 - \bar{\chi}_i) \Phi
    \end{align*}
    and  
    \begin{align*}
        \Phi (\bar{\chi}_2 - \bar{\chi}_3) 
        &= \begin{pmatrix}
            \tau(\chi_2 - \chi_3) & {} \\ 
            & 0
        \end{pmatrix} \\ 
        &= \begin{pmatrix}
            (\chi'_3 - \chi'_2)\tau & {} \\
            & 0
        \end{pmatrix} \\ 
        &= -(\bar{\chi}_2 - \bar{\chi}_3) \Phi.
    \end{align*}
    Similarly, one can show that
    \begin{align*}
        \Phi (\overline{\chi_1 \chi_2} + \overline{\chi_1 \chi_3}) 
        &= \begin{pmatrix}
            \tau \chi_1(\chi_2 + \chi_3) he & -\tau \chi_1(\chi_2 + \chi_3) f \\ 
            & 0
        \end{pmatrix} \\ 
        &= \begin{pmatrix}
            -\chi'_1(\chi'_2 + \chi'_3) h'e' \tau & -\chi'_1(\chi'_2 + \chi'_3) f' \tau \\
            & 0
        \end{pmatrix} \\ 
        &= -(\overline{\chi_1 \chi_2}' + \overline{\chi_1 \chi_3}') \Phi
    \end{align*}
    and 
    \begin{align*}
        \Phi \overline{\chi_2 \chi_3} 
        &= \begin{pmatrix}
            \tau \chi_2\chi_3 & {} \\ 
            & 0
        \end{pmatrix} \\ 
        &= \begin{pmatrix}
            \chi'_3 \chi'_2 \tau & {} \\
            & 0
        \end{pmatrix} \\ 
        &= - \overline{\chi_2 \chi_3}' \Phi.
    \end{align*}
    Together with equations \eqref{eq:R3-xi}, \eqref{eq:R3-xi'}, \eqref{eqn:R3-u}, \eqref{eqn:R3-u'}, one sees that $\Phi$ is $\mcA^w_n$-linear.

    \bigskip 
    \noindent \textbf{Step 7.} Define $\Psi, \Psi'$ and verify the conditions. 
    \medskip 
    
    Finally, define maps $\Psi$ and $\Psi'$ by the compositions of the horizontal arrows
    \[
        \begin{tikzcd}
        {[T]} 
            \arrow[r, "\tilde{G}", shift left] 
            \arrow["H"', loop, distance=2em, in=305, out=235] 
        & {E} 
            \arrow[r, "\Phi", shift left] 
            \arrow[l, "F", shift left] 
        & {E'} 
            \arrow[r, "F'", shift left] 
            \arrow[l, "\Phi^{-1}", shift left] 
        & {[T'].} 
            \arrow[l, "\tilde{G}'", shift left] 
            \arrow["H'"', loop, distance=2em, in=305, out=235]
        \end{tikzcd}
    \]
    Since $F, F', \Phi$ are strict and $\tilde{G}, \tilde{G}'$ are strongly 1-truncated, the composite maps $\Psi, \Psi'$ are strongly 1-truncated. Furthermore, 
    \[
        \Psi' \Psi = (F \Phi^{-1} \tilde{G}') (F' \Phi \tilde{G}) = F \tilde{G}
    \]
    is strongly 1-truncated, and is homotopic to $I$ by the strongly $1$-truncated homotopy $(H, 0, 0)$. The exact same argument applies to $\Psi \Psi'$.
\end{proof}

\begin{rem}
    Passing to $y$-ifications, the strong deformation retraction
    \[
    \begin{tikzcd}
        {y([T])} 
            \arrow[r, "y(\tilde{G})", shift left] 
            \arrow["H"', loop, distance=2em, in=305, out=235] 
        & {y(E)} 
            \arrow[l, "F", shift left]
    \end{tikzcd}
    \]
    is described as follows. First, the differential of $y([T])$ is given by 
    \begin{align*}
        D &= d + \sum_i \xi_i y_i \\ 
        &= d + \epsilon (\chi_1 - \chi_2) y_1 - \epsilon (\chi_1 - \chi_3) y_2 + \epsilon (\chi_2 - \chi_3) y_3 \\ 
        &= \begin{pmatrix}
            d_0 \\ 
            e & -d_1
        \end{pmatrix}
        + 
        \epsilon 
        \begin{pmatrix}
            \chi_2(y_3 - y_1) + \chi_3(y_2 - y_3) & \chi_1 (y_1 - y_2)\\ 
            & -\chi_2(y_3 - y_1) - \chi_3(y_2 - y_3) 
        \end{pmatrix}
    \end{align*}
    and the differential of $y(E)$ is given by 
    \begin{align*}
        D &= d + \sum_i \bar{\xi_i} y_i \\
            &= d + \epsilon (\bar{\chi}_1 - \bar{\chi}_2) y_1 - \epsilon (\bar{\chi}_1 - \bar{\chi}_3) y_2 + \epsilon (\bar{\chi}_2 - \bar{\chi}_3) y_3 \\
        &= \begin{pmatrix}
            d_0 \\ 
            ge & -d_1
        \end{pmatrix}
        + 
        \epsilon 
        \begin{pmatrix}
            -\chi_1 h e (y_1 - y_2) + \chi_2 (y_3 - y_1) + \chi_3 (y_2 - y_3) & \chi_1 f (y_1 - y_2) \\ 
            & 0
        \end{pmatrix}.
    \end{align*}
    Maps $F, H$ are unchanged
    \[
        F = \begin{pmatrix}
            I \\ 
            -he & f
        \end{pmatrix}, \quad 
        H = \begin{pmatrix}
            0 \\ 
             & -h
        \end{pmatrix},
    \]
    and the map $y(\tilde{G})$ is given by 
    \begin{align*}
        y(\tilde{G}) &= G + \sum_i G_i y_i \\
        &= \begin{pmatrix}
            I \\ 
            & g
        \end{pmatrix}
        + 
        \epsilon 
        \begin{pmatrix}
            0 & -\chi_1 h(y_1 - y_2) \\ 
            & 0
        \end{pmatrix}.
    \end{align*}
    Setting $y_1 = y_2 = y_3$ recovers the original complex given in \cite[Theorem 1]{BarNatan:2005}.
\end{rem}

\subsection
[y-ified Khovanov homology for links]
{$y$-ified Khovanov homology for links}
\label{subsec:y-ified-khovanov-links}

Next, we consider the closure $\bbeta = \cl(\beta) \in \Diag(\emptyset)$ of a braid diagram $\beta \in \BDiag(n)$. Let $l$ be the number of components of $\bbeta$. For each $i$, the actions $x_i, x'_i$ become equal on $[\bbeta]$, so the $\mcA^w_n$-action on $[\beta]$ becomes central. By choosing a cycle decomposition of $w$,
\[
    w = (i_{1, 1} \cdots i_{1, n_1}) \cdots (i_{l, 1} \cdots i_{l, n_l}),
\]
we may apply \Cref{prop:reduce-central-dga-action} so that $[\bbeta]$ admits a $\mcCA_l$-module structure. The actions are described explicitly as 
\begin{equation}
\label{eqn:A-for-link}
    \bar{x}_k := x_{k, 1},\quad
    \bar{\xi}_k := \sum_{j = 1}^{n_k} \xi_{k, j},\quad 
    \bar{u} := u + \sum_k \tilde{u}_k, \quad
    \tilde{u}_k := \sum_{j < j'} \xi_{k, j} \xi_{k, j'}.
\end{equation}

\begin{defn}
    The \textit{$y$-ified Khovanov complex} of $\bbeta$, denoted $y[\bbeta]$, is the $y$-complex associated to the $\mcCA_l$-module $[\bbeta]$. 
\end{defn}

The differential $D$ of $y[\bbeta]$ is given by 
\[
    D = d + \sum_k \bar{\xi}_k y_k
\]
with $D^2 = 0$. The associated $\sle$-operator on $y[\bbeta]$ is given by 
\[
    \sle = \bar{u} + \sum_k 2 \bar{x}_k \ddel{}{y_k}
\]
Recall from \Cref{prop:reduced-y-ification} that the $R_l$-actions and the $\sle$-operator are independent up to homotopy of the choice of the cycle decomposition of $w$. 

\begin{thm}
\label{thm:closed-y-invariance}
    The $\sle$-equivariant chain homotopy type of $y[\bbeta]$ is invariant under the braid relations and the Markov moves. 
\end{thm}

The proof will be given in the following subsections. \Cref{thm:closed-y-invariance} justifies the following definition. 

\begin{defn}
    Let $L$ be a link represented by a braid diagram $\beta$. The \textit{$y$-ified $U(2)$-equivariant Khovanov complex} of $\bbeta$, denoted $y\CKh_{U(2)}(\bbeta)$, is the bigraded complex over the polynomial ring $R_{U(2)}[y_1, \ldots, y_l]$ obtained by applying the TQFT $\mcF$ to the complex $y[\bbeta]$. Its homology, denoted $y\Kh_{U(2)}(L)$, is called the \textit{$y$-ified $U(2)$-equivariant Khovanov homology} of $L$. By specializing $h, t \in R_{U(2)}$ before taking the homology, we obtain the $y$-ified homology for other variants. In particular, for $h = t = 0$, the corresponding complex $y\CKh(\bbeta)$ is called the \textit{$y$-ified Khovanov complex} of $\bbeta$, and its homology $y\Kh(L)$ the \textit{$y$-ified Khovanov homology} of $L$.
\end{defn}

\subsubsection{Invariance under the braid relations.}
\label{subsubsec:inv-closed-braid}

We state the closed versions of \Cref{prop:y-cpx-braid-rel1,prop:y-cpx-braid-rel2,prop:y-cpx-braid-rel3}. In each of the cases, it is assumed that a cycle decomposition \eqref{eqn:w-cycle-decomp} of the underlying permutation $w$ of the braid is fixed. 

\begin{prop}
\label{prop:y-cpx-cl-braid-rel1}
    Let $\beta, \beta'$ be braid diagrams represented as 
    \[
        \beta = \beta_2 \beta_1, \quad \beta' = \beta_2 (\sigma^{-1}_k \sigma_k) \beta_1
    \]
    with $1 \leq k \leq n - 1$. Then $[\bbeta]$ is a strong deformation retract of $[\bbeta']$ as $\mcCA_l$-modules.
\end{prop}

\begin{prop}
\label{prop:y-cpx-cl-braid-rel2}
    Let $\beta, \beta'$ be braid diagrams represented as 
    \[
        \beta = \beta_2 (\sigma_i \sigma_j) \beta_1, \quad \beta' = \beta_2 (\sigma_j \sigma_i) \beta_1
    \]
    where $|i - j| > 2$. Then $[\bbeta]$ and $[\bbeta']$ are isomorphic as $\mcCA_l$-modules. 
\end{prop}

\begin{prop}
\label{prop:y-cpx-cl-braid-rel3}
    Let $\beta, \beta'$ be braids that are represented as 
    \[
        \beta = \beta_2 (\sigma_i \sigma_{i+1} \sigma_i) \beta_1, \quad 
        \beta' = \beta_2 (\sigma_{i+1} \sigma_i \sigma_{i+1}) \beta_1
    \]
    where $1 \leq i \leq n - 2$. Then $[\bbeta]$ and $[\bbeta']$ are homotopy equivalent in the category $\Kob_{h\mcCA_l}(\emptyset)$.
\end{prop}

For braid relations I and II (\Cref{prop:y-cpx-cl-braid-rel1,prop:y-cpx-cl-braid-rel2}), the result is immediate from \Cref{prop:y-cpx-braid-rel1,prop:y-cpx-braid-rel2} combined with \Cref{prop:reduce-functor}. Indeed, for each case, we have $\mcA^w_n$-linear homotopy equivalences
\[
    \begin{tikzcd}
        {[\beta]} 
            \arrow[r, "f", shift left] 
        & {[\beta'].} 
            \arrow[l, "g", shift left] 
    \end{tikzcd}
\]
The closure functor $\cl$ preserves the homotopy equivalence, and give $\mcA^w_n$-morphisms between the central $\mcA^w_n$-modules
\[
    \begin{tikzcd}
        {[\bbeta]} 
            \arrow[r, "\bar{f}", shift left] 
        & {[\bbeta'].} 
            \arrow[l, "\bar{g}", shift left] 
    \end{tikzcd}
\]
Then, from \Cref{prop:reduce-functor}, by applying the reduction functor, these maps can be regarded as $\mcCA_l$-equivalences between the $\mcCA_l$-modules.

For braid relation III, we only have homotopy coherent equivalences instead of strict ones, and the reduction functor of \Cref{prop:reduce-functor} cannot be applied directly. Nonetheless, the following ``partial functoriality" is available from the additional properties stated in the latter part of \Cref{prop:y-cpx-braid-rel3}, and then \Cref{prop:y-cpx-cl-braid-rel3} can be deduced from an analogous argument.

\begin{lem}
\label{prop:reduce-1-truncated-mor}
    Let $C, C'$ be $\mcCA^w_n$-modules, and $\bar{C}, \bar{C}'$ be the reduction of $C, C'$ as $\mcCA_l$-modules, obtained from a cycle decomposition \eqref{eqn:w-cycle-decomp} of $w$. Let $W^1(C, C')$ denote the set of all strongly $1$-truncated $h\mcCA^w_n$-morphisms from $C$ to $C'$, and similarly $W^1(\bar{C}, \bar{C}')$ the corresponding set for $(\bar{C}, \bar{C}')$. We say $F, F'$ are homotopic in $W^1(C, C')$ if $F, F'$ is homotopic by a strongly $1$-truncated homotopy $H$. Now, for each pair $(C, C')$ of $\mcCA^w_n$-modules, there is a map
    \[
        W^1(C, C') \to W^1(\bar{C}, \bar{C}'), \quad
        F \mapsto \bar{F} 
    \]
    satisfying the following properties: 
    \begin{enumerate}
        \item If $F \in W^1(C, C'), G \in W^1(C', C'')$ and $GF \in W^1(C, C'')$, then $\bar{G} \bar{F} \htpy \overline{GF}$ in $W^1(\bar{C}, \bar{C}'')$.
        \item If $F \htpy F'$ in $W^1(C, C')$, then $\bar{F} \htpy \bar{F}'$ in $W^1(\bar{C}, \bar{C}')$.
    \end{enumerate}
    In particular, if $F \in W^1(C, C')$, $G \in W^1(C', C)$, $GF \in W^1(C, C)$ and $GF \htpy 1_C$ in $W^1(C, C)$, then $\bar{G} \bar{F} \htpy 1_{\bar{C}}$ in $W^1(\bar{C}, \bar{C})$. 
\end{lem}

\begin{proof}
    Let $F = (f, f_i, f_u)$ be a strongly $1$-truncated $h\mcA^w_n$-morphism from $C$ to $C'$. The corresponding map $\bar{F} = (\bar{f}, \bar{f}_k, \bar{f}_u)$ is given by 
    \[
        \bar{f} := f,\quad
        \bar{f}_k := \sum_j f_{k, j},\quad 
        \bar{f}_u := f_u + \sum_k \tilde{f}_k
    \]
    where we put
    \[
        \tilde{f}_k := -\sum_{j < j'} (\xi_{k, j} f_{k, j'} - f_{k, j} \xi_{k, j'}).
    \]
    First, we verify that $\bar{F}$ is indeed a strongly $1$-truncated $h\mcA^w_n$-morphism. The required condition \eqref{eqn:f-extension-lowest1} for $\bar{f}_i$ holds from
    \[
        [d, \bar{f}_k] + [\bar{\xi}_k, f] = \sum_j [d, f_{k, j}] + \sum_j [\xi_{k, j}, f] = 0.
    \]
    The required condition \eqref{eqn:f-extension-lowest2} for $\bar{f}_u$ is
    \[
        [d, \bar{f}_u] + [\bar{u}, \bar{f}] + 2 \sum_k \bar{x}_k \bar{f}_k = 0.
    \]
    From the assumption, we have 
    \[
        [d, f_u] + [u, f] + \sum_i (x_i + x_{w(i)}) f_i = 0
    \]
    so it suffices to show that for each $k$,
    \[
        [d, \tilde{f}_k] + [\tilde{u}_k, f] + \sum_j ( 2\bar{x}_k - (x_{k, j} + x_{k, j+1})) f_{k, j} = 0 
    \]
    Indeed, we have 
    \begin{align*}
        [\tilde{u}_k, f] 
            &= \sum_{j < j'} [\xi_{k, j} \xi_{k, j'}, f] \\
            &= \sum_{j < j'} \left( \xi_{k, j}[\xi_{k, j'}, f] + [\xi_{k, j}, f] \xi_{k, j'} \right) \\ 
            &= -\sum_{j < j'} \left( \xi_{k, j}[d, f_{k, j'}] + [d, f_{k, j}] \xi_{k, j'} \right)
    \end{align*}
    and 
    \[
        \bar{x}_k - x_{k, j} = \sum_{j' < j} [d, \xi_{k, j'}],\quad 
        \bar{x}_k - x_{k, j + 1} = -\sum_{j < j'} [d, \xi_{k, j'}].
    \]
    Therefore, 
    \begin{align*}
        &[\tilde{u}_k, f] + \sum_j ( 2\bar{x}_k - (x_{k, j} + x_{k, j+1})) f_{k, j}\\
        &= -\sum_{j < j'} \left(\xi_{k, j}[d, f_{k, j'}] + [d, f_{k, j}] \xi_{k, j'} \right) 
            + \sum_{j < j'} \left ( [d, \xi_{k, j}] f_{k, j'} -  f_{k, j} [d, \xi_{k, j'}] \right ) \\ 
        &= \sum_{j < j'} \left( [d, \xi_{k, j} f_{k, j'}] - [d, f_{k, j} \xi_{k, j'}] \right) \\ 
        &= -[d, \tilde{f}_k].
    \end{align*}
    Furthermore, $f_i$ and $f_u$ are assumed to be $\mcA^w_n$-linear, so $\bar{f}_k$ and $\bar{f}_u$ both commute with $\xi_{k, j}$ and $u$, and hence are $\mcCA_l$-linear. Therefore, $\bar{F}$ is a strongly $1$-truncated $h\mcCA_l$-morphism. 
    
    Next, we verify the two properties. For $F = (f, f_i, f_u)$ and $G = (g, g_i, g_u)$, we have $GF = (gf, g_i f + g f_i, g_u f + g f_u)$. The corresponding maps are, $\bar{F} = (\bar{f}, \bar{f}_k, \bar{f}_u)$, $\bar{G} = (\bar{g}, \bar{g}_k, \bar{g}_u)$, and 
    \begin{align*}
        \overline{GF}_0 &= gf = \bar{g} \bar{f} \\ 
        \overline{GF}_k &= \sum_j (g_{k, j} f + g f_{k, j}) = \bar{g}_k \bar{f} + \bar{g} \bar{f}_k, \\ 
        \overline{GF}_u &= (g_u f + g f_u) - \sum_k \sum_{j < j'} (\xi_{k, j} (g_{k, j'} f + g f_{k, j'}) - (g_{k, j} f + g f_{k, j}) \xi_{k, j'}) \\ 
            &= (\bar{g}_u \bar{f} + \bar{g} \bar{f}_u) + \sum_k \sum_{j < j'} ([d, g_{k, j}] f_{k, j'} + g_{k, j} [d, f_{k, j'}]) \\
            &= (\bar{g}_u \bar{f} + \bar{g} \bar{f}_u) + \sum_k \sum_{j < j'} [d, g_{k, j} f_{k, j'}].
    \end{align*}
    Therefore, 
    \[
        H = (0, 0, \sum_k \sum_{j < j'} g_{k, j} f_{k, j'})
    \]
    gives $\bar{G} \bar{F} \htpy \overline{GF}$. 
    
    For the second property, suppose we have a strongly $1$-truncated null-homotopy $H = (h, h_i, h_u)$ for $F \htpy 0$. Consider the correspondence $H \mapsto \bar{H} = (\bar{h}, \bar{h}_i, \bar{h}_u)$ is given by 
    \[
        \bar{h} := h,\quad
        \bar{h}_k := \sum_j f_{k, j},\quad 
        \bar{h}_u := h_u + \sum_k \tilde{h}_k
    \]
    where 
    \[
        \tilde{h}_k := \sum_{j < j'} \left( \xi_{k, j} h_{k, j'} + h_{k, j} \xi_{k, j'} \right).
    \]
    We claim that $\bar{H}$ is a strongly $1$-truncated null homotopy for $\bar{F} \htpy 0$. First, \eqref{eqn:homotopy-hi} is easy to verify. For \eqref{eqn:homotopy-hu}, we must have 
    \[
        \bar{f}_u = [d, \bar{h}_u] + [\bar{u}, \bar{h}] + \sum_k 2\bar{x}_k \bar{h}_k
    \]
    which follows from, 
    \[
        \tilde{f}_k = [d, \tilde{h}_k] + [\tilde{u}_k, h] + \sum_j (2\bar{x}_k - (x_{k, j} + x_{k, j + 1}))h_{k, j} .
    \]
    The verification is similar to the one we showed for $[d, \tilde{f}_k]$. It is also easy to see that $\bar{h}_k, \bar{h}_u$ are $\mcCA_l$-linear, and hence $\bar{H}$ is strongly $1$-truncated. 
\end{proof}

\subsubsection{Invariance under Markov move I.}

\begin{prop}
\label{prop:markov-inv1}
    Let $\beta, \beta'$ be braids that are represented as 
    \[
        \beta = \beta_2 \beta_1,\quad 
        \beta' = \beta_1 \beta_2.
    \]
    Then $y[\bbeta]$ and $y[\bbeta']$ are isomorphic as $\sle$-equivariant $y$-complexes. 
\end{prop}

\begin{proof}
    First, observe that $\bbeta$ and $\bbeta'$ are planar isotopic, so we may naturally identify the complexes $[\bbeta]$ and $[\bbeta']$. 
    For $a = 1, 2$, let $w_a$ be the underlying permutation of $\beta_a$. Then $w = w_2 w_1$ is the underlying permutation of $\beta$, and $w' = w_1 w_2$ is the underlying permutation of $\beta'$. Take a cycle decomposition of $w$ 
    \[
        w = (i_{1, 1} \cdots i_{1, n_1}) \cdots (i_{l, 1} \cdots i_{l, n_l})
    \]
    and one of $w'$
    \[
        w' = (i'_{1, 1} \cdots i'_{1, n_1}) \cdots (i'_{l, 1} \cdots i'_{l, n_l})
    \]
    so that $i'_{k, j} = w_1(i_{k, j})$ for every $(k, j)$. 
    
    Let $\xi^a_i$ denote the actions of $\xi_i$ on $[\beta_a]$, and $\xi_i, \xi'_i$ denote those for $[\beta], [\beta']$ respectively. Then we have 
    \begin{align*}
        \xi_{i_{k, j}} 
        &= \xi^1_{i_{k, j}} + \xi^2_{w_1(i_{k, j})}
        = \xi^1_{k, j} + \xi^2_{k, j},\\
        \xi'_{i'_{k, j}}
        &= \xi^2_{i'_{k, j}} + \xi^1_{w_2(i'_{k, j})}
        = \xi^2_{k, j} + \xi^1_{k, j + 1}.
    \end{align*}
    Therefore, 
    \[
        \bar{\xi}_k 
        = \sum_j \xi_{k, j}
        = \sum_j \xi'_{k, j}
        = \bar{\xi}'_k
    \]
    on $[\bar{\beta}] = [\bar{\beta}']$. Thus, the $y$-complexes $y[\bbeta]$ and $[\bbeta']$ can also be identified. 
    Next, the $\sle$-operators for $y[\bbeta]$ and $y[\bbeta']$ are given by 
    \begin{align*}
        \sle &= 2 \sum_k \bar{x}_k \ddel{}{y_k} + u + \sum_k \tilde{u}_k, \quad
        \tilde{u}_k = \sum_{j < j'} \xi_{k, j} \xi_{k, j'}, \\
        \sle' &= 2 \sum_k \bar{x}'_k \ddel{}{y_k} + u' + \sum_k \tilde{u}'_k, \quad
        \tilde{u}'_k = \sum_{j < j'} \xi'_{k, j} \xi'_{k, j'}.
    \end{align*}
    We claim that 
    \[
        H := \sum_k H_k, \quad H_k := 2 \xi^1_{k, 1} \ddel{}{y_k}
    \]
    gives 
    \[
        \sle - \sle' = [D, H].
    \]
    This can be verified directly using the given equations and \eqref{eqn:u-by-braid}.  
\end{proof}

\subsubsection{Invariance under Markov move II.}

\begin{prop}
\label{prop:markov-inv2}
    Let $\beta$ be a braid on $n$ strands, and $\beta'$ a braid on $n + 1$ strands, satisfying
    \[
        \beta' = \beta \sigma_n.
    \]
    Then $[\beta]$ is a strong deformation retract of $[\beta']$ as $\mcCA_l$-modules. 
\end{prop}

\begin{proof}
    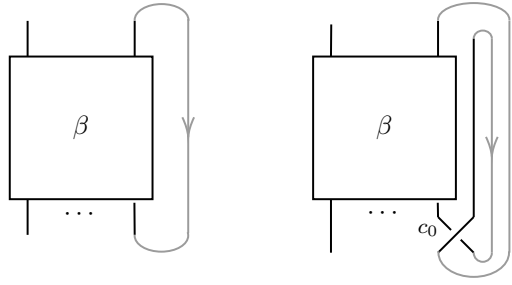
\begin{figure}[t]
        \centering
        \input{tikzpictures/R1}
        \caption{Markov move II}
        \label{fig:R1}
    \end{figure}
    \Cref{fig:R1} depicts the two braids $\beta$ and $\beta'$. Let $c_0$ denote the additional crossing of $\beta'$ that corresponds to $\sigma_n$. Since their closures $\bbeta$ and $\bbeta'$ are related by an R1-move, from \cite[Theorem 1]{BarNatan:2005}, there is a strong deformation retraction between the complexes 
    \[
    \begin{tikzcd}
        {[\bbeta]} \arrow[r, "f", shift left] & {[\bbeta'].} \arrow[l, "g", shift left] \arrow["h"', loop, distance=2em, in=305, out=235]
    \end{tikzcd}
    \]
    
    First, consider the $\mcA^w_n$-action on $[\beta]$ given by endomorphisms $\xi_i, u$ $(i = 1, \ldots, n)$, and the $\mcA^w_{n + 1}$-action on $[\beta']$ given by $\xi'_i, u'$ $(i = 1, \ldots, n + 1)$. Obviously $\xi_i = \xi'_i$ for $i < n$ (regarded locally, away from $c_0$), and 
    \[
        \xi'_n + \xi'_{n + 1} = \hchi_{c_0} - \hchi_{c_0} + \xi_n = \xi_n. 
    \]
    Therefore, after taking the closure and passing to their $\mcCA_l$-reductions, we have $\bar{\xi}_k = \bar{\xi}'_k$ for all $k$, and since all $\xi_i$ commute with $f$ and $g$, we have 
    \[
        g \bar{\xi}'_k f = g f \bar{\xi}_k = \bar{\xi}_k.
    \]
    Similarly, we easily see from the formula \eqref{eqn:u-by-chi-crossing} that $g \bar{u}' f = \bar{u}$. Therefore, the $\mcCA_l$-action on $\beta$ is transferred from $\beta'$. Furthermore, it is obvious that $[\bar{\xi}'_k, h] = [\bar{\xi}_k, h] = 0$ and $[\bar{u}', h] = [\bar{u}, h] = 0$. Therefore, the result follows from \Cref{prop:strong-defr-induced-dga2}. 
\end{proof}

\subsection{Reduced versions}
\label{subsec:reduced}

The construction of the $y$-ified Khovanov homology extends verbatim to the reduced setting. First, let us recall the definition of \textit{reduced Khovanov homology}, originally introduced in \cite{Khovanov:2003}. 

Consider a pointed link $(L, p)$ represented as a braid $\beta$ with a regular point $p \in \bbeta$. Let us consider the $U(1)$-equivariant theory, given by the ring $R_{U(1)} = \ZZ[h]$ and the Frobenius algebra $A_{U(1)} = R_{U(1)}/(X^2 - hX)$. The \textit{reduced $U(1)$-Khovanov complex} $\rCKh_{U(1)}(\bbeta, p)$ of the pointed link diagram $(\bbeta, p)$ is defined by the image of the endomorphism $x_p$,
\[
    \rCKh_{U(1)}(\bbeta, p) := \Ima(x_p \colon \CKh_{U(1)}(\bbeta) \to \CKh_{U(1)}(\bbeta)),
\]
regarded as a subcomplex of $\CKh_{U(1)}(\bbeta)$. One may alternatively define $\rCKh(\bbeta, p)$ by the quotient
\[
    \CKh_{U(1)}(\bbeta) / \Ima(x_p),
\]
or by using the conjugate map $h - x_p$ instead of $x_p$; all give isomorphic complexes (see \cite[Propositions 3.12, 3.14]{Sano-Sato:2023}). The \textit{reduced $U(1)$-Khovanov homology} $\rKh_{U(1)}(\bbeta, p)$ is defined by the homology of $\rCKh_{U(1)}(\bbeta, p)$, which is an invariant of the pointed link $(L, p)$, depending only on $L$ and the component to which $p$ belongs (see \cite[Propositions 3.16, 3.18]{Sano-Sato:2023}). In particular, for a knot $K$, the reduced homology is independent of the choice of $p$, and hence we write $\rKh_{U(1)}(K)$ and call it the \textit{reduced $U(1)$-Khovanov homology} of $K$. Setting $h = 0$ before taking the homology recovers the original reduced theory. 

\begin{prop}
    The $\mcCA_l$-action on $\CKh_{U(1)}(\bbeta)$ restricts to the subcomplex $\rCKh_{U(1)}(\bbeta, p)$.
\end{prop}

\begin{proof}
    The action of $R^e_n$ on $\CKh_{U(1)}(\bbeta)$ restricts to $\rCKh_{U(1)}(\bbeta, p)$. Since the dot-sliding homotopies also restrict to $\rCKh_{U(1)}(\bbeta, p)$, so do the endomorphisms $\xi_i$ and $u$.
\end{proof}

Therefore, the \textit{$y$-ified reduced $U(1)$-Khovanov complex} $y\rCKh_{U(1)}(\bbeta, p)$ is defined as the $y$-complex associated to the $\mcCA_l$-module $\rCKh_{U(1)}(\bbeta, p)$, and is regarded as a subcomplex of the unreduced $y$-ification $y\CKh_{U(1)}(\bbeta)$. 

\begin{prop}
    The $\sle$-equivariant chain homotopy type of $y\rCKh_{U(1)}(\bbeta, p)$ is an invariant of the pointed link.
\end{prop}

\begin{proof}
    After sliding the basepoint away from the local moves, all of the maps used to prove \Cref{thm:closed-y-invariance} restrict to the reduced complex $\rCKh_{U(1)}(\bbeta, p)$.
\end{proof}

We call the homology of $y\rCKh_{U(1)}(\bbeta, p)$ the \textit{$y$-ified reduced $U(1)$-Khovanov homology} of the pointed link $(L, p)$, and denote it by $y\rKh_{U(1)}(L, p)$. Setting $h = 0$ before taking the homology gives the \textit{$y$-ified reduced Khovanov homology} $y\rKh(L, p)$. 

\begin{rem}
    Similarly, we may consider the reduced version for the \textit{$U(1) \times U(1)$-equivariant theory}, given by the Frobenius extension 
    \[
        R_{U(1) \times U(1)} = \ZZ[\alpha_1, \alpha_2],\quad 
        A_{U(1) \times U(1)} = R_{U(1) \times U(1)}[X]/(X - \alpha_1)(X - \alpha_2).
    \]
    By replacing the definition of the endomorphism $x_p$ by 
    \[
        x_p = \begin{cases}
            X_p - \alpha_1 & \text{ if $\theta(p) = a$, } \\
            -(X_p - \alpha_2) & \text{ if $\theta(p) = b$},
        \end{cases}
    \]
    we similarly obtain a $y$-ification equipped with the $\sle$-operator for the reduced homology. 
\end{rem}

\subsection{Extension to tangles}
\label{subsec:ext-to-tangles}

The dg-module structures on $[\beta]$ and $[\bbeta]$ can be naturally generalized to dg-module structures  on tangles. 
Let $B \subset \partial D^2$ be a finite set of signed points with equal numbers of negative points $P_1, \ldots, P_n$ and positive points $P'_1, \ldots, P'_n$. Let $T \in \Diag(B)$ be an oriented tangle diagram with $n$ arcs and $l$ loops. Let $w \in \mfS_n$ denote the permutation giving $P_i \sim P'_{w(i)}$ for each index $i$. For each loop $C_k$ of $T$, choose one regular point $\bar{P}_k$ on $C_k$. Fix an admissible $ab$-coloring $\theta$ on $T$. Then $[T]$ admits an $R^e_n \otimes R_l$-module structure given by 
\begin{equation}
\label{eqn:A-for-tangle}
    x_i := x_{P_i}, \quad
    x'_i := x_{P'_i}, \quad 
    \bar{x}_k := x_{\bar{P}_k}\quad 
    (i = 1, \ldots, n;\ k = 1, \ldots, r).
\end{equation}
We also have dot-sliding homotopies $\chi_c$ for each crossing $c$ of $T$. Let $\gamma_i$ be the path from $P_i$ to $P'_i$, and $\bar{\gamma}_k$ the loop that starts from $\bar{P}_k$, runs along $C_k$ in the direction given by the orientation and returns back to $\bar{P}_k$. Define degree $-1$ homotopies
\[
    \xi_i := \xi_{\gamma_i},\quad 
    \bar{\xi}_k := \xi_{\bar{\gamma}_k}
\]
which give 
\[
    [d, \xi_i] = x_i - x'_{w(i)},\quad 
    [d, \bar{\xi}_k] = 0.
\]
Furthermore, as in the proof of \Cref{prop:homotopy-u}, define degree $-2$ homotopies
\[
    u_i := u_{\gamma_i},\quad   
    \bar{u}_k := u_{\bar{\gamma}_k}
\]
and
\[
    u := \sum_i u_i + \sum_k \bar{u}_k.
\]
Then, by an argument similar to the proof of \Cref{prop:homotopy-u},
\[
    [d, u] = \sum_i (x_i + x'_{w(i)}) \xi_i + 2\sum_k \bar{x}_k \bar{\xi}_k.
\]
Define
\[
    \mcA^w_{n, l} := \mcA^w_n \otimes_R \mcCA_l
\]
which is a dga over $R^e_n \otimes_R R_l$. The above constructions show that $[T]$ admits an $\mcA^w_{n, l}$-module structure. Since $\mcA^w_{n, l}$ naturally embeds into $\mcA^w_{n + l}$, we obtain the associated $y$-complex by
\[
    y[T] = [T] \otimes_R R[y_1, \ldots, y_n; \bar{y}_1, \ldots, \bar{y}_l]
\]
with differential 
\[
    D := d + \sum_i \xi_i y_i + \sum_k \bar{\xi}_k \bar{y}_k
\]
and the $\sle$-operator 
\[
    \sle = u + \sum_i (x_i + x'_{w(i)}) \ddel{}{y_i} + 2\sum_k \bar{x}_k \ddel{}{\bar{y}_k}. 
\]

\begin{prop}
\label{thm:tangle-y-invariance}
    The $\sle$-equivariant chain homotopy type of $y[T]$ is independent of the choice of the base points on the loops, and is invariant under the Reidemeister moves. 
\end{prop}

\begin{proof}
    Given a tangle diagram $T \in \Diag(\sfB_n)$ with $l$ loops, we obtain a $2(n + l)$-end tangle diagram $T^\times$ by cutting $T$ open at each point $\bar{P}_k$. Although the endpoints may lie in the interior of $D^2$, we may define the complex $[T^\times]$ and give it an $\mcA^w_{n + l}$-action in the obvious way. Now, suppose two tangle diagrams $T$ and $T'$ are related by a Reidemeister move. Then by repeating the proofs of \Cref{prop:markov-inv2,prop:y-cpx-cl-braid-rel1,prop:y-cpx-cl-braid-rel3}, we can show that $[T^\times]$ and $[T'^\times]$ are homotopy equivalent as $\mcA^w_{n + l}$- or $h\mcA^w_{n + l}$-modules. Passing from $[T^\times], [T'^\times]$ back to $[T], [T']$ is realized by a correspondence similar to the reduction functor \Cref{prop:reduce-functor}, and we may repeat the argument of \Cref{subsubsec:inv-closed-braid} to show that $[T]$ and $[T']$ are homotopy equivalent. Further applying the $y$-ification functor gives the invariance of the corresponding $y$-complexes. Finally, the independence of the choice of the base points can be proved by repeating the proof of \Cref{prop:markov-inv1}. 
\end{proof}

Next, we extend the horizontal composition by planar arc diagrams to the category of dg-modules. Let $D$ be a 2-input planar arc diagram with inner boundary points $B_1, B_2$ and outer boundary points $B$. Put $n = |B|/2$ and $n_a = |B_a|/2$ ($a = 1, 2$). For $a = 1, 2$, let $X_a$ be an object of $\Kob(B_a)$ equipped with an $A^{w_a}_{n_a, l_a}$-action given by endomorphisms $x^a_i$, $x'^a_i$, $\bar{x}^a_k$, $\xi^a_i$, $\bar{\xi}^a_k$, $u^a$. Let $X = D(X_1, X_2)$ be the composed complex as an object of $\Kob(B)$. We say $X_1, X_2$ have \textit{compatible actions} with respect to $D$ if, for any pair of boundary points $p_1 \in B_1$ and $p_2 \in B_2$ that are connected by an arc in $D$, the actions $x^1_{p_1}$ from $X_1$ and $x^2_{p_2}$ from $X_2$ coincide. In such case, $X$ is given an $\mcA^w_{n, l}$-module structure as follows. 


We regard each homotopy $\xi^a_i$ as a \textit{virtual arc} that connects $P^a_i$ and $P'^a_{w_a(i)}$ within the input disk $D_a$. A \textit{virtual path} $\gamma$ in $X$ is a sequence,
\[
    \gamma\colon 
    p_1 \xrightarrow{\zeta_1} p_2 \xrightarrow{\zeta_2} \cdots \xrightarrow{\zeta_r} p_{r + 1}
\]
where each $p_i$ is a point on an arc of $D$, and $\zeta_i$ is one of the virtual arcs $\{\xi^a_i\}$. $\gamma$ is called a \textit{virtual strand} if it starts from a negative endpoint of $B$ and ends with a positive endpoint of $B$. $\gamma$ is called a \textit{virtual loop} if $p_1 = p_{r + 1}$ and it does not contain the same virtual arc twice. Let $l'$ be the number of virtual loops, and for each virtual loop of $X$, choose a point on an arc of $D$. Put $l = l_1 + l_2 + l'$. The permutation $w \in \mfS_n$ is defined by the correspondence of the virtual strands, and $X$ is given an $R^e_n \otimes_R R_l$-action from the polynomial actions of $X_1$ and $X_2$ in the obvious way. Furthermore, for the above virtual path $\gamma$, we define
\begin{equation}
\label{eqn:virtual-path-homotopy}
    \xi_{\gamma} := \sum_j \zeta_j. 
\end{equation}
We define homotopies $\xi_i$ ($i = 1, \ldots, n$) for the virtual strands $\gamma_i$ and homotopies $\bar{\xi}_k$ ($k = 1, \ldots, l'$) for the virtual loops $\bar{\gamma}_k$ by \eqref{eqn:virtual-path-homotopy}. Combining the endomorphisms $\bar{\xi}^a_k$ from $X_a$, we get endomorphisms $\bar{\xi}_k$ ($k = 1, \ldots, l$). Furthermore, define
\[
    u := u^1 + u^2 + \sum_{\gamma} \sum_{j < j'} \zeta_j \zeta_{j'}
\]
where $\gamma$ runs over all virtual strands and virtual loops, starting from the chosen base points. 

\begin{prop}
    The above defined maps $x_i, x'_i, \bar{x}_k, \xi_i, \bar{\xi}_k, u$ give $X$ an $\mcA^w_{n, l}$-module structure.
\end{prop}

\begin{proof}
    Formally repeat the proof of \Cref{prop:homotopy-u}. 
\end{proof}

\begin{prop}
    For tangles $T_a \in \Diag(B_a)$ $(a = 1, 2)$, we have $[D(T_1, T_2)] = D([T_1], [T_2])$ as $\mcA^w_{n, l}$-modules. Here, the dg-module structures on $[T_1], [T_2]$ and $[D(T_1, T_2)]$ are given by choosing an admissible $ab$-coloring $\theta$ for $D(T_1, T_2)$ and restricting it to $T_1, T_2$. 
\end{prop}

\begin{proof}
    Obvious from the construction. 
\end{proof}

The above action of a $2$-input planar arc diagram can be generalized to a planar arc diagram $D$ with arbitrary number of inputs. In particular, for the $1$-input planar arc diagram $D$ for the braid closure operation, we have
\[
    D([\beta]) = [D(\beta)] = [\bbeta],
\]
thereby specializing to the reduction of an $\mcA^w_n$-module to a $\mcCA_l$-module described in \Cref{subsec:reduce-dga}. 

\begin{prop}
    For $a = 1, 2$, let $f_a\colon X_a \to Y_a$ be an $\mcA^{w_a}_{n_a, l_a}$-module homomorphism. The horizontal composition
    \[
        D(f_1, f_2) \colon D(X_1, X_2) \to D(Y_1, Y_2)
    \]
    is an $\mcA^w_{n, l}$-module homomorphism.  
\end{prop}

\subsection{Restriction to knots}
\label{subsec:restr-to-knots}

\begin{prop}
\label{prop:y-complex-knot}
    If $\beta$ closes up to a knot, then on the complex $[\bbeta] \in \Kob_{\mcCA_1}(\emptyset)$, we have 
    \[
        \bar{\xi} = 0, \quad 
        [d, \bar{u}] = 0.
    \]
\end{prop}

\begin{proof}
    In the sum of $\bar{\xi}$, since $\chi_c$ appears exactly twice with opposite signs, the sum
    \[
        \bar{\xi} = \sum_{i = 1}^n \xi_i
    \]
    vanishes. The second equation follows from the defining property of $\bar{u}$. 
\end{proof}

\begin{cor}
    If $\beta$ closes up to a knot, then the $y$-ified complex $y[\bbeta] = [\bbeta] \otimes R[y]$ has differential
    \[
        D = d.
    \]
    Therefore, $([\bbeta], d)$ is a subcomplex of $(y[\bbeta], D)$ and $\sle$ restricts to $\bar{u}$ on $[\bbeta]$. 
\end{cor}

\begin{thm}
\label{thm:e-inv-knot}
    If $\beta$ closes to a knot, then the $\bar{u}$-equivariant chain homotopy type of $[\bbeta]$ is an invariant of the knot. 
\end{thm}

\Cref{thm:e-inv-knot} follows from \Cref{prop:u-invariance-braid,prop:u-invariance-markov}. 

\begin{prop}
\label{prop:u-invariance-braid}
    Let $\beta, \beta'$ be two braids that are related by one of the braid relations. The $\sle$-equivariant homotopy equivalences between the corresponding $y$-complexes $y[\bbeta], y[\bbeta']$ given in \Cref{prop:y-cpx-cl-braid-rel1,prop:y-cpx-cl-braid-rel2,prop:y-cpx-braid-rel3} restrict to $\bar{u}$-equivariant homotopy equivalences between the non-$y$-ified complexes.
\end{prop}

\begin{proof}
    For the first and the second braid relations, recall from \Cref{prop:y-cpx-cl-braid-rel1,prop:y-cpx-cl-braid-rel2} that the homotopy equivalences 
    \[
    \begin{tikzcd}
        {[\bbeta]} 
            \arrow[r, "f", shift left] 
        & {[\bbeta']}
            \arrow[l, "g", shift left]
    \end{tikzcd}    
    \]
    are $\mcC^w_1$-linear, thus in particular strictly $\bar{u}$-equivariant. 
    
    For the third braid relation, recall from the proof of \Cref{prop:y-cpx-cl-braid-rel3} that the maps $F, H$ are $\mcA^w_n$-linear, but $G$ is not. There, we had
    \[
        G_1 = -G_2, \quad G_3 = 0
    \]
    so the corresponding $y$-chain map is given by $G$. Therefore, $\Psi, \Psi'$ directly give the homotopy equivalences between the $y$-complexes (and it coincides with the maps for the R3 move given in \cite[Theorem 1]{BarNatan:2005}.) Furthermore, we know that there is a homotopy $\tilde{\Psi}$, represented as a polynomial in $y$, giving
    \[
        [\sle, \Psi] = [D, \tilde{\Psi}].
    \]
    However, 
    \[
        [\sle, \Psi] = [\bar{u} + 2x \ddel{}{y}, \Psi] = [\bar{u}, \Psi]
    \]
    and
    \[
        [D, \tilde{\Psi}] = [d, \tilde{\Psi}]
    \]
    so the constant term of $\tilde{\Psi}$ gives a homotopy in the non-$y$-ified category.
\end{proof}

\begin{prop}
\label{prop:u-invariance-markov}
    Let $\beta, \beta'$ be two braids that are related by a Markov move. The $\sle$-equivariant homotopy equivalences between the corresponding $y$-complexes given in \Cref{prop:markov-inv1,prop:markov-inv2} restrict to $\bar{u}$-equivariant homotopy equivalences between the non-$y$-ified complexes. 
\end{prop}

\begin{proof}
    Obvious from the proofs of \Cref{prop:markov-inv1,prop:markov-inv2}. 
\end{proof}

The following is also worth noting, which is immediate from \Cref{prop:ab-coloring-dependence}.

\begin{prop}
    On $[\bbeta]$, $\bar{u}$ is independent of the choice of an admissible coloring $\theta$.
\end{prop}

The above results justify the following definition.

\begin{defn}
    For a knot $K$ represented by a braid $\beta$, we rewrite $\sle$ for the endomorphism $\bar{u}$ on $[\bbeta]$, and call it the \textit{$\sle$-operator} on $[\bbeta]$. By applying the TQFT $\mcF$, we obtain the \textit{$\sle$-operator} on $\CKh_{U(2)}(\bbeta)$, and the \textit{$\sle$-operator} on $\Kh_{U(2)}(K)$. By specializing $h, t \in R_{U(2)}$ before taking the homology, we obtain the $\sle$-operator for other variants. 
\end{defn}

Now, let us focus on the ordinary Khovanov homology $\Kh(K; k)$ of a knot $K$ over a field $k$. Since $\sle$ has degree $(-2, 4)$, it preserves the \textit{$\delta$-grading}
\[
    \gr_\delta = 2\gr_h + \gr_q.
\]
Recall that a \textit{$\Kh$-thin} knot $K$ has unreduced Khovanov homology $\Kh(K)$ supported on two diagonals $\delta = -s(K) \pm 1$, and reduced Khovanov homology $\rKh(K)$ supported on a single diagonal $\delta = -s(K)$, where $s(K)$ is the \textit{Rasmussen invariant} of $K$.\footnote{The minus sign is due to the convention of the quantum grading.} Here, we use $(\gr_\delta, \gr_q)$ in place of $(\gr_h, \gr_q)$ for the bigrading of $\Kh(K)$. 

Let $\Kh^{\delta, *}(K)$ denote a $\delta$-slice of $\Kh(K)$. There is a bounded sequence of the form
\[
\begin{tikzcd}
    \cdots \arrow[rr, "\sle", bend left] & 
    {\Kh^{\delta, j-4}(K)} \arrow[rr, "\sle"', bend right] & 
    {\Kh^{\delta, j-2}(K)} \arrow[rr, "\sle", bend left] & 
    {\Kh^{\delta, j}(K)} \arrow[rr, "\sle"', bend right] & 
    {\Kh^{\delta, j+2}(K)} \arrow[rr, "\sle", bend left]& 
    {\Kh^{\delta, j+4}(K)} & 
    \cdots 
\end{tikzcd}
\]
We regard $\Kh^{\delta, *}(K)$ as a finitely generated graded $k[e]$-module, where $e$ acts as the nilpotent endomorphism $\sle$. Then, from the structure theorem of finitely generated graded $k[e]$-modules, there is a decomposition
\begin{equation}
\label{eqn:kh-e-string-decomp}
    \Kh^{\delta, *}(K; k) \isom \bigoplus_i q^{j_i} (k[e] / (e^{l_i}))
\end{equation}
which we call an \textit{$e$-string decomposition} of $\Kh(K)$. Each summand $q^{j_i} k[e] / (e^{l_i})$ corresponds to a string of dimension $l_i$ starting from $q$-grading $j_i$:
\[
\begin{tikzcd}
    q^{j_i} k \arrow[r, "\sle"] & q^{j_i + 4} k \arrow[r, "\sle"] & \cdots \arrow[r, "\sle"] & q^{j_i + 4(l_i - 1)} k.
\end{tikzcd}
\]

Now suppose $k$ is a field of characteristic $0$. By a standard argument, each $e$-string can be given a unique $\sl_2$-action extending the action of $e$. However, since the decomposition itself is not unique in general, $\sle$ does not uniquely extend to an $\sl_2$-action on $\Kh(K)$. Nonetheless, in simple cases, such unique extension is possible. 

\begin{prop}
\label{prop:e-extension}
    Let $K$ be a knot and $k$ a field of characteristic $0$. If $\Kh^{i, j}(K; k)$ is at most $1$-dimensional for each bidegree $(i, j)$, the $\sle$-operator uniquely extends to an $\sl_2$-action on $\Kh(K; k)$.
\end{prop}

We shall see in \Cref{subsec:torus-knots} that the $\QQ$-Khovanov homology of the $(2, 2k+1)$-torus knot gives an example that satisfies the condition of \Cref{prop:e-extension}. 
In \Cref{sec:direct-computation}, we shall give an explicit algorithm that computes the $k[e]$-module structure of $\Kh(K)$. 

\subsection{Mirror duality}

Here, we state the duality between the $\sle$-operators for a knot $K$ and its mirror $K^*$. First, let us recall the chain level description of the canonical isomorphism 
\[
    \CKh_{U(2)}(D^*) \isom \CKh_{U(2)}(D)^*.
\]

For brevity, we write $R = R_{U(2)}$ and $A = A_{U(2)}$. The Frobenius algebra operations of $A$ are denoted $(m, \iota, \Delta, \epsilon)$ as in \Cref{subsubsec:preliminary-khovanov}. Consider the non-degenerate pairing $B$ on $A$ defined by
\[
    B := \epsilon \circ m\colon A \otimes A \to R.
\]
Let $\dual$ denote the associated isomorphism
\[
    \dual\colon A \to A^*,\quad 
    x \mapsto B(x, -)
\]
where $A^* := \Hom_{R}(A, R)$ is the dual $R$-module of $A$. Put
\[
    \check{1} := \dual(1),\quad
    \check{X} := \dual(X)
\]
and call $\{\check{1}, \check{X}\}$ the \textit{standard basis} of $A^*$. These elements can be described as cobordisms:
\begin{center}
    \input{tikzpictures/duality-map}
\end{center}
With $Y = X - h$, one can see (pictorially or algebraically) that $\{1, X\}$ and $\{Y, 1\}$ are mutually dual with respect to the non-degenerate pairing $B$. Thus the standard basis $\{\check{1}, \check{X}\}$ for $ A^*$ is precisely the algebraic dual of $\{Y, 1\}$ with respect to $B$. 

The dual module $A^*$ admits the \textit{dual Frobenius algebra structure} with multiplication $\Delta^*$, unit $\epsilon^*$, comultiplication $m^*$, and counit $\iota^*$. 
One sees that $\dual$ is a Frobenius algebra isomorphism, since the operations of $A^*$ described with the basis $\{\check{1}, \check{X} \}$ are nothing but those of $A$ described with the basis $\{1, X\}$ turned around.

\begin{prop}
\label{prop:A-adj}
    Let $f\colon A \to A$ be an endomorphism given by a cobordism, and $\overline{f}$ be the reversal of $f$. Then $f$ and $\overline{f}$ are adjoint with respect to the pairing $B$. Equivalently, $\dual \circ f = \bar{f}^* \circ \dual$. 
\end{prop}

\begin{proof}
    Obvious from the following diagram, which shows $B(fx, y) = B(x, \bar{f} y)$.
    \vspace{1em}
    \begin{center}
        \input{tikzpictures/dual-adj}
    \end{center}
\end{proof}

Let $D$ be a link diagram, and $D^*$ its mirror. Here, we orient $D^*$ oppositely to $D$, and also order the crossings of $D^*$ oppositely to $D$. Furthermore, take base points of $D$ and $D^*$ that they correspond accordingly. The duality isomorphism $\dual$ extends to a chain isomorphism
\[
    \dual\colon \CKh_{U(2)}(D^*) \isom \CKh_{U(2)}(D)^*,\quad 
    x_1 \otimes \cdots \otimes x_r \mapsto 
    \dual(x_1) \otimes \cdots \otimes \dual(x_r).
\]
where $\CKh_{U(2)}(D)^*$ is the dual chain complex of $\CKh_{U(2)}(D)$, whose bigrading (over $\ZZ$) is given by 
\[
    (\CKh_{U(2)}(D)^*)^{i, j} = (\CKh^{-i, -j}_{h, t}(D))^*.  
\]
This gives a duality pairing 
\[
    \langle \cdot, \cdot \rangle \colon 
    \CKh_{U(2)}(D) \otimes \CKh_{U(2)}(D^*) \to R
\]
by 
\[
    \langle x, y \rangle := \langle x, \dual(y) \rangle_{\text{alg}}
\]
where $\langle \cdot, \cdot \rangle_{\text{alg}}$ denotes the standard pairing between $\CKh_{U(2)}(D)$ and its dual module. 

\begin{prop}
\label{prop:mirror-adjoint}
    Let $f$ be an endomorphism of $\CKh_{U(2)}(D)$ given by a matrix of cobordisms, and $\overline{f}$ be the endomorphism of $\CKh_{U(2)}(D^*)$ given by the entry-wise reversal of $f$. Then $f$ and $\overline{f}$ are adjoint with respect to the pairing $\langle \cdot, \cdot \rangle$. 
\end{prop}

\begin{proof}
    Immediate from \Cref{prop:A-adj}. 
\end{proof}

From \Cref{prop:mirror-adjoint}, one easily sees that
\[
    \langle dx, y \rangle = \langle x, dy \rangle
\]
holds. A similar adjunction formula holds for the $\sle$-operator. 

\begin{prop}
\label{prop:e-adjoint}
\[
    \langle \sle x, y \rangle = \langle x, \sle y \rangle.
\]    
\end{prop}

\begin{proof}
    Write 
    \[
        \sle(D) = \sum_{j < j'} \chi_j \chi_{j'}
    \]
    for $D$. Then
    \[
        \sle(D^*) = \sum_{j < j'} \overline{\chi}_{j'} \overline{\chi}_j
    \]
    for $D^*$, where $\overline{\chi}_j$ is the reversal of $\chi_j$. The result is immediate from \Cref{prop:mirror-adjoint}. 
\end{proof}

\begin{cor}
    Let $K$ be a knot and $k$ a field. Suppose we have an $e$-string decomposition
    \[
        \Kh(K; k) \isom \bigoplus_{(i, j, l) \in I} \delta^i q^j (k[e] / (e^l))
    \]
    for some multiset $I$. Then, an $e$-string decomposition of $\Kh(K^*; k)$ is given by 
    \[
        \Kh(K^*; k) \isom \bigoplus_{(i, j, l) \in I} \delta^{-i} q^{-j - 4(l - 1)} (k[e] / (e^l)).
    \]
    Furthermore, suppose we choose a homogeneous generating set $\{v_{i, j, l}\}$ of $\Kh(K; k)$ corresponding to the above $e$-string decomposition. Then there is a unique homogeneous generating set $\{w_{i, j, l}\}$ of $\Kh(K^*; k)$ such that 
    \begin{align*}
        \langle e^a v_{i, j, l},\ e^b w_{i, j, l} \rangle 
        &= \delta_{a + b, l - 1},\\
        \langle e^a v_{i, j, l},\ e^b w_{i', j', l'} \rangle 
        &= 0 \ ((i, j, l) \neq (i', j', l')). 
    \end{align*}
\end{cor}


%% file: tikzpictures/loc-relation-dot.tex
\tikzset{every picture/.style={line width=0.75pt}} 

\begin{tikzpicture}[x=0.75pt,y=0.75pt,yscale=-1,xscale=1]

\draw  [draw opacity=0] (157.31,28.4) .. controls (156.57,31) and (148.66,33.04) .. (139.01,33.04) .. controls (128.87,33.04) and (120.65,30.78) .. (120.65,28) .. controls (120.65,27.83) and (120.68,27.66) .. (120.74,27.5) -- (139.01,28) -- cycle ; \draw  [color={rgb, 255:red, 128; green, 128; blue, 128 }  ,draw opacity=1 ] (157.31,28.4) .. controls (156.57,31) and (148.66,33.04) .. (139.01,33.04) .. controls (128.87,33.04) and (120.65,30.78) .. (120.65,28) .. controls (120.65,27.83) and (120.68,27.66) .. (120.74,27.5) ;  
\draw  [draw opacity=0][dash pattern={on 0.84pt off 2.51pt}] (121.1,26.75) .. controls (122.57,24.37) and (130.12,22.56) .. (139.21,22.56) .. controls (149.2,22.56) and (157.33,24.75) .. (157.56,27.48) -- (139.21,27.6) -- cycle ; \draw  [color={rgb, 255:red, 128; green, 128; blue, 128 }  ,draw opacity=1 ][dash pattern={on 0.84pt off 2.51pt}] (121.1,26.75) .. controls (122.57,24.37) and (130.12,22.56) .. (139.21,22.56) .. controls (149.2,22.56) and (157.33,24.75) .. (157.56,27.48) ;  
\draw   (120.51,28) .. controls (120.51,17.78) and (128.79,9.5) .. (139.01,9.5) .. controls (149.23,9.5) and (157.51,17.78) .. (157.51,28) .. controls (157.51,38.22) and (149.23,46.5) .. (139.01,46.5) .. controls (128.79,46.5) and (120.51,38.22) .. (120.51,28) -- cycle ;

\draw  [draw opacity=0] (342.81,29.4) .. controls (342.07,32) and (334.16,34.04) .. (324.51,34.04) .. controls (314.37,34.04) and (306.15,31.78) .. (306.15,29) .. controls (306.15,28.83) and (306.18,28.66) .. (306.24,28.5) -- (324.51,29) -- cycle ; \draw  [color={rgb, 255:red, 128; green, 128; blue, 128 }  ,draw opacity=1 ] (342.81,29.4) .. controls (342.07,32) and (334.16,34.04) .. (324.51,34.04) .. controls (314.37,34.04) and (306.15,31.78) .. (306.15,29) .. controls (306.15,28.83) and (306.18,28.66) .. (306.24,28.5) ;  
\draw  [draw opacity=0][dash pattern={on 0.84pt off 2.51pt}] (306.6,27.75) .. controls (308.07,25.37) and (315.62,23.56) .. (324.71,23.56) .. controls (334.7,23.56) and (342.83,25.75) .. (343.06,28.48) -- (324.71,28.6) -- cycle ; \draw  [color={rgb, 255:red, 128; green, 128; blue, 128 }  ,draw opacity=1 ][dash pattern={on 0.84pt off 2.51pt}] (306.6,27.75) .. controls (308.07,25.37) and (315.62,23.56) .. (324.71,23.56) .. controls (334.7,23.56) and (342.83,25.75) .. (343.06,28.48) ;  
\draw   (306.01,29) .. controls (306.01,18.78) and (314.29,10.5) .. (324.51,10.5) .. controls (334.73,10.5) and (343.01,18.78) .. (343.01,29) .. controls (343.01,39.22) and (334.73,47.5) .. (324.51,47.5) .. controls (314.29,47.5) and (306.01,39.22) .. (306.01,29) -- cycle ;
\draw  [fill={rgb, 255:red, 0; green, 0; blue, 0 }  ,fill opacity=1 ] (322.34,19.72) .. controls (322.34,18.41) and (323.4,17.35) .. (324.71,17.35) .. controls (326.02,17.35) and (327.08,18.41) .. (327.08,19.72) .. controls (327.08,21.03) and (326.02,22.09) .. (324.71,22.09) .. controls (323.4,22.09) and (322.34,21.03) .. (322.34,19.72) -- cycle ;

\draw  [color={rgb, 255:red, 0; green, 0; blue, 0 }  ,draw opacity=1 ] (49.43,84.19) .. controls (49.43,76.48) and (52,70.23) .. (55.18,70.23) .. controls (58.36,70.23) and (60.93,76.48) .. (60.93,84.19) .. controls (60.93,91.89) and (58.36,98.14) .. (55.18,98.14) .. controls (52,98.14) and (49.43,91.89) .. (49.43,84.19) -- cycle ;
\draw [color={rgb, 255:red, 0; green, 0; blue, 0 }  ,draw opacity=1 ]   (55.18,98.14) -- (114.44,98.14) ;
\draw [color={rgb, 255:red, 0; green, 0; blue, 0 }  ,draw opacity=1 ]   (55.18,70.23) -- (114.44,70.23) ;
\draw  [draw opacity=0][dash pattern={on 0.84pt off 2.51pt}] (112.15,97.93) .. controls (109.93,96.27) and (108.31,90.84) .. (108.31,84.39) .. controls (108.31,77.95) and (109.93,72.52) .. (112.15,70.86) -- (113.58,84.39) -- cycle ; \draw  [dash pattern={on 0.84pt off 2.51pt}] (112.15,97.93) .. controls (109.93,96.27) and (108.31,90.84) .. (108.31,84.39) .. controls (108.31,77.95) and (109.93,72.52) .. (112.15,70.86) ;  
\draw  [draw opacity=0] (113,70.41) .. controls (113.19,70.36) and (113.38,70.33) .. (113.58,70.33) .. controls (116.49,70.33) and (118.85,76.63) .. (118.85,84.39) .. controls (118.85,92.15) and (116.49,98.45) .. (113.58,98.45) .. controls (113.38,98.45) and (113.19,98.42) .. (113,98.37) -- (113.58,84.39) -- cycle ; \draw   (113,70.41) .. controls (113.19,70.36) and (113.38,70.33) .. (113.58,70.33) .. controls (116.49,70.33) and (118.85,76.63) .. (118.85,84.39) .. controls (118.85,92.15) and (116.49,98.45) .. (113.58,98.45) .. controls (113.38,98.45) and (113.19,98.42) .. (113,98.37) ;  

\draw  [color={rgb, 255:red, 0; green, 0; blue, 0 }  ,draw opacity=1 ] (156.6,84.61) .. controls (156.6,76.9) and (159.18,70.66) .. (162.35,70.66) .. controls (165.53,70.66) and (168.11,76.9) .. (168.11,84.61) .. controls (168.11,92.31) and (165.53,98.56) .. (162.35,98.56) .. controls (159.18,98.56) and (156.6,92.31) .. (156.6,84.61) -- cycle ;
\draw  [draw opacity=0] (161.77,70.56) .. controls (161.97,70.55) and (162.16,70.55) .. (162.35,70.55) .. controls (171.48,70.55) and (178.88,76.84) .. (178.88,84.61) .. controls (178.88,92.37) and (171.48,98.67) .. (162.35,98.67) .. controls (162.16,98.67) and (161.97,98.66) .. (161.77,98.66) -- (162.35,84.61) -- cycle ; \draw   (161.77,70.56) .. controls (161.97,70.55) and (162.16,70.55) .. (162.35,70.55) .. controls (171.48,70.55) and (178.88,76.84) .. (178.88,84.61) .. controls (178.88,92.37) and (171.48,98.67) .. (162.35,98.67) .. controls (162.16,98.67) and (161.97,98.66) .. (161.77,98.66) ;  
\draw  [draw opacity=0][dash pattern={on 0.84pt off 2.51pt}] (214.17,97.61) .. controls (211.95,95.96) and (210.32,90.52) .. (210.32,84.07) .. controls (210.32,77.63) and (211.95,72.2) .. (214.16,70.54) -- (215.59,84.07) -- cycle ; \draw  [dash pattern={on 0.84pt off 2.51pt}] (214.17,97.61) .. controls (211.95,95.96) and (210.32,90.52) .. (210.32,84.07) .. controls (210.32,77.63) and (211.95,72.2) .. (214.16,70.54) ;  
\draw  [draw opacity=0] (215.02,70.1) .. controls (215.21,70.04) and (215.4,70.02) .. (215.59,70.02) .. controls (218.5,70.02) and (220.86,76.31) .. (220.86,84.07) .. controls (220.86,91.84) and (218.5,98.13) .. (215.59,98.13) .. controls (215.4,98.13) and (215.21,98.11) .. (215.02,98.05) -- (215.59,84.07) -- cycle ; \draw   (215.02,70.1) .. controls (215.21,70.04) and (215.4,70.02) .. (215.59,70.02) .. controls (218.5,70.02) and (220.86,76.31) .. (220.86,84.07) .. controls (220.86,91.84) and (218.5,98.13) .. (215.59,98.13) .. controls (215.4,98.13) and (215.21,98.11) .. (215.02,98.05) ;  
\draw  [draw opacity=0] (216.17,70.02) .. controls (215.98,70.02) and (215.79,70.02) .. (215.59,70.02) .. controls (206.47,70.02) and (199.07,76.31) .. (199.07,84.07) .. controls (199.07,91.84) and (206.47,98.13) .. (215.59,98.13) .. controls (215.79,98.13) and (215.98,98.13) .. (216.17,98.12) -- (215.59,84.07) -- cycle ; \draw   (216.17,70.02) .. controls (215.98,70.02) and (215.79,70.02) .. (215.59,70.02) .. controls (206.47,70.02) and (199.07,76.31) .. (199.07,84.07) .. controls (199.07,91.84) and (206.47,98.13) .. (215.59,98.13) .. controls (215.79,98.13) and (215.98,98.13) .. (216.17,98.12) ;  

\draw  [color={rgb, 255:red, 0; green, 0; blue, 0 }  ,draw opacity=1 ] (255.13,84.61) .. controls (255.13,76.9) and (257.71,70.66) .. (260.89,70.66) .. controls (264.07,70.66) and (266.64,76.9) .. (266.64,84.61) .. controls (266.64,92.31) and (264.07,98.56) .. (260.89,98.56) .. controls (257.71,98.56) and (255.13,92.31) .. (255.13,84.61) -- cycle ;
\draw  [draw opacity=0] (260.31,70.56) .. controls (260.5,70.55) and (260.69,70.55) .. (260.89,70.55) .. controls (270.01,70.55) and (277.41,76.84) .. (277.41,84.61) .. controls (277.41,92.37) and (270.01,98.67) .. (260.89,98.67) .. controls (260.69,98.67) and (260.5,98.66) .. (260.31,98.66) -- (260.89,84.61) -- cycle ; \draw   (260.31,70.56) .. controls (260.5,70.55) and (260.69,70.55) .. (260.89,70.55) .. controls (270.01,70.55) and (277.41,76.84) .. (277.41,84.61) .. controls (277.41,92.37) and (270.01,98.67) .. (260.89,98.67) .. controls (260.69,98.67) and (260.5,98.66) .. (260.31,98.66) ;  
\draw  [draw opacity=0][dash pattern={on 0.84pt off 2.51pt}] (312.7,97.61) .. controls (310.49,95.96) and (308.86,90.52) .. (308.86,84.07) .. controls (308.86,77.63) and (310.48,72.2) .. (312.7,70.54) -- (314.13,84.07) -- cycle ; \draw  [dash pattern={on 0.84pt off 2.51pt}] (312.7,97.61) .. controls (310.49,95.96) and (308.86,90.52) .. (308.86,84.07) .. controls (308.86,77.63) and (310.48,72.2) .. (312.7,70.54) ;  
\draw  [draw opacity=0] (313.55,70.1) .. controls (313.74,70.04) and (313.93,70.02) .. (314.13,70.02) .. controls (317.04,70.02) and (319.4,76.31) .. (319.4,84.07) .. controls (319.4,91.84) and (317.04,98.13) .. (314.13,98.13) .. controls (313.93,98.13) and (313.74,98.11) .. (313.55,98.05) -- (314.13,84.07) -- cycle ; \draw   (313.55,70.1) .. controls (313.74,70.04) and (313.93,70.02) .. (314.13,70.02) .. controls (317.04,70.02) and (319.4,76.31) .. (319.4,84.07) .. controls (319.4,91.84) and (317.04,98.13) .. (314.13,98.13) .. controls (313.93,98.13) and (313.74,98.11) .. (313.55,98.05) ;  
\draw  [draw opacity=0] (314.71,70.02) .. controls (314.51,70.02) and (314.32,70.02) .. (314.13,70.02) .. controls (305,70.02) and (297.61,76.31) .. (297.61,84.07) .. controls (297.61,91.84) and (305,98.13) .. (314.13,98.13) .. controls (314.32,98.13) and (314.51,98.13) .. (314.71,98.12) -- (314.13,84.07) -- cycle ; \draw   (314.71,70.02) .. controls (314.51,70.02) and (314.32,70.02) .. (314.13,70.02) .. controls (305,70.02) and (297.61,76.31) .. (297.61,84.07) .. controls (297.61,91.84) and (305,98.13) .. (314.13,98.13) .. controls (314.32,98.13) and (314.51,98.13) .. (314.71,98.12) ;  

\draw  [color={rgb, 255:red, 0; green, 0; blue, 0 }  ,draw opacity=1 ] (353.66,84.61) .. controls (353.66,76.9) and (356.24,70.66) .. (359.42,70.66) .. controls (362.59,70.66) and (365.17,76.9) .. (365.17,84.61) .. controls (365.17,92.31) and (362.59,98.56) .. (359.42,98.56) .. controls (356.24,98.56) and (353.66,92.31) .. (353.66,84.61) -- cycle ;
\draw  [draw opacity=0] (358.84,70.57) .. controls (359.03,70.56) and (359.22,70.55) .. (359.42,70.55) .. controls (366.1,70.55) and (371.52,76.84) .. (371.52,84.61) .. controls (371.52,92.37) and (366.1,98.67) .. (359.42,98.67) .. controls (359.22,98.67) and (359.03,98.66) .. (358.84,98.65) -- (359.42,84.61) -- cycle ; \draw   (358.84,70.57) .. controls (359.03,70.56) and (359.22,70.55) .. (359.42,70.55) .. controls (366.1,70.55) and (371.52,76.84) .. (371.52,84.61) .. controls (371.52,92.37) and (366.1,98.67) .. (359.42,98.67) .. controls (359.22,98.67) and (359.03,98.66) .. (358.84,98.65) ;  
\draw  [draw opacity=0][dash pattern={on 0.84pt off 2.51pt}] (411.23,97.61) .. controls (409.01,95.96) and (407.39,90.52) .. (407.39,84.07) .. controls (407.39,77.63) and (409.01,72.2) .. (411.23,70.54) -- (412.66,84.07) -- cycle ; \draw  [dash pattern={on 0.84pt off 2.51pt}] (411.23,97.61) .. controls (409.01,95.96) and (407.39,90.52) .. (407.39,84.07) .. controls (407.39,77.63) and (409.01,72.2) .. (411.23,70.54) ;  
\draw  [draw opacity=0] (412.08,70.1) .. controls (412.27,70.04) and (412.46,70.02) .. (412.66,70.02) .. controls (415.57,70.02) and (417.93,76.31) .. (417.93,84.07) .. controls (417.93,91.84) and (415.57,98.13) .. (412.66,98.13) .. controls (412.46,98.13) and (412.27,98.11) .. (412.08,98.05) -- (412.66,84.07) -- cycle ; \draw   (412.08,70.1) .. controls (412.27,70.04) and (412.46,70.02) .. (412.66,70.02) .. controls (415.57,70.02) and (417.93,76.31) .. (417.93,84.07) .. controls (417.93,91.84) and (415.57,98.13) .. (412.66,98.13) .. controls (412.46,98.13) and (412.27,98.11) .. (412.08,98.05) ;  
\draw  [draw opacity=0] (413.23,70.03) .. controls (413.04,70.02) and (412.85,70.02) .. (412.66,70.02) .. controls (406.21,70.02) and (400.98,76.31) .. (400.98,84.07) .. controls (400.98,91.84) and (406.21,98.13) .. (412.66,98.13) .. controls (412.85,98.13) and (413.04,98.13) .. (413.23,98.12) -- (412.66,84.07) -- cycle ; \draw   (413.23,70.03) .. controls (413.04,70.02) and (412.85,70.02) .. (412.66,70.02) .. controls (406.21,70.02) and (400.98,76.31) .. (400.98,84.07) .. controls (400.98,91.84) and (406.21,98.13) .. (412.66,98.13) .. controls (412.85,98.13) and (413.04,98.13) .. (413.23,98.12) ;  
\draw  [draw opacity=0] (397.87,84.26) .. controls (397.38,85.96) and (392.21,87.3) .. (385.89,87.3) .. controls (379.26,87.3) and (373.88,85.82) .. (373.88,84) .. controls (373.88,83.89) and (373.9,83.78) .. (373.94,83.67) -- (385.89,84) -- cycle ; \draw  [color={rgb, 255:red, 128; green, 128; blue, 128 }  ,draw opacity=1 ] (397.87,84.26) .. controls (397.38,85.96) and (392.21,87.3) .. (385.89,87.3) .. controls (379.26,87.3) and (373.88,85.82) .. (373.88,84) .. controls (373.88,83.89) and (373.9,83.78) .. (373.94,83.67) ;  
\draw  [draw opacity=0][dash pattern={on 0.84pt off 2.51pt}] (374.17,83.18) .. controls (375.13,81.62) and (380.07,80.44) .. (386.02,80.44) .. controls (392.56,80.44) and (397.88,81.87) .. (398.03,83.66) -- (386.02,83.74) -- cycle ; \draw  [color={rgb, 255:red, 128; green, 128; blue, 128 }  ,draw opacity=1 ][dash pattern={on 0.84pt off 2.51pt}] (374.17,83.18) .. controls (375.13,81.62) and (380.07,80.44) .. (386.02,80.44) .. controls (392.56,80.44) and (397.88,81.87) .. (398.03,83.66) ;  
\draw   (373.79,84) .. controls (373.79,77.31) and (379.21,71.89) .. (385.89,71.89) .. controls (392.58,71.89) and (398,77.31) .. (398,84) .. controls (398,90.69) and (392.58,96.11) .. (385.89,96.11) .. controls (379.21,96.11) and (373.79,90.69) .. (373.79,84) -- cycle ;
\draw  [fill={rgb, 255:red, 0; green, 0; blue, 0 }  ,fill opacity=1 ] (381.52,77.93) .. controls (381.52,77.07) and (382.22,76.37) .. (383.07,76.37) .. controls (383.93,76.37) and (384.63,77.07) .. (384.63,77.93) .. controls (384.63,78.78) and (383.93,79.48) .. (383.07,79.48) .. controls (382.22,79.48) and (381.52,78.78) .. (381.52,77.93) -- cycle ;
\draw  [fill={rgb, 255:red, 0; green, 0; blue, 0 }  ,fill opacity=1 ] (387.42,77.93) .. controls (387.42,77.07) and (388.11,76.37) .. (388.97,76.37) .. controls (389.82,76.37) and (390.52,77.07) .. (390.52,77.93) .. controls (390.52,78.78) and (389.82,79.48) .. (388.97,79.48) .. controls (388.11,79.48) and (387.42,78.78) .. (387.42,77.93) -- cycle ;
\draw  [fill={rgb, 255:red, 0; green, 0; blue, 0 }  ,fill opacity=1 ] (172.32,84.56) .. controls (172.32,83.7) and (173.01,83) .. (173.87,83) .. controls (174.72,83) and (175.42,83.7) .. (175.42,84.56) .. controls (175.42,85.41) and (174.72,86.11) .. (173.87,86.11) .. controls (173.01,86.11) and (172.32,85.41) .. (172.32,84.56) -- cycle ;
\draw  [fill={rgb, 255:red, 0; green, 0; blue, 0 }  ,fill opacity=1 ] (300.49,83.82) .. controls (300.49,82.96) and (301.19,82.27) .. (302.04,82.27) .. controls (302.9,82.27) and (303.6,82.96) .. (303.6,83.82) .. controls (303.6,84.68) and (302.9,85.37) .. (302.04,85.37) .. controls (301.19,85.37) and (300.49,84.68) .. (300.49,83.82) -- cycle ;

\draw (76,18) node [anchor=north west][inner sep=0.75pt]   [align=left] {(S)};
\draw (164,18.4) node [anchor=north west][inner sep=0.75pt]    {$=0$};
\draw (3.5,72.5) node [anchor=north west][inner sep=0.75pt]   [align=left] {(NC)};
\draw (261.5,19) node [anchor=north west][inner sep=0.75pt]   [align=left] {(S$\displaystyle _{\bullet }$)};
\draw (349.5,19.4) node [anchor=north west][inner sep=0.75pt]    {$=1$};
\draw (125.28,76.42) node [anchor=north west][inner sep=0.75pt]    {$=$};
\draw (232.07,76.42) node [anchor=north west][inner sep=0.75pt]    {$+$};
\draw (330.78,76.42) node [anchor=north west][inner sep=0.75pt]    {$-$};

\end{tikzpicture}

%% file: tikzpictures/cob-ht.tex
\tikzset{every picture/.style={line width=0.75pt}} 

\begin{tikzpicture}[x=0.75pt,y=0.75pt,yscale=-1,xscale=1]

\draw  [draw opacity=0] (89.94,40.75) .. controls (89.29,43) and (82.45,44.76) .. (74.11,44.76) .. controls (65.35,44.76) and (58.24,42.81) .. (58.24,40.41) .. controls (58.24,40.26) and (58.27,40.11) .. (58.32,39.97) -- (74.11,40.41) -- cycle ; \draw  [color={rgb, 255:red, 128; green, 128; blue, 128 }  ,draw opacity=1 ] (89.94,40.75) .. controls (89.29,43) and (82.45,44.76) .. (74.11,44.76) .. controls (65.35,44.76) and (58.24,42.81) .. (58.24,40.41) .. controls (58.24,40.26) and (58.27,40.11) .. (58.32,39.97) ;  
\draw  [draw opacity=0][dash pattern={on 0.84pt off 2.51pt}] (58.63,39.33) .. controls (59.9,37.27) and (66.43,35.7) .. (74.28,35.7) .. controls (82.93,35.7) and (89.96,37.6) .. (90.15,39.96) -- (74.28,40.06) -- cycle ; \draw  [color={rgb, 255:red, 128; green, 128; blue, 128 }  ,draw opacity=1 ][dash pattern={on 0.84pt off 2.51pt}] (58.63,39.33) .. controls (59.9,37.27) and (66.43,35.7) .. (74.28,35.7) .. controls (82.93,35.7) and (89.96,37.6) .. (90.15,39.96) ;  
\draw   (58.12,40.41) .. controls (58.12,31.57) and (65.28,24.41) .. (74.11,24.41) .. controls (82.95,24.41) and (90.11,31.57) .. (90.11,40.41) .. controls (90.11,49.24) and (82.95,56.4) .. (74.11,56.4) .. controls (65.28,56.4) and (58.12,49.24) .. (58.12,40.41) -- cycle ;
\draw  [fill={rgb, 255:red, 0; green, 0; blue, 0 }  ,fill opacity=1 ] (67.63,31.68) .. controls (67.63,30.55) and (68.55,29.63) .. (69.68,29.63) .. controls (70.82,29.63) and (71.73,30.55) .. (71.73,31.68) .. controls (71.73,32.81) and (70.82,33.73) .. (69.68,33.73) .. controls (68.55,33.73) and (67.63,32.81) .. (67.63,31.68) -- cycle ;
\draw  [fill={rgb, 255:red, 0; green, 0; blue, 0 }  ,fill opacity=1 ] (77.23,31.68) .. controls (77.23,30.55) and (78.15,29.63) .. (79.28,29.63) .. controls (80.42,29.63) and (81.33,30.55) .. (81.33,31.68) .. controls (81.33,32.81) and (80.42,33.73) .. (79.28,33.73) .. controls (78.15,33.73) and (77.23,32.81) .. (77.23,31.68) -- cycle ;
\draw  [draw opacity=0] (249.56,23.31) .. controls (249,25.25) and (243.11,26.77) .. (235.92,26.77) .. controls (228.36,26.77) and (222.23,25.09) .. (222.23,23.01) .. controls (222.23,22.88) and (222.26,22.76) .. (222.3,22.64) -- (235.92,23.01) -- cycle ; \draw  [color={rgb, 255:red, 128; green, 128; blue, 128 }  ,draw opacity=1 ] (249.56,23.31) .. controls (249,25.25) and (243.11,26.77) .. (235.92,26.77) .. controls (228.36,26.77) and (222.23,25.09) .. (222.23,23.01) .. controls (222.23,22.88) and (222.26,22.76) .. (222.3,22.64) ;  
\draw  [draw opacity=0][dash pattern={on 0.84pt off 2.51pt}] (222.57,22.08) .. controls (223.67,20.31) and (229.29,18.95) .. (236.06,18.95) .. controls (243.51,18.95) and (249.57,20.59) .. (249.75,22.62) -- (236.06,22.71) -- cycle ; \draw  [color={rgb, 255:red, 128; green, 128; blue, 128 }  ,draw opacity=1 ][dash pattern={on 0.84pt off 2.51pt}] (222.57,22.08) .. controls (223.67,20.31) and (229.29,18.95) .. (236.06,18.95) .. controls (243.51,18.95) and (249.57,20.59) .. (249.75,22.62) ;  
\draw   (222.13,23.01) .. controls (222.13,15.4) and (228.3,9.22) .. (235.92,9.22) .. controls (243.53,9.22) and (249.71,15.4) .. (249.71,23.01) .. controls (249.71,30.63) and (243.53,36.8) .. (235.92,36.8) .. controls (228.3,36.8) and (222.13,30.63) .. (222.13,23.01) -- cycle ;
\draw  [fill={rgb, 255:red, 0; green, 0; blue, 0 }  ,fill opacity=1 ] (230.33,15.49) .. controls (230.33,14.51) and (231.12,13.72) .. (232.1,13.72) .. controls (233.07,13.72) and (233.87,14.51) .. (233.87,15.49) .. controls (233.87,16.47) and (233.07,17.26) .. (232.1,17.26) .. controls (231.12,17.26) and (230.33,16.47) .. (230.33,15.49) -- cycle ;
\draw  [fill={rgb, 255:red, 0; green, 0; blue, 0 }  ,fill opacity=1 ] (238.61,15.49) .. controls (238.61,14.51) and (239.4,13.72) .. (240.37,13.72) .. controls (241.35,13.72) and (242.14,14.51) .. (242.14,15.49) .. controls (242.14,16.47) and (241.35,17.26) .. (240.37,17.26) .. controls (239.4,17.26) and (238.61,16.47) .. (238.61,15.49) -- cycle ;

\draw  [draw opacity=0] (249.56,60.11) .. controls (249,62.05) and (243.11,63.57) .. (235.92,63.57) .. controls (228.36,63.57) and (222.23,61.89) .. (222.23,59.81) .. controls (222.23,59.68) and (222.26,59.56) .. (222.3,59.44) -- (235.92,59.81) -- cycle ; \draw  [color={rgb, 255:red, 128; green, 128; blue, 128 }  ,draw opacity=1 ] (249.56,60.11) .. controls (249,62.05) and (243.11,63.57) .. (235.92,63.57) .. controls (228.36,63.57) and (222.23,61.89) .. (222.23,59.81) .. controls (222.23,59.68) and (222.26,59.56) .. (222.3,59.44) ;  
\draw  [draw opacity=0][dash pattern={on 0.84pt off 2.51pt}] (222.57,58.88) .. controls (223.67,57.11) and (229.29,55.75) .. (236.06,55.75) .. controls (243.51,55.75) and (249.57,57.39) .. (249.75,59.42) -- (236.06,59.51) -- cycle ; \draw  [color={rgb, 255:red, 128; green, 128; blue, 128 }  ,draw opacity=1 ][dash pattern={on 0.84pt off 2.51pt}] (222.57,58.88) .. controls (223.67,57.11) and (229.29,55.75) .. (236.06,55.75) .. controls (243.51,55.75) and (249.57,57.39) .. (249.75,59.42) ;  
\draw   (222.13,59.81) .. controls (222.13,52.2) and (228.3,46.02) .. (235.92,46.02) .. controls (243.53,46.02) and (249.71,52.2) .. (249.71,59.81) .. controls (249.71,67.43) and (243.53,73.6) .. (235.92,73.6) .. controls (228.3,73.6) and (222.13,67.43) .. (222.13,59.81) -- cycle ;
\draw  [fill={rgb, 255:red, 0; green, 0; blue, 0 }  ,fill opacity=1 ] (230.33,52.29) .. controls (230.33,51.31) and (231.12,50.52) .. (232.1,50.52) .. controls (233.07,50.52) and (233.87,51.31) .. (233.87,52.29) .. controls (233.87,53.27) and (233.07,54.06) .. (232.1,54.06) .. controls (231.12,54.06) and (230.33,53.27) .. (230.33,52.29) -- cycle ;
\draw  [fill={rgb, 255:red, 0; green, 0; blue, 0 }  ,fill opacity=1 ] (238.61,52.29) .. controls (238.61,51.31) and (239.4,50.52) .. (240.37,50.52) .. controls (241.35,50.52) and (242.14,51.31) .. (242.14,52.29) .. controls (242.14,53.27) and (241.35,54.06) .. (240.37,54.06) .. controls (239.4,54.06) and (238.61,53.27) .. (238.61,52.29) -- cycle ;

\draw  [draw opacity=0] (195.94,39.76) .. controls (195.29,42) and (188.45,43.77) .. (180.11,43.77) .. controls (171.35,43.77) and (164.24,41.82) .. (164.24,39.41) .. controls (164.24,39.26) and (164.27,39.12) .. (164.32,38.98) -- (180.11,39.41) -- cycle ; \draw  [color={rgb, 255:red, 128; green, 128; blue, 128 }  ,draw opacity=1 ] (195.94,39.76) .. controls (195.29,42) and (188.45,43.77) .. (180.11,43.77) .. controls (171.35,43.77) and (164.24,41.82) .. (164.24,39.41) .. controls (164.24,39.26) and (164.27,39.12) .. (164.32,38.98) ;  
\draw  [draw opacity=0][dash pattern={on 0.84pt off 2.51pt}] (164.63,38.33) .. controls (165.9,36.27) and (172.43,34.71) .. (180.28,34.71) .. controls (188.93,34.71) and (195.96,36.6) .. (196.15,38.96) -- (180.28,39.06) -- cycle ; \draw  [color={rgb, 255:red, 128; green, 128; blue, 128 }  ,draw opacity=1 ][dash pattern={on 0.84pt off 2.51pt}] (164.63,38.33) .. controls (165.9,36.27) and (172.43,34.71) .. (180.28,34.71) .. controls (188.93,34.71) and (195.96,36.6) .. (196.15,38.96) ;  
\draw   (164.12,39.41) .. controls (164.12,30.58) and (171.28,23.42) .. (180.11,23.42) .. controls (188.95,23.42) and (196.11,30.58) .. (196.11,39.41) .. controls (196.11,48.25) and (188.95,55.41) .. (180.11,55.41) .. controls (171.28,55.41) and (164.12,48.25) .. (164.12,39.41) -- cycle ;
\draw  [fill={rgb, 255:red, 0; green, 0; blue, 0 }  ,fill opacity=1 ] (170.63,30.49) .. controls (170.63,29.35) and (171.55,28.44) .. (172.68,28.44) .. controls (173.82,28.44) and (174.73,29.35) .. (174.73,30.49) .. controls (174.73,31.62) and (173.82,32.54) .. (172.68,32.54) .. controls (171.55,32.54) and (170.63,31.62) .. (170.63,30.49) -- cycle ;
\draw  [fill={rgb, 255:red, 0; green, 0; blue, 0 }  ,fill opacity=1 ] (178.23,30.49) .. controls (178.23,29.35) and (179.15,28.44) .. (180.28,28.44) .. controls (181.42,28.44) and (182.33,29.35) .. (182.33,30.49) .. controls (182.33,31.62) and (181.42,32.54) .. (180.28,32.54) .. controls (179.15,32.54) and (178.23,31.62) .. (178.23,30.49) -- cycle ;
\draw  [fill={rgb, 255:red, 0; green, 0; blue, 0 }  ,fill opacity=1 ] (185.83,30.49) .. controls (185.83,29.35) and (186.75,28.44) .. (187.88,28.44) .. controls (189.02,28.44) and (189.93,29.35) .. (189.93,30.49) .. controls (189.93,31.62) and (189.02,32.54) .. (187.88,32.54) .. controls (186.75,32.54) and (185.83,31.62) .. (185.83,30.49) -- cycle ;

\draw (204.5,30.9) node [anchor=north west][inner sep=0.75pt]    {$-$};
\draw (21,31.61) node [anchor=north west][inner sep=0.75pt]    {$h\ =\ $};
\draw (127,31.61) node [anchor=north west][inner sep=0.75pt]    {$t\ =\ $};
\draw (97,31.61) node [anchor=north west][inner sep=0.75pt]    {$,$};

\end{tikzpicture}

%% file: tikzpictures/1X-cob.tex
\tikzset{every picture/.style={line width=0.75pt}} 

\begin{tikzpicture}[x=0.75pt,y=0.75pt,yscale=-1,xscale=1]

\draw  [draw opacity=0][dash pattern={on 0.84pt off 2.51pt}] (25.7,33.61) .. controls (23.49,31.96) and (21.86,26.52) .. (21.86,20.07) .. controls (21.86,13.63) and (23.48,8.2) .. (25.7,6.54) -- (27.13,20.07) -- cycle ; \draw  [dash pattern={on 0.84pt off 2.51pt}] (25.7,33.61) .. controls (23.49,31.96) and (21.86,26.52) .. (21.86,20.07) .. controls (21.86,13.63) and (23.48,8.2) .. (25.7,6.54) ;  
\draw  [draw opacity=0] (26.55,6.1) .. controls (26.74,6.04) and (26.93,6.02) .. (27.13,6.02) .. controls (30.04,6.02) and (32.4,12.31) .. (32.4,20.07) .. controls (32.4,27.84) and (30.04,34.13) .. (27.13,34.13) .. controls (26.93,34.13) and (26.74,34.11) .. (26.55,34.05) -- (27.13,20.07) -- cycle ; \draw   (26.55,6.1) .. controls (26.74,6.04) and (26.93,6.02) .. (27.13,6.02) .. controls (30.04,6.02) and (32.4,12.31) .. (32.4,20.07) .. controls (32.4,27.84) and (30.04,34.13) .. (27.13,34.13) .. controls (26.93,34.13) and (26.74,34.11) .. (26.55,34.05) ;  
\draw  [draw opacity=0] (27.71,6.02) .. controls (27.51,6.02) and (27.32,6.02) .. (27.13,6.02) .. controls (18,6.02) and (10.61,12.31) .. (10.61,20.07) .. controls (10.61,27.84) and (18,34.13) .. (27.13,34.13) .. controls (27.32,34.13) and (27.51,34.13) .. (27.71,34.12) -- (27.13,20.07) -- cycle ; \draw   (27.71,6.02) .. controls (27.51,6.02) and (27.32,6.02) .. (27.13,6.02) .. controls (18,6.02) and (10.61,12.31) .. (10.61,20.07) .. controls (10.61,27.84) and (18,34.13) .. (27.13,34.13) .. controls (27.32,34.13) and (27.51,34.13) .. (27.71,34.12) ;  

\draw  [draw opacity=0][dash pattern={on 0.84pt off 2.51pt}] (117.7,34.61) .. controls (115.49,32.96) and (113.86,27.52) .. (113.86,21.07) .. controls (113.86,14.63) and (115.48,9.2) .. (117.7,7.54) -- (119.13,21.07) -- cycle ; \draw  [dash pattern={on 0.84pt off 2.51pt}] (117.7,34.61) .. controls (115.49,32.96) and (113.86,27.52) .. (113.86,21.07) .. controls (113.86,14.63) and (115.48,9.2) .. (117.7,7.54) ;  
\draw  [draw opacity=0] (118.55,7.1) .. controls (118.74,7.04) and (118.93,7.02) .. (119.13,7.02) .. controls (122.04,7.02) and (124.4,13.31) .. (124.4,21.07) .. controls (124.4,28.84) and (122.04,35.13) .. (119.13,35.13) .. controls (118.93,35.13) and (118.74,35.11) .. (118.55,35.05) -- (119.13,21.07) -- cycle ; \draw   (118.55,7.1) .. controls (118.74,7.04) and (118.93,7.02) .. (119.13,7.02) .. controls (122.04,7.02) and (124.4,13.31) .. (124.4,21.07) .. controls (124.4,28.84) and (122.04,35.13) .. (119.13,35.13) .. controls (118.93,35.13) and (118.74,35.11) .. (118.55,35.05) ;  
\draw  [draw opacity=0] (119.71,7.02) .. controls (119.51,7.02) and (119.32,7.02) .. (119.13,7.02) .. controls (110,7.02) and (102.61,13.31) .. (102.61,21.07) .. controls (102.61,28.84) and (110,35.13) .. (119.13,35.13) .. controls (119.32,35.13) and (119.51,35.13) .. (119.71,35.12) -- (119.13,21.07) -- cycle ; \draw   (119.71,7.02) .. controls (119.51,7.02) and (119.32,7.02) .. (119.13,7.02) .. controls (110,7.02) and (102.61,13.31) .. (102.61,21.07) .. controls (102.61,28.84) and (110,35.13) .. (119.13,35.13) .. controls (119.32,35.13) and (119.51,35.13) .. (119.71,35.12) ;  
\draw  [fill={rgb, 255:red, 0; green, 0; blue, 0 }  ,fill opacity=1 ] (105.49,20.82) .. controls (105.49,19.96) and (106.19,19.27) .. (107.04,19.27) .. controls (107.9,19.27) and (108.6,19.96) .. (108.6,20.82) .. controls (108.6,21.68) and (107.9,22.37) .. (107.04,22.37) .. controls (106.19,22.37) and (105.49,21.68) .. (105.49,20.82) -- cycle ;

\draw (29.13,20.07) node [anchor=west] [inner sep=0.75pt]    {$\ \leftrightarrow \ 1,$};
\draw (121.13,21.07) node [anchor=west] [inner sep=0.75pt]    {$\ \leftrightarrow \ X.$};

\end{tikzpicture}

%% file: tikzpictures/frobenius-ops-cob.tex
\tikzset{every picture/.style={line width=0.75pt}} 

\begin{tikzpicture}[x=0.75pt,y=0.75pt,yscale=-1,xscale=1]

\draw  [draw opacity=0][dash pattern={on 0.84pt off 2.51pt}] (22.18,44.98) .. controls (20.25,43.54) and (18.84,38.82) .. (18.84,33.22) .. controls (18.84,27.63) and (20.25,22.91) .. (22.17,21.47) -- (23.42,33.22) -- cycle ; \draw  [dash pattern={on 0.84pt off 2.51pt}] (22.18,44.98) .. controls (20.25,43.54) and (18.84,38.82) .. (18.84,33.22) .. controls (18.84,27.63) and (20.25,22.91) .. (22.17,21.47) ;  
\draw  [draw opacity=0] (22.92,21.09) .. controls (23.08,21.04) and (23.25,21.02) .. (23.42,21.02) .. controls (25.94,21.02) and (27.99,26.48) .. (27.99,33.22) .. controls (27.99,39.96) and (25.94,45.43) .. (23.42,45.43) .. controls (23.25,45.43) and (23.08,45.4) .. (22.92,45.36) -- (23.42,33.22) -- cycle ; \draw   (22.92,21.09) .. controls (23.08,21.04) and (23.25,21.02) .. (23.42,21.02) .. controls (25.94,21.02) and (27.99,26.48) .. (27.99,33.22) .. controls (27.99,39.96) and (25.94,45.43) .. (23.42,45.43) .. controls (23.25,45.43) and (23.08,45.4) .. (22.92,45.36) ;  
\draw  [draw opacity=0] (23.92,21.02) .. controls (23.75,21.02) and (23.58,21.02) .. (23.42,21.02) .. controls (15.49,21.02) and (9.07,26.48) .. (9.07,33.22) .. controls (9.07,39.96) and (15.49,45.43) .. (23.42,45.43) .. controls (23.58,45.43) and (23.75,45.43) .. (23.92,45.42) -- (23.42,33.22) -- cycle ; \draw   (23.92,21.02) .. controls (23.75,21.02) and (23.58,21.02) .. (23.42,21.02) .. controls (15.49,21.02) and (9.07,26.48) .. (9.07,33.22) .. controls (9.07,39.96) and (15.49,45.43) .. (23.42,45.43) .. controls (23.58,45.43) and (23.75,45.43) .. (23.92,45.42) ;  

\draw  [color={rgb, 255:red, 0; green, 0; blue, 0 }  ,draw opacity=1 ] (49.93,33.69) .. controls (49.93,27) and (52.17,21.57) .. (54.93,21.57) .. controls (57.69,21.57) and (59.92,27) .. (59.92,33.69) .. controls (59.92,40.38) and (57.69,45.8) .. (54.93,45.8) .. controls (52.17,45.8) and (49.93,40.38) .. (49.93,33.69) -- cycle ;
\draw  [draw opacity=0] (54.42,21.49) .. controls (54.59,21.48) and (54.76,21.48) .. (54.93,21.48) .. controls (62.85,21.48) and (69.27,26.94) .. (69.27,33.69) .. controls (69.27,40.43) and (62.85,45.89) .. (54.93,45.89) .. controls (54.76,45.89) and (54.59,45.89) .. (54.42,45.88) -- (54.93,33.69) -- cycle ; \draw   (54.42,21.49) .. controls (54.59,21.48) and (54.76,21.48) .. (54.93,21.48) .. controls (62.85,21.48) and (69.27,26.94) .. (69.27,33.69) .. controls (69.27,40.43) and (62.85,45.89) .. (54.93,45.89) .. controls (54.76,45.89) and (54.59,45.89) .. (54.42,45.88) ;  
\draw  [color={rgb, 255:red, 0; green, 0; blue, 0 }  ,draw opacity=1 ] (108.62,17.29) .. controls (108.62,11.03) and (110.71,5.95) .. (113.3,5.95) .. controls (115.88,5.95) and (117.98,11.03) .. (117.98,17.29) .. controls (117.98,23.56) and (115.88,28.64) .. (113.3,28.64) .. controls (110.71,28.64) and (108.62,23.56) .. (108.62,17.29) -- cycle ;
\draw  [color={rgb, 255:red, 0; green, 0; blue, 0 }  ,draw opacity=1 ] (108.62,49.01) .. controls (108.62,42.74) and (110.71,37.66) .. (113.3,37.66) .. controls (115.88,37.66) and (117.98,42.74) .. (117.98,49.01) .. controls (117.98,55.27) and (115.88,60.35) .. (113.3,60.35) .. controls (110.71,60.35) and (108.62,55.27) .. (108.62,49.01) -- cycle ;
\draw  [draw opacity=0] (113.44,37.5) .. controls (122.36,37.47) and (129.56,35.5) .. (129.56,33.07) .. controls (129.56,30.62) and (122.28,28.64) .. (113.3,28.64) -- (113.3,33.07) -- cycle ; \draw   (113.44,37.5) .. controls (122.36,37.47) and (129.56,35.5) .. (129.56,33.07) .. controls (129.56,30.62) and (122.28,28.64) .. (113.3,28.64) ;  
\draw    (113.3,5.95) .. controls (136.59,5.33) and (146.51,22.76) .. (162.25,21.76) ;
\draw    (113.3,60.35) .. controls (136.59,59.73) and (146.51,45.49) .. (162.25,44.49) ;
\draw  [draw opacity=0][dash pattern={on 0.84pt off 2.51pt}] (161.56,44.13) .. controls (159.75,42.79) and (158.43,38.37) .. (158.43,33.12) .. controls (158.43,27.89) and (159.75,23.47) .. (161.55,22.12) -- (162.71,33.12) -- cycle ; \draw  [dash pattern={on 0.84pt off 2.51pt}] (161.56,44.13) .. controls (159.75,42.79) and (158.43,38.37) .. (158.43,33.12) .. controls (158.43,27.89) and (159.75,23.47) .. (161.55,22.12) ;  
\draw  [draw opacity=0] (162.25,21.76) .. controls (162.4,21.72) and (162.56,21.69) .. (162.71,21.69) .. controls (165.08,21.69) and (167,26.81) .. (167,33.12) .. controls (167,39.44) and (165.08,44.55) .. (162.71,44.55) .. controls (162.56,44.55) and (162.4,44.53) .. (162.25,44.49) -- (162.71,33.12) -- cycle ; \draw   (162.25,21.76) .. controls (162.4,21.72) and (162.56,21.69) .. (162.71,21.69) .. controls (165.08,21.69) and (167,26.81) .. (167,33.12) .. controls (167,39.44) and (165.08,44.55) .. (162.71,44.55) .. controls (162.56,44.55) and (162.4,44.53) .. (162.25,44.49) ;  

\draw  [color={rgb, 255:red, 0; green, 0; blue, 0 }  ,draw opacity=1 ] (215.47,33.12) .. controls (215.47,26.86) and (213.37,21.78) .. (210.79,21.78) .. controls (208.2,21.78) and (206.11,26.86) .. (206.11,33.12) .. controls (206.11,39.39) and (208.2,44.47) .. (210.79,44.47) .. controls (213.37,44.47) and (215.47,39.39) .. (215.47,33.12) -- cycle ;
\draw  [draw opacity=0] (259.59,37.51) .. controls (250.68,37.49) and (243.47,35.52) .. (243.47,33.09) .. controls (243.47,30.64) and (250.76,28.66) .. (259.74,28.66) -- (259.74,33.09) -- cycle ; \draw   (259.59,37.51) .. controls (250.68,37.49) and (243.47,35.52) .. (243.47,33.09) .. controls (243.47,30.64) and (250.76,28.66) .. (259.74,28.66) ;  
\draw    (259.74,5.97) .. controls (236.44,5.35) and (226.53,22.78) .. (210.79,21.78) ;
\draw    (259.74,60.37) .. controls (236.44,59.75) and (226.53,45.51) .. (210.79,44.51) ;
\draw  [draw opacity=0][dash pattern={on 0.84pt off 2.51pt}] (259.16,28.15) .. controls (257.36,26.81) and (256.03,22.39) .. (256.03,17.14) .. controls (256.03,11.9) and (257.36,7.49) .. (259.16,6.14) -- (260.32,17.14) -- cycle ; \draw  [dash pattern={on 0.84pt off 2.51pt}] (259.16,28.15) .. controls (257.36,26.81) and (256.03,22.39) .. (256.03,17.14) .. controls (256.03,11.9) and (257.36,7.49) .. (259.16,6.14) ;  
\draw  [draw opacity=0] (259.85,5.78) .. controls (260.01,5.73) and (260.16,5.71) .. (260.32,5.71) .. controls (262.69,5.71) and (264.61,10.83) .. (264.61,17.14) .. controls (264.61,23.46) and (262.69,28.57) .. (260.32,28.57) .. controls (260.16,28.57) and (260.01,28.55) .. (259.85,28.51) -- (260.32,17.14) -- cycle ; \draw   (259.85,5.78) .. controls (260.01,5.73) and (260.16,5.71) .. (260.32,5.71) .. controls (262.69,5.71) and (264.61,10.83) .. (264.61,17.14) .. controls (264.61,23.46) and (262.69,28.57) .. (260.32,28.57) .. controls (260.16,28.57) and (260.01,28.55) .. (259.85,28.51) ;  

\draw  [draw opacity=0][dash pattern={on 0.84pt off 2.51pt}] (258.16,60.15) .. controls (256.36,58.81) and (255.03,54.39) .. (255.03,49.14) .. controls (255.03,43.9) and (256.36,39.49) .. (258.16,38.14) -- (259.32,49.14) -- cycle ; \draw  [dash pattern={on 0.84pt off 2.51pt}] (258.16,60.15) .. controls (256.36,58.81) and (255.03,54.39) .. (255.03,49.14) .. controls (255.03,43.9) and (256.36,39.49) .. (258.16,38.14) ;  
\draw  [draw opacity=0] (258.85,37.78) .. controls (259.01,37.73) and (259.16,37.71) .. (259.32,37.71) .. controls (261.69,37.71) and (263.61,42.83) .. (263.61,49.14) .. controls (263.61,55.46) and (261.69,60.57) .. (259.32,60.57) .. controls (259.16,60.57) and (259.01,60.55) .. (258.85,60.51) -- (259.32,49.14) -- cycle ; \draw   (258.85,37.78) .. controls (259.01,37.73) and (259.16,37.71) .. (259.32,37.71) .. controls (261.69,37.71) and (263.61,42.83) .. (263.61,49.14) .. controls (263.61,55.46) and (261.69,60.57) .. (259.32,60.57) .. controls (259.16,60.57) and (259.01,60.55) .. (258.85,60.51) ;

\draw (14,63.8) node [anchor=north west][inner sep=0.75pt]    {$\iota $};
\draw (55.42,63.8) node [anchor=north west][inner sep=0.75pt]    {$\epsilon $};
\draw (132.94,63.8) node [anchor=north west][inner sep=0.75pt]    {$m$};
\draw (229.6,63.8) node [anchor=north west][inner sep=0.75pt]    {$\Delta $};

\end{tikzpicture}

%% file: tikzpictures/planar-diagram.tex
\tikzset{every picture/.style={line width=0.75pt}} 

\begin{tikzpicture}[x=0.75pt,y=0.75pt,yscale=-.9,xscale=.9]

\draw  [fill={rgb, 255:red, 241; green, 241; blue, 241 }  ,fill opacity=1 ][dash pattern={on 4.5pt off 4.5pt}] (13.38,70.46) .. controls (13.38,37.4) and (40.18,10.6) .. (73.24,10.6) .. controls (106.31,10.6) and (133.11,37.4) .. (133.11,70.46) .. controls (133.11,103.52) and (106.31,130.33) .. (73.24,130.33) .. controls (40.18,130.33) and (13.38,103.52) .. (13.38,70.46) -- cycle ;
\draw  [fill={rgb, 255:red, 255; green, 255; blue, 255 }  ,fill opacity=1 ][dash pattern={on 4.5pt off 4.5pt}] (28.91,70.46) .. controls (28.91,60.82) and (36.72,53) .. (46.37,53) .. controls (56.01,53) and (63.83,60.82) .. (63.83,70.46) .. controls (63.83,80.11) and (56.01,87.92) .. (46.37,87.92) .. controls (36.72,87.92) and (28.91,80.11) .. (28.91,70.46) -- cycle ;
\draw  [fill={rgb, 255:red, 255; green, 255; blue, 255 }  ,fill opacity=1 ][dash pattern={on 4.5pt off 4.5pt}] (84.1,70.46) .. controls (84.1,62.7) and (90.39,56.41) .. (98.15,56.41) .. controls (105.91,56.41) and (112.2,62.7) .. (112.2,70.46) .. controls (112.2,78.22) and (105.91,84.51) .. (98.15,84.51) .. controls (90.39,84.51) and (84.1,78.22) .. (84.1,70.46) -- cycle ;
\draw  [fill={rgb, 255:red, 0; green, 0; blue, 0 }  ,fill opacity=1 ] (50.76,14.14) .. controls (50.76,12.87) and (51.78,11.85) .. (53.05,11.85) .. controls (54.32,11.85) and (55.34,12.87) .. (55.34,14.14) .. controls (55.34,15.41) and (54.32,16.43) .. (53.05,16.43) .. controls (51.78,16.43) and (50.76,15.41) .. (50.76,14.14) -- cycle ;
\draw  [fill={rgb, 255:red, 0; green, 0; blue, 0 }  ,fill opacity=1 ] (92.6,14.43) .. controls (92.6,13.16) and (93.63,12.13) .. (94.89,12.13) .. controls (96.16,12.13) and (97.18,13.16) .. (97.18,14.43) .. controls (97.18,15.69) and (96.16,16.72) .. (94.89,16.72) .. controls (93.63,16.72) and (92.6,15.69) .. (92.6,14.43) -- cycle ;
\draw  [fill={rgb, 255:red, 0; green, 0; blue, 0 }  ,fill opacity=1 ] (130.81,70.46) .. controls (130.81,69.2) and (131.84,68.17) .. (133.11,68.17) .. controls (134.37,68.17) and (135.4,69.2) .. (135.4,70.46) .. controls (135.4,71.73) and (134.37,72.76) .. (133.11,72.76) .. controls (131.84,72.76) and (130.81,71.73) .. (130.81,70.46) -- cycle ;
\draw  [fill={rgb, 255:red, 0; green, 0; blue, 0 }  ,fill opacity=1 ] (107.5,117.41) .. controls (107.5,116.14) and (108.53,115.12) .. (109.79,115.12) .. controls (111.06,115.12) and (112.09,116.14) .. (112.09,117.41) .. controls (112.09,118.68) and (111.06,119.7) .. (109.79,119.7) .. controls (108.53,119.7) and (107.5,118.68) .. (107.5,117.41) -- cycle ;
\draw  [fill={rgb, 255:red, 0; green, 0; blue, 0 }  ,fill opacity=1 ] (11.21,70.46) .. controls (11.21,69.2) and (12.24,68.17) .. (13.5,68.17) .. controls (14.77,68.17) and (15.8,69.2) .. (15.8,70.46) .. controls (15.8,71.73) and (14.77,72.76) .. (13.5,72.76) .. controls (12.24,72.76) and (11.21,71.73) .. (11.21,70.46) -- cycle ;
\draw  [fill={rgb, 255:red, 0; green, 0; blue, 0 }  ,fill opacity=1 ] (34.71,117.31) .. controls (34.71,116.04) and (35.74,115.02) .. (37,115.02) .. controls (38.27,115.02) and (39.3,116.04) .. (39.3,117.31) .. controls (39.3,118.57) and (38.27,119.6) .. (37,119.6) .. controls (35.74,119.6) and (34.71,118.57) .. (34.71,117.31) -- cycle ;
\draw  [fill={rgb, 255:red, 0; green, 0; blue, 0 }  ,fill opacity=1 ] (32.3,56.54) .. controls (32.3,55.28) and (33.33,54.25) .. (34.59,54.25) .. controls (35.86,54.25) and (36.89,55.28) .. (36.89,56.54) .. controls (36.89,57.81) and (35.86,58.84) .. (34.59,58.84) .. controls (33.33,58.84) and (32.3,57.81) .. (32.3,56.54) -- cycle ;
\draw  [fill={rgb, 255:red, 0; green, 0; blue, 0 }  ,fill opacity=1 ] (53.75,56.04) .. controls (53.75,54.78) and (54.78,53.75) .. (56.04,53.75) .. controls (57.31,53.75) and (58.34,54.78) .. (58.34,56.04) .. controls (58.34,57.31) and (57.31,58.34) .. (56.04,58.34) .. controls (54.78,58.34) and (53.75,57.31) .. (53.75,56.04) -- cycle ;
\draw  [fill={rgb, 255:red, 0; green, 0; blue, 0 }  ,fill opacity=1 ] (33.3,83.48) .. controls (33.3,82.21) and (34.32,81.19) .. (35.59,81.19) .. controls (36.86,81.19) and (37.88,82.21) .. (37.88,83.48) .. controls (37.88,84.75) and (36.86,85.77) .. (35.59,85.77) .. controls (34.32,85.77) and (33.3,84.75) .. (33.3,83.48) -- cycle ;
\draw  [fill={rgb, 255:red, 0; green, 0; blue, 0 }  ,fill opacity=1 ] (54.75,82.98) .. controls (54.75,81.72) and (55.78,80.69) .. (57.04,80.69) .. controls (58.31,80.69) and (59.33,81.72) .. (59.33,82.98) .. controls (59.33,84.25) and (58.31,85.27) .. (57.04,85.27) .. controls (55.78,85.27) and (54.75,84.25) .. (54.75,82.98) -- cycle ;
\draw  [fill={rgb, 255:red, 0; green, 0; blue, 0 }  ,fill opacity=1 ] (95.46,56.41) .. controls (95.46,55.15) and (96.48,54.12) .. (97.75,54.12) .. controls (99.02,54.12) and (100.04,55.15) .. (100.04,56.41) .. controls (100.04,57.68) and (99.02,58.71) .. (97.75,58.71) .. controls (96.48,58.71) and (95.46,57.68) .. (95.46,56.41) -- cycle ;
\draw  [fill={rgb, 255:red, 0; green, 0; blue, 0 }  ,fill opacity=1 ] (95.46,84.51) .. controls (95.46,83.25) and (96.48,82.22) .. (97.75,82.22) .. controls (99.02,82.22) and (100.04,83.25) .. (100.04,84.51) .. controls (100.04,85.78) and (99.02,86.81) .. (97.75,86.81) .. controls (96.48,86.81) and (95.46,85.78) .. (95.46,84.51) -- cycle ;
\draw    (56.04,56.04) .. controls (69.31,43.28) and (90.57,48.43) .. (97.75,56.41) ;
\draw    (59.33,82.98) .. controls (70.81,90.17) and (86.28,92.25) .. (97.75,84.51) ;
\draw    (34.59,56.54) .. controls (36.39,34.79) and (50.36,34.5) .. (53.05,14.14) ;
\draw    (130.81,70.46) .. controls (126.46,71.09) and (92.26,25.32) .. (94.89,14.43) ;
\draw    (107.5,117.41) .. controls (102.24,103.64) and (48.86,101.64) .. (37,117.31) ;
\draw    (13.5,70.46) .. controls (21.92,68.72) and (25.91,90.67) .. (35.59,85.77) ;

\draw (98.15,70.46) node  [font=\scriptsize]  {$D_{2}$};
\draw (46.37,70.46) node  [font=\scriptsize]  {$D_{1}$};

\end{tikzpicture}

%% file: tikzpictures/crossing-endpoints.tex
\tikzset{every picture/.style={line width=0.75pt}} 

\begin{tikzpicture}[x=0.75pt,y=0.75pt,yscale=-.75,xscale=.75]

\draw [color={rgb, 255:red, 0; green, 0; blue, 0 }  ,draw opacity=1 ][line width=1.5]    (51.53,35.12) -- (106.05,89.65) ;
\draw [shift={(49.4,33)}, rotate = 45] [color={rgb, 255:red, 0; green, 0; blue, 0 }  ,draw opacity=1 ][line width=1.5]    (8.53,-2.57) .. controls (5.42,-1.09) and (2.58,-0.23) .. (0,0) .. controls (2.58,0.23) and (5.42,1.09) .. (8.53,2.57)   ;
\draw [color={rgb, 255:red, 0; green, 0; blue, 0 }  ,draw opacity=1 ][line width=1.5]    (104.32,35.11) -- (49,89.65) ;
\draw [shift={(106.46,33)}, rotate = 135.41] [color={rgb, 255:red, 0; green, 0; blue, 0 }  ,draw opacity=1 ][line width=1.5]    (8.53,-2.57) .. controls (5.42,-1.09) and (2.58,-0.23) .. (0,0) .. controls (2.58,0.23) and (5.42,1.09) .. (8.53,2.57)   ;

\draw (47,93.05) node [anchor=north east] [inner sep=0.75pt]    {$p_{1}( c)$};
\draw (108.05,93.05) node [anchor=north west][inner sep=0.75pt]    {$p_{2}( c)$};
\draw (47.4,29.6) node [anchor=south east] [inner sep=0.75pt]    {$p'_{1}( c)$};
\draw (108.46,29.6) node [anchor=south west] [inner sep=0.75pt]    {$p'_{2}( c)$};

\end{tikzpicture}

%% file: tikzpictures/crossing-admissible.tex
\tikzset{every picture/.style={line width=0.75pt}} 

\begin{tikzpicture}[x=0.75pt,y=0.75pt,yscale=-.75,xscale=.75]

\draw [color={rgb, 255:red, 208; green, 2; blue, 27 }  ,draw opacity=1 ][line width=1.5]    (52.12,32.12) -- (80,60) ;
\draw [shift={(50,30)}, rotate = 45] [color={rgb, 255:red, 208; green, 2; blue, 27 }  ,draw opacity=1 ][line width=1.5]    (8.53,-2.57) .. controls (5.42,-1.09) and (2.58,-0.23) .. (0,0) .. controls (2.58,0.23) and (5.42,1.09) .. (8.53,2.57)   ;
\draw [color={rgb, 255:red, 208; green, 2; blue, 27 }  ,draw opacity=1 ][line width=1.5]    (80,60) -- (49,89.65) ;
\draw [color={rgb, 255:red, 74; green, 144; blue, 226 }  ,draw opacity=1 ][line width=1.5]    (80,60) -- (110,90) ;
\draw [color={rgb, 255:red, 74; green, 144; blue, 226 }  ,draw opacity=1 ][line width=1.5]    (107.54,34.08) -- (80,60) ;
\draw [shift={(109.73,32.03)}, rotate = 136.74] [color={rgb, 255:red, 74; green, 144; blue, 226 }  ,draw opacity=1 ][line width=1.5]    (8.53,-2.57) .. controls (5.42,-1.09) and (2.58,-0.23) .. (0,0) .. controls (2.58,0.23) and (5.42,1.09) .. (8.53,2.57)   ;

\draw (47,93.05) node [anchor=north east] [inner sep=0.75pt]    {$p_{1}( c)$};
\draw (112,93.4) node [anchor=north west][inner sep=0.75pt]    {$p_{2}( c)$};
\draw (48,26.6) node [anchor=south east] [inner sep=0.75pt]    {$p'_{1}( c)$};
\draw (111.73,28.63) node [anchor=south west] [inner sep=0.75pt]    {$p'_{2}( c)$};

\end{tikzpicture}

%% file: tikzpictures/a-coeffs.tex
\tikzset{every picture/.style={line width=0.75pt}} 

\begin{tikzpicture}[x=0.75pt,y=0.75pt,yscale=-.9,xscale=.9]

\draw    (80,80) -- (110,110) ;
\draw  [color={rgb, 255:red, 255; green, 255; blue, 255 }  ,draw opacity=1 ][fill={rgb, 255:red, 255; green, 255; blue, 255 }  ,fill opacity=1 ] (90,95) .. controls (90,92.24) and (92.24,90) .. (95,90) .. controls (97.76,90) and (100,92.24) .. (100,95) .. controls (100,97.76) and (97.76,100) .. (95,100) .. controls (92.24,100) and (90,97.76) .. (90,95) -- cycle ;
\draw    (110,80) -- (80,110) ;

\draw    (80,10) -- (110,40) ;
\draw  [color={rgb, 255:red, 255; green, 255; blue, 255 }  ,draw opacity=1 ][fill={rgb, 255:red, 255; green, 255; blue, 255 }  ,fill opacity=1 ] (90,25) .. controls (90,22.24) and (92.24,20) .. (95,20) .. controls (97.76,20) and (100,22.24) .. (100,25) .. controls (100,27.76) and (97.76,30) .. (95,30) .. controls (92.24,30) and (90,27.76) .. (90,25) -- cycle ;
\draw    (110,10) -- (80,40) ;

\draw [color={rgb, 255:red, 0; green, 0; blue, 0 }  ,draw opacity=1 ] [dash pattern={on 0.84pt off 2.51pt}]  (80,80) -- (110,40) ;
\draw [color={rgb, 255:red, 0; green, 0; blue, 0 }  ,draw opacity=1 ] [dash pattern={on 0.84pt off 2.51pt}]  (80,40) -- (110,80) ;
\draw    (140,80) -- (170,110) ;
\draw  [color={rgb, 255:red, 255; green, 255; blue, 255 }  ,draw opacity=1 ][fill={rgb, 255:red, 255; green, 255; blue, 255 }  ,fill opacity=1 ] (150,95) .. controls (150,92.24) and (152.24,90) .. (155,90) .. controls (157.76,90) and (160,92.24) .. (160,95) .. controls (160,97.76) and (157.76,100) .. (155,100) .. controls (152.24,100) and (150,97.76) .. (150,95) -- cycle ;
\draw    (170,80) -- (140,110) ;

\draw    (140,10) -- (170,40) ;
\draw  [color={rgb, 255:red, 255; green, 255; blue, 255 }  ,draw opacity=1 ][fill={rgb, 255:red, 255; green, 255; blue, 255 }  ,fill opacity=1 ] (150,25) .. controls (150,22.24) and (152.24,20) .. (155,20) .. controls (157.76,20) and (160,22.24) .. (160,25) .. controls (160,27.76) and (157.76,30) .. (155,30) .. controls (152.24,30) and (150,27.76) .. (150,25) -- cycle ;
\draw    (170,10) -- (140,40) ;

\draw [color={rgb, 255:red, 0; green, 0; blue, 0 }  ,draw opacity=1 ] [dash pattern={on 0.84pt off 2.51pt}]  (140,80) -- (170,40) ;
\draw    (200,80) -- (230,110) ;
\draw  [color={rgb, 255:red, 255; green, 255; blue, 255 }  ,draw opacity=1 ][fill={rgb, 255:red, 255; green, 255; blue, 255 }  ,fill opacity=1 ] (210,95) .. controls (210,92.24) and (212.24,90) .. (215,90) .. controls (217.76,90) and (220,92.24) .. (220,95) .. controls (220,97.76) and (217.76,100) .. (215,100) .. controls (212.24,100) and (210,97.76) .. (210,95) -- cycle ;
\draw    (230,80) -- (200,110) ;

\draw    (200,10) -- (230,40) ;
\draw  [color={rgb, 255:red, 255; green, 255; blue, 255 }  ,draw opacity=1 ][fill={rgb, 255:red, 255; green, 255; blue, 255 }  ,fill opacity=1 ] (210,25) .. controls (210,22.24) and (212.24,20) .. (215,20) .. controls (217.76,20) and (220,22.24) .. (220,25) .. controls (220,27.76) and (217.76,30) .. (215,30) .. controls (212.24,30) and (210,27.76) .. (210,25) -- cycle ;
\draw    (230,10) -- (200,40) ;

\draw    (260,80) -- (290,110) ;
\draw  [color={rgb, 255:red, 255; green, 255; blue, 255 }  ,draw opacity=1 ][fill={rgb, 255:red, 255; green, 255; blue, 255 }  ,fill opacity=1 ] (270,95) .. controls (270,92.24) and (272.24,90) .. (275,90) .. controls (277.76,90) and (280,92.24) .. (280,95) .. controls (280,97.76) and (277.76,100) .. (275,100) .. controls (272.24,100) and (270,97.76) .. (270,95) -- cycle ;
\draw    (290,80) -- (260,110) ;

\draw    (260,10) -- (290,40) ;
\draw  [color={rgb, 255:red, 255; green, 255; blue, 255 }  ,draw opacity=1 ][fill={rgb, 255:red, 255; green, 255; blue, 255 }  ,fill opacity=1 ] (270,25) .. controls (270,22.24) and (272.24,20) .. (275,20) .. controls (277.76,20) and (280,22.24) .. (280,25) .. controls (280,27.76) and (277.76,30) .. (275,30) .. controls (272.24,30) and (270,27.76) .. (270,25) -- cycle ;
\draw    (290,10) -- (260,40) ;

\draw    (320,80) -- (350,110) ;
\draw  [color={rgb, 255:red, 255; green, 255; blue, 255 }  ,draw opacity=1 ][fill={rgb, 255:red, 255; green, 255; blue, 255 }  ,fill opacity=1 ] (330,95) .. controls (330,92.24) and (332.24,90) .. (335,90) .. controls (337.76,90) and (340,92.24) .. (340,95) .. controls (340,97.76) and (337.76,100) .. (335,100) .. controls (332.24,100) and (330,97.76) .. (330,95) -- cycle ;
\draw    (350,80) -- (320,110) ;

\draw    (320,10) -- (350,40) ;
\draw  [color={rgb, 255:red, 255; green, 255; blue, 255 }  ,draw opacity=1 ][fill={rgb, 255:red, 255; green, 255; blue, 255 }  ,fill opacity=1 ] (330,25) .. controls (330,22.24) and (332.24,20) .. (335,20) .. controls (337.76,20) and (340,22.24) .. (340,25) .. controls (340,27.76) and (337.76,30) .. (335,30) .. controls (332.24,30) and (330,27.76) .. (330,25) -- cycle ;
\draw    (350,10) -- (320,40) ;

\draw [color={rgb, 255:red, 0; green, 0; blue, 0 }  ,draw opacity=1 ] [dash pattern={on 0.84pt off 2.51pt}]  (320,40) -- (320,80) ;
\draw    (380,80) -- (410,110) ;
\draw  [color={rgb, 255:red, 255; green, 255; blue, 255 }  ,draw opacity=1 ][fill={rgb, 255:red, 255; green, 255; blue, 255 }  ,fill opacity=1 ] (390,95) .. controls (390,92.24) and (392.24,90) .. (395,90) .. controls (397.76,90) and (400,92.24) .. (400,95) .. controls (400,97.76) and (397.76,100) .. (395,100) .. controls (392.24,100) and (390,97.76) .. (390,95) -- cycle ;
\draw    (410,80) -- (380,110) ;

\draw    (380,10) -- (410,40) ;
\draw  [color={rgb, 255:red, 255; green, 255; blue, 255 }  ,draw opacity=1 ][fill={rgb, 255:red, 255; green, 255; blue, 255 }  ,fill opacity=1 ] (390,25) .. controls (390,22.24) and (392.24,20) .. (395,20) .. controls (397.76,20) and (400,22.24) .. (400,25) .. controls (400,27.76) and (397.76,30) .. (395,30) .. controls (392.24,30) and (390,27.76) .. (390,25) -- cycle ;
\draw    (410,10) -- (380,40) ;

\draw [color={rgb, 255:red, 0; green, 0; blue, 0 }  ,draw opacity=1 ] [dash pattern={on 0.84pt off 2.51pt}]  (200,40) -- (230,80) ;
\draw    (440,80) -- (470,110) ;
\draw  [color={rgb, 255:red, 255; green, 255; blue, 255 }  ,draw opacity=1 ][fill={rgb, 255:red, 255; green, 255; blue, 255 }  ,fill opacity=1 ] (450,95) .. controls (450,92.24) and (452.24,90) .. (455,90) .. controls (457.76,90) and (460,92.24) .. (460,95) .. controls (460,97.76) and (457.76,100) .. (455,100) .. controls (452.24,100) and (450,97.76) .. (450,95) -- cycle ;
\draw    (470,80) -- (440,110) ;

\draw    (440,10) -- (470,40) ;
\draw  [color={rgb, 255:red, 255; green, 255; blue, 255 }  ,draw opacity=1 ][fill={rgb, 255:red, 255; green, 255; blue, 255 }  ,fill opacity=1 ] (450,25) .. controls (450,22.24) and (452.24,20) .. (455,20) .. controls (457.76,20) and (460,22.24) .. (460,25) .. controls (460,27.76) and (457.76,30) .. (455,30) .. controls (452.24,30) and (450,27.76) .. (450,25) -- cycle ;
\draw    (470,10) -- (440,40) ;

\draw [color={rgb, 255:red, 0; green, 0; blue, 0 }  ,draw opacity=1 ] [dash pattern={on 0.84pt off 2.51pt}]  (410,40) -- (410,80) ;
\draw [color={rgb, 255:red, 0; green, 0; blue, 0 }  ,draw opacity=1 ] [dash pattern={on 0.84pt off 2.51pt}]  (440,40) -- (440,80) ;
\draw [color={rgb, 255:red, 0; green, 0; blue, 0 }  ,draw opacity=1 ] [dash pattern={on 0.84pt off 2.51pt}]  (470,40) -- (470,80) ;

\draw (11,132.5) node [anchor=west] [inner sep=0.75pt]    {$a_{c,\ c'} \ =\ \ $};
\draw (96.5,132.5) node    {$2$};
\draw (153.5,132.5) node    {$1$};
\draw (213.5,132.5) node    {$1$};
\draw (276.5,132.5) node    {$0$};
\draw (336.5,132.5) node    {$-1$};
\draw (399,132.5) node    {$-1$};
\draw (459,132.5) node    {$-2$};
\draw (64.22,92.5) node    {$c$};
\draw (66.5,20.5) node    {$c'$};

\end{tikzpicture}

%% file: tikzpictures/R2.tex
\tikzset{every picture/.style={line width=0.75pt}} 

\begin{tikzpicture}[x=0.75pt,y=0.75pt,yscale=-1,xscale=1]

\draw    (138.59,11.41) -- (110,40) ;
\draw [shift={(140,10)}, rotate = 135] [color={rgb, 255:red, 0; green, 0; blue, 0 }  ][line width=0.75]    (6.56,-1.97) .. controls (4.17,-0.84) and (1.99,-0.18) .. (0,0) .. controls (1.99,0.18) and (4.17,0.84) .. (6.56,1.97)   ;
\draw  [color={rgb, 255:red, 255; green, 255; blue, 255 }  ,draw opacity=1 ][fill={rgb, 255:red, 255; green, 255; blue, 255 }  ,fill opacity=1 ] (120,25) .. controls (120,22.24) and (122.24,20) .. (125,20) .. controls (127.76,20) and (130,22.24) .. (130,25) .. controls (130,27.76) and (127.76,30) .. (125,30) .. controls (122.24,30) and (120,27.76) .. (120,25) -- cycle ;
\draw    (110,40) -- (140,70) ;
\draw  [color={rgb, 255:red, 255; green, 255; blue, 255 }  ,draw opacity=1 ][fill={rgb, 255:red, 255; green, 255; blue, 255 }  ,fill opacity=1 ] (120,55) .. controls (120,52.24) and (122.24,50) .. (125,50) .. controls (127.76,50) and (130,52.24) .. (130,55) .. controls (130,57.76) and (127.76,60) .. (125,60) .. controls (122.24,60) and (120,57.76) .. (120,55) -- cycle ;
\draw    (140,40) -- (110,70) ;

\draw    (10,12) -- (10,70) ;
\draw [shift={(10,10)}, rotate = 90] [color={rgb, 255:red, 0; green, 0; blue, 0 }  ][line width=0.75]    (6.56,-1.97) .. controls (4.17,-0.84) and (1.99,-0.18) .. (0,0) .. controls (1.99,0.18) and (4.17,0.84) .. (6.56,1.97)   ;
\draw    (40,12) -- (40,70) ;
\draw [shift={(40,10)}, rotate = 90] [color={rgb, 255:red, 0; green, 0; blue, 0 }  ][line width=0.75]    (6.56,-1.97) .. controls (4.17,-0.84) and (1.99,-0.18) .. (0,0) .. controls (1.99,0.18) and (4.17,0.84) .. (6.56,1.97)   ;

\draw    (111.41,11.41) -- (140,40) ;
\draw [shift={(110,10)}, rotate = 45] [color={rgb, 255:red, 0; green, 0; blue, 0 }  ][line width=0.75]    (6.56,-1.97) .. controls (4.17,-0.84) and (1.99,-0.18) .. (0,0) .. controls (1.99,0.18) and (4.17,0.84) .. (6.56,1.97)   ;

\draw (146,52) node [anchor=west] [inner sep=0.75pt]  [font=\footnotesize]  {$c_{1}$};
\draw (146,22) node [anchor=west] [inner sep=0.75pt]  [font=\footnotesize]  {$c_{2}$};
\draw (26.5,82.4) node [anchor=north] [inner sep=0.75pt]    {$T$};
\draw (126.5,82.4) node [anchor=north] [inner sep=0.75pt]    {$T'$};

\end{tikzpicture}

%% file: tikzpictures/R2-maps.tex
\tikzset{every picture/.style={line width=0.75pt}} 

\begin{tikzpicture}[x=0.75pt,y=0.75pt,yscale=-.85,xscale=.85]

\draw    (310,190) -- (280,220) ;
\draw    (280,220) -- (310,250) ;
\draw    (310,220) -- (280,250) ;
\draw    (20,100) -- (20,160) ;
\draw    (50,100) -- (50,160) ;
\draw    (280,190) -- (310,220) ;
\draw    (310,11) -- (280,41) ;
\draw    (280,41) -- (310,71) ;
\draw    (310,41) -- (280,71) ;
\draw    (280,11) -- (310,41) ;
\draw    (200,100) -- (170,130) ;
\draw    (170,130) -- (200,160) ;
\draw    (200,130) -- (170,160) ;
\draw    (170,100) -- (200,130) ;
\draw    (420,100) -- (390,130) ;
\draw    (390,130) -- (420,160) ;
\draw    (420,130) -- (390,160) ;
\draw    (390,100) -- (420,130) ;
\draw  [color={rgb, 255:red, 255; green, 255; blue, 255 }  ,draw opacity=1 ][fill={rgb, 255:red, 255; green, 255; blue, 255 }  ,fill opacity=1 ] (290,51) -- (300,51) -- (300,61) -- (290,61) -- cycle ;
\draw    (290,51) -- (290,61) ;
\draw    (300,51) -- (300,61) ;

\draw  [color={rgb, 255:red, 255; green, 255; blue, 255 }  ,draw opacity=1 ][fill={rgb, 255:red, 255; green, 255; blue, 255 }  ,fill opacity=1 ] (290,21) -- (300,21) -- (300,31) -- (290,31) -- cycle ;
\draw    (290,21) -- (290,31) ;
\draw    (300,21) -- (300,31) ;

\draw  [color={rgb, 255:red, 255; green, 255; blue, 255 }  ,draw opacity=1 ][fill={rgb, 255:red, 255; green, 255; blue, 255 }  ,fill opacity=1 ] (410,140) -- (410,150) -- (400,150) -- (400,140) -- cycle ;
\draw    (410,140) -- (400,140) ;
\draw    (410,150) -- (400,150) ;

\draw  [color={rgb, 255:red, 255; green, 255; blue, 255 }  ,draw opacity=1 ][fill={rgb, 255:red, 255; green, 255; blue, 255 }  ,fill opacity=1 ] (400,110) -- (410,110) -- (410,120) -- (400,120) -- cycle ;
\draw    (400,110) -- (400,120) ;
\draw    (410,110) -- (410,120) ;

\draw  [color={rgb, 255:red, 255; green, 255; blue, 255 }  ,draw opacity=1 ][fill={rgb, 255:red, 255; green, 255; blue, 255 }  ,fill opacity=1 ] (190,110) -- (190,120) -- (180,120) -- (180,110) -- cycle ;
\draw    (190,110) -- (180,110) ;
\draw    (190,120) -- (180,120) ;

\draw  [color={rgb, 255:red, 255; green, 255; blue, 255 }  ,draw opacity=1 ][fill={rgb, 255:red, 255; green, 255; blue, 255 }  ,fill opacity=1 ] (180,140) -- (190,140) -- (190,150) -- (180,150) -- cycle ;
\draw    (180,140) -- (180,150) ;
\draw    (190,140) -- (190,150) ;

\draw  [color={rgb, 255:red, 255; green, 255; blue, 255 }  ,draw opacity=1 ][fill={rgb, 255:red, 255; green, 255; blue, 255 }  ,fill opacity=1 ] (300,200) -- (300,210) -- (290,210) -- (290,200) -- cycle ;
\draw    (300,200) -- (290,200) ;
\draw    (300,210) -- (290,210) ;

\draw  [color={rgb, 255:red, 255; green, 255; blue, 255 }  ,draw opacity=1 ][fill={rgb, 255:red, 255; green, 255; blue, 255 }  ,fill opacity=1 ] (300,230) -- (300,240) -- (290,240) -- (290,230) -- cycle ;
\draw    (300,230) -- (290,230) ;
\draw    (300,240) -- (290,240) ;

\draw    (70,298) -- (258,298) ;
\draw [shift={(260,298)}, rotate = 180] [color={rgb, 255:red, 0; green, 0; blue, 0 }  ][line width=0.75]    (10.93,-3.29) .. controls (6.95,-1.4) and (3.31,-0.3) .. (0,0) .. controls (3.31,0.3) and (6.95,1.4) .. (10.93,3.29)   ;
\draw    (72,308) -- (260,308) ;
\draw [shift={(70,308)}, rotate = 0] [color={rgb, 255:red, 0; green, 0; blue, 0 }  ][line width=0.75]    (10.93,-3.29) .. controls (6.95,-1.4) and (3.31,-0.3) .. (0,0) .. controls (3.31,0.3) and (6.95,1.4) .. (10.93,3.29)   ;
\draw [color={rgb, 255:red, 155; green, 155; blue, 155 }  ,draw opacity=1 ]   (70.72,142.07) .. controls (94.48,210.25) and (118.93,250.43) .. (260,250) ;
\draw [shift={(260,250)}, rotate = 179.83] [color={rgb, 255:red, 155; green, 155; blue, 155 }  ,draw opacity=1 ][line width=0.75]    (10.93,-3.29) .. controls (6.95,-1.4) and (3.31,-0.3) .. (0,0) .. controls (3.31,0.3) and (6.95,1.4) .. (10.93,3.29)   ;
\draw [shift={(70,140)}, rotate = 70.93] [color={rgb, 255:red, 155; green, 155; blue, 155 }  ,draw opacity=1 ][line width=0.75]    (10.93,-3.29) .. controls (6.95,-1.4) and (3.31,-0.3) .. (0,0) .. controls (3.31,0.3) and (6.95,1.4) .. (10.93,3.29)   ;
\draw    (220,90) -- (258.59,51.41) ;
\draw [shift={(260,50)}, rotate = 135] [color={rgb, 255:red, 0; green, 0; blue, 0 }  ][line width=0.75]    (10.93,-3.29) .. controls (6.95,-1.4) and (3.31,-0.3) .. (0,0) .. controls (3.31,0.3) and (6.95,1.4) .. (10.93,3.29)   ;
\draw    (330,50) -- (368.59,88.59) ;
\draw [shift={(370,90)}, rotate = 225] [color={rgb, 255:red, 0; green, 0; blue, 0 }  ][line width=0.75]    (10.93,-3.29) .. controls (6.95,-1.4) and (3.31,-0.3) .. (0,0) .. controls (3.31,0.3) and (6.95,1.4) .. (10.93,3.29)   ;
\draw    (330,220) -- (368.59,181.41) ;
\draw [shift={(370,180)}, rotate = 135] [color={rgb, 255:red, 0; green, 0; blue, 0 }  ][line width=0.75]    (10.93,-3.29) .. controls (6.95,-1.4) and (3.31,-0.3) .. (0,0) .. controls (3.31,0.3) and (6.95,1.4) .. (10.93,3.29)   ;
\draw    (220,180) -- (258.59,218.59) ;
\draw [shift={(260,220)}, rotate = 225] [color={rgb, 255:red, 0; green, 0; blue, 0 }  ][line width=0.75]    (10.93,-3.29) .. controls (6.95,-1.4) and (3.31,-0.3) .. (0,0) .. controls (3.31,0.3) and (6.95,1.4) .. (10.93,3.29)   ;
\draw  [dash pattern={on 4.5pt off 4.5pt}]  (221.41,98.59) -- (260,60) ;
\draw [shift={(220,100)}, rotate = 315] [color={rgb, 255:red, 0; green, 0; blue, 0 }  ][line width=0.75]    (10.93,-3.29) .. controls (6.95,-1.4) and (3.31,-0.3) .. (0,0) .. controls (3.31,0.3) and (6.95,1.4) .. (10.93,3.29)   ;
\draw  [dash pattern={on 4.5pt off 4.5pt}]  (331.41,61.41) -- (370,100) ;
\draw [shift={(330,60)}, rotate = 45] [color={rgb, 255:red, 0; green, 0; blue, 0 }  ][line width=0.75]    (10.93,-3.29) .. controls (6.95,-1.4) and (3.31,-0.3) .. (0,0) .. controls (3.31,0.3) and (6.95,1.4) .. (10.93,3.29)   ;
\draw  [dash pattern={on 4.5pt off 4.5pt}]  (331.41,208.59) -- (370,170) ;
\draw [shift={(330,210)}, rotate = 315] [color={rgb, 255:red, 0; green, 0; blue, 0 }  ][line width=0.75]    (10.93,-3.29) .. controls (6.95,-1.4) and (3.31,-0.3) .. (0,0) .. controls (3.31,0.3) and (6.95,1.4) .. (10.93,3.29)   ;
\draw  [dash pattern={on 4.5pt off 4.5pt}]  (221.41,171.41) -- (260,210) ;
\draw [shift={(220,170)}, rotate = 45] [color={rgb, 255:red, 0; green, 0; blue, 0 }  ][line width=0.75]    (10.93,-3.29) .. controls (6.95,-1.4) and (3.31,-0.3) .. (0,0) .. controls (3.31,0.3) and (6.95,1.4) .. (10.93,3.29)   ;
\draw [color={rgb, 255:red, 0; green, 0; blue, 0 }  ,draw opacity=1 ] [dash pattern={on 4.5pt off 4.5pt}]  (401,180) .. controls (400.67,219.08) and (373.5,238.43) .. (332.87,239.94) ;
\draw [shift={(331,240)}, rotate = 358.64] [color={rgb, 255:red, 0; green, 0; blue, 0 }  ,draw opacity=1 ][line width=0.75]    (10.93,-3.29) .. controls (6.95,-1.4) and (3.31,-0.3) .. (0,0) .. controls (3.31,0.3) and (6.95,1.4) .. (10.93,3.29)   ;
\draw [color={rgb, 255:red, 0; green, 0; blue, 0 }  ,draw opacity=1 ] [dash pattern={on 4.5pt off 4.5pt}]  (260,240) .. controls (226.84,240.34) and (199.17,233.23) .. (199.94,191.91) ;
\draw [shift={(200,190)}, rotate = 92.22] [color={rgb, 255:red, 0; green, 0; blue, 0 }  ,draw opacity=1 ][line width=0.75]    (10.93,-3.29) .. controls (6.95,-1.4) and (3.31,-0.3) .. (0,0) .. controls (3.31,0.3) and (6.95,1.4) .. (10.93,3.29)   ;
\draw [color={rgb, 255:red, 155; green, 155; blue, 155 }  ,draw opacity=1 ]   (70.72,117.93) .. controls (94.48,49.95) and (118.93,19.57) .. (260,20) ;
\draw [shift={(260,20)}, rotate = 180.17] [color={rgb, 255:red, 155; green, 155; blue, 155 }  ,draw opacity=1 ][line width=0.75]    (10.93,-3.29) .. controls (6.95,-1.4) and (3.31,-0.3) .. (0,0) .. controls (3.31,0.3) and (6.95,1.4) .. (10.93,3.29)   ;
\draw [shift={(70,120)}, rotate = 289.07] [color={rgb, 255:red, 155; green, 155; blue, 155 }  ,draw opacity=1 ][line width=0.75]    (10.93,-3.29) .. controls (6.95,-1.4) and (3.31,-0.3) .. (0,0) .. controls (3.31,0.3) and (6.95,1.4) .. (10.93,3.29)   ;
\draw  [draw opacity=0] (320.88,293.16) .. controls (323.61,289.42) and (328.02,287) .. (333,287) .. controls (341.28,287) and (348,293.72) .. (348,302) .. controls (348,310.28) and (341.28,317) .. (333,317) .. controls (328.44,317) and (324.35,314.96) .. (321.6,311.75) -- (333,302) -- cycle ; \draw   (320.88,293.16) .. controls (323.61,289.42) and (328.02,287) .. (333,287) .. controls (341.28,287) and (348,293.72) .. (348,302) .. controls (348,310.28) and (341.28,317) .. (333,317) .. controls (328.44,317) and (324.35,314.96) .. (321.6,311.75) ;  
\draw    (320.88,293.16) -- (319.2,295.4) ;
\draw [shift={(318,297)}, rotate = 306.87] [color={rgb, 255:red, 0; green, 0; blue, 0 }  ][line width=0.75]    (6.56,-1.97) .. controls (4.17,-0.84) and (1.99,-0.18) .. (0,0) .. controls (1.99,0.18) and (4.17,0.84) .. (6.56,1.97)   ;

\draw (32.5,292.4) node [anchor=north] [inner sep=0.75pt]    {$[ T]$};
\draw (295.5,292.4) node [anchor=north] [inner sep=0.75pt]    {$[ T']$};
\draw (186,168.4) node [anchor=north] [inner sep=0.75pt]  [font=\scriptsize]  {$00$};
\draw (405.5,169.4) node [anchor=north] [inner sep=0.75pt]  [font=\scriptsize]  {$11$};
\draw (295.5,258.4) node [anchor=north] [inner sep=0.75pt]  [font=\scriptsize]  {$10$};
\draw (295.5,75.4) node [anchor=north] [inner sep=0.75pt]  [font=\scriptsize]  {$01$};
\draw (238,66.6) node [anchor=south east] [inner sep=0.75pt]  [font=\scriptsize]  {$e_{2}$};
\draw (352,66.6) node [anchor=south west] [inner sep=0.75pt]  [font=\scriptsize]  {$e_{1}$};
\draw (352,203.4) node [anchor=north west][inner sep=0.75pt]  [font=\scriptsize]  {$-e_{2}$};
\draw (238,203.4) node [anchor=north east] [inner sep=0.75pt]  [font=\scriptsize]  {$e_{1}$};
\draw (242,83.4) node [anchor=north west][inner sep=0.75pt]  [font=\scriptsize]  {$\chi _{2}$};
\draw (348,83.4) node [anchor=north east] [inner sep=0.75pt]  [font=\scriptsize]  {$\chi _{1}$};
\draw (348,186.6) node [anchor=south east] [inner sep=0.75pt]  [font=\scriptsize]  {$-\chi _{2}$};
\draw (242,186.6) node [anchor=south west] [inner sep=0.75pt]  [font=\scriptsize]  {$\chi _{1}$};
\draw (392,216.4) node [anchor=north west][inner sep=0.75pt]  [font=\footnotesize,color={rgb, 255:red, 0; green, 0; blue, 0 }  ,opacity=1 ]  {$h=\iota $};
\draw (199,216.4) node [anchor=north east] [inner sep=0.75pt]  [font=\footnotesize,color={rgb, 255:red, 0; green, 0; blue, 0 }  ,opacity=1 ]  {$h\ =\ -\epsilon $};
\draw (165,294.6) node [anchor=south] [inner sep=0.75pt]    {$f$};
\draw (165,311.4) node [anchor=north] [inner sep=0.75pt]    {$g$};
\draw (351,303) node [anchor=west] [inner sep=0.75pt]    {$h$};

\end{tikzpicture}

%% file: tikzpictures/R3.tex
\tikzset{every picture/.style={line width=0.75pt}} 

\begin{tikzpicture}[x=0.75pt,y=0.75pt,yscale=-1,xscale=1]

\draw    (21.41,21.41) -- (50,50) ;
\draw [shift={(20,20)}, rotate = 45] [color={rgb, 255:red, 0; green, 0; blue, 0 }  ][line width=0.75]    (6.56,-1.97) .. controls (4.17,-0.84) and (1.99,-0.18) .. (0,0) .. controls (1.99,0.18) and (4.17,0.84) .. (6.56,1.97)   ;
\draw  [color={rgb, 255:red, 255; green, 255; blue, 255 }  ,draw opacity=1 ][fill={rgb, 255:red, 255; green, 255; blue, 255 }  ,fill opacity=1 ] (30,35) .. controls (30,32.24) and (32.24,30) .. (35,30) .. controls (37.76,30) and (40,32.24) .. (40,35) .. controls (40,37.76) and (37.76,40) .. (35,40) .. controls (32.24,40) and (30,37.76) .. (30,35) -- cycle ;
\draw    (48.59,21.41) -- (20,50) ;
\draw [shift={(50,20)}, rotate = 135] [color={rgb, 255:red, 0; green, 0; blue, 0 }  ][line width=0.75]    (6.56,-1.97) .. controls (4.17,-0.84) and (1.99,-0.18) .. (0,0) .. controls (1.99,0.18) and (4.17,0.84) .. (6.56,1.97)   ;

\draw    (20,80) -- (50,110) ;
\draw  [color={rgb, 255:red, 255; green, 255; blue, 255 }  ,draw opacity=1 ][fill={rgb, 255:red, 255; green, 255; blue, 255 }  ,fill opacity=1 ] (30,95) .. controls (30,92.24) and (32.24,90) .. (35,90) .. controls (37.76,90) and (40,92.24) .. (40,95) .. controls (40,97.76) and (37.76,100) .. (35,100) .. controls (32.24,100) and (30,97.76) .. (30,95) -- cycle ;
\draw    (50,80) -- (20,110) ;

\draw    (50,50) -- (80,80) ;
\draw  [color={rgb, 255:red, 255; green, 255; blue, 255 }  ,draw opacity=1 ][fill={rgb, 255:red, 255; green, 255; blue, 255 }  ,fill opacity=1 ] (60,65) .. controls (60,62.24) and (62.24,60) .. (65,60) .. controls (67.76,60) and (70,62.24) .. (70,65) .. controls (70,67.76) and (67.76,70) .. (65,70) .. controls (62.24,70) and (60,67.76) .. (60,65) -- cycle ;
\draw    (80,50) -- (50,80) ;

\draw    (20,50) -- (20,80) ;
\draw    (80,80) -- (80,110) ;
\draw    (80,22) -- (80,50) ;
\draw [shift={(80,20)}, rotate = 90] [color={rgb, 255:red, 0; green, 0; blue, 0 }  ][line width=0.75]    (6.56,-1.97) .. controls (4.17,-0.84) and (1.99,-0.18) .. (0,0) .. controls (1.99,0.18) and (4.17,0.84) .. (6.56,1.97)   ;
\draw    (205,80) -- (235,110) ;
\draw  [color={rgb, 255:red, 255; green, 255; blue, 255 }  ,draw opacity=1 ][fill={rgb, 255:red, 255; green, 255; blue, 255 }  ,fill opacity=1 ] (215,95) .. controls (215,92.24) and (217.24,90) .. (220,90) .. controls (222.76,90) and (225,92.24) .. (225,95) .. controls (225,97.76) and (222.76,100) .. (220,100) .. controls (217.24,100) and (215,97.76) .. (215,95) -- cycle ;
\draw    (235,80) -- (205,110) ;

\draw    (175,50) -- (205,80) ;
\draw  [color={rgb, 255:red, 255; green, 255; blue, 255 }  ,draw opacity=1 ][fill={rgb, 255:red, 255; green, 255; blue, 255 }  ,fill opacity=1 ] (185,65) .. controls (185,62.24) and (187.24,60) .. (190,60) .. controls (192.76,60) and (195,62.24) .. (195,65) .. controls (195,67.76) and (192.76,70) .. (190,70) .. controls (187.24,70) and (185,67.76) .. (185,65) -- cycle ;
\draw    (205,50) -- (175,80) ;

\draw    (235,50) -- (235,80) ;
\draw    (175,80) -- (175,110) ;
\draw    (175,22) -- (175,50) ;
\draw [shift={(175,20)}, rotate = 90] [color={rgb, 255:red, 0; green, 0; blue, 0 }  ][line width=0.75]    (6.56,-1.97) .. controls (4.17,-0.84) and (1.99,-0.18) .. (0,0) .. controls (1.99,0.18) and (4.17,0.84) .. (6.56,1.97)   ;
\draw    (206.41,21.41) -- (235,50) ;
\draw [shift={(205,20)}, rotate = 45] [color={rgb, 255:red, 0; green, 0; blue, 0 }  ][line width=0.75]    (6.56,-1.97) .. controls (4.17,-0.84) and (1.99,-0.18) .. (0,0) .. controls (1.99,0.18) and (4.17,0.84) .. (6.56,1.97)   ;
\draw  [color={rgb, 255:red, 255; green, 255; blue, 255 }  ,draw opacity=1 ][fill={rgb, 255:red, 255; green, 255; blue, 255 }  ,fill opacity=1 ] (215,35) .. controls (215,32.24) and (217.24,30) .. (220,30) .. controls (222.76,30) and (225,32.24) .. (225,35) .. controls (225,37.76) and (222.76,40) .. (220,40) .. controls (217.24,40) and (215,37.76) .. (215,35) -- cycle ;
\draw    (233.59,21.41) -- (205,50) ;
\draw [shift={(235,20)}, rotate = 135] [color={rgb, 255:red, 0; green, 0; blue, 0 }  ][line width=0.75]    (6.56,-1.97) .. controls (4.17,-0.84) and (1.99,-0.18) .. (0,0) .. controls (1.99,0.18) and (4.17,0.84) .. (6.56,1.97)   ;

\draw (19.55,92) node [anchor=east] [inner sep=0.75pt]  [font=\footnotesize]  {$c_{1}$};
\draw (81,62) node [anchor=west] [inner sep=0.75pt]  [font=\footnotesize]  {$c_{2}$};
\draw (50.5,122.9) node [anchor=north] [inner sep=0.75pt]    {$T$};
\draw (206,122.9) node [anchor=north] [inner sep=0.75pt]    {$T'$};
\draw (19.55,32) node [anchor=east] [inner sep=0.75pt]  [font=\footnotesize]  {$c_{3}$};
\draw (235.45,32) node [anchor=west] [inner sep=0.75pt]  [font=\footnotesize]  {$c'_{1}$};
\draw (174,62) node [anchor=east] [inner sep=0.75pt]  [font=\footnotesize]  {$c'_{2}$};
\draw (235.45,98) node [anchor=west] [inner sep=0.75pt]  [font=\footnotesize]  {$c'_{3}$};

\end{tikzpicture}

%% file: tikzpictures/R3-maps.tex
\tikzset{every picture/.style={line width=0.75pt}} 

\begin{tikzpicture}[x=0.75pt,y=0.75pt,yscale=-.7,xscale=.7]

\draw    (103.4,10) -- (133.4,40) ;
\draw  [color={rgb, 255:red, 255; green, 255; blue, 255 }  ,draw opacity=1 ][fill={rgb, 255:red, 255; green, 255; blue, 255 }  ,fill opacity=1 ] (113.4,25) .. controls (113.4,22.24) and (115.64,20) .. (118.4,20) .. controls (121.16,20) and (123.4,22.24) .. (123.4,25) .. controls (123.4,27.76) and (121.16,30) .. (118.4,30) .. controls (115.64,30) and (113.4,27.76) .. (113.4,25) -- cycle ;
\draw    (133.4,10) -- (103.4,40) ;
\draw    (103.4,70) -- (133.4,100) ;
\draw  [color={rgb, 255:red, 255; green, 255; blue, 255 }  ,draw opacity=1 ][fill={rgb, 255:red, 255; green, 255; blue, 255 }  ,fill opacity=1 ] (113.4,85) .. controls (113.4,82.24) and (115.64,80) .. (118.4,80) .. controls (121.16,80) and (123.4,82.24) .. (123.4,85) .. controls (123.4,87.76) and (121.16,90) .. (118.4,90) .. controls (115.64,90) and (113.4,87.76) .. (113.4,85) -- cycle ;
\draw    (133.4,70) -- (103.4,100) ;

\draw    (133.4,40) -- (163.4,70) ;
\draw  [color={rgb, 255:red, 255; green, 255; blue, 255 }  ,draw opacity=1 ][fill={rgb, 255:red, 255; green, 255; blue, 255 }  ,fill opacity=1 ] (143.4,55) .. controls (143.4,52.24) and (145.64,50) .. (148.4,50) .. controls (151.16,50) and (153.4,52.24) .. (153.4,55) .. controls (153.4,57.76) and (151.16,60) .. (148.4,60) .. controls (145.64,60) and (143.4,57.76) .. (143.4,55) -- cycle ;
\draw    (163.4,40) -- (133.4,70) ;

\draw    (103.4,40) -- (103.4,70) ;
\draw    (163.4,70) -- (163.4,100) ;
\draw    (163.4,10) -- (163.4,40) ;
\draw  [color={rgb, 255:red, 255; green, 255; blue, 255 }  ,draw opacity=1 ][fill={rgb, 255:red, 255; green, 255; blue, 255 }  ,fill opacity=1 ] (113.4,80) -- (123.4,80) -- (123.4,90) -- (113.4,90) -- cycle ;
\draw    (113.4,80) -- (113.4,90) ;
\draw    (123.4,80) -- (123.4,90) ;

\draw    (310,10) -- (340,40) ;
\draw  [color={rgb, 255:red, 255; green, 255; blue, 255 }  ,draw opacity=1 ][fill={rgb, 255:red, 255; green, 255; blue, 255 }  ,fill opacity=1 ] (320,25) .. controls (320,22.24) and (322.24,20) .. (325,20) .. controls (327.76,20) and (330,22.24) .. (330,25) .. controls (330,27.76) and (327.76,30) .. (325,30) .. controls (322.24,30) and (320,27.76) .. (320,25) -- cycle ;
\draw    (340,10) -- (310,40) ;
\draw    (310,70) -- (340,100) ;
\draw  [color={rgb, 255:red, 255; green, 255; blue, 255 }  ,draw opacity=1 ][fill={rgb, 255:red, 255; green, 255; blue, 255 }  ,fill opacity=1 ] (320,85) .. controls (320,82.24) and (322.24,80) .. (325,80) .. controls (327.76,80) and (330,82.24) .. (330,85) .. controls (330,87.76) and (327.76,90) .. (325,90) .. controls (322.24,90) and (320,87.76) .. (320,85) -- cycle ;
\draw    (340,70) -- (310,100) ;

\draw    (340,40) -- (370,70) ;
\draw  [color={rgb, 255:red, 255; green, 255; blue, 255 }  ,draw opacity=1 ][fill={rgb, 255:red, 255; green, 255; blue, 255 }  ,fill opacity=1 ] (350,55) .. controls (350,52.24) and (352.24,50) .. (355,50) .. controls (357.76,50) and (360,52.24) .. (360,55) .. controls (360,57.76) and (357.76,60) .. (355,60) .. controls (352.24,60) and (350,57.76) .. (350,55) -- cycle ;
\draw    (370,40) -- (340,70) ;

\draw    (310,40) -- (310,70) ;
\draw    (370,70) -- (370,100) ;
\draw    (370,10) -- (370,40) ;
\draw  [color={rgb, 255:red, 255; green, 255; blue, 255 }  ,draw opacity=1 ][fill={rgb, 255:red, 255; green, 255; blue, 255 }  ,fill opacity=1 ] (330,80) -- (330,90) -- (320,90) -- (320,80) -- cycle ;
\draw    (330,80) -- (320,80) ;
\draw    (330,90) -- (320,90) ;

\draw    (103.4,190) -- (163.4,270) ;
\draw    (310,190) -- (370,270) ;
\draw    (340,190) -- (340,210) ;
\draw    (370,190) -- (370,210) ;
\draw    (351,170) -- (351,122) ;
\draw [shift={(351,120)}, rotate = 90] [color={rgb, 255:red, 0; green, 0; blue, 0 }  ][line width=0.75]    (10.93,-3.29) .. controls (6.95,-1.4) and (3.31,-0.3) .. (0,0) .. controls (3.31,0.3) and (6.95,1.4) .. (10.93,3.29)   ;
\draw    (341,168) -- (341,120) ;
\draw [shift={(341,170)}, rotate = 270] [color={rgb, 255:red, 0; green, 0; blue, 0 }  ][line width=0.75]    (10.93,-3.29) .. controls (6.95,-1.4) and (3.31,-0.3) .. (0,0) .. controls (3.31,0.3) and (6.95,1.4) .. (10.93,3.29)   ;
\draw  [color={rgb, 255:red, 255; green, 255; blue, 255 }  ,draw opacity=1 ][fill={rgb, 255:red, 255; green, 255; blue, 255 }  ,fill opacity=1 ] (118.4,217) .. controls (118.4,214.24) and (120.64,212) .. (123.4,212) .. controls (126.16,212) and (128.4,214.24) .. (128.4,217) .. controls (128.4,219.76) and (126.16,222) .. (123.4,222) .. controls (120.64,222) and (118.4,219.76) .. (118.4,217) -- cycle ;
\draw  [color={rgb, 255:red, 255; green, 255; blue, 255 }  ,draw opacity=1 ][fill={rgb, 255:red, 255; green, 255; blue, 255 }  ,fill opacity=1 ] (138.4,243) .. controls (138.4,240.24) and (140.64,238) .. (143.4,238) .. controls (146.16,238) and (148.4,240.24) .. (148.4,243) .. controls (148.4,245.76) and (146.16,248) .. (143.4,248) .. controls (140.64,248) and (138.4,245.76) .. (138.4,243) -- cycle ;
\draw    (133.4,190) -- (103.4,270) ;
\draw    (163.4,190) -- (133.4,270) ;
\draw    (310,250) -- (310,270) ;
\draw    (340,250) -- (340,270) ;
\draw    (330,240) -- (340,250) ;
\draw    (340,210) -- (350,220) ;
\draw    (360,220) -- (370,210) ;
\draw    (310,250) -- (320,240) ;
\draw    (320,240) -- (330,240) ;
\draw    (350,220) -- (360,220) ;
\draw    (131.4,120) -- (131.4,170) ;
\draw    (127.4,120) -- (127.4,170) ;

\draw    (200,50) -- (268,50) ;
\draw [shift={(270,50)}, rotate = 180] [color={rgb, 255:red, 0; green, 0; blue, 0 }  ][line width=0.75]    (6.56,-1.97) .. controls (4.17,-0.84) and (1.99,-0.18) .. (0,0) .. controls (1.99,0.18) and (4.17,0.84) .. (6.56,1.97)   ;
\draw  [dash pattern={on 4.5pt off 4.5pt}]  (202,60) -- (270,60) ;
\draw [shift={(200,60)}, rotate = 0] [color={rgb, 255:red, 0; green, 0; blue, 0 }  ][line width=0.75]    (6.56,-1.97) .. controls (4.17,-0.84) and (1.99,-0.18) .. (0,0) .. controls (1.99,0.18) and (4.17,0.84) .. (6.56,1.97)   ;
\draw    (200,229) -- (268,229) ;
\draw [shift={(270,229)}, rotate = 180] [color={rgb, 255:red, 0; green, 0; blue, 0 }  ][line width=0.75]    (6.56,-1.97) .. controls (4.17,-0.84) and (1.99,-0.18) .. (0,0) .. controls (1.99,0.18) and (4.17,0.84) .. (6.56,1.97)   ;
\draw  [dash pattern={on 4.5pt off 4.5pt}]  (202,239) -- (270,239) ;
\draw [shift={(200,239)}, rotate = 0] [color={rgb, 255:red, 0; green, 0; blue, 0 }  ][line width=0.75]    (6.56,-1.97) .. controls (4.17,-0.84) and (1.99,-0.18) .. (0,0) .. controls (1.99,0.18) and (4.17,0.84) .. (6.56,1.97)   ;
\draw  [draw opacity=0] (382.88,46.16) .. controls (385.61,42.42) and (390.02,40) .. (395,40) .. controls (403.28,40) and (410,46.72) .. (410,55) .. controls (410,63.28) and (403.28,70) .. (395,70) .. controls (390.44,70) and (386.35,67.96) .. (383.6,64.75) -- (395,55) -- cycle ; \draw   (382.88,46.16) .. controls (385.61,42.42) and (390.02,40) .. (395,40) .. controls (403.28,40) and (410,46.72) .. (410,55) .. controls (410,63.28) and (403.28,70) .. (395,70) .. controls (390.44,70) and (386.35,67.96) .. (383.6,64.75) ;  
\draw    (382.88,46.16) -- (381.2,48.4) ;
\draw [shift={(380,50)}, rotate = 306.87] [color={rgb, 255:red, 0; green, 0; blue, 0 }  ][line width=0.75]    (6.56,-1.97) .. controls (4.17,-0.84) and (1.99,-0.18) .. (0,0) .. controls (1.99,0.18) and (4.17,0.84) .. (6.56,1.97)   ;
\draw    (31.5,80) -- (31.5,198) ;
\draw [shift={(31.5,200)}, rotate = 270] [color={rgb, 255:red, 0; green, 0; blue, 0 }  ][line width=0.75]    (10.93,-3.29) .. controls (6.95,-1.4) and (3.31,-0.3) .. (0,0) .. controls (3.31,0.3) and (6.95,1.4) .. (10.93,3.29)   ;
\draw    (40,82) -- (40,200) ;
\draw [shift={(40,80)}, rotate = 90] [color={rgb, 255:red, 0; green, 0; blue, 0 }  ][line width=0.75]    (10.93,-3.29) .. controls (6.95,-1.4) and (3.31,-0.3) .. (0,0) .. controls (3.31,0.3) and (6.95,1.4) .. (10.93,3.29)   ;

\draw    (196.6,180) -- (278.56,101.38) ;
\draw [shift={(280,100)}, rotate = 136.19] [color={rgb, 255:red, 0; green, 0; blue, 0 }  ][line width=0.75]    (6.56,-1.97) .. controls (4.17,-0.84) and (1.99,-0.18) .. (0,0) .. controls (1.99,0.18) and (4.17,0.84) .. (6.56,1.97)   ;

\draw (353,145) node [anchor=west] [inner sep=0.75pt]    {$f$};
\draw (339,145) node [anchor=east] [inner sep=0.75pt]    {$g$};
\draw (233.5,46.6) node [anchor=south] [inner sep=0.75pt]    {$e$};
\draw (233.5,62.4) node [anchor=north] [inner sep=0.75pt]    {$\chi _{1}$};
\draw (233.5,225.6) node [anchor=south] [inner sep=0.75pt]    {$ge$};
\draw (233.5,241.4) node [anchor=north] [inner sep=0.75pt]    {$\chi _{1} f$};
\draw (411,50) node [anchor=west] [inner sep=0.75pt]    {$-h$};
\draw (37,42.4) node [anchor=north] [inner sep=0.75pt]    {$[ T]$};
\draw (37,222.4) node [anchor=north] [inner sep=0.75pt]    {$E$};
\draw (29.5,140) node [anchor=east] [inner sep=0.75pt]    {$G$};
\draw (42,140) node [anchor=west] [inner sep=0.75pt]    {$F$};
\draw (236.3,136.6) node [anchor=south east] [inner sep=0.75pt]    {$-he$};

\end{tikzpicture}

%% file: tikzpictures/R3-tau.tex
\tikzset{every picture/.style={line width=0.75pt}} 

\begin{tikzpicture}[x=0.75pt,y=0.75pt,yscale=-.55,xscale=.55]

\draw [color={rgb, 255:red, 155; green, 155; blue, 155 }  ,draw opacity=1 ]   (430,130) -- (430,248) ;
\draw [shift={(430,250)}, rotate = 270] [color={rgb, 255:red, 155; green, 155; blue, 155 }  ,draw opacity=1 ][line width=0.75]    (10.93,-3.29) .. controls (6.95,-1.4) and (3.31,-0.3) .. (0,0) .. controls (3.31,0.3) and (6.95,1.4) .. (10.93,3.29)   ;
\draw    (40,10) -- (70,40) ;
\draw  [color={rgb, 255:red, 255; green, 255; blue, 255 }  ,draw opacity=1 ][fill={rgb, 255:red, 255; green, 255; blue, 255 }  ,fill opacity=1 ] (50,25) .. controls (50,22.24) and (52.24,20) .. (55,20) .. controls (57.76,20) and (60,22.24) .. (60,25) .. controls (60,27.76) and (57.76,30) .. (55,30) .. controls (52.24,30) and (50,27.76) .. (50,25) -- cycle ;
\draw    (70,10) -- (40,40) ;
\draw    (40,70) -- (70,100) ;
\draw  [color={rgb, 255:red, 255; green, 255; blue, 255 }  ,draw opacity=1 ][fill={rgb, 255:red, 255; green, 255; blue, 255 }  ,fill opacity=1 ] (50,85) .. controls (50,82.24) and (52.24,80) .. (55,80) .. controls (57.76,80) and (60,82.24) .. (60,85) .. controls (60,87.76) and (57.76,90) .. (55,90) .. controls (52.24,90) and (50,87.76) .. (50,85) -- cycle ;
\draw    (70,70) -- (40,100) ;

\draw    (70,40) -- (100,70) ;
\draw  [color={rgb, 255:red, 255; green, 255; blue, 255 }  ,draw opacity=1 ][fill={rgb, 255:red, 255; green, 255; blue, 255 }  ,fill opacity=1 ] (80,55) .. controls (80,52.24) and (82.24,50) .. (85,50) .. controls (87.76,50) and (90,52.24) .. (90,55) .. controls (90,57.76) and (87.76,60) .. (85,60) .. controls (82.24,60) and (80,57.76) .. (80,55) -- cycle ;
\draw    (100,40) -- (70,70) ;

\draw    (40,40) -- (40,70) ;
\draw    (100,70) -- (100,100) ;
\draw    (100,10) -- (100,40) ;
\draw  [color={rgb, 255:red, 255; green, 255; blue, 255 }  ,draw opacity=1 ][fill={rgb, 255:red, 255; green, 255; blue, 255 }  ,fill opacity=1 ] (60,90) -- (50,90) -- (50,80) -- (60,80) -- cycle ;
\draw    (60,90) -- (60,80) ;
\draw    (50,90) -- (50,80) ;

\draw    (100,360) -- (70,330) ;
\draw  [color={rgb, 255:red, 255; green, 255; blue, 255 }  ,draw opacity=1 ][fill={rgb, 255:red, 255; green, 255; blue, 255 }  ,fill opacity=1 ] (90,345) .. controls (90,347.76) and (87.76,350) .. (85,350) .. controls (82.24,350) and (80,347.76) .. (80,345) .. controls (80,342.24) and (82.24,340) .. (85,340) .. controls (87.76,340) and (90,342.24) .. (90,345) -- cycle ;
\draw    (70,360) -- (100,330) ;
\draw    (100,300) -- (70,270) ;
\draw  [color={rgb, 255:red, 255; green, 255; blue, 255 }  ,draw opacity=1 ][fill={rgb, 255:red, 255; green, 255; blue, 255 }  ,fill opacity=1 ] (90,285) .. controls (90,287.76) and (87.76,290) .. (85,290) .. controls (82.24,290) and (80,287.76) .. (80,285) .. controls (80,282.24) and (82.24,280) .. (85,280) .. controls (87.76,280) and (90,282.24) .. (90,285) -- cycle ;
\draw    (70,300) -- (100,270) ;

\draw    (70,330) -- (40,300) ;
\draw  [color={rgb, 255:red, 255; green, 255; blue, 255 }  ,draw opacity=1 ][fill={rgb, 255:red, 255; green, 255; blue, 255 }  ,fill opacity=1 ] (60,315) .. controls (60,317.76) and (57.76,320) .. (55,320) .. controls (52.24,320) and (50,317.76) .. (50,315) .. controls (50,312.24) and (52.24,310) .. (55,310) .. controls (57.76,310) and (60,312.24) .. (60,315) -- cycle ;
\draw    (40,330) -- (70,300) ;

\draw    (100,330) -- (100,300) ;
\draw    (40,300) -- (40,270) ;
\draw    (40,360) -- (40,330) ;
\draw  [color={rgb, 255:red, 255; green, 255; blue, 255 }  ,draw opacity=1 ][fill={rgb, 255:red, 255; green, 255; blue, 255 }  ,fill opacity=1 ] (80,280) -- (90,280) -- (90,290) -- (80,290) -- cycle ;
\draw    (80,280) -- (80,290) ;
\draw    (90,280) -- (90,290) ;

\draw    (190,10) -- (220,40) ;
\draw  [color={rgb, 255:red, 255; green, 255; blue, 255 }  ,draw opacity=1 ][fill={rgb, 255:red, 255; green, 255; blue, 255 }  ,fill opacity=1 ] (200,25) .. controls (200,22.24) and (202.24,20) .. (205,20) .. controls (207.76,20) and (210,22.24) .. (210,25) .. controls (210,27.76) and (207.76,30) .. (205,30) .. controls (202.24,30) and (200,27.76) .. (200,25) -- cycle ;
\draw    (220,10) -- (190,40) ;
\draw    (190,70) -- (220,100) ;
\draw  [color={rgb, 255:red, 255; green, 255; blue, 255 }  ,draw opacity=1 ][fill={rgb, 255:red, 255; green, 255; blue, 255 }  ,fill opacity=1 ] (200,85) .. controls (200,82.24) and (202.24,80) .. (205,80) .. controls (207.76,80) and (210,82.24) .. (210,85) .. controls (210,87.76) and (207.76,90) .. (205,90) .. controls (202.24,90) and (200,87.76) .. (200,85) -- cycle ;
\draw    (220,70) -- (190,100) ;

\draw    (220,40) -- (250,70) ;
\draw  [color={rgb, 255:red, 255; green, 255; blue, 255 }  ,draw opacity=1 ][fill={rgb, 255:red, 255; green, 255; blue, 255 }  ,fill opacity=1 ] (230,55) .. controls (230,52.24) and (232.24,50) .. (235,50) .. controls (237.76,50) and (240,52.24) .. (240,55) .. controls (240,57.76) and (237.76,60) .. (235,60) .. controls (232.24,60) and (230,57.76) .. (230,55) -- cycle ;
\draw    (250,40) -- (220,70) ;

\draw    (190,40) -- (190,70) ;
\draw    (250,70) -- (250,100) ;
\draw    (250,10) -- (250,40) ;
\draw  [color={rgb, 255:red, 255; green, 255; blue, 255 }  ,draw opacity=1 ][fill={rgb, 255:red, 255; green, 255; blue, 255 }  ,fill opacity=1 ] (210,90) -- (200,90) -- (200,80) -- (210,80) -- cycle ;
\draw    (210,90) -- (210,80) ;
\draw    (200,90) -- (200,80) ;

\draw  [color={rgb, 255:red, 255; green, 255; blue, 255 }  ,draw opacity=1 ][fill={rgb, 255:red, 255; green, 255; blue, 255 }  ,fill opacity=1 ] (230,50) -- (240,50) -- (240,60) -- (230,60) -- cycle ;
\draw    (230,50) -- (230,60) ;
\draw    (240,50) -- (240,60) ;

\draw  [color={rgb, 255:red, 255; green, 255; blue, 255 }  ,draw opacity=1 ][fill={rgb, 255:red, 255; green, 255; blue, 255 }  ,fill opacity=1 ] (200,20) -- (210,20) -- (210,30) -- (200,30) -- cycle ;
\draw    (200,20) -- (200,30) ;
\draw    (210,20) -- (210,30) ;

\draw    (250,360) -- (220,330) ;
\draw  [color={rgb, 255:red, 255; green, 255; blue, 255 }  ,draw opacity=1 ][fill={rgb, 255:red, 255; green, 255; blue, 255 }  ,fill opacity=1 ] (240,345) .. controls (240,347.76) and (237.76,350) .. (235,350) .. controls (232.24,350) and (230,347.76) .. (230,345) .. controls (230,342.24) and (232.24,340) .. (235,340) .. controls (237.76,340) and (240,342.24) .. (240,345) -- cycle ;
\draw    (220,360) -- (250,330) ;
\draw    (250,300) -- (220,270) ;
\draw  [color={rgb, 255:red, 255; green, 255; blue, 255 }  ,draw opacity=1 ][fill={rgb, 255:red, 255; green, 255; blue, 255 }  ,fill opacity=1 ] (240,285) .. controls (240,287.76) and (237.76,290) .. (235,290) .. controls (232.24,290) and (230,287.76) .. (230,285) .. controls (230,282.24) and (232.24,280) .. (235,280) .. controls (237.76,280) and (240,282.24) .. (240,285) -- cycle ;
\draw    (220,300) -- (250,270) ;

\draw    (220,330) -- (190,300) ;
\draw  [color={rgb, 255:red, 255; green, 255; blue, 255 }  ,draw opacity=1 ][fill={rgb, 255:red, 255; green, 255; blue, 255 }  ,fill opacity=1 ] (210,315) .. controls (210,317.76) and (207.76,320) .. (205,320) .. controls (202.24,320) and (200,317.76) .. (200,315) .. controls (200,312.24) and (202.24,310) .. (205,310) .. controls (207.76,310) and (210,312.24) .. (210,315) -- cycle ;
\draw    (190,330) -- (220,300) ;

\draw    (250,330) -- (250,300) ;
\draw    (190,300) -- (190,270) ;
\draw    (190,360) -- (190,330) ;
\draw  [color={rgb, 255:red, 255; green, 255; blue, 255 }  ,draw opacity=1 ][fill={rgb, 255:red, 255; green, 255; blue, 255 }  ,fill opacity=1 ] (230,280) -- (240,280) -- (240,290) -- (230,290) -- cycle ;
\draw    (230,280) -- (230,290) ;
\draw    (240,280) -- (240,290) ;

\draw  [color={rgb, 255:red, 255; green, 255; blue, 255 }  ,draw opacity=1 ][fill={rgb, 255:red, 255; green, 255; blue, 255 }  ,fill opacity=1 ] (230,340) -- (240,340) -- (240,350) -- (230,350) -- cycle ;
\draw    (230,340) -- (230,350) ;
\draw    (240,340) -- (240,350) ;

\draw  [color={rgb, 255:red, 255; green, 255; blue, 255 }  ,draw opacity=1 ][fill={rgb, 255:red, 255; green, 255; blue, 255 }  ,fill opacity=1 ] (200,310) -- (210,310) -- (210,320) -- (200,320) -- cycle ;
\draw    (200,310) -- (200,320) ;
\draw    (210,310) -- (210,320) ;

\draw    (400,10) -- (430,40) ;
\draw  [color={rgb, 255:red, 255; green, 255; blue, 255 }  ,draw opacity=1 ][fill={rgb, 255:red, 255; green, 255; blue, 255 }  ,fill opacity=1 ] (410,25) .. controls (410,22.24) and (412.24,20) .. (415,20) .. controls (417.76,20) and (420,22.24) .. (420,25) .. controls (420,27.76) and (417.76,30) .. (415,30) .. controls (412.24,30) and (410,27.76) .. (410,25) -- cycle ;
\draw    (430,10) -- (400,40) ;
\draw    (400,70) -- (430,100) ;
\draw  [color={rgb, 255:red, 255; green, 255; blue, 255 }  ,draw opacity=1 ][fill={rgb, 255:red, 255; green, 255; blue, 255 }  ,fill opacity=1 ] (410,85) .. controls (410,82.24) and (412.24,80) .. (415,80) .. controls (417.76,80) and (420,82.24) .. (420,85) .. controls (420,87.76) and (417.76,90) .. (415,90) .. controls (412.24,90) and (410,87.76) .. (410,85) -- cycle ;
\draw    (430,70) -- (400,100) ;

\draw    (430,40) -- (460,70) ;
\draw  [color={rgb, 255:red, 255; green, 255; blue, 255 }  ,draw opacity=1 ][fill={rgb, 255:red, 255; green, 255; blue, 255 }  ,fill opacity=1 ] (440,55) .. controls (440,52.24) and (442.24,50) .. (445,50) .. controls (447.76,50) and (450,52.24) .. (450,55) .. controls (450,57.76) and (447.76,60) .. (445,60) .. controls (442.24,60) and (440,57.76) .. (440,55) -- cycle ;
\draw    (460,40) -- (430,70) ;

\draw    (460,70) -- (460,100) ;
\draw    (460,10) -- (460,40) ;
\draw  [color={rgb, 255:red, 255; green, 255; blue, 255 }  ,draw opacity=1 ][fill={rgb, 255:red, 255; green, 255; blue, 255 }  ,fill opacity=1 ] (420,90) -- (410,90) -- (410,80) -- (420,80) -- cycle ;
\draw    (420,90) -- (420,80) ;
\draw    (410,90) -- (410,80) ;

\draw  [color={rgb, 255:red, 255; green, 255; blue, 255 }  ,draw opacity=1 ][fill={rgb, 255:red, 255; green, 255; blue, 255 }  ,fill opacity=1 ] (450,50) -- (450,60) -- (440,60) -- (440,50) -- cycle ;
\draw    (450,50) -- (440,50) ;
\draw    (450,60) -- (440,60) ;

\draw  [color={rgb, 255:red, 255; green, 255; blue, 255 }  ,draw opacity=1 ][fill={rgb, 255:red, 255; green, 255; blue, 255 }  ,fill opacity=1 ] (410,20) -- (420,20) -- (420,30) -- (410,30) -- cycle ;
\draw    (410,20) -- (410,30) ;
\draw    (420,20) -- (420,30) ;

\draw    (400,40) -- (400,70) ;

\draw    (460,360) -- (430,330) ;
\draw  [color={rgb, 255:red, 255; green, 255; blue, 255 }  ,draw opacity=1 ][fill={rgb, 255:red, 255; green, 255; blue, 255 }  ,fill opacity=1 ] (450,345) .. controls (450,347.76) and (447.76,350) .. (445,350) .. controls (442.24,350) and (440,347.76) .. (440,345) .. controls (440,342.24) and (442.24,340) .. (445,340) .. controls (447.76,340) and (450,342.24) .. (450,345) -- cycle ;
\draw    (430,360) -- (460,330) ;
\draw    (460,300) -- (430,270) ;
\draw  [color={rgb, 255:red, 255; green, 255; blue, 255 }  ,draw opacity=1 ][fill={rgb, 255:red, 255; green, 255; blue, 255 }  ,fill opacity=1 ] (450,285) .. controls (450,287.76) and (447.76,290) .. (445,290) .. controls (442.24,290) and (440,287.76) .. (440,285) .. controls (440,282.24) and (442.24,280) .. (445,280) .. controls (447.76,280) and (450,282.24) .. (450,285) -- cycle ;
\draw    (430,300) -- (460,270) ;

\draw    (430,330) -- (400,300) ;
\draw  [color={rgb, 255:red, 255; green, 255; blue, 255 }  ,draw opacity=1 ][fill={rgb, 255:red, 255; green, 255; blue, 255 }  ,fill opacity=1 ] (420,315) .. controls (420,317.76) and (417.76,320) .. (415,320) .. controls (412.24,320) and (410,317.76) .. (410,315) .. controls (410,312.24) and (412.24,310) .. (415,310) .. controls (417.76,310) and (420,312.24) .. (420,315) -- cycle ;
\draw    (400,330) -- (430,300) ;

\draw    (460,330) -- (460,300) ;
\draw    (400,300) -- (400,270) ;
\draw    (400,360) -- (400,330) ;
\draw  [color={rgb, 255:red, 255; green, 255; blue, 255 }  ,draw opacity=1 ][fill={rgb, 255:red, 255; green, 255; blue, 255 }  ,fill opacity=1 ] (440,280) -- (450,280) -- (450,290) -- (440,290) -- cycle ;
\draw    (440,280) -- (440,290) ;
\draw    (450,280) -- (450,290) ;

\draw  [color={rgb, 255:red, 255; green, 255; blue, 255 }  ,draw opacity=1 ][fill={rgb, 255:red, 255; green, 255; blue, 255 }  ,fill opacity=1 ] (450,340) -- (450,350) -- (440,350) -- (440,340) -- cycle ;
\draw    (450,340) -- (440,340) ;
\draw    (450,350) -- (440,350) ;

\draw  [color={rgb, 255:red, 255; green, 255; blue, 255 }  ,draw opacity=1 ][fill={rgb, 255:red, 255; green, 255; blue, 255 }  ,fill opacity=1 ] (410,310) -- (420,310) -- (420,320) -- (410,320) -- cycle ;
\draw    (410,310) -- (410,320) ;
\draw    (420,310) -- (420,320) ;

\draw    (270,140) -- (300,170) ;
\draw  [color={rgb, 255:red, 255; green, 255; blue, 255 }  ,draw opacity=1 ][fill={rgb, 255:red, 255; green, 255; blue, 255 }  ,fill opacity=1 ] (280,155) .. controls (280,152.24) and (282.24,150) .. (285,150) .. controls (287.76,150) and (290,152.24) .. (290,155) .. controls (290,157.76) and (287.76,160) .. (285,160) .. controls (282.24,160) and (280,157.76) .. (280,155) -- cycle ;
\draw    (300,140) -- (270,170) ;
\draw    (270,200) -- (300,230) ;
\draw  [color={rgb, 255:red, 255; green, 255; blue, 255 }  ,draw opacity=1 ][fill={rgb, 255:red, 255; green, 255; blue, 255 }  ,fill opacity=1 ] (280,215) .. controls (280,212.24) and (282.24,210) .. (285,210) .. controls (287.76,210) and (290,212.24) .. (290,215) .. controls (290,217.76) and (287.76,220) .. (285,220) .. controls (282.24,220) and (280,217.76) .. (280,215) -- cycle ;
\draw    (300,200) -- (270,230) ;

\draw    (300,170) -- (330,200) ;
\draw  [color={rgb, 255:red, 255; green, 255; blue, 255 }  ,draw opacity=1 ][fill={rgb, 255:red, 255; green, 255; blue, 255 }  ,fill opacity=1 ] (310,185) .. controls (310,182.24) and (312.24,180) .. (315,180) .. controls (317.76,180) and (320,182.24) .. (320,185) .. controls (320,187.76) and (317.76,190) .. (315,190) .. controls (312.24,190) and (310,187.76) .. (310,185) -- cycle ;
\draw    (330,170) -- (300,200) ;

\draw    (270,170) -- (270,200) ;
\draw    (330,200) -- (330,230) ;
\draw    (330,140) -- (330,170) ;
\draw  [color={rgb, 255:red, 255; green, 255; blue, 255 }  ,draw opacity=1 ][fill={rgb, 255:red, 255; green, 255; blue, 255 }  ,fill opacity=1 ] (290,220) -- (280,220) -- (280,210) -- (290,210) -- cycle ;
\draw    (290,220) -- (290,210) ;
\draw    (280,220) -- (280,210) ;

\draw  [color={rgb, 255:red, 255; green, 255; blue, 255 }  ,draw opacity=1 ][fill={rgb, 255:red, 255; green, 255; blue, 255 }  ,fill opacity=1 ] (310,180) -- (320,180) -- (320,190) -- (310,190) -- cycle ;
\draw    (310,180) -- (310,190) ;
\draw    (320,180) -- (320,190) ;

\draw  [color={rgb, 255:red, 255; green, 255; blue, 255 }  ,draw opacity=1 ][fill={rgb, 255:red, 255; green, 255; blue, 255 }  ,fill opacity=1 ] (290,150) -- (290,160) -- (280,160) -- (280,150) -- cycle ;
\draw    (290,150) -- (280,150) ;
\draw    (290,160) -- (280,160) ;

\draw    (330,490) -- (300,460) ;
\draw  [color={rgb, 255:red, 255; green, 255; blue, 255 }  ,draw opacity=1 ][fill={rgb, 255:red, 255; green, 255; blue, 255 }  ,fill opacity=1 ] (320,475) .. controls (320,477.76) and (317.76,480) .. (315,480) .. controls (312.24,480) and (310,477.76) .. (310,475) .. controls (310,472.24) and (312.24,470) .. (315,470) .. controls (317.76,470) and (320,472.24) .. (320,475) -- cycle ;
\draw    (300,490) -- (330,460) ;
\draw    (330,430) -- (300,400) ;
\draw  [color={rgb, 255:red, 255; green, 255; blue, 255 }  ,draw opacity=1 ][fill={rgb, 255:red, 255; green, 255; blue, 255 }  ,fill opacity=1 ] (320,415) .. controls (320,417.76) and (317.76,420) .. (315,420) .. controls (312.24,420) and (310,417.76) .. (310,415) .. controls (310,412.24) and (312.24,410) .. (315,410) .. controls (317.76,410) and (320,412.24) .. (320,415) -- cycle ;
\draw    (300,430) -- (330,400) ;

\draw    (300,460) -- (270,430) ;
\draw  [color={rgb, 255:red, 255; green, 255; blue, 255 }  ,draw opacity=1 ][fill={rgb, 255:red, 255; green, 255; blue, 255 }  ,fill opacity=1 ] (290,445) .. controls (290,447.76) and (287.76,450) .. (285,450) .. controls (282.24,450) and (280,447.76) .. (280,445) .. controls (280,442.24) and (282.24,440) .. (285,440) .. controls (287.76,440) and (290,442.24) .. (290,445) -- cycle ;
\draw    (270,460) -- (300,430) ;

\draw    (330,460) -- (330,430) ;
\draw  [color={rgb, 255:red, 255; green, 255; blue, 255 }  ,draw opacity=1 ][fill={rgb, 255:red, 255; green, 255; blue, 255 }  ,fill opacity=1 ] (310,410) -- (320,410) -- (320,420) -- (310,420) -- cycle ;
\draw    (310,410) -- (310,420) ;
\draw    (320,410) -- (320,420) ;

\draw  [color={rgb, 255:red, 255; green, 255; blue, 255 }  ,draw opacity=1 ][fill={rgb, 255:red, 255; green, 255; blue, 255 }  ,fill opacity=1 ] (310,470) -- (320,470) -- (320,480) -- (310,480) -- cycle ;
\draw    (310,470) -- (310,480) ;
\draw    (320,470) -- (320,480) ;

\draw  [color={rgb, 255:red, 255; green, 255; blue, 255 }  ,draw opacity=1 ][fill={rgb, 255:red, 255; green, 255; blue, 255 }  ,fill opacity=1 ] (290,440) -- (290,450) -- (280,450) -- (280,440) -- cycle ;
\draw    (290,440) -- (280,440) ;
\draw    (290,450) -- (280,450) ;

\draw    (270,430) -- (270,400) ;
\draw    (270,490) -- (270,460) ;
\draw    (480,140) -- (510,170) ;
\draw  [color={rgb, 255:red, 255; green, 255; blue, 255 }  ,draw opacity=1 ][fill={rgb, 255:red, 255; green, 255; blue, 255 }  ,fill opacity=1 ] (490,155) .. controls (490,152.24) and (492.24,150) .. (495,150) .. controls (497.76,150) and (500,152.24) .. (500,155) .. controls (500,157.76) and (497.76,160) .. (495,160) .. controls (492.24,160) and (490,157.76) .. (490,155) -- cycle ;
\draw    (510,140) -- (480,170) ;
\draw    (480,200) -- (510,230) ;
\draw  [color={rgb, 255:red, 255; green, 255; blue, 255 }  ,draw opacity=1 ][fill={rgb, 255:red, 255; green, 255; blue, 255 }  ,fill opacity=1 ] (490,215) .. controls (490,212.24) and (492.24,210) .. (495,210) .. controls (497.76,210) and (500,212.24) .. (500,215) .. controls (500,217.76) and (497.76,220) .. (495,220) .. controls (492.24,220) and (490,217.76) .. (490,215) -- cycle ;
\draw    (510,200) -- (480,230) ;

\draw    (510,170) -- (540,200) ;
\draw  [color={rgb, 255:red, 255; green, 255; blue, 255 }  ,draw opacity=1 ][fill={rgb, 255:red, 255; green, 255; blue, 255 }  ,fill opacity=1 ] (520,185) .. controls (520,182.24) and (522.24,180) .. (525,180) .. controls (527.76,180) and (530,182.24) .. (530,185) .. controls (530,187.76) and (527.76,190) .. (525,190) .. controls (522.24,190) and (520,187.76) .. (520,185) -- cycle ;
\draw    (540,170) -- (510,200) ;

\draw    (540,200) -- (540,230) ;
\draw    (540,140) -- (540,170) ;
\draw  [color={rgb, 255:red, 255; green, 255; blue, 255 }  ,draw opacity=1 ][fill={rgb, 255:red, 255; green, 255; blue, 255 }  ,fill opacity=1 ] (500,220) -- (490,220) -- (490,210) -- (500,210) -- cycle ;
\draw    (500,220) -- (500,210) ;
\draw    (490,220) -- (490,210) ;

\draw  [color={rgb, 255:red, 255; green, 255; blue, 255 }  ,draw opacity=1 ][fill={rgb, 255:red, 255; green, 255; blue, 255 }  ,fill opacity=1 ] (530,180) -- (530,190) -- (520,190) -- (520,180) -- cycle ;
\draw    (530,180) -- (520,180) ;
\draw    (530,190) -- (520,190) ;

\draw  [color={rgb, 255:red, 255; green, 255; blue, 255 }  ,draw opacity=1 ][fill={rgb, 255:red, 255; green, 255; blue, 255 }  ,fill opacity=1 ] (500,150) -- (500,160) -- (490,160) -- (490,150) -- cycle ;
\draw    (500,150) -- (490,150) ;
\draw    (500,160) -- (490,160) ;

\draw    (480,170) -- (480,200) ;

\draw    (540,490) -- (510,460) ;
\draw  [color={rgb, 255:red, 255; green, 255; blue, 255 }  ,draw opacity=1 ][fill={rgb, 255:red, 255; green, 255; blue, 255 }  ,fill opacity=1 ] (530,475) .. controls (530,477.76) and (527.76,480) .. (525,480) .. controls (522.24,480) and (520,477.76) .. (520,475) .. controls (520,472.24) and (522.24,470) .. (525,470) .. controls (527.76,470) and (530,472.24) .. (530,475) -- cycle ;
\draw    (510,490) -- (540,460) ;
\draw    (540,430) -- (510,400) ;
\draw  [color={rgb, 255:red, 255; green, 255; blue, 255 }  ,draw opacity=1 ][fill={rgb, 255:red, 255; green, 255; blue, 255 }  ,fill opacity=1 ] (530,415) .. controls (530,417.76) and (527.76,420) .. (525,420) .. controls (522.24,420) and (520,417.76) .. (520,415) .. controls (520,412.24) and (522.24,410) .. (525,410) .. controls (527.76,410) and (530,412.24) .. (530,415) -- cycle ;
\draw    (510,430) -- (540,400) ;

\draw    (510,460) -- (480,430) ;
\draw  [color={rgb, 255:red, 255; green, 255; blue, 255 }  ,draw opacity=1 ][fill={rgb, 255:red, 255; green, 255; blue, 255 }  ,fill opacity=1 ] (500,445) .. controls (500,447.76) and (497.76,450) .. (495,450) .. controls (492.24,450) and (490,447.76) .. (490,445) .. controls (490,442.24) and (492.24,440) .. (495,440) .. controls (497.76,440) and (500,442.24) .. (500,445) -- cycle ;
\draw    (480,460) -- (510,430) ;

\draw    (540,460) -- (540,430) ;
\draw    (480,430) -- (480,400) ;
\draw    (480,490) -- (480,460) ;
\draw  [color={rgb, 255:red, 255; green, 255; blue, 255 }  ,draw opacity=1 ][fill={rgb, 255:red, 255; green, 255; blue, 255 }  ,fill opacity=1 ] (520,410) -- (530,410) -- (530,420) -- (520,420) -- cycle ;
\draw    (520,410) -- (520,420) ;
\draw    (530,410) -- (530,420) ;

\draw  [color={rgb, 255:red, 255; green, 255; blue, 255 }  ,draw opacity=1 ][fill={rgb, 255:red, 255; green, 255; blue, 255 }  ,fill opacity=1 ] (530,470) -- (530,480) -- (520,480) -- (520,470) -- cycle ;
\draw    (530,470) -- (520,470) ;
\draw    (530,480) -- (520,480) ;

\draw  [color={rgb, 255:red, 255; green, 255; blue, 255 }  ,draw opacity=1 ][fill={rgb, 255:red, 255; green, 255; blue, 255 }  ,fill opacity=1 ] (500,440) -- (500,450) -- (490,450) -- (490,440) -- cycle ;
\draw    (500,440) -- (490,440) ;
\draw    (500,450) -- (490,450) ;

\draw    (70,130) -- (70,238) ;
\draw [shift={(70,240)}, rotate = 270] [color={rgb, 255:red, 0; green, 0; blue, 0 }  ][line width=0.75]    (10.93,-3.29) .. controls (6.95,-1.4) and (3.31,-0.3) .. (0,0) .. controls (3.31,0.3) and (6.95,1.4) .. (10.93,3.29)   ;
\draw    (280,60) -- (368,60) ;
\draw [shift={(370,60)}, rotate = 180] [color={rgb, 255:red, 0; green, 0; blue, 0 }  ][line width=0.75]    (10.93,-3.29) .. controls (6.95,-1.4) and (3.31,-0.3) .. (0,0) .. controls (3.31,0.3) and (6.95,1.4) .. (10.93,3.29)   ;
\draw    (360,190) -- (448,190) ;
\draw [shift={(450,190)}, rotate = 180] [color={rgb, 255:red, 0; green, 0; blue, 0 }  ][line width=0.75]    (10.93,-3.29) .. controls (6.95,-1.4) and (3.31,-0.3) .. (0,0) .. controls (3.31,0.3) and (6.95,1.4) .. (10.93,3.29)   ;
\draw    (280,320) -- (368,320) ;
\draw [shift={(370,320)}, rotate = 180] [color={rgb, 255:red, 0; green, 0; blue, 0 }  ][line width=0.75]    (10.93,-3.29) .. controls (6.95,-1.4) and (3.31,-0.3) .. (0,0) .. controls (3.31,0.3) and (6.95,1.4) .. (10.93,3.29)   ;
\draw    (360,450) -- (448,450) ;
\draw [shift={(450,450)}, rotate = 180] [color={rgb, 255:red, 0; green, 0; blue, 0 }  ][line width=0.75]    (10.93,-3.29) .. controls (6.95,-1.4) and (3.31,-0.3) .. (0,0) .. controls (3.31,0.3) and (6.95,1.4) .. (10.93,3.29)   ;
\draw    (230,130) -- (258.59,158.59) ;
\draw [shift={(260,160)}, rotate = 225] [color={rgb, 255:red, 0; green, 0; blue, 0 }  ][line width=0.75]    (10.93,-3.29) .. controls (6.95,-1.4) and (3.31,-0.3) .. (0,0) .. controls (3.31,0.3) and (6.95,1.4) .. (10.93,3.29)   ;
\draw    (440,130) -- (468.59,158.59) ;
\draw [shift={(470,160)}, rotate = 225] [color={rgb, 255:red, 0; green, 0; blue, 0 }  ][line width=0.75]    (10.93,-3.29) .. controls (6.95,-1.4) and (3.31,-0.3) .. (0,0) .. controls (3.31,0.3) and (6.95,1.4) .. (10.93,3.29)   ;
\draw    (230,390) -- (258.59,418.59) ;
\draw [shift={(260,420)}, rotate = 225] [color={rgb, 255:red, 0; green, 0; blue, 0 }  ][line width=0.75]    (10.93,-3.29) .. controls (6.95,-1.4) and (3.31,-0.3) .. (0,0) .. controls (3.31,0.3) and (6.95,1.4) .. (10.93,3.29)   ;
\draw    (440,390) -- (468.59,418.59) ;
\draw [shift={(470,420)}, rotate = 225] [color={rgb, 255:red, 0; green, 0; blue, 0 }  ][line width=0.75]    (10.93,-3.29) .. controls (6.95,-1.4) and (3.31,-0.3) .. (0,0) .. controls (3.31,0.3) and (6.95,1.4) .. (10.93,3.29)   ;
\draw [color={rgb, 255:red, 155; green, 155; blue, 155 }  ,draw opacity=1 ]   (220,130) -- (220,248) ;
\draw [shift={(220,250)}, rotate = 270] [color={rgb, 255:red, 155; green, 155; blue, 155 }  ,draw opacity=1 ][line width=0.75]    (10.93,-3.29) .. controls (6.95,-1.4) and (3.31,-0.3) .. (0,0) .. controls (3.31,0.3) and (6.95,1.4) .. (10.93,3.29)   ;
\draw [color={rgb, 255:red, 155; green, 155; blue, 155 }  ,draw opacity=1 ]   (510,260) -- (510,378) ;
\draw [shift={(510,380)}, rotate = 270] [color={rgb, 255:red, 155; green, 155; blue, 155 }  ,draw opacity=1 ][line width=0.75]    (10.93,-3.29) .. controls (6.95,-1.4) and (3.31,-0.3) .. (0,0) .. controls (3.31,0.3) and (6.95,1.4) .. (10.93,3.29)   ;
\draw [color={rgb, 255:red, 155; green, 155; blue, 155 }  ,draw opacity=1 ]   (300,260) -- (300,378) ;
\draw [shift={(300,380)}, rotate = 270] [color={rgb, 255:red, 155; green, 155; blue, 155 }  ,draw opacity=1 ][line width=0.75]    (10.93,-3.29) .. controls (6.95,-1.4) and (3.31,-0.3) .. (0,0) .. controls (3.31,0.3) and (6.95,1.4) .. (10.93,3.29)   ;

\draw (221.5,107.4) node [anchor=north] [inner sep=0.75pt]  [font=\footnotesize]  {$00$};
\draw (431.5,107.4) node [anchor=north] [inner sep=0.75pt]  [font=\footnotesize]  {$10$};
\draw (301.5,237.4) node [anchor=north] [inner sep=0.75pt]  [font=\footnotesize]  {$01$};
\draw (511.5,237.4) node [anchor=north] [inner sep=0.75pt]  [font=\footnotesize]  {$11$};
\draw (215.5,368.4) node [anchor=north] [inner sep=0.75pt]  [font=\footnotesize]  {$00$};
\draw (425.5,368.4) node [anchor=north] [inner sep=0.75pt]  [font=\footnotesize]  {$01$};
\draw (295.5,497.4) node [anchor=north] [inner sep=0.75pt]  [font=\footnotesize]  {$10$};
\draw (509.5,497.4) node [anchor=north] [inner sep=0.75pt]  [font=\footnotesize]  {$11$};
\draw (68,185) node [anchor=east] [inner sep=0.75pt]    {$\tau $};
\draw (457,141.6) node [anchor=south west] [inner sep=0.75pt]    {$-$};
\draw (405,446.6) node [anchor=south] [inner sep=0.75pt]    {$-$};
\draw (512,320) node [anchor=west] [inner sep=0.75pt]  [color={rgb, 255:red, 155; green, 155; blue, 155 }  ,opacity=1 ]  {$-$};
\draw (101,54.5) node [anchor=west] [inner sep=0.75pt]  [font=\footnotesize]  {$c_{2}$};
\draw (39,22) node [anchor=east] [inner sep=0.75pt]  [font=\footnotesize]  {$c_{3}$};
\draw (39,310.5) node [anchor=east] [inner sep=0.75pt]  [font=\footnotesize]  {$c'_{2}$};
\draw (100.45,344.5) node [anchor=west] [inner sep=0.75pt]  [font=\footnotesize]  {$c'_{3}$};
\draw (325,56.6) node [anchor=south] [inner sep=0.75pt]    {$e_{2}$};
\draw (325,316.6) node [anchor=south] [inner sep=0.75pt]    {$e'_{3}$};
\draw (244,138.6) node [anchor=south west] [inner sep=0.75pt]    {$e_{3}$};
\draw (243,408.4) node [anchor=north east] [inner sep=0.75pt]    {$e'_{2}$};

\end{tikzpicture}

%% file: tikzpictures/R1.tex
\tikzset{every picture/.style={line width=0.75pt}} 

\begin{tikzpicture}[x=0.75pt,y=0.75pt,yscale=-.9,xscale=.9]

\draw   (10,36.43) -- (90,36.43) -- (90,116.43) -- (10,116.43) -- cycle ;
\draw    (20,16.43) -- (20,36.43) ;
\draw    (80,16.43) -- (80,36.43) ;
\draw  [draw opacity=0] (80,16.43) .. controls (79.99,16.29) and (79.98,16.15) .. (79.98,16) .. controls (79.98,11.03) and (86.73,7) .. (95.06,7) .. controls (103.39,7) and (110.14,11.03) .. (110.14,16) .. controls (110.14,16.13) and (110.13,16.26) .. (110.12,16.39) -- (95.06,16) -- cycle ; \draw  [color={rgb, 255:red, 155; green, 155; blue, 155 }  ,draw opacity=1 ] (80,16.43) .. controls (79.99,16.29) and (79.98,16.15) .. (79.98,16) .. controls (79.98,11.03) and (86.73,7) .. (95.06,7) .. controls (103.39,7) and (110.14,11.03) .. (110.14,16) .. controls (110.14,16.13) and (110.13,16.26) .. (110.12,16.39) ;  
\draw [color={rgb, 255:red, 155; green, 155; blue, 155 }  ,draw opacity=1 ]   (110,15.43) -- (110,135.43) ;
\draw [shift={(110,81.43)}, rotate = 270] [color={rgb, 255:red, 155; green, 155; blue, 155 }  ,draw opacity=1 ][line width=0.75]    (10.93,-3.29) .. controls (6.95,-1.4) and (3.31,-0.3) .. (0,0) .. controls (3.31,0.3) and (6.95,1.4) .. (10.93,3.29)   ;
\draw   (180,36.43) -- (260,36.43) -- (260,116.43) -- (180,116.43) -- cycle ;
\draw    (190.15,18.43) -- (190,36.43) ;
\draw    (250,16.43) -- (250,36.43) ;
\draw  [draw opacity=0] (270,26.43) .. controls (270,26.39) and (270,26.34) .. (270,26.29) .. controls (270,23.84) and (272.27,21.86) .. (275.07,21.86) .. controls (277.88,21.86) and (280.15,23.84) .. (280.15,26.29) .. controls (280.15,26.33) and (280.15,26.37) .. (280.15,26.42) -- (275.07,26.29) -- cycle ; \draw  [color={rgb, 255:red, 155; green, 155; blue, 155 }  ,draw opacity=1 ] (270,26.43) .. controls (270,26.39) and (270,26.34) .. (270,26.29) .. controls (270,23.84) and (272.27,21.86) .. (275.07,21.86) .. controls (277.88,21.86) and (280.15,23.84) .. (280.15,26.29) .. controls (280.15,26.33) and (280.15,26.37) .. (280.15,26.42) ;  
\draw [color={rgb, 255:red, 155; green, 155; blue, 155 }  ,draw opacity=1 ]   (290,16.43) -- (290,146.43) ;
\draw    (270,26.43) -- (270,126.43) ;
\draw  [draw opacity=0] (250.05,16.31) .. controls (250.02,16.12) and (250.01,15.93) .. (250.01,15.73) .. controls (250.01,10.6) and (258.96,6.43) .. (270.01,6.43) .. controls (281.05,6.43) and (290,10.6) .. (290,15.73) .. controls (290,15.9) and (289.99,16.07) .. (289.97,16.24) -- (270.01,15.73) -- cycle ; \draw  [color={rgb, 255:red, 155; green, 155; blue, 155 }  ,draw opacity=1 ] (250.05,16.31) .. controls (250.02,16.12) and (250.01,15.93) .. (250.01,15.73) .. controls (250.01,10.6) and (258.96,6.43) .. (270.01,6.43) .. controls (281.05,6.43) and (290,10.6) .. (290,15.73) .. controls (290,15.9) and (289.99,16.07) .. (289.97,16.24) ;  
\draw [color={rgb, 255:red, 155; green, 155; blue, 155 }  ,draw opacity=1 ]   (280.15,26.42) -- (280.15,146.42) ;
\draw [shift={(280.15,92.42)}, rotate = 270] [color={rgb, 255:red, 155; green, 155; blue, 155 }  ,draw opacity=1 ][line width=0.75]    (10.93,-3.29) .. controls (6.95,-1.4) and (3.31,-0.3) .. (0,0) .. controls (3.31,0.3) and (6.95,1.4) .. (10.93,3.29)   ;
\draw    (20,116.43) -- (20,136.43) ;
\draw    (79.85,118.43) -- (80,136.43) ;
\draw  [draw opacity=0] (80,136.44) .. controls (79.99,136.58) and (79.98,136.72) .. (79.98,136.86) .. controls (79.98,141.84) and (86.7,145.86) .. (94.99,145.86) .. controls (103.28,145.86) and (110,141.84) .. (110,136.86) .. controls (110,136.22) and (109.89,135.59) .. (109.67,134.98) -- (94.99,136.86) -- cycle ; \draw  [color={rgb, 255:red, 155; green, 155; blue, 155 }  ,draw opacity=1 ] (80,136.44) .. controls (79.99,136.58) and (79.98,136.72) .. (79.98,136.86) .. controls (79.98,141.84) and (86.7,145.86) .. (94.99,145.86) .. controls (103.28,145.86) and (110,141.84) .. (110,136.86) .. controls (110,136.22) and (109.89,135.59) .. (109.67,134.98) ;  
\draw    (190,116.43) -- (190,146.43) ;
\draw    (249.85,118.43) -- (250,126.43) ;
\draw    (250,126.43) -- (270,146.43) ;
\draw  [color={rgb, 255:red, 255; green, 255; blue, 255 }  ,draw opacity=1 ][fill={rgb, 255:red, 255; green, 255; blue, 255 }  ,fill opacity=1 ] (256.67,136.43) .. controls (256.67,134.59) and (258.16,133.1) .. (260,133.1) .. controls (261.84,133.1) and (263.33,134.59) .. (263.33,136.43) .. controls (263.33,138.27) and (261.84,139.77) .. (260,139.77) .. controls (258.16,139.77) and (256.67,138.27) .. (256.67,136.43) -- cycle ;
\draw    (270,126.43) -- (250,146.43) ;

\draw  [draw opacity=0] (270,146.4) .. controls (270,146.45) and (270,146.5) .. (270,146.55) .. controls (270,149.31) and (272.27,151.55) .. (275.07,151.55) .. controls (277.88,151.55) and (280.15,149.31) .. (280.15,146.55) .. controls (280.15,146.5) and (280.15,146.46) .. (280.15,146.42) -- (275.07,146.55) -- cycle ; \draw  [color={rgb, 255:red, 155; green, 155; blue, 155 }  ,draw opacity=1 ] (270,146.4) .. controls (270,146.45) and (270,146.5) .. (270,146.55) .. controls (270,149.31) and (272.27,151.55) .. (275.07,151.55) .. controls (277.88,151.55) and (280.15,149.31) .. (280.15,146.55) .. controls (280.15,146.5) and (280.15,146.46) .. (280.15,146.42) ;  
\draw  [draw opacity=0] (250.26,145.63) .. controls (250.24,145.89) and (250.23,146.16) .. (250.23,146.43) .. controls (250.23,153.93) and (259.08,160) .. (270,160) .. controls (280.92,160) and (289.77,153.93) .. (289.77,146.43) .. controls (289.77,146.17) and (289.76,145.91) .. (289.74,145.66) -- (270,146.43) -- cycle ; \draw  [color={rgb, 255:red, 155; green, 155; blue, 155 }  ,draw opacity=1 ] (250.26,145.63) .. controls (250.24,145.89) and (250.23,146.16) .. (250.23,146.43) .. controls (250.23,153.93) and (259.08,160) .. (270,160) .. controls (280.92,160) and (289.77,153.93) .. (289.77,146.43) .. controls (289.77,146.17) and (289.76,145.91) .. (289.74,145.66) ;  

\draw (50.2,120.03) node [anchor=north] [inner sep=0.75pt]    {$\cdots $};
\draw (220.2,119.23) node [anchor=north] [inner sep=0.75pt]    {$\cdots $};
\draw (50,76.43) node    {$\beta $};
\draw (220,76.43) node    {$\beta $};
\draw (252,133.93) node [anchor=east] [inner sep=0.75pt]  [font=\footnotesize]  {$c_{0}$};

\end{tikzpicture}

%% file: tikzpictures/duality-map.tex
\tikzset{every picture/.style={line width=0.75pt}} 

\begin{tikzpicture}[x=0.75pt,y=0.75pt,yscale=-1,xscale=1]

\draw  [draw opacity=0][dash pattern={on 0.84pt off 2.51pt}] (23.7,36.61) .. controls (21.49,34.96) and (19.86,29.52) .. (19.86,23.07) .. controls (19.86,16.63) and (21.48,11.2) .. (23.7,9.54) -- (25.13,23.07) -- cycle ; \draw  [dash pattern={on 0.84pt off 2.51pt}] (23.7,36.61) .. controls (21.49,34.96) and (19.86,29.52) .. (19.86,23.07) .. controls (19.86,16.63) and (21.48,11.2) .. (23.7,9.54) ;  
\draw  [draw opacity=0] (24.55,9.1) .. controls (24.74,9.04) and (24.93,9.02) .. (25.13,9.02) .. controls (28.04,9.02) and (30.4,15.31) .. (30.4,23.07) .. controls (30.4,30.84) and (28.04,37.13) .. (25.13,37.13) .. controls (24.93,37.13) and (24.74,37.11) .. (24.55,37.05) -- (25.13,23.07) -- cycle ; \draw   (24.55,9.1) .. controls (24.74,9.04) and (24.93,9.02) .. (25.13,9.02) .. controls (28.04,9.02) and (30.4,15.31) .. (30.4,23.07) .. controls (30.4,30.84) and (28.04,37.13) .. (25.13,37.13) .. controls (24.93,37.13) and (24.74,37.11) .. (24.55,37.05) ;  
\draw  [draw opacity=0] (25.71,9.02) .. controls (25.51,9.02) and (25.32,9.02) .. (25.13,9.02) .. controls (16,9.02) and (8.61,15.31) .. (8.61,23.07) .. controls (8.61,30.84) and (16,37.13) .. (25.13,37.13) .. controls (25.32,37.13) and (25.51,37.13) .. (25.71,37.12) -- (25.13,23.07) -- cycle ; \draw   (25.71,9.02) .. controls (25.51,9.02) and (25.32,9.02) .. (25.13,9.02) .. controls (16,9.02) and (8.61,15.31) .. (8.61,23.07) .. controls (8.61,30.84) and (16,37.13) .. (25.13,37.13) .. controls (25.32,37.13) and (25.51,37.13) .. (25.71,37.12) ;  

\draw  [color={rgb, 255:red, 0; green, 0; blue, 0 }  ,draw opacity=1 ] (160,56.69) .. controls (160,49.42) and (162.43,43.53) .. (165.43,43.53) .. controls (168.42,43.53) and (170.85,49.42) .. (170.85,56.69) .. controls (170.85,63.95) and (168.42,69.84) .. (165.43,69.84) .. controls (162.43,69.84) and (160,63.95) .. (160,56.69) -- cycle ;
\draw  [draw opacity=0] (164.88,43.44) .. controls (165.06,43.44) and (165.24,43.43) .. (165.43,43.43) .. controls (174.03,43.43) and (181,49.37) .. (181,56.69) .. controls (181,64) and (174.03,69.94) .. (165.43,69.94) .. controls (165.24,69.94) and (165.06,69.93) .. (164.88,69.93) -- (165.43,56.69) -- cycle ; \draw   (164.88,43.44) .. controls (165.06,43.44) and (165.24,43.43) .. (165.43,43.43) .. controls (174.03,43.43) and (181,49.37) .. (181,56.69) .. controls (181,64) and (174.03,69.94) .. (165.43,69.94) .. controls (165.24,69.94) and (165.06,69.93) .. (164.88,69.93) ;  
\draw  [fill={rgb, 255:red, 0; green, 0; blue, 0 }  ,fill opacity=1 ] (227.45,20.29) .. controls (227.45,18.86) and (228.61,17.7) .. (230.04,17.7) .. controls (231.47,17.7) and (232.63,18.86) .. (232.63,20.29) .. controls (232.63,21.72) and (231.47,22.88) .. (230.04,22.88) .. controls (228.61,22.88) and (227.45,21.72) .. (227.45,20.29) -- cycle ;
\draw  [fill={rgb, 255:red, 0; green, 0; blue, 0 }  ,fill opacity=1 ] (296.35,19.19) .. controls (296.35,17.75) and (297.51,16.59) .. (298.94,16.59) .. controls (300.37,16.59) and (301.53,17.75) .. (301.53,19.19) .. controls (301.53,20.62) and (300.37,21.78) .. (298.94,21.78) .. controls (297.51,21.78) and (296.35,20.62) .. (296.35,19.19) -- cycle ;
\draw  [color={rgb, 255:red, 0; green, 0; blue, 0 }  ,draw opacity=1 ] (88.74,56.91) .. controls (88.74,50.09) and (91.02,44.56) .. (93.83,44.56) .. controls (96.64,44.56) and (98.92,50.09) .. (98.92,56.91) .. controls (98.92,63.72) and (96.64,69.25) .. (93.83,69.25) .. controls (91.02,69.25) and (88.74,63.72) .. (88.74,56.91) -- cycle ;
\draw  [draw opacity=0] (93.83,44.56) .. controls (98.64,44.51) and (102.52,42.34) .. (102.52,39.65) .. controls (102.52,36.94) and (98.56,34.75) .. (93.67,34.75) -- (93.67,39.65) -- cycle ; \draw   (93.83,44.56) .. controls (98.64,44.51) and (102.52,42.34) .. (102.52,39.65) .. controls (102.52,36.94) and (98.56,34.75) .. (93.67,34.75) ;  
\draw    (93.83,10.06) .. controls (112.27,10.38) and (122.07,23.31) .. (122.07,38.13) .. controls (122.06,52.95) and (116.08,68.6) .. (93.83,69.25) ;
\draw  [draw opacity=0][dash pattern={on 0.84pt off 2.51pt}] (90.24,33.87) .. controls (88.32,32.44) and (86.91,27.74) .. (86.91,22.16) .. controls (86.91,16.59) and (88.31,11.89) .. (90.23,10.45) -- (91.47,22.16) -- cycle ; \draw  [dash pattern={on 0.84pt off 2.51pt}] (90.24,33.87) .. controls (88.32,32.44) and (86.91,27.74) .. (86.91,22.16) .. controls (86.91,16.59) and (88.31,11.89) .. (90.23,10.45) ;  
\draw  [draw opacity=0] (90.97,10.07) .. controls (91.13,10.03) and (91.3,10) .. (91.47,10) .. controls (93.99,10) and (96.03,15.45) .. (96.03,22.16) .. controls (96.03,28.88) and (93.99,34.32) .. (91.47,34.32) .. controls (91.3,34.32) and (91.13,34.3) .. (90.97,34.25) -- (91.47,22.16) -- cycle ; \draw   (90.97,10.07) .. controls (91.13,10.03) and (91.3,10) .. (91.47,10) .. controls (93.99,10) and (96.03,15.45) .. (96.03,22.16) .. controls (96.03,28.88) and (93.99,34.32) .. (91.47,34.32) .. controls (91.3,34.32) and (91.13,34.3) .. (90.97,34.25) ;  
\draw  [draw opacity=0] (91.97,10.01) .. controls (91.8,10) and (91.64,10) .. (91.47,10) .. controls (83.58,10) and (77.18,15.45) .. (77.18,22.16) .. controls (77.18,28.88) and (83.58,34.32) .. (91.47,34.32) .. controls (91.64,34.32) and (91.8,34.32) .. (91.97,34.31) -- (91.47,22.16) -- cycle ; \draw   (91.97,10.01) .. controls (91.8,10) and (91.64,10) .. (91.47,10) .. controls (83.58,10) and (77.18,15.45) .. (77.18,22.16) .. controls (77.18,28.88) and (83.58,34.32) .. (91.47,34.32) .. controls (91.64,34.32) and (91.8,34.32) .. (91.97,34.31) ;  

\draw  [draw opacity=0][dash pattern={on 0.84pt off 2.51pt}] (239.7,34.11) .. controls (237.49,32.46) and (235.86,27.02) .. (235.86,20.57) .. controls (235.86,14.13) and (237.48,8.7) .. (239.7,7.04) -- (241.13,20.57) -- cycle ; \draw  [dash pattern={on 0.84pt off 2.51pt}] (239.7,34.11) .. controls (237.49,32.46) and (235.86,27.02) .. (235.86,20.57) .. controls (235.86,14.13) and (237.48,8.7) .. (239.7,7.04) ;  
\draw  [draw opacity=0] (240.55,6.6) .. controls (240.74,6.54) and (240.93,6.52) .. (241.13,6.52) .. controls (244.04,6.52) and (246.4,12.81) .. (246.4,20.57) .. controls (246.4,28.34) and (244.04,34.63) .. (241.13,34.63) .. controls (240.93,34.63) and (240.74,34.61) .. (240.55,34.55) -- (241.13,20.57) -- cycle ; \draw   (240.55,6.6) .. controls (240.74,6.54) and (240.93,6.52) .. (241.13,6.52) .. controls (244.04,6.52) and (246.4,12.81) .. (246.4,20.57) .. controls (246.4,28.34) and (244.04,34.63) .. (241.13,34.63) .. controls (240.93,34.63) and (240.74,34.61) .. (240.55,34.55) ;  
\draw  [draw opacity=0] (241.71,6.52) .. controls (241.51,6.52) and (241.32,6.52) .. (241.13,6.52) .. controls (232,6.52) and (224.61,12.81) .. (224.61,20.57) .. controls (224.61,28.34) and (232,34.63) .. (241.13,34.63) .. controls (241.32,34.63) and (241.51,34.63) .. (241.71,34.62) -- (241.13,20.57) -- cycle ; \draw   (241.71,6.52) .. controls (241.51,6.52) and (241.32,6.52) .. (241.13,6.52) .. controls (232,6.52) and (224.61,12.81) .. (224.61,20.57) .. controls (224.61,28.34) and (232,34.63) .. (241.13,34.63) .. controls (241.32,34.63) and (241.51,34.63) .. (241.71,34.62) ;  

\draw  [color={rgb, 255:red, 0; green, 0; blue, 0 }  ,draw opacity=1 ] (376,54.19) .. controls (376,46.92) and (378.43,41.03) .. (381.43,41.03) .. controls (384.42,41.03) and (386.85,46.92) .. (386.85,54.19) .. controls (386.85,61.45) and (384.42,67.34) .. (381.43,67.34) .. controls (378.43,67.34) and (376,61.45) .. (376,54.19) -- cycle ;
\draw  [draw opacity=0] (380.88,40.94) .. controls (381.06,40.94) and (381.24,40.93) .. (381.43,40.93) .. controls (390.03,40.93) and (397,46.87) .. (397,54.19) .. controls (397,61.5) and (390.03,67.44) .. (381.43,67.44) .. controls (381.24,67.44) and (381.06,67.43) .. (380.88,67.43) -- (381.43,54.19) -- cycle ; \draw   (380.88,40.94) .. controls (381.06,40.94) and (381.24,40.93) .. (381.43,40.93) .. controls (390.03,40.93) and (397,46.87) .. (397,54.19) .. controls (397,61.5) and (390.03,67.44) .. (381.43,67.44) .. controls (381.24,67.44) and (381.06,67.43) .. (380.88,67.43) ;  
\draw  [color={rgb, 255:red, 0; green, 0; blue, 0 }  ,draw opacity=1 ] (304.74,54.41) .. controls (304.74,47.59) and (307.02,42.06) .. (309.83,42.06) .. controls (312.64,42.06) and (314.92,47.59) .. (314.92,54.41) .. controls (314.92,61.22) and (312.64,66.75) .. (309.83,66.75) .. controls (307.02,66.75) and (304.74,61.22) .. (304.74,54.41) -- cycle ;
\draw  [draw opacity=0] (309.83,42.06) .. controls (314.64,42.01) and (318.52,39.84) .. (318.52,37.15) .. controls (318.52,34.44) and (314.56,32.25) .. (309.67,32.25) -- (309.67,37.15) -- cycle ; \draw   (309.83,42.06) .. controls (314.64,42.01) and (318.52,39.84) .. (318.52,37.15) .. controls (318.52,34.44) and (314.56,32.25) .. (309.67,32.25) ;  
\draw    (309.83,7.56) .. controls (328.27,7.88) and (338.07,20.81) .. (338.07,35.63) .. controls (338.06,50.45) and (332.08,66.1) .. (309.83,66.75) ;
\draw  [draw opacity=0][dash pattern={on 0.84pt off 2.51pt}] (306.24,31.37) .. controls (304.32,29.94) and (302.91,25.24) .. (302.91,19.66) .. controls (302.91,14.09) and (304.31,9.39) .. (306.23,7.95) -- (307.47,19.66) -- cycle ; \draw  [dash pattern={on 0.84pt off 2.51pt}] (306.24,31.37) .. controls (304.32,29.94) and (302.91,25.24) .. (302.91,19.66) .. controls (302.91,14.09) and (304.31,9.39) .. (306.23,7.95) ;  
\draw  [draw opacity=0] (306.97,7.57) .. controls (307.13,7.53) and (307.3,7.5) .. (307.47,7.5) .. controls (309.99,7.5) and (312.03,12.95) .. (312.03,19.66) .. controls (312.03,26.38) and (309.99,31.82) .. (307.47,31.82) .. controls (307.3,31.82) and (307.13,31.8) .. (306.97,31.75) -- (307.47,19.66) -- cycle ; \draw   (306.97,7.57) .. controls (307.13,7.53) and (307.3,7.5) .. (307.47,7.5) .. controls (309.99,7.5) and (312.03,12.95) .. (312.03,19.66) .. controls (312.03,26.38) and (309.99,31.82) .. (307.47,31.82) .. controls (307.3,31.82) and (307.13,31.8) .. (306.97,31.75) ;  
\draw  [draw opacity=0] (307.97,7.51) .. controls (307.8,7.5) and (307.64,7.5) .. (307.47,7.5) .. controls (299.58,7.5) and (293.18,12.95) .. (293.18,19.66) .. controls (293.18,26.38) and (299.58,31.82) .. (307.47,31.82) .. controls (307.64,31.82) and (307.8,31.82) .. (307.97,31.81) -- (307.47,19.66) -- cycle ; \draw   (307.97,7.51) .. controls (307.8,7.5) and (307.64,7.5) .. (307.47,7.5) .. controls (299.58,7.5) and (293.18,12.95) .. (293.18,19.66) .. controls (293.18,26.38) and (299.58,31.82) .. (307.47,31.82) .. controls (307.64,31.82) and (307.8,31.82) .. (307.97,31.81) ;  

\draw  [fill={rgb, 255:red, 0; green, 0; blue, 0 }  ,fill opacity=1 ] (389.35,54.69) .. controls (389.35,53.25) and (390.51,52.09) .. (391.94,52.09) .. controls (393.37,52.09) and (394.53,53.25) .. (394.53,54.69) .. controls (394.53,56.12) and (393.37,57.28) .. (391.94,57.28) .. controls (390.51,57.28) and (389.35,56.12) .. (389.35,54.69) -- cycle ;

\draw (27.13,23.07) node [anchor=west] [inner sep=0.75pt]    {$\ \ \ \mapsto $};
\draw (129,45.9) node [anchor=north west][inner sep=0.75pt]    {$=$};
\draw (243.13,20.57) node [anchor=west] [inner sep=0.75pt]    {$\ \ \ \mapsto $};
\draw (345,43.4) node [anchor=north west][inner sep=0.75pt]    {$=$};

\end{tikzpicture}

%% file: tikzpictures/dual-adj.tex
\tikzset{every picture/.style={line width=0.75pt}} 

\begin{tikzpicture}[x=0.75pt,y=0.75pt,yscale=-1,xscale=1]

\draw    (209.75,52.9) -- (173.35,52.9) ;
\draw    (209.76,71.31) -- (173.36,71.31) ;

\draw  [draw opacity=0] (211.03,52.55) .. controls (216.47,52.35) and (220.83,47.55) .. (220.83,41.66) .. controls (220.83,35.64) and (216.29,30.77) .. (210.68,30.77) -- (210.68,41.66) -- cycle ; \draw   (211.03,52.55) .. controls (216.47,52.35) and (220.83,47.55) .. (220.83,41.66) .. controls (220.83,35.64) and (216.29,30.77) .. (210.68,30.77) ;  
\draw    (209.76,12.5) .. controls (228.21,12.82) and (238.01,25.66) .. (238,40.39) .. controls (237.99,55.12) and (232.02,70.66) .. (209.76,71.31) ;

\draw  [dash pattern={on 0.84pt off 2.51pt}]  (209.75,52.9) -- (209.76,71.31) ;
\draw    (70.26,11.97) -- (33.86,11.97) ;
\draw    (70.27,30.38) -- (33.86,30.38) ;

\draw  [draw opacity=0] (71.74,52.02) .. controls (77.18,51.82) and (81.53,47.02) .. (81.53,41.13) .. controls (81.53,35.11) and (76.99,30.23) .. (71.39,30.23) -- (71.39,41.13) -- cycle ; \draw   (71.74,52.02) .. controls (77.18,51.82) and (81.53,47.02) .. (81.53,41.13) .. controls (81.53,35.11) and (76.99,30.23) .. (71.39,30.23) ;  
\draw    (70.47,11.97) .. controls (88.91,12.29) and (98.71,25.13) .. (98.71,39.86) .. controls (98.7,54.58) and (92.72,70.13) .. (70.47,70.78) ;

\draw  [dash pattern={on 0.84pt off 2.51pt}]  (70.47,11.97) -- (70.47,30.37) ;
\draw   (186.83,10) -- (209.83,10) -- (209.83,33) -- (186.83,33) -- cycle ;

\draw   (10,10) -- (33,10) -- (33,33) -- (10,33) -- cycle ;

\draw   (46.4,49.6) -- (69.4,49.6) -- (69.4,72.6) -- (46.4,72.6) -- cycle ;

\draw   (150,50) -- (173,50) -- (173,73) -- (150,73) -- cycle ;

\draw (127.5,40.5) node    {$=$};
\draw (191.56,68.91) node [anchor=south] [inner sep=0.75pt]  [font=\small]  {$\overline{f}$};
\draw (52.07,26.98) node [anchor=south] [inner sep=0.75pt]  [font=\small]  {$f$};
\draw (90.65,39.09) node  [font=\scriptsize]  {$B$};
\draw (229.95,39.42) node  [font=\scriptsize]  {$B$};
\draw (198.33,21.5) node    {$x$};
\draw (21.5,21.5) node    {$x$};
\draw (57.9,61.1) node    {$y$};
\draw (161.5,61.5) node    {$y$};

\end{tikzpicture}

%% file: 4.tex
\section
[Relation with y-ified HOMFLY--PT homology]
{Relation with $y$-ified HOMFLY--PT homology}
\label{sec:homfly}

\subsection{Overview}
\label{subsec:homfly-overview}
In this section, we work over $R = \QQ$, and omit $R$ from the notation.

The \textit{reduced HOMFLY--PT homology} $\mcH^{*,*,*}(L)$ of a link $L$, introduced by Khovanov and Rozansky in \cite{KR:2008b}, is a triply graded homology group whose graded Euler characteristic is the \textit{HOMFLY--PT polynomial}
\[
    P(L)(q, a) = \sum_{i,j,k} (-1)^{(k-j)/2} q^i a^j \dim_\QQ \mcH^{i,\,j,\,k}(L).
\]
Here, the HOMFLY--PT polynomial is a bivariate polynomial in $(q, a)$ defined by the skein relation
\[
    a P(
    \begin{tikzpicture}[baseline=1]
        \draw[->] (0.3,0)--(0,0.3);
        \draw (0,0)--(0.1,0.1);
        \draw[->] (0.2,0.2)--(0.3,0.3);
    \end{tikzpicture}
    )
    -a^{-1}P(
    \begin{tikzpicture}[baseline=1]
        \draw[->] (0,0)--(0.3,0.3);
        \draw (0.3,0)--(0.2,0.1);
        \draw[->] (0.1,0.2)--(0,0.3);
    \end{tikzpicture}
    )
    =(q - q^{-1})P(
    \begin{tikzpicture}[baseline=1]
        \draw [bend right = 45,->] (0,0) to (0,0.3);
        \draw [bend left = 45,->] (0.3,0) to (0.3,0.3);
    \end{tikzpicture}
    )
\]
and the normalization $P(\bigcirc) = 1$. For the triple grading of $\mcH^{*,*,*}(L)$, we follow the one given in \cite{Rasmussen:2015}. The three gradings of $\mcH^{*,*,*}$ are respectively called the \textit{$q$-grading}, the \textit{$a$-grading}, and the \textit{$v$-grading}, and their grading functions are denoted $\gr_q, \gr_a, \gr_v$.\footnote{%
    There seems to be no standard naming for the third grading of $\mcH^{*,*,*}$, so we call it the \textit{$v$-grading}, since it corresponds to the \textit{vertical} direction, whereas $a$ corresponds to the \textit{horizontal} direction (see \cite[Definition 2.7]{Rasmussen:2015}). 
} 
All three gradings of $\mcH$ have the opposite parity as the number of components of $L$, and hence $k - j$ is always even (see Remark after \cite[Definition 2.14]{Rasmussen:2015}). 

Next, for each $N > 0$, the \textit{reduced $\sl_N$-homology} $\overline{H}_N(L)$ of a link $L$, introduced by Khovanov and Rozansky in \cite{KR:2008a}, is a bigraded homology group, whose graded Euler characteristic gives the \textit{$\sl_N$-polynomial}
\[
    P_N(L)(q) = P(L)(q, q^N) = \sum_{i,j} (-1)^{i} q^j \dim_\QQ \overline{H}^{i, j}_N(L).
\]
In particular, when $N = 2$, the reduced $\sl_2$-homology $\overline{H}_2(L)$ is isomorphic to the reduced Khovanov homology of the mirror $L^*$ of $L$,\footnote{
    Since we follow the $q$-grading convention opposite to that of \cite{Khovanov:2000}, there is no minus sign on $j$ which is present in the expression given in \cite[p.\ 11]{KR:2008a}.
}
\[
    \overline{H}{}^{i, j}_2(L) \isom \rKh^{i, j}(L^*).
\]

In \cite{Rasmussen:2015}, Rasmussen showed that, for each $N \geq 1$, there is a spectral sequence from the reduced HOMFLY--PT homology to the reduced $\sl_N$-polynomial, thereby lifting the equation on the level of polynomials. Before stating the result, we define two additional gradings on $\mcH$, namely, the \textit{$t$-grading} by\footnote{
    Our $t$-grading is the negative of the $t$-grading defined in \cite{CG:structure-in-homfly2024}. 
}
\[
    \gr_t := \frac{\gr_v - \gr_a}{2},
\]
and for each $N > 0$, the \textit{$\sl_N$-modified $q$-grading} by 
\[
    \gr_{q_N} := \gr_q + N \gr_a.
\]

\begin{prop}[{\cite[Theorem 5.1]{Rasmussen:2015}}]
\label{prop:homfly-to-slN}
    For each $N \geq 1$, there is a spectral sequence $E = \{E_k\}_{k \geq 1}$ of triply graded $\QQ$-modules, starting from $E_1 \isom \mcH(L)$ and converging to $E_\infty \isom \overline{H}_N(L)$ as bigraded $\QQ$-modules with bigrading $(\gr_t, \gr_{q_N})$. The $k$-th differential $d^k_N$ is homogeneous of degree $(2Nk, -2k, 2 - 2k)$. 
\end{prop}

We also make use of the \textit{$\Delta$-grading} on $\mcH$, introduced in \cite{Rasmussen:2015,CG:structure-in-homfly2024}, defined as
\[
    \gr_\Delta := \gr_q + \gr_a + \gr_v.
\]
For a knot $K$, the \textit{$\Delta$-thickness} of a knot $K$ is defined by 
\[
    |\Delta|(K) = \frac{1}{2} \left( \Delta_{\max} - \Delta_{\min} \right) + 1
\]
where $\Delta_{\max}, \Delta_{\min}$ denote the maximal and the minimal $\Delta$-degrees of $\mcH(K)$ respectively. We say $K$ is \textit{KR-thin} if $|\Delta|(K) = 1$, i.e. if $\mcH(K)$ is supported on a single $\Delta$-slice. For $N = 2$, we have 
\[
    \gr_\Delta = 2 \gr_t + \gr_{q_2}
\]
and it coincides with the $\delta$-grading on $\widetilde{\Kh}(K^*)$. 
The $\Delta$-degree of the differential $d^k_N$ is given by 
\[
    \deg_\Delta(d^k_N) = 2(N - 2)k + 2.
\]
This immediately implies the following result. 

\begin{prop}[{\cite[Corollary 2.14]{CG:structure-in-homfly2024}}]
\label{prop:kr-thin}
    For any $N > |\Delta|(K)$, the $N$-th spectral sequence collapses at the $E_1$ page and give $\mcH(K) \isom \overline{H}_{N}(K)$. In particular, if $K$ is KR-thin, then $\mcH(K) \isom \overline{H}_{N}(K)$ for all $N \geq 2$. 
\end{prop}

Next, we briefly review the $y$-ification of HOMFLY--PT homology and the $\sl_2$-action, given in \cite{GH:y-ification2022,GHM:symmety-kr2024,CG:structure-in-homfly2024}. In \cite{GH:y-ification2022}, Gorsky and Hogancamp gave a $y$-ification $y\mcH(L)$ of the reduced HOMFLY homology $\mcH(L)$.\footnote{
    In \cite{GH:y-ification2022}, the HOMFLY--PT homology is denoted $H_{KR}$ and its $y$-ification is denoted $HY$. Here, we changed the notation for the sake of consistency throughout the paper. 
} 
Furthermore, in \cite{CG:structure-in-homfly2024}, Chandler and Gorsky defines the (reduced) \textit{$y$-ified $\sl_N$-homology} $y\overline{H}_N$ for each $N \geq 1$, extends Rasmussen's spectral sequence to the $y$-ified version. Namely,

\begin{prop}[{\cite[Theorem 3.19]{CG:structure-in-homfly2024}}]
\label{prop:y-homfly-to-slN}
    For each $N \geq 1$, there is a spectral sequence $E = \{E_k\}_{k \geq 1}$ of triply graded $\QQ[\bfy]$-modules, starting from $E_1 \isom y\mcH(L)$ and converging to $E_\infty \isom y\overline{H}_N(L)$. 
\end{prop}

In \cite{GHM:symmety-kr2024}, Gorsky, Hogancamp and Mellit define a sequence of operators $F_k$ on $y\mcH(L)$. In particular, it is proved that the second operator $\sle = F_2$ of degree $(4, 0, -4)$ satisfies a ``hard Lefshetz" property:

\begin{prop}[{\cite[Theorem 1.11, Corollary 1.12]{GHM:symmety-kr2024}}]
\label{prop:hard-lefshetz}
    For each $i \geq 0$, the $i$-th power $\sle^i$ of $\sle$ gives an isomorphism
    \[
        \sle^i\colon 
        y\mcH^{-2i, j, k}(L) \xrightarrow{\ \isom\ } y\mcH^{2i, j, k - 4i}(L).
    \]
    In particular, for a knot $K$, the induced operator $\sle$ on $\mcH(K)$ gives an isomorphism
    \[
        \sle^i\colon 
        \mcH^{-2i, j, k}(K) \xrightarrow{\ \isom\ } \mcH^{2i, j, k - 4i}(K).
    \]
\end{prop}

This implies that the operator $\sle$ with $\slh = \frac12 \gr_q$ extend to an $\sl_2$-action on $y\mcH(L)$, and when $K$ is a knot, it induces an $\sl_2$-action on $\mcH(K)$. Note that $\sle$ preserves the $\Delta$-gradings, and thus exhibits a symmetry on each $\Delta$-slice by the reflection $(q, a) \leftrightarrow (-q, a)$ (see \Cref{fig:H-7_1,fig:H-conway} in \Cref{sec:intro}). Furthermore, in \cite{CG:structure-in-homfly2024}, it is proved that $\sle$ extends over the the $y$-ified spectral sequence. Namely, 

\begin{prop}[{\cite[Theorem 3.27]{CG:structure-in-homfly2024}}]
\label{prop:sl2-action-spectral-seq}
    For each $N \geq 1$, the operator $\sle$ on $y\mcH(K)$ extends over the $y$-ified spectral sequence of \Cref{prop:y-homfly-to-slN}, so that it commutes with all the differentials of the spectral sequence and induces an operator $\sle$ on $y\overline{H}_N(L)$. In particular, when $K$ is a knot, the operator $\sle$ on $\mcH(K)$ extends over Rasmussen's spectral sequence and induces an operator $\sle$ on $\overline{H}_N(K)$. 
\end{prop}

We finally remark that Rasmussen's spectral sequence and Chandler--Gorsky's $y$-ified extension both hold for the \textit{deformed $\sl_N$-homology} $\overline{H}_{\del W}$, given by an arbitrary \textit{potential} $W \in \QQ[X]$ of degree $N + 1$; the special case $\del W = X^N$ recovers the original $\sl_N$-homology. In general, the spectral sequence does not preserve the $q$-grading and $\bar{H}_{\del W}$ instead admits a $q$-filtration. 

\subsection
[Proof of Theorem \ref{mainthm:homfly}]
{Proof of \Cref{mainthm:homfly}}
\label{subsec:y-kh-proof}

We restate \Cref{mainthm:homfly} given in \Cref{sec:intro} with a generalization to the deformed setting. Here, for a potential $W$ of degree $3$, we let $\CKh_{\del W}$ denote the complex obtained from the Frobenius algebra $A_{\del W} := \QQ[X]/(\del W)$. 

\begin{thm}
\label{thm:y-compatibility-deformed}
    Our constructions are compatible with the constructions of \cite{GH:y-ification2022,GHM:symmety-kr2024} under the $y$-ified spectral sequence of \cite{CG:structure-in-homfly2024}. Namely, for any degree $3$ potential $W \in \QQ[X]$, our reduced $y$-ified Khovanov homology $y\rKh_{\del W}(L^*)$ is isomorphic to the reduced $y$-ified $\sl_2$-homology $y\overline{H}_{\del W}(L)$, and the associated $\sle$-operator coincides with the induced action of $e$ on $y\overline{H}_{\del W}(L)$. In particular, for a knot $K$, our $\sle$-operator on $\rKh_{\del W}(K^*)$ coincides with the induced action of $e$ on $\overline{H}_{\del W}(K)$.
\end{thm}

The overall proof of \Cref{thm:y-compatibility-deformed} is straightforward and proceeds by unraveling the relations between various formulations of the relevant theories. First, we prove that the Khovanov complex and the $\sl_2$-complex are isomorphic as $\mcCA_l$-modules. We follow Hughes' approach given in \cite{Hughes:2014}, where they give an explicit isomorphism between the corresponding chain complexes,\footnote{
    In \cite{Hughes:2014}, the construction of $C_2(D)$ is mirrored so that the mirroring is implicit in the expression. 
}
\[
    C_2(D) \isom \CKh(D^*).
\]
We shall also see why we need an admissible $ab$-coloring only on the Khovanov side. 

\begin{figure}[t]
    \centering
    \input{tikzpictures/gamma-resolution}
    \caption{Resolutions $\Gamma_0, \Gamma_1$ of a crossing and maps $\chi_0, \chi_1$.}
    \label{fig:gamma-resolution}
    \vspace{1em}
    \input{tikzpictures/gamma-thick-edge-removal}
    \caption{Removal of a thick edge}
    \label{fig:thick-edge-removal}
\end{figure}

First, take a degree $3$ \textit{potential} $W \in \QQ[X]$ and put
\[
    \del W = X^2 - hX - t.
\]
For the $\sl_2$-complex, each \textit{state space} of a trivalent graph $\Gamma$ that arise in the cube of resolutions is given by the tensor product of \textit{matrix factorizations}, one for each resolution of a crossing of $D$ which is either $\Gamma_0$ or $\Gamma_1$ depicted in \Cref{fig:gamma-resolution}. From \cite[Section 6]{KR:2008a} (and \cite[Section 3.4]{Rasmussen:2015} for a general potential), the matrix factorizations $M_{\del W}(\Gamma_0)$, $M_{\del W}(\Gamma_1)$ are both given over the polynomial ring $\QQ[z_1, z_2, z'_1, z'_2]$. There are defined by
\[
    M_{\del W}(\Gamma_0) = 
    \begin{pmatrix}
        \pi(z_1, z'_1) & z'_1 - z_1\\ 
        \pi(z_2, z'_2) & z'_2 - z_2
    \end{pmatrix},
\]
where $\pi(z, z')$ denotes the polynomial
\[
    \pi(z, z') := \frac{W(z) - W(z')}{z - z'},
\]
and 
\[
    M_{\del W}(\Gamma_1) = 
    \begin{pmatrix}
        u_1 & z'_1 + z'_2 - z_1 - z_2 \\ 
        u_2 & z'_1 z'_2 - z_1 z_2
    \end{pmatrix}
\]
where $u_1, u_2$ are polynomials satisfying 
\[
    W = u_1 \cdot (z'_1 + z'_2 - z_1 - z_2) + u_2 \cdot (z'_1 z'_2 - z_1 z_2).  
\]
In particular, for our potential $W$, the polynomial $u_2$ is computed as
\[
    u_2 = h - (z_1 + z_2).
\]
Also note that 
\[
    \pi(z, z) = \del W(z).
\]

From \cite[Proposition 9]{KR:2008a}, we may reduce the number of variables using the linear terms. $M_{\del W}(\Gamma_0)$ reduces to $\QQ[z_1, z_2]$ by the identifications
\[
    z'_1 = z_1,\quad 
    z'_2 = z_2
\]
and $M_{\del W}(\Gamma_1)$ reduces to the polynomial ring $\QQ[z_1, z'_1]$ by the identification
\[
    z_2 = h - z_1, \quad 
    z'_2 = h - z_1'.
\]
The reduction for $M_{\del W}(\Gamma_1)$ can be regarded as \textit{removing the thick edge} (see \Cref{fig:thick-edge-removal}). Thus, for the whole trivalent graph $\Gamma$, after removing all the thick edges, each component transforms into a circle, and the associated variables can be reduced to a single variable $z$, subject to the rule that as $z$ passes an endpoint of a removed thick edge, it must be flipped by the involution
\[
    \sigma\colon z \mapsto h - z.
\]

Now, consider the (deformed) Khovanov complex $\CKh_{\del W}(D)$ of $D$ and the (deformed) $\sl_2$-complex $C_{\del W}(D^*)$ of $D^*$ corresponding to the potential $W$ (note that we mirrored the $C_2$ side, as in \cite{Hughes:2014}). For each vertex $v \in \{0, 1\}^n$, the vertex-module of $\CKh_{\del W}(D)$ at $v$ is given by 
\[
    V_v = \bigotimes_r A_{\del W}
\]
where $r$ denotes the number of circles in the resolved diagram $D_v$. On the other hand, the vertex-module of $C_{\del W}(D^*)$ at $v$ is given by 
\[
    V'_v \isom \QQ[z_1, \cdots, z_r] / (\del W(z_1), \ldots, \del W(z_r))
\]
(see \cite[Lemma 3]{Hughes:2014}). Here, each variable $z_i$ is associated to a point $p_i$ chosen on each component $C_i$ of $D_v$. The above descriptions shows that $V_v \isom V'_v$ as $\QQ$-modules. It remains to describe the explicit vertex-wise isomorphism that commutes with the edge maps of the two complexes. 

In \cite[Section 5.2]{Hughes:2014}, an auxiliary signing function $\tau$ is given, by arbitrarily choosing a set of base points $\{b_j\}$ on $D$ such that each connected component $G_j$ of the underlying 4-valent graph $G(D)$ of $D$ contains a unique base point $b_j$. Observe that the base points $\{b_j\}$ gives a unique admissible $ab$-coloring $\theta$, by requiring that the edges containing the base points are always colored $a$. Then we define, for each regular point $p$, 
\[
    \tau_p(z) := \begin{cases}
        z & \text{if $\theta(p) = a$}, \\
        \sigma(z) & \text{if $\theta(p) = b$}. \\
    \end{cases}
\]
One easily verifies that setting $h = 0$ recovers the original signing function $\tau$. We simply write $\tau(z)$ for $\tau_p(z)$ when $z$ is a variable associated to a point $p$. 

Using this function $\tau$, we may choose a generator $g_v \in V'_v$ for each vertex $v$, so that if $D_v \to D_{v'}$ is a merge, the edge map sends 
\[
    g_v \mapsto g_{v'},
\]
and if $D_v \to D_{v'}$ is a split, 
\[
    g_v \mapsto (\tau(z_i) + \tau(z_j) - h) g_{v'}
\]
where $z_i, z_j$ are the variables corresponding to the newly produced circles (see \cite[Section 5.3]{Hughes:2014}). Note that these relations correspond to the operations of $A_{\del W}$, where
\[
    1 \otimes 1 \xmapsto{m} 1 
\]
and 
\[
    1 \xmapsto{\Delta} X \otimes 1 + 1 \otimes X - h(1 \otimes 1).
\]
As stated in \cite[Section 5.4]{Hughes:2014}, the desired chain isomorphism is given by the following vertex-wise correspondence
\[
    V_v \to V'_v,\quad 
    X^{\eta_1} \otimes \cdots \otimes X^{\eta_r} \mapsto (\tau(z_1)^{\eta_1} \cdots \tau(z_r)^{\eta_r}) g_v.
\]
Here, $\eta_i \in \{0, 1\}$ and $X^0 = 1$. Multiplying $\tau_{p_i}(X)$ to the $i$-th tensor factor of $V(D_v)$ is exactly the action of $x_{p_i}$, so the above correspondence can be rewritten as 
\[
    (x_{p_1}^{\eta_1} \cdots x_{p_r}^{\eta_r})(1 \otimes \cdots \otimes 1) \mapsto (z_1^{\eta_1} \cdots z_r^{\eta_r}) g_v.
\]
Thus the identification $C_2(D) \isom \CKh(D^*)$ is established. Furthermore, we see from the construction that,

\begin{prop}
\label{prop:C2-and-CKh}
    The isomorphism $\CKh_{\del W}(D) \isom C_{\del W}(D^*)$ intertwines the actions $x_p$ and $z_p$. 
\end{prop}


Dot-sliding homotopies are defined similarly on $C_{\del W}(D^*)$. Recall from \cite[Section 6]{KR:2008a} that, local edge maps $\chi_0, \chi_1$ are defined between the two matrix factorizations $M_2(\Gamma_0)$, $M_2(\Gamma_1)$ associated to the elementary trivalent graphs $\Gamma_0$, $\Gamma_1$,
\[
\begin{tikzcd}
    M_2(\Gamma_0)
        \arrow[r, "\chi_0", shift left] 
    & M_2(\Gamma_1).
        \arrow[l, "\chi_1", shift left]
\end{tikzcd}
\]
As stated in \cite[Equation (20)]{KR:2008a}, they satisfy
\[
    \chi_1 \chi_0 = \chi_0 \chi_1 = z_1 - z'_2 = -(z_2 - z'_1).
\]
The dot-sliding homotopy $\chi_c$ for $c$ is given by the reversal of the edge map for $c$, i.e.\ if $\chi_i$ is the edge map, then $\chi_c = \chi_{1 - i}$. The above equation shows that
\[
    [d, \chi_c] = z_{p_1(c)} - z_{p'_2(c)} = -(z_{p_2(c)} - z_{p'_1(c)}).
\]

\begin{prop}
\label{prop:C2-and-CKh-dot-sliding}
    The isomorphism $\CKh_{\del W}(D) \isom C_{\del W}(D^*)$ intertwines the dot-sliding homotopies. 
\end{prop}


Now, $C_{\del W}(D^*)$ is given a $\mcCA_l$-action from the polynomial action and the dot-sliding homotopies, by formally repeating the argument of \Cref{subsec:dg-module-on-kh-cpx}. Therefore, we have 

\begin{prop}
\label{prop:C2-and-CKh-dg-module}
    The isomorphism $\CKh_{\del W}(D) \isom C_{\del W}(D^*)$ extends to an $\mcCA_l$-module isomorphism. 
\end{prop}

\begin{proof}[Proof of \Cref{thm:y-compatibility-deformed} (sketch)]
    The constructions of the $y$-ified HOMFLY-PT homology and its associated $\sle$-operator given in \cite{GH:y-ification2022,GHM:symmety-kr2024} are based on \textit{Rouquier complexes} \cite{Rouquier:2004,Khovanov:2007}. In \cite[Lemma 3.9]{CG:structure-in-homfly2024}, Chandler and Gorsky prove that their description of $\mcH$ coincides with Rasmussen's formulation, and in \cite[Lemma 3.13]{CG:structure-in-homfly2024}, they prove that the dot-sliding homotopies for the Rouquier complex commute with the horizontal differentials that give rise to the spectral sequence from $\mcH$ to $\overline{H}_{\del W}$. Therefore, the spectral sequence extends verbatim to the $y$-ified setting, and since $\sle$ is defined using the dot-sliding homotopies, it extends over the $y$-ified spectral sequence (\cite[Theorems 3.19, 3.27]{CG:structure-in-homfly2024}). Since our constructions are formally identical to the constructions of \cite{GH:y-ification2022,GHM:symmety-kr2024}, the claim is immediate from the reduced version of \Cref{prop:C2-and-CKh-dg-module}. 
\end{proof}

\subsection{Relation with Batson--Seed's link splitting spectral sequence}
\label{subsec:batson-seed}

In \cite{Batson-Seed:2015}, Batson and Seed give a \textit{link-splitting spectral sequence} from a deformation of Khovanov homology. We shall prove that their construction is obtained from a specialization of the $y$-ified Khovanov complex. Here, we may work over any commutative ring $R$. 

\begin{prop}
    Let $L$ be an $l$-component link with diagram $D$, and assume that $D$ has no crossingless circle components. Let $\nu = \{\nu_k\}$ be a set of \textit{weights} assigned to each component of $L$, and $s$ a \textit{sign assignment} in the sense of \cite[(1)]{Batson-Seed:2015}. Let $d'$ be the endomorphism of the Khovanov complex $(\CKh(D), d)$, given in \cite{Batson-Seed:2015}, defined from $\nu$ and $s$. Then there is an admissible $ab$-coloring $\theta$ on $D$, such that $d'$ can be written as 
    \[
        d' = \sum_k \bar{\xi}_k \nu_k.
    \]
\end{prop}

\begin{proof}[Proof (sketch)]
    From the explicit description of $d'$ given in \cite{Batson-Seed:2015}, it can be written using our notation as 
    \[
        d' = \sum_c s(c) \sgn(c) (\nu^c_1 - \nu^c_2) \chi_{c} 
    \]
    where $c$ ranges over all crossings of $D$, $s(c)$ is the chosen sign assignment, $\sgn(c)$ is the sign of the crossing, $\nu^{c}_1, \nu^{c}_2$ are the weights corresponding to the components containing $p_1(c), p_2(c)$, and $\chi_c$ is the (unsigned) dot-sliding homotopy for $c$. We claim that, there is a unique admissible coloring $\theta$ such that for each crossing $c$ of $D$, 
    \begin{equation}
    \label{eqn:sign-assign-cond}
        s(c) = \sgn(c) \epsilon_{p_1(c)}.
    \end{equation}
    Then the claim immediately follow from the description of $D$ given in \eqref{eqn:y-ify-Kh-alt}. To prove the equality, let us recall the condition that $s$ is required to satisfy. 
    
    Consider the CW decomposition $X$ of $S^2$ given by the underlying 4-valent graph $G(D) \subset S^2$ of $D$. The $0$-cells of $X$ correspond to the crossings of $D$, and the $1$-cells correspond to the edges of $D$. A $1$-cochain $g \in C^1(X; \ZZ/2)$ is defined as follows: for each edge $e$ of $D$, 
    \[
        g(e) := \begin{cases}
            -1 & \text{if the two ends of $e$ are both over strands or both under strands}, \\
            1 & \text{otherwise}.
        \end{cases}
    \]
    It is proved in \cite[Proposition 2.1]{Batson-Seed:2015} that $g$ is in fact a cocycle. A sign assignment $s$ is a $0$-cochain such that $\delta s = g$. Such $s$ always exists since $H^1(X; \ZZ/2) = 0$. Moreover, for each component $C$ of $G(D)$, there are exactly two choices of $s$ that give $\delta s = g$. 
    
    Having chosen one such $s$, we may associate an admissible coloring $\theta$ as follows. For each component $C$ of $G(D)$, choose one crossing $c$ and define $\theta$ over $C$ by requiring that 
    \[
        \theta(p_1(c)) = \begin{cases}
            a & \text{if $s(c) = \sgn(c)$,} \\ 
            b & \text{otherwise}. 
        \end{cases}
    \]
    Therefore, \eqref{eqn:sign-assign-cond} holds for the chosen crossing $c$, and it suffices to prove that the coboundary of the right-hand side gives $g$. This follows from an easy combinatorial argument. 
\end{proof}

With this identification, the assumption that $D$ has no crossingless components is unnecessary, since any crossingless component does not contribute to $d'$. 

\begin{prop}
\label{prop:bs-complex-from-y}
    The filtered complex $(\CKh(D), d + d')$ given in \cite{Batson-Seed:2015} is isomorphic to the specialization of the $y$-ified Khovanov complex
    \[
        y \CKh(D; R_\nu) := y \CKh(D) \otimes_{R[y_1, \ldots, y_l]} R_\nu
    \]
    where $R_\nu$ denotes $R$ regarded as a $R[y_1, \ldots, y_l]$-module with $y_k$ acted on as $\nu_k$. 
\end{prop}

Combining the results obtained thus far, we obtain the following proposition claimed in \Cref{sec:intro}. First, recall from \cite{Rasmussen:2015} that there is a triply graded complex $\mathbf{C}(D)$ equipped with differentials $d_+, d_v$ such that the HOMFLY--PT homology is given by 
\[
    \mcH(D) = H(H(\mathbf{C}(D), d_+), d_v).
\]
We let $\mcC(D) = H(\mathbf{C}(D), d_+)$ and call it the \textit{HOMFLY--PT complex}. 

\begin{prop}
\label{prop:4-ss}
    The HOMFLY--PT complex $(\mcC(D), d_v)$ admits a double filtration with respect to differentials $d_-$, $d'_v$ that give rise to the following four spectral sequences: (i) Rasmussen's spectral sequence \cite{Rasmussen:2015} for $N = 2$, (ii) its $y$-ified extension \cite{CG:structure-in-homfly2024} specialized by $y_k = \nu_k$, (iii) Batson--Seed's link-splitting spectral sequence~\cite{Batson-Seed:2015} for reduced Khovanov homology, and (iv) its HOMFLY--PT analogue~\cite{GH:y-ification2022}.
    \[
        \begin{tikzcd}[row sep=4em, column sep=3em]
        \mcH(L) = H(\mcC(D), d_v)
            \arrow[r, "\text{(ii)}", Rightarrow] 
            \arrow[d, "\text{(i)}"', Rightarrow] 
        & y \mcH(L; \QQ_\nu) = H(\mcC(D), d_v + d_v')
            \arrow[d, "\text{(iv)}", Rightarrow] \\
        \rKh(L^*) = H(\mcC(D), d_v + d_-)
            \arrow[r, "\text{(iii)}", Rightarrow] 
        & y\rKh(L^*; \QQ_\nu) = H(\mcC(D), d_v + d_- + d_v').                      
        \end{tikzcd}
    \]
\end{prop}

\begin{proof}[Proof (sketch)]
    The differential $d_-$ is given in \cite{Rasmussen:2015}, such that 
    \[
        \overline{C}_2(D) \isom H(\mathbf{C}(D), d_+ + d_-)
    \]
    and 
    \[
        \overline{H}_2(D) \isom H(H(\mathbf{C}(D), d_+ + d_-), d_v).
    \]
    The differential $d'_v$ is the specialization $d'_v := (D_v - d_v) |_{y_k = \nu_k}$ of the linear term of the differential $D_v$ of the $y$-ified HOMFLY--PT complex (see \cite[Section 3.4]{CG:structure-in-homfly2024}). These maps have degrees
    \[
        \deg(d_+) = (2, 2, 0),\quad
        \deg(d_-) = (4, -2, 0),\quad
        \deg(d_v) = (0, 0, 2),\quad
        \deg(d'_v) = (2, 0, -2).
    \]
    Note the following: (i) maps $d_-, d_v, d'_v$ are homogeneous of degree $1$ with respect to $\gr_t + \gr_{q_2} = \gr_q + \gr_a + $, (ii) maps $d_v, d'_v$ preserve $\gr_a$ while $d_-$ decreases it by $2$, and (iii) maps $d_-, d_v$ preserve $\gr_{q_2}$ while $d'_v$ increases it by $2$. Therefore, we obtain a doubly filtered complex 
    \[
        (\mcC(D) , d_v + d_- + d'_v)
    \]
    with respect to the triple grading 
    \[
        (\gr_t + \gr_{q_2},\ \gr_a,\ \gr_{q_2}).
    \]
    One sees that the four spectral sequences are indeed obtained from this double filtration on $\mcC(D)$, together with the identifications of \Cref{mainthm:homfly} and \Cref{prop:bs-complex-from-y}. 
\end{proof}

%% file: tikzpictures/gamma-resolution.tex
\tikzset{every picture/.style={line width=0.75pt}} 

\begin{tikzpicture}[x=0.75pt,y=0.75pt,yscale=-1,xscale=1]

\draw    (29,31) -- (29,89) ;
\draw [shift={(29,29)}, rotate = 90] [color={rgb, 255:red, 0; green, 0; blue, 0 }  ][line width=0.75]    (6.56,-1.97) .. controls (4.17,-0.84) and (1.99,-0.18) .. (0,0) .. controls (1.99,0.18) and (4.17,0.84) .. (6.56,1.97)   ;
\draw    (59,31) -- (59,89) ;
\draw [shift={(59,29)}, rotate = 90] [color={rgb, 255:red, 0; green, 0; blue, 0 }  ][line width=0.75]    (6.56,-1.97) .. controls (4.17,-0.84) and (1.99,-0.18) .. (0,0) .. controls (1.99,0.18) and (4.17,0.84) .. (6.56,1.97)   ;

\draw    (161.31,31.51) -- (173,45) ;
\draw [shift={(160,30)}, rotate = 49.09] [color={rgb, 255:red, 0; green, 0; blue, 0 }  ][line width=0.75]    (6.56,-1.97) .. controls (4.17,-0.84) and (1.99,-0.18) .. (0,0) .. controls (1.99,0.18) and (4.17,0.84) .. (6.56,1.97)   ;
\draw    (188.75,31.56) -- (178,45) ;
\draw [shift={(190,30)}, rotate = 128.66] [color={rgb, 255:red, 0; green, 0; blue, 0 }  ][line width=0.75]    (6.56,-1.97) .. controls (4.17,-0.84) and (1.99,-0.18) .. (0,0) .. controls (1.99,0.18) and (4.17,0.84) .. (6.56,1.97)   ;
\draw  [fill={rgb, 255:red, 0; green, 0; blue, 0 }  ,fill opacity=1 ] (173,45) -- (178,45) -- (178,75) -- (173,75) -- cycle ;
\draw    (171.69,76.51) -- (160,90) ;
\draw [shift={(173,75)}, rotate = 130.91] [color={rgb, 255:red, 0; green, 0; blue, 0 }  ][line width=0.75]    (6.56,-1.97) .. controls (4.17,-0.84) and (1.99,-0.18) .. (0,0) .. controls (1.99,0.18) and (4.17,0.84) .. (6.56,1.97)   ;
\draw    (179.25,76.56) -- (190,90) ;
\draw [shift={(178,75)}, rotate = 51.34] [color={rgb, 255:red, 0; green, 0; blue, 0 }  ][line width=0.75]    (6.56,-1.97) .. controls (4.17,-0.84) and (1.99,-0.18) .. (0,0) .. controls (1.99,0.18) and (4.17,0.84) .. (6.56,1.97)   ;
\draw    (90,50) -- (128,50) ;
\draw [shift={(130,50)}, rotate = 180] [color={rgb, 255:red, 0; green, 0; blue, 0 }  ][line width=0.75]    (10.93,-3.29) .. controls (6.95,-1.4) and (3.31,-0.3) .. (0,0) .. controls (3.31,0.3) and (6.95,1.4) .. (10.93,3.29)   ;
\draw    (130,60) -- (92,60) ;
\draw [shift={(90,60)}, rotate = 360] [color={rgb, 255:red, 0; green, 0; blue, 0 }  ][line width=0.75]    (10.93,-3.29) .. controls (6.95,-1.4) and (3.31,-0.3) .. (0,0) .. controls (3.31,0.3) and (6.95,1.4) .. (10.93,3.29)   ;

\draw (110,46.6) node [anchor=south] [inner sep=0.75pt]  [font=\small]  {$\chi _{0}$};
\draw (110,63.4) node [anchor=north] [inner sep=0.75pt]  [font=\small]  {$\chi _{1}$};
\draw (29,92.4) node [anchor=north] [inner sep=0.75pt]  [font=\small]  {$z_{1}$};
\draw (59,92.4) node [anchor=north] [inner sep=0.75pt]  [font=\small]  {$z_{2}$};
\draw (29,25.6) node [anchor=south] [inner sep=0.75pt]  [font=\small]  {$z'_{1}$};
\draw (59,25.6) node [anchor=south] [inner sep=0.75pt]  [font=\small]  {$z'_{2}$};
\draw (160,26.6) node [anchor=south] [inner sep=0.75pt]  [font=\small]  {$z'_{1}$};
\draw (190,26.6) node [anchor=south] [inner sep=0.75pt]  [font=\small]  {$z'_{2}$};
\draw (160,93.4) node [anchor=north] [inner sep=0.75pt]  [font=\small]  {$z_{1}$};
\draw (190,93.4) node [anchor=north] [inner sep=0.75pt]  [font=\small]  {$z_{2}$};
\draw (43.47,131.05) node    {$\Gamma _{0}$};
\draw (175.47,130.65) node    {$\Gamma _{1}$};

\end{tikzpicture}

%% file: tikzpictures/gamma-thick-edge-removal.tex
\tikzset{every picture/.style={line width=0.75pt}} 

\begin{tikzpicture}[x=0.75pt,y=0.75pt,yscale=-1,xscale=1]

\draw    (22.81,32.51) -- (34.5,46) ;
\draw [shift={(21.5,31)}, rotate = 49.09] [color={rgb, 255:red, 0; green, 0; blue, 0 }  ][line width=0.75]    (6.56,-1.97) .. controls (4.17,-0.84) and (1.99,-0.18) .. (0,0) .. controls (1.99,0.18) and (4.17,0.84) .. (6.56,1.97)   ;
\draw    (50.25,32.56) -- (39.5,46) ;
\draw [shift={(51.5,31)}, rotate = 128.66] [color={rgb, 255:red, 0; green, 0; blue, 0 }  ][line width=0.75]    (6.56,-1.97) .. controls (4.17,-0.84) and (1.99,-0.18) .. (0,0) .. controls (1.99,0.18) and (4.17,0.84) .. (6.56,1.97)   ;
\draw  [fill={rgb, 255:red, 0; green, 0; blue, 0 }  ,fill opacity=1 ] (34.5,46) -- (39.5,46) -- (39.5,76) -- (34.5,76) -- cycle ;
\draw    (33.19,77.51) -- (21.5,91) ;
\draw [shift={(34.5,76)}, rotate = 130.91] [color={rgb, 255:red, 0; green, 0; blue, 0 }  ][line width=0.75]    (6.56,-1.97) .. controls (4.17,-0.84) and (1.99,-0.18) .. (0,0) .. controls (1.99,0.18) and (4.17,0.84) .. (6.56,1.97)   ;
\draw    (40.75,77.56) -- (51.5,91) ;
\draw [shift={(39.5,76)}, rotate = 51.34] [color={rgb, 255:red, 0; green, 0; blue, 0 }  ][line width=0.75]    (6.56,-1.97) .. controls (4.17,-0.84) and (1.99,-0.18) .. (0,0) .. controls (1.99,0.18) and (4.17,0.84) .. (6.56,1.97)   ;
\draw    (141.81,31.71) -- (153.5,45.2) ;
\draw [shift={(140.5,30.2)}, rotate = 49.09] [color={rgb, 255:red, 0; green, 0; blue, 0 }  ][line width=0.75]    (6.56,-1.97) .. controls (4.17,-0.84) and (1.99,-0.18) .. (0,0) .. controls (1.99,0.18) and (4.17,0.84) .. (6.56,1.97)   ;
\draw    (169.25,31.76) -- (158.5,45.2) ;
\draw [shift={(170.5,30.2)}, rotate = 128.66] [color={rgb, 255:red, 0; green, 0; blue, 0 }  ][line width=0.75]    (6.56,-1.97) .. controls (4.17,-0.84) and (1.99,-0.18) .. (0,0) .. controls (1.99,0.18) and (4.17,0.84) .. (6.56,1.97)   ;
\draw    (152.19,76.71) -- (140.5,90.2) ;
\draw [shift={(153.5,75.2)}, rotate = 130.91] [color={rgb, 255:red, 0; green, 0; blue, 0 }  ][line width=0.75]    (6.56,-1.97) .. controls (4.17,-0.84) and (1.99,-0.18) .. (0,0) .. controls (1.99,0.18) and (4.17,0.84) .. (6.56,1.97)   ;
\draw    (159.75,76.76) -- (170.5,90.2) ;
\draw [shift={(158.5,75.2)}, rotate = 51.34] [color={rgb, 255:red, 0; green, 0; blue, 0 }  ][line width=0.75]    (6.56,-1.97) .. controls (4.17,-0.84) and (1.99,-0.18) .. (0,0) .. controls (1.99,0.18) and (4.17,0.84) .. (6.56,1.97)   ;
\draw    (70,59.59) .. controls (71.68,57.94) and (73.35,57.95) .. (75,59.63) .. controls (76.65,61.31) and (78.32,61.32) .. (80,59.67) .. controls (81.68,58.02) and (83.35,58.04) .. (85,59.72) .. controls (86.65,61.4) and (88.32,61.41) .. (90,59.76) .. controls (91.68,58.11) and (93.35,58.12) .. (95,59.8) .. controls (96.65,61.48) and (98.32,61.49) .. (100,59.84) .. controls (101.68,58.19) and (103.35,58.2) .. (105,59.88) .. controls (106.65,61.56) and (108.32,61.57) .. (110,59.92) -- (110,59.92) -- (118,59.98) ;
\draw [shift={(120,60)}, rotate = 180.47] [color={rgb, 255:red, 0; green, 0; blue, 0 }  ][line width=0.75]    (10.93,-3.29) .. controls (6.95,-1.4) and (3.31,-0.3) .. (0,0) .. controls (3.31,0.3) and (6.95,1.4) .. (10.93,3.29)   ;
\draw  [fill={rgb, 255:red, 0; green, 0; blue, 0 }  ,fill opacity=1 ] (157.85,46.6) .. controls (157.85,47.7) and (156.95,48.6) .. (155.85,48.6) .. controls (154.75,48.6) and (153.85,47.7) .. (153.85,46.6) .. controls (153.85,45.5) and (154.75,44.6) .. (155.85,44.6) .. controls (156.95,44.6) and (157.85,45.5) .. (157.85,46.6) -- cycle ;
\draw  [fill={rgb, 255:red, 0; green, 0; blue, 0 }  ,fill opacity=1 ] (157.85,74.2) .. controls (157.85,75.3) and (156.95,76.2) .. (155.85,76.2) .. controls (154.75,76.2) and (153.85,75.3) .. (153.85,74.2) .. controls (153.85,73.1) and (154.75,72.2) .. (155.85,72.2) .. controls (156.95,72.2) and (157.85,73.1) .. (157.85,74.2) -- cycle ;

\draw (21.5,27.6) node [anchor=south] [inner sep=0.75pt]  [font=\small]  {$z'_{1}$};
\draw (51.5,27.6) node [anchor=south] [inner sep=0.75pt]  [font=\small]  {$z'_{2}$};
\draw (21.5,94.4) node [anchor=north] [inner sep=0.75pt]  [font=\small]  {$z_{1}$};
\draw (51.5,94.4) node [anchor=north] [inner sep=0.75pt]  [font=\small]  {$z_{2}$};
\draw (140.5,26.8) node [anchor=south] [inner sep=0.75pt]  [font=\small]  {$z'_{1}$};
\draw (172.5,26.8) node [anchor=south west] [inner sep=0.75pt]  [font=\small]  {$z'_{2} \ =\ h\ -\ z'_{1}$};
\draw (140.5,93.6) node [anchor=north] [inner sep=0.75pt]  [font=\small]  {$z_{1}$};
\draw (172.5,93.6) node [anchor=north west][inner sep=0.75pt]  [font=\small]  {$z_{2} \ =\ h\ -\ z'_{1}$};

\end{tikzpicture}

%% file: 5.tex
\section{Diagrammatic computations}
\label{sec:diagrammatic-computation}

In the remaining two sections, we show that our constructions of $y$-ified Khovanov homology and the $\sle$-operator are well-suited for actual computations. In this section, we provide a \textit{diagrammatic computation method} in the spirit of Bar-Natan's computation of Khovanov homology, given in \cite{BarNatan:2007}. 

\subsection{General methods}

Here, we collect methods that are useful in the diagrammatic computation. First, we introduce the notation of a \textit{hollow dot}%
\footnote{
    Our definition of the hollow dot $\hdot$ differs by an overall sign from the one defined in \cite{BHP:2023} and in \cite{Khovanov:2004}.
} 
\begin{center}
    \input{tikzpictures/hollow-dot}
\end{center}
which corresponds to the element $Y = X - h$ in $A_{U(2)}$. Using the hollow dot, the neck-cutting relation (NC) can be written as:
\begin{center}
    \input{tikzpictures/hollow-dot-NC}
\end{center}
The following relations are also useful:
\begin{center}
    \input{tikzpictures/hollow-dot-relations}
\end{center}

The \textit{delooping isomorphism}, originally introduced in \cite{BarNatan:2007}, is an immediate consequence of the neck-cutting relation.\footnote{
    Again, the $q$-grading shifts are negated from the original \cite[Lemma 4.1]{BarNatan:2007}, due to the change of convention. 
}

\begin{prop}[Delooping]
\label{prop:delooping}
    The following matrix of morphisms give mutually inverse isomorphisms in the category $\Mat(\Cob_\bullet(B))$
    \begin{center}
        \input{tikzpictures/delooping}
    \end{center}
\end{prop}

The following proposition appears in \cite[Lemma 4.2]{BarNatan:2007}, which allows reduction of chain complexes in any additive category. 

\begin{prop}[Gaussian elimination]
\label{prop:gauss-elim}
    Suppose the top row of the following diagram is a part of a chain complex $\Omega$ such that the morphism $A \xrightarrow{a} A'$ is invertible. Then there is a strong deformation retract $\Omega'$ of $\Omega$ described by the bottom row. 
    \[
\begin{tikzcd}
[row sep=6em, column sep=4.5em, ampersand replacement=\&]
B_1 
    \arrow[r, "{\begin{pmatrix} \alpha \\ \beta \end{pmatrix}}"] 
    \arrow[d, equal] 
\& \begin{pmatrix} A \\ B_2 \end{pmatrix} 
    \arrow[rr, "{\begin{pmatrix} a & b \\ c & d \end{pmatrix}}", shift left] 
    \arrow[d, "{\begin{pmatrix} 0 & 1 \end{pmatrix}}"', shift right] 
\& \& \begin{pmatrix} A' \\ B_3 \end{pmatrix} 
    \arrow[r, "{\begin{pmatrix} \gamma & \delta \end{pmatrix}}"] 
    \arrow[d, "{\begin{pmatrix} -ca^{-1} & 1 \end{pmatrix}}"', shift right] 
    \arrow[ll, "{\begin{pmatrix} -a & 0 \\ 0 & 0 \end{pmatrix}}", dashed, shift left]
\& B_4 
    \arrow[d, equal] 
\\ B_1 
    \arrow[r, "\beta"] 
\& B_2 
    \arrow[rr, "d - c a^{-1} b"] 
    \arrow[u, "{\begin{pmatrix} -a^{-1}b \\ 1 \end{pmatrix}}"', dashed, shift right]
\& \& B_3 
    \arrow[r, "\delta"] 
    \arrow[u, "{\begin{pmatrix} 0 \\ 1 \end{pmatrix}}"', dashed, shift right]
\& B_4
\end{tikzcd}
    \]
    The strong deformation retraction $g$ is indicated by the downward vertical arrows, the inclusion $f$ by the upward dashed arrows, and the homotopy $h$ by the leftward dashed arrow. 
\end{prop}

The morphism $d - c a^{-1} b$ appearing in the bottom row of \Cref{prop:gauss-elim} is the \textit{Schur complement} of the matrix
$
    \begin{psmallmatrix} a & b \\ c & d \end{psmallmatrix}
$
by the invertible entry $a$. One may view that the subtracted term $c a^{-1} b$ is produced by the reversal of the arrow $a$
\[
\begin{tikzcd}[row sep=3em, column sep=5em]
A \arrow[rd, "c" description, pos=.8] & A' \arrow[l, "a^{-1}"', dashed] \\
B_2 \arrow[ru, "b" description, pos=.8] & B_3.
\end{tikzcd}
\]

In \Cref{subsec:pushforwad-by-strong-def-retract}, we discussed the problem of transferring an $\mcA$-action by a strong deformation retraction. The following proposition describes the transferred maps when the strong deformation retraction is given by the maps of \Cref{prop:gauss-elim}. 

\begin{prop}
\label{prop:gauss-elim-map}
    Retain notations from \Cref{prop:gauss-elim}. Suppose we have maps in $\Omega$, 
    \begin{align*}
        u = \begin{pmatrix}
            u_1 & u_2
        \end{pmatrix} \colon 
        \begin{pmatrix}
            A \\ B_2
        \end{pmatrix} \to X
    \end{align*}
    and 
    \begin{align*}
        v = \begin{pmatrix}
            v_1 \\ v_2
        \end{pmatrix} \colon 
        Y \to \begin{pmatrix}
            A' \\ B_3
        \end{pmatrix}
    \end{align*}
    where $X, Y$ are objects other than $A, A'$. The maps $u' := guf$, $v' := gvf$ in $\Omega'$ are given by 
    \begin{align*}
        u' = u_2 - u_1 a^{-1} b &\colon B_2 \to X, \\
        v' = v_2 - c a^{-1} v_1 &\colon Y \to B_3.
    \end{align*}
\end{prop}

Again, the descriptions given in \Cref{prop:gauss-elim-map} can be viewed diagrammatically as 
\[
\begin{tikzcd}[row sep=3.5em, column sep=4em]
& A \arrow[rd, "c" description, pos=.8] \arrow[ld, "u_1"'] & A' \arrow[l, "a^{-1}"', dashed] & \\
X & B_2 \arrow[ru, "b" description, pos=.8] \arrow[l, "u_2"] & B_3 & Y. \arrow[lu, "v_1"'] \arrow[l, "v_2"]
\end{tikzcd}
\]
Here, $u, v$ are drawn as homological degree decreasing maps, which is indeed the case for $\xi_i$ and $u$.

Now, let $C$ be an $\mcA^w_n$-module, and suppose there is a strong deformation retract $C'$ of $C$ as $R$-chain complexes, given by maps
\[
    \begin{tikzcd}
    C 
        \arrow[r, "g", shift left] 
        \arrow["h"', loop, distance=2em, in=305, out=235] & 
    C'. 
        \arrow[l, "f", shift left]
    \end{tikzcd}
\]
From \Cref{prop:strong-defr-induced-dga}, the maps $\xi'_i := g \xi_i f$, $u' := g u f$ do not necessarily give a strict $\mcA^w_n$-action on $C'$, and even if it does, from \Cref{prop:strong-defr-induced-dga2}, the maps $f, g, h$ are not necessarily (strictly or homotopy coherently) $\mcA^w_n$-linear. In what follows, we aim to describe the $\sle$-operator on ordinary Khovanov homology for specific classes of knots $K$, by decomposing $K$ into tangle pieces whose complexes admit strong deformation retractions over $R$. We do not claim (strict or homotopy coherent) $\mcA^w_n$-linearity of the intermediate maps. This lazy approach is justified by the following proposition. 

\begin{prop} 
\label{prop:lazy}
    Let $D$ be a 2-input planar arc diagram. For $a = 1, 2$, let $X_a$ be a complex in $\Kob(B_a)$ equipped with an $\mcA^{w_a}_{n_a}$-action, such that $X_1$ and $X_2$ are composable by $D$ and the actions are compatible. Further suppose that (i) $X = D(X_1, X_2)$ consists of a single virtual loop $\gamma$, (ii) the action of $\xi_1$ on $X$ is trivial, and (iii) there is a strong deformation retraction $(g^a, f^a, h^a)\colon X_a \to Y_a$ over $R$ that transfers the $\mcA^{w_a}_{n_a}$-action to $Y_a$. Put $f = D(f^1, f^2)$, $g = D(g^1, g^2)$, and let $u, u'$ denote the induced endomorphisms on the $\mcCA_1$-modules $X = D(X_1, X_2)$ and $Y = D(Y_1, Y_2)$. Then $f$, $g$ are $u$-equivariant, i.e.\ $[u, f] \htpy [u, g] \htpy 0$ over $R$. 
\end{prop}

\begin{proof}
    For $a = 1, 2$, let $\xi^a_i, u^a$ denote the endomorphisms for the $\mcA^{w_a}_{n_a}$-action of $X_a$. From \Cref{prop:strong-defr-induced-dga2}, there are homotopies $f^a_i := h^a \xi^a_i f^a$ and $f^a_u := h^a u^a f^a$ giving
    \[
        [d, f^a_i] + [\xi^a_i, f^a] = 0
    \]
    and 
    \[
        [d, f^a_u] + [u^a, f^a] + \sum_i(x^a_i + x'^a_{w_a(i)})f^a_i = 0.
    \]
    Let 
    \[
        \gamma\colon 
        p_1 \xrightarrow{\zeta_1} p_2 \xrightarrow{\zeta_2} \cdots \xrightarrow{\zeta_r} p_{r + 1} = p_1
    \]
    be the unique virtual loop of $X$, consisting of virtual arcs $\{\zeta_j\} = \{\xi^a_i\}$. We reorder $\{f_j\} = \{f^a_i\}$ so that $[\zeta_j, f] = [d, f_j]$ for each $j$. From \Cref{subsec:ext-to-tangles}, the endomorphism $u$ on $D(X_1, X_2)$ is defined by 
    \[
        u = u^1 + u^2 + \sum_{j < j'} \zeta_j \zeta_{j'}.
    \]
    We also have a similar description for $u'$. 
    By a computation similar to the one given in the proof of \Cref{prop:homotopy-u}, we obtain
    \[
        [u, f] \htpy -2 x_{p_1} \sum_{j} f_j
    \]
    However, from assumption (ii), we have
    \[
        \sum_j f_j = h (\sum_j \zeta_j ) f = 0.
    \]
    Analogous argument holds for $g$. 
\end{proof}

\begin{cor}
\label{cor:lazy-computation}
    Let $K$ be a knot diagram decomposed as $K = D(T_1, T_2)$. Suppose that, for $a = 1, 2$, the complex $[T_a]$ has a strongly deformation retract $E_a$ over $R$, and the $\mcA^{w_a}_{n_a}$-action on $[T_a]$ is transferred to $E_a$. Then $[K] = D([T_1], [T_2])$ and $E = D(E_1, E_2)$ are $\sle$-equivariantly homotopy equivalent. In particular, $\Kh(K)$ is $\sle$-equivariantly isomorphic to the homology of $\mcF_0(E)$.
\end{cor}


\subsection{Twist tangles}

\begin{figure}[t]
    \centering
    \input{tikzpictures/T_q}
    \caption{Twist tangle $T_q$}
    \label{fig:T_q}
    \vspace{2em}
    \input{tikzpictures/T_q-ori}
    \caption{Various orientations on $T_q$}
    \label{fig:T_q-ori}
\end{figure}

The chain homotopy type of the Khovanov complex of twist tangles has been studied previously in \cite[Section 6.2]{Khovanov:2000}, \cite[Proposition 4.1]{Thompson:2018}, \cite[Proposition 5.1]{Schuetz:2021} and in \cite[Proposition 4.3]{Kim-Sano:2025}. Here, we further study how the $\mcA^w_2$-actions are transformed by the strong deformation retractions. 

For $q \in \ZZ$, let $T_q$ denote the unoriented tangle diagram obtained by adding $|q|$ half-twists to a pair of crossingless vertical strands, either positively or negatively, depending on the sign of $q$ (see \Cref{fig:T_q}). 
The orientations on $T_q$ will be essential in the following arguments, so we fix the notation. $T^{\up\up}_q$ denotes the tangle diagram equipped with the orientation so that both strands are oriented upwards. Crossings are ordered from bottom to top, and $P_1, P_2$ are the bottom-left and the bottom-right endpoints of $T_q$. This coincides with the convention for the oriented braid diagram $\sigma_1^q$. $T^{\dn\dn}_q$ denotes the diagram obtained by rotating $T^{\up\up}_q$ 180 degrees about the axis perpendicular to the plane, so that the crossings are ordered from top to bottom, and $P_1, P_2$ are top-right and top-left endpoints of $T_q$. $T^{\up\dn}_q$ denotes the diagram equipped with the orientation so that the strand containing the bottom-left endpoint $P_1$ is oriented upwards, while the other strand is oriented downwards. When $q$ is even, the top-right endpoint is $P_2$, and when $q$ is odd, the top-left endpoint is $P_2$. $T^{\dn\up}_q$ denotes the diagram obtained by rotating $T^{\up\dn}_q$ 180 degrees about the axis parallel to the plane, so that $P_1$ is the bottom-right endpoint. For all of the diagrams, we set $P'_1, P'_2$ to be the endpoints such that $P_1 \sim P'_1$ and $P_2 \sim P'_2$. See \Cref{fig:T_q-ori} for the case $q = 2$. For all of these tangle diagrams, we choose the admissible $ab$-coloring $\theta$ such that $\theta(P_1) = a$. Note that the orientation preserving resolutions of $T^{\up\up}_q$ and $T^{\up\dn}_q$ are not isotopic. With this admissible $ab$-coloring, the homotopies $\xi_1$, $\xi_2$ are defined by the paths starting from $P_1$ and $P_2$. Hereafter, we prefer to study $T^{\dn\dn}_q$ rather than $T^{\up\up}_q$, since it better aligns with results of the previous works mentioned above. 

Let $\sfE_0$ and $\sfE_1$ denote the following unoriented $4$-end crossingless tangle diagrams, 
\begin{center}
    \input{tikzpictures/E01}
\end{center}
Define degree $-2$ endomorphisms $a, b$ on $\sfE_1$ by 
\begin{align*}
    a &= u_X + l_Y = u_Y + l_X, \\
    b &= u_X - l_X = u_Y - l_Y
\end{align*}
where $u_X, l_X$ are endomorphisms defined by
\begin{center}
    \input{tikzpictures/E1_uxlx}
\end{center}
and $u_Y, l_Y$ are similar endomorphisms defined using the hollow dot $\circ$. Let $s$ denote the genus $0$ saddle cobordism from $\sfE_0$ to $\sfE_1$ and also from $\sfE_1$ to $\sfE_0$. Note that we have $ab = ba = 0$, $sb = 0$ and $bs = 0$. 

\begin{lem}
\label{lem:Tq-lem1}
    Consider the following tangle diagrams $T, T'$ and morphisms $m, \Delta, a, b$ given by 
    \begin{center}
        \input{tikzpictures/endo-ab}
    \end{center}
    Here, it is assumed that these morphisms are $q$-grading preserving, by appropriately shifting the gradings of the objects. By delooping the circle appearing in tangle $T'$, morphisms $\Delta, m$ transform into
    \[
\begin{tikzcd}[row sep=.3em] 
 & & T & & q^{2}T \arrow[rrd, "X"] & & \\
T \arrow[rrd, "Y"'] \arrow[rru, "I"] & & \oplus & & & & T \\
 & & q^{-2}T & & T \arrow[uu, "\oplus", phantom] \arrow[rru, "I"'] & & 
\end{tikzcd} 
    \]
    and endomorphisms $a, b$ transform into
    \[
\begin{tikzcd}[row sep=2em]
q^{2}T \arrow[rrr, "X"] \arrow[rrrd, "t" description, pos=.75] & & & T & q^{2}T \arrow[rrr, "-Y"] \arrow[rrrd, "t" description, pos=.75] & & & T. \\
T \arrow[u, "\oplus", phantom] \arrow[rrr, "Y"'] \arrow[rrru, "I" description, pos=.75] & & & q^{-2}T \arrow[u, "\oplus", phantom] & T \arrow[u, "\oplus", phantom] \arrow[rrru, "I" description, pos=.75] \arrow[rrr, "-X"'] & & & q^{-2}T \arrow[u, "\oplus", phantom]
\end{tikzcd}
    \]
\end{lem}

The following examples are special cases of \Cref{prop:twist-tangle-retract,prop:twist-tangle-retract-op}, which illustrate the ideas of the general proof.

\begin{ex}
\label{ex:T1}
    The complex $[T^{\dn\dn}_1]$ is given by  
    \[
\begin{tikzcd}[column sep=3em]
\underline{\sfE_0} \arrow[r, "s", shift left] & \sfE_1 \ \arrow[l, "s", dashed, shift left]
\end{tikzcd}
    \]
    The solid arrow describes the differential, and the dashed arrow describes the homotopy $\xi_1$. We have $\xi_2 = -\xi_1$ and $u = 0$. 
\end{ex}

\begin{ex}
\label{ex:T2-para}
    The complex $[T^{\dn\dn}_2]$ is given by 
    \[
\begin{tikzcd}[row sep=4em, column sep=4em]
\underline{\sfE_0} \arrow[r, "s", shift left] \arrow[d, "s"', shift right]     & \sfE_1 \arrow[l, "-s", dashed, shift left] \arrow[d, "\Delta_1"', shift right]          \\
\sfE_1 \arrow[u, "s"', dashed, shift right] \arrow[r, "-\Delta_2", shift left] & \sfE_{11} \arrow[u, "m_1"', dashed, shift right] \arrow[l, "m_2", dashed, shift left] \arrow[ul, Rightarrow, "u"']
\end{tikzcd}
    \]
    Here, $\slE_{11}$ denotes the diagram obtained by vertically stacking two copies of $\slE_1$. The solid arrows describe the differential, where the vertical arrows correspond to the edge maps of the first crossing, and the horizontal arrow correspond to those of the second. The subscripts indicate the indices of the crossings on which the morphisms are applied.  The dashed arrows describe the homotopy $\xi_1$, which is given by 
    \[
        \xi_1 = \chi_1 - \chi_2.
    \]
    Obviously, we have $\xi_2 = -\xi_1$. The double arrow describes the homotopy $u$, which is given by 
    \[
        u = \chi_1 (-\chi_2) + (-\chi_1)\chi_2 = -2 \chi_1 \chi_2 = 2 s_1 s_2
    \]
    where $s_1 s_2$ is a cobordism with two saddles that connects the center circle with the top and bottom arcs of $\slE_{11}$. 

    Now, apply delooping to the center circle of $\slE_{11}$. From \Cref{lem:Tq-lem1}, the complex transforms into 
    \[
\begin{tikzcd}[row sep=3em, column sep=2em]
\underline{\sfE_0} \arrow[rr, "s", shift left] \arrow[d, "s"', shift right] &  & \sfE_1 \arrow[ll, "-s", dashed, shift left] \arrow[rd, "Y_1"', bend left, shift right] \arrow[d, "I"'] & \\
\sfE_1 \arrow[u, "s"', dashed, shift right] \arrow[rrr, "-Y_2", bend right, shift left] \arrow[rr, "-I", shift left] &  & \sfE^+_1 \arrow[u, dashed, shift right=2] \arrow[ll, dashed, shift left] & \sfE^-_1 \arrow[lu, "I"', dashed, bend right, shift right] \arrow[lll, dashed, bend left, shift left] \arrow[l, "\oplus", phantom]
\end{tikzcd}
    \]
    The double arrows are suppressed for the moment. We then apply Gaussian elimination to the bottom-most horizontal arrow labeled $-I$, and get 
    \[
\begin{tikzcd}[row sep=3em, column sep=4em]
\underline{\sfE_0} \arrow[r, "s", shift left] & \sfE_1 \arrow[l, "0", dashed, shift left] \arrow[d, "Y_1 - Y_2"', shift right] \\
& \sfE^-_1. \arrow[u, "I"', dashed, shift right]
\end{tikzcd}    
    \]
    Note that the solid arrow mapping into the bottom right $\sfE^+$ and the dashed arrow $\xi_1$ mapping out from the top right $\sfE$ have been modified according to \Cref{prop:gauss-elim,prop:gauss-elim-map}. 
    
    Next, the delooping transforms $u$ into
    \[
\begin{tikzcd}[row sep=3em, column sep=2em]
\underline{\sfE_0} \arrow[rr] \arrow[d]       &  & \sfE_1 \arrow[rd] \arrow[d]      &                                                                                  \\
\slE \arrow[rrr, bend right] \arrow[rr, "-I"] &  & \sfE^+_1 \arrow[llu, Rightarrow] & \sfE^-_1 \arrow[l, "\oplus", phantom] \arrow[lllu, "2s" description, Rightarrow, bend right=10]
\end{tikzcd}
    \]
    The elimination has the effect of simply removing the map starting from $\sfE^+_1$. Thus we obtain a simplified complex equipped with the transferred $\mcA^w_2$-action
    \[
\begin{tikzcd}[column sep=3em]
\underline{\sfE_0} \arrow[r, "s"] & \sfE_1 \arrow[r, "b", shift left] & \sfE^-_1. \arrow[l, "I", dashed, shift left] \arrow[ll, "2s"', Rightarrow, bend right, shift right]
\end{tikzcd}
    \]
\end{ex}

\begin{ex}
    Without the dga-action, the complex $[T^{\up\dn}_2]$ is identical to $[T^{\dn\dn}_2]$ up to bigrading shift. Let $P'_2$ denote the top-left endpoint of $T^{\up\dn}_2$. From the $ab$-coloring on $T^{\up\dn}_2$, one sees that $\xi_1 = \xi_{P_1} = \xi_{P'_2} = -\xi_{P_2} = -\xi_2$. As for $u$, we have 
    \[
        u = \chi_1 (-\chi_2) + \chi_2 (-\chi_1) = 0.
    \]
    Under the transformation of \Cref{ex:T2-para}, the transferred $\mcA^w_2$-action is given by
    \[
\begin{tikzcd}[column sep=3em]
\underline{\sfE_0} \arrow[r, "s"] & \sfE_1 \arrow[r, "b", shift left] & \sfE^-_1 \arrow[l, "I", dashed, shift left]
\end{tikzcd}
    \]
    with $u = 0$. 
\end{ex}

\begin{prop}
\label{prop:twist-tangle-retract}
    For each $k \geq 0$, the complex $[T^{\dn\dn}_{2k}]$ strongly deformation retracts (over $R$) onto a complex $E^{\dn\dn}_{2k}$ that admits the transferred $\mcA^w_2$-action:
    \[
\begin{tikzcd}[column sep=2.5em, row sep=4em]
\underline{\sfE_0} \arrow[r, "s"] & q^{-1}\sfE_1 \arrow[r, "b", shift left] & q^{-3}\sfE_1 \arrow[r, "a"] \arrow[l, "I", dashed, shift left] \arrow[ll, "2ks"', Rightarrow, bend right, shift right] & q^{-5}\sfE_1 \arrow[r, "b", shift left] \arrow[ll, "2(k - 1)I", Rightarrow, bend left, shift left] & \cdots \arrow[l, "I", dashed, shift left] \arrow[ll, "2(k - 1)I"', Rightarrow, bend right, shift right] & \\
& & \cdots \arrow[r, "b", shift left] & q^{-4k + 5}\sfE_1 \arrow[r, "a"] \arrow[l, "I", dashed, shift left] & q^{-4k + 3}\sfE_1 \arrow[r, "b", shift left] \arrow[ll, "2I", Rightarrow, bend left, shift left] & q^{-4k + 1}\sfE_1. \arrow[l, "I", dashed, shift left] \arrow[ll, "2I"', Rightarrow, bend right]
\end{tikzcd}
    \]
    Similarly, the complex $[T^{\dn\dn}_{2k+1}]$ strongly deformation retracts (over $R$) onto a complex $E^{\dn\dn}_{2k}$ that admits the transferred $\mcA^w_2$-action:
    \[
\begin{tikzcd}[column sep=2.5em, row sep=4em]
\underline{\sfE_0} \arrow[r, "s", shift left] & q^{-1}\sfE_1 \arrow[r, "b"] \arrow[l, "s", dashed, shift left] & q^{-3}\sfE_1 \arrow[r, "a", shift left] \arrow[ll, "2ks"', Rightarrow, bend right, shift right] & q^{-5}\sfE_1 \arrow[r, "b"] \arrow[ll, "2kI", Rightarrow, bend left, shift left] \arrow[l, "I", dashed, shift left] & \cdots \arrow[ll, "2(k - 1)I"', Rightarrow, bend right, shift right] & \\
 & & \cdots \arrow[r, "a", shift left] & q^{-4k + 3}\sfE_1 \arrow[r, "b"] \arrow[l, "I", dashed, shift left] & q^{-4k + 1}\sfE_1 \arrow[r, "a", shift left] \arrow[ll, "2I"', Rightarrow, bend right, shift right] & q^{-4k - 1}\sfE_1. \arrow[ll, "2I", Rightarrow, bend left, shift left] \arrow[l, "I", dashed, shift left]
\end{tikzcd}
    \]
    Here, the bigradings are relative with respect to the underlined object $\underline{\sfE_0}$. Dashed arrows indicate the action of $\xi_1 = -\xi_2$ and the double arrows indicate the action of $u$. The case $k < 0$ is obtained by reversing all of the arrows, so that the object $\underline{\sfE}_0$ is placed at the right end. 
\end{prop}

\begin{proof}
    The statement is obtained by repeating the proof of \cite[Proposition 4.3]{Kim-Sano:2025}, while also keeping track of homotopies $\xi_1$ and $u$ using \Cref{prop:gauss-elim-map}. We have seen that the statement holds for $k = 0$. Now, assume that the statement is true up to $2k$. 

    Let $D$ be a 2-input planar arc diagram giving $T_{2k + 1} = D(T_1, T_{2k})$. Then from the induction hypothesis, the complex $[T_{2k+1}] = D([T_1], [T_{2k}])$ deformation retracts to $D(E_1, E_{2k})$, described by 
    \[
\begin{tikzcd}[column sep=3em, row sep=3em]
\underline{\sfE_0} \arrow[d, "s"', shift right] \arrow[r, "s"]     & \sfE_1 \arrow[d, "\Delta_1"', shift right] \arrow[r, "b", shift left]          & \sfE_1 \arrow[r, "a"] \arrow[d, "\Delta_1"', shift right] \arrow[l, "-I", dashed, shift left]          & \cdots \arrow[r, "a"]    & \sfE_1 \arrow[d, "\Delta_1"', shift right] \arrow[r, "b", shift left]          & \sfE_1 \arrow[d, "\Delta_1"', shift right] \arrow[l, "-I", dashed, shift left]        \\
\sfE_1 \arrow[r, "-\Delta_2"] \arrow[u, "s"', dashed, shift right] & \sfE_{11} \arrow[r, "-b_2", shift left] \arrow[u, "m_1"', dashed, shift right] & \sfE_{11} \arrow[r, "-a_2"] \arrow[u, "m_1"', dashed, shift right] \arrow[l, "I", dashed, shift left] & \cdots \arrow[r, "-a_2"] & \sfE_{11} \arrow[r, "-b_2", shift left] \arrow[u, "m_1"', dashed, shift right] & \sfE_{11}. \arrow[u, "m_1"', dashed, shift right] \arrow[l, "I", dashed, shift left]
\end{tikzcd}
    \]
    Here, the solid vertical arrows correspond to the edge maps of $E_1$, and the solid horizontal arrows correspond to those of $E_{2k}$. The dashed arrow describes the homotopy $\xi_1$, given by 
    \[
        \xi_1 = \xi_1(E_1) + \xi_2(E_{2k}) = \xi_1(E_1) - \xi_1(E_{2k}).
    \]
    The double arrows are suppressed for the moment.

    For each of the objects in the bottom row, the delooping isomorphism gives $\sfE_{11} \isom \sfE^+_1 \oplus \sfE^-_1$, where we write $\sfE^\pm_1 := q^{\pm1} \sfE_1$. From \Cref{lem:Tq-lem1}, the bottom row can be written as 
    \[
\begin{tikzcd} [row sep=.5em] 
 & \sfE^+_1 \arrow[r] & \sfE^+_1 \arrow[r] & \sfE^+_1 \arrow[r] & \cdots \arrow[r] & \sfE^+_1. \\
\sfE_1 \arrow[rd] \arrow[ru, "-I" description] & & & & & \\
 & \sfE^-_1 \arrow[r] \arrow[ruu, "-I" description] & \sfE^-_1 \arrow[ruu, "-I" description] \arrow[r] & \sfE^-_1 \arrow[ruu, "-I" description] \arrow[r] & \cdots \arrow[r] \arrow[ruu, "-I" description] & \sfE^-_1 
\end{tikzcd} 
    \]
    By eliminating the diagonal arrows labeled $-I$ from the left end, one sees that $D(E_1, E_{2k})$ transforms into the complex $E_{2k+1}$, 
    \[
\begin{tikzcd}[column sep=3em, row sep=3em]
\underline{\sfE_0} \arrow[r, "s", shift left] & \sfE_1 \arrow[r, "b", shift left] \arrow[l, "s", dashed, shift left] & \sfE_1 \arrow[r, "a", shift left] \arrow[l, "0", dashed, shift left] & \cdots \arrow[r, "a", shift left] \arrow[l, "I", dashed, shift left] & \sfE_1 \arrow[l, "I", dashed, shift left] \arrow[r, "b", shift left] & \sfE_1 \arrow[l, "0", dashed, shift left] \arrow[d, "Y_1 + X_2"', shift right] \\
 & & & & & \sfE_1. \arrow[u, "I"', dashed, shift right] \end{tikzcd} 
    \]
    Note the effect of \Cref{prop:gauss-elim-map} on the dashed arrows. 
    
    Next, let us observe how $u$ is transformed. In $D(E_1, E_{2k})$ we have:
    \[
\begin{tikzcd}[column sep=3em, row sep=3em]
\underline{\sfE_0} \arrow[r, "s"] \arrow[d, "s"'] & \sfE_1 \arrow[r, "b"] \arrow[d, "\Delta_1"'] & \sfE_1 \arrow[d, "\Delta_1"'] \arrow[ll, "2ks"', Rightarrow, bend right] \arrow[r, "a"] & \sfE_1 \arrow[ll, "2(k - 1)I"', Rightarrow, bend right] \arrow[d, "\Delta_1"'] \arrow[r, "b"] & \cdots \arrow[ll, "2(k - 1)I"', Rightarrow, bend right] \arrow[r, "a"] & \sfE_1 \arrow[r, "b"] \arrow[d, "\Delta_1"'] & \sfE_1 \arrow[d, "\Delta_1"] \arrow[ll, "2I"', Rightarrow, bend right] \\
\sfE_{1} \arrow[r, "-a_2"] & \sfE_{11} \arrow[r, "-b_2"] & \sfE_{11} \arrow[ll, "2ks", Rightarrow, bend left] \arrow[lu, "2m_1" description, Rightarrow] \arrow[r, "-a_2"] & \sfE_{11} \arrow[ll, "2(k - 1)I", Rightarrow, bend left] \arrow[r, "-b_2"] & \cdots \arrow[ll, "2(k-1)I", Rightarrow, bend left] \arrow[lu, "2m_1" description, Rightarrow] \arrow[r, "-a_2"] & \sfE_1 \arrow[r, "-b_2"] & \sfE_{11} \arrow[ll, "2I", Rightarrow, bend left] \arrow[lu, "2m_1" description, Rightarrow]
\end{tikzcd}
    \]
    Note the up-left pointing diagonal arrows, which comes from the cross terms of $u$, 
    \[
        \xi_1(E_1) \xi_2(E_{2k}) + \xi_2(E_1) \xi_1(E_{2k}) = 2 \xi_1(E_1) \xi_2(E_{2k}).
    \]
    After performing delooping, we see that there are parts of the following pattern
    \[
\begin{tikzcd}
\sfE_1 \arrow[r, "b"] & \sfE_1 \arrow[r, "a"] \arrow[d] & \sfE_1 \arrow[ll, "(2i)I"', Rightarrow, bend right] \arrow[d, "I"] \\
& \sfE^-_1 \arrow[lu, "2I", Rightarrow] \arrow[r, "-I"] & \sfE^+_1.
\end{tikzcd}
    \]
    By eliminating the bottom horizontal arrow, the factor coming from the lower route adds up to the existing factor, resulting in 
    \[
\begin{tikzcd}
\sfE_1 \arrow[r, "b"] & \sfE_1 \arrow[r, "a"] & \sfE_1. \arrow[ll, "(2i+1)I"', Rightarrow, bend right]
\end{tikzcd}
    \]
    Therefore, after the Gaussian elimination, the complex transforms into
    \[
\begin{tikzcd}[column sep=3em, row sep=3em]
\underline{\sfE_0} \arrow[r, "s"] & \sfE_1 \arrow[r, "b"] & \sfE_1 \arrow[ll, "2ks"', Rightarrow, bend right] \arrow[r, "a"] & \sfE_1 \arrow[ll, "2kI"', Rightarrow, bend right] \arrow[r, "b"] & \cdots \arrow[ll, "2(k - 1)I"', Rightarrow, bend right] \arrow[r, "a"] & \sfE_1 \arrow[r, "b"] & \sfE_1 \arrow[d, "a"] \arrow[ll, "2I"', Rightarrow, bend right] \\
 & & & & & & \sfE^-_1 \arrow[lu, "2I" description, Rightarrow] \end{tikzcd}
 \]
    and obtain the complex $E_{2k+1}$ as expected. A similar argument shows that the complex $E_{2k+2}$ can be obtained from $D(E_1, E_{2k+1})$. 
\end{proof}


\begin{prop}
\label{prop:twist-tangle-retract-op}
    For each $k \geq 0$, the complex $[T^{\up\dn}_{2k}]$ strongly deformation retracts (over $R$) onto to a complex $E^{\up\dn}_{2k}$, which is similar (up to bigrading shift) to the complex $E^{\dn\dn}_{2k}$ of \Cref{prop:twist-tangle-retract}, except that double arrows (representing $u$) are all removed. Similar statement holds for the complex $[T^{\up\dn}_{2k+1}]$.
\end{prop}

\begin{proof}
    For both $l = 2k, 2k+1$, the strong deformation retraction $[T_l] \to E_l$ is defined independently of the orientations, so the descriptions of the complex $E^{\up\dn}_l$ and the homotopy $\xi_1$ are identical to those of \Cref{prop:twist-tangle-retract}. For $u$, we see that for both $l = 2k, 2k+1$, 
    \begin{align*}
        u_1 &= \sum_{j < j'} (-1)^{j - 1}(-1)^{j - 1} \chi_{j} \chi_{j'} \\ 
        u_2 &= \sum_{j < j'} (-1)^{j}(-1)^{j'} \chi_{j'} \chi_{j}
    \end{align*}
    and thus $u = u_1 + u_2 = 0$. 
\end{proof}

\subsection
[(2, q)-torus knots]
{$(2, q)$-torus knots}
\label{subsec:torus-knots}

For odd $q > 0$, let $T_{2, q}$ denote the diagram of the \textit{$(2, q)$-torus knot}, obtained from the closure of $T^{\dn\dn}_{q}$. For $l \geq 1$, let $\sfU_l$ denote the $l$-component unnested crossingless unlink diagram. A morphism $H\colon \sfU_1 \to \sfU_1$ is defined by $H = X + Y = 2X - h$. 

\begin{prop}
\label{prop:torus-knot-str}
    For $k \geq 0$, the complex $[T_{2, 2k + 1}]$ strongly deformation retracts (over $R$) onto the following complex that admits the transferred $\mcCA_1$-action:
    \[
\begin{tikzcd}[column sep=2.5em]
\underline{\sfU_1} \arrow[r, dotted] & 0 \arrow[r, dotted] & \sfU_1 \arrow[r, "H"] \arrow[ll, "2kI"', Rightarrow, bend right, shift right] & \sfU_1 \arrow[r, dotted] & \sfU_1 \arrow[ll, "2(k - 1)I"', Rightarrow, bend right, shift right] \arrow[r, "H"] & \cdots \arrow[ll, "2(k - 1)I", Rightarrow, bend left] \arrow[r, "H"] & \sfU_1 \arrow[r, dotted] & \sfU_1 \arrow[r, "H"] \arrow[ll, "2I"', Rightarrow, bend right, shift right] & \sfU_1. \arrow[ll, "2I", Rightarrow, bend left, shift left]
\end{tikzcd}
    \]
    Here, the bigrading is relative to the underlined object of bigrading $(0, -2k)$. For each arrow, the $q$-grading of the target object is shifted by $-2$ relative to the source. 
\end{prop}

\begin{proof}
    By applying the closure functor $\cl$ to \Cref{prop:twist-tangle-retract}, we get 
    \[
        \cl(\sfE_0) = \sfU'_2,\quad 
        \cl(\sfE_1) = \sfU_1
    \]
    where $\sfU'_2$ denotes the nested $2$-component unlink diagram, and 
    \[
        \cl(a) = H,\quad \cl(b) = 0.
    \]
    The left-end of the diagram is given by 
    \[
        \underline{\sfU'_2} \xrightarrow{\ m\ } q^{-1} \sfU_1
    \]
    which reduces to 
    \[
        q^{1}\underline{\sfU_1} \xrightarrow{\ 0\ } 0.
    \]

    Now, let us observe the $\mcCA_1$-action on the reduced complex. First, we have 
    \[
        \bar{\xi}_1 = \xi_1 + \xi_2 = 0.
    \]
    (Recall from \Cref{prop:y-complex-knot} that this holds for any knot diagram.) Furthermore, the correction term of $\bar{u}$ is given by 
    \[
        \xi_1 \xi_2 = 0
    \]
    so $\bar{u} = u$. 
\end{proof}

\begin{cor}
\label{cor:kh-2q-torus}
    For any $k \geq 0$, the unreduced homology $\Kh(T_{2, 2k+1}, \QQ)$ is supported in two $\delta$-gradings $\delta = -2k \pm 1$, and its $\QQ[e]$-module structure is given by
    \[
 \begin{tikzcd}[row sep=1em, column sep=1.5em]
\delta = -2k + 1 \arrow[r, ":", phantom] 
& \underline{\QQ} & & \QQ \arrow[ll, "2k"'] & & \cdots \arrow[ll, "2(k - 1)"'] & \cdots & & \QQ \arrow[ll, "2"'] & \\
\delta = -2k - 1 \arrow[r, ":", phantom] 
& \underline{\QQ} \arrow["0"', loop, distance=2em, in=305, out=235] & 0 & 0 & \QQ & & \cdots \arrow[ll, "2(k-1)"'] & \QQ & & \QQ \arrow[ll, "2"']
\end{tikzcd}
    \]
    Similarly, the reduced homology $\rKh(T_{2, 2k+1}, \QQ)$ is supported in a single $\delta$-grading, and its $\QQ[e]$-module structure is given by
    \[
\begin{tikzcd}[row sep=1em, column sep=1.5em]
\delta = -2k \arrow[r, ":", phantom] 
& \underline{\QQ} & & \QQ \arrow[ll, "2k"'] & & \cdots \arrow[ll, "2(k - 1)"'] & \cdots & & \QQ \arrow[ll, "2"'] & \\
& \underline{0} & 0 & 0 & \QQ & & \cdots \arrow[ll, "2(k-1)"'] & \QQ & & \QQ \arrow[ll, "2"']
\end{tikzcd}
    \]
    In both cases, the $\sle$-operators uniquely extend to $\sl_2$-actions on the homology groups. 
\end{cor}

\begin{proof}
    By applying the TQFT $\mcF$ to the complex of \Cref{prop:torus-knot-str}, the segment 
    \[
        \sfU_1 \xrightarrow{\ H\ } q^{-2} \sfU_1
    \]
    transforms into 
    \[
    \begin{tikzcd}[ampersand replacement=\&, column sep=6em]
    {q^{-1}\begin{pmatrix}R\{X\} \\ R\{1\}\end{pmatrix}} 
    \arrow[r, "{\begin{pmatrix} h & 2 \\ 2t & -h \end{pmatrix}}"]
    \& {q^{-3}\begin{pmatrix}R\{X\} \\ R\{1\}\end{pmatrix}}
    \end{tikzcd}
    \]
    In particular, when $R = \QQ$ and $h = t = 0$, the segment is homotopic to 
    \[
    \begin{tikzcd}[column sep=4em]
    {q^{-1} \QQ\{X\}} 
    \arrow[r, "0", dotted] & 
    {q^{-3} \QQ\{1\}}
    \end{tikzcd}
    \]
    In the reduced case, the segment transforms into
    \[
        q^{-2} \bar{A}  \xrightarrow{\ X\ } q^{-4} \bar{A}
    \]
    where $\bar{A} = A / (X)$. Thus, when $h = t = 0$, the map becomes trivial. The final statement is immediate from \Cref{prop:e-extension}.
\end{proof}

\subsection{Twist knots}

\begin{figure}[t]
    \centering
    \input{tikzpictures/twist-knot}
    \caption{Twist knot $K_p$}
    \label{fig:twist-knot}
\end{figure}

For $p \in \ZZ$, let $K_p$ denote the diagram of a twist knot with $-p$ half-twists, as depicted in \Cref{fig:twist-knot}. In our convention, the knots $K_p$ with smaller $p$ correspond to 
\[
    K_0 = 0_1,\quad
    K_1 = 3_1^*,\quad
    K_2 = 4_1^*,\quad
    K_3 = 5_2^*,\quad 
    K_4 = 6_1^*,\quad
    \ldots
\]
in Rolfsen's knot table. Also note that, for $p > 0$,  we have $K_{-p} = K^*_{p - 1}$, so we may assume that $p > 0$. 
With the 2-input planar arc diagram $D$ depicted in \Cref{fig:twist-knot}, we have $K_p = D(T_2, T_{-p})$. When $p$ is odd, the tangles are oriented as $D(T^{\up\up}_2, T^{\dn\up}_{-p})$, and when $p$ is even, they are oriented as $D(T^{\up\dn}_2, T^{\dn\up}_{-p})$. Note that the signs of the crossings of $T_{-p}$ are positive independent of $p$, whereas those of $T_2$ are positive when $p$ is odd, and negative otherwise. 

\begin{prop}
\label{prop:twist-knot-reduced}
    For $k \geq 0$, the complex $[K_{2k}]$ strongly deformation retracts (over $R$) onto the following complex that admits the transferred $\mcCA_1$-action:
    \[
\begin{tikzcd}
\sfU_1 \arrow[d, "H"'] & & & & & & & & \\
\sfU_1 \arrow[r, dotted] \arrow[d, dotted] & \sfU_1 \arrow[d, dotted] \arrow[r, "-H"] & \sfU_1 \arrow[d, dotted] \arrow[r, dotted] & \sfU_1 \arrow[r, "-H"] \arrow[d, dotted] & \sfU_1 \arrow[d, dotted] \arrow[r, dotted] & \cdots \arrow[r, "-H"] & \sfU_1 \arrow[d, dotted] \arrow[r, dotted] & 0 \arrow[d, dotted] \arrow[r, dotted] & 0 \arrow[d, dotted] \\
\underline{\sfU_1} \arrow[r, dotted] & \sfU_1 \arrow[r, "H"] \arrow[lu, "2I" description, Rightarrow] & \sfU_1 \arrow[r, dotted] & \sfU_1 \arrow[r, "H"] \arrow[lu, "2I" description, Rightarrow] & \sfU_1 \arrow[r, dotted] & \cdots \arrow[r, "H"] \arrow[lu, "2I" description, Rightarrow] & \sfU_1 \arrow[r, dotted] & \sfU_1 \arrow[r, "H"] \arrow[lu, Rightarrow] & \sfU_1 
\end{tikzcd}
    \]
    The bigrading is relative to the underlined object of bigrading $(0, 0)$. For each horizontal and vertical arrow, the $q$-grading of the target object is shifted by $-2$ relative to the source.

    Similarly, the complex $[K_{2k + 1}]$ strongly deformation retracts (over $R$) onto the following complex that admits the transferred $\mcCA_1$-action:
    \[
\begin{tikzcd}
\underline{\sfU_1} \arrow[d, dotted] & & & & & & & & \\
\sfU_1 \arrow[r, "-H"] \arrow[d, dotted] & \sfU_1 \arrow[d, dotted] \arrow[r, dotted] & \sfU_1 \arrow[r, "-H"] \arrow[d, dotted] & \sfU_1 \arrow[r, dotted] \arrow[d, dotted] & \sfU_1 \arrow[d, dotted] \arrow[r, "-H"] & \cdots \arrow[r, "-H"] & \sfU_1 \arrow[d, dotted] \arrow[r, dotted] & 0 \arrow[d, dotted] \arrow[r, dotted] & 0 \arrow[d, dotted] \\
\sfU_1 \arrow[r, "H"] \arrow[uu, "2I", Rightarrow, bend left] & \sfU_1 \arrow[r, dotted] & \sfU_1 \arrow[r, "H"] \arrow[lu, "-2I" description, Rightarrow] & \sfU_1 \arrow[r, dotted] & \sfU_1 \arrow[r, "H"] \arrow[lu, "-2I" description, Rightarrow] & \cdots \arrow[r, "H"] & \sfU_1 \arrow[r, dotted] & \sfU_1 \arrow[r, "H"] \arrow[lu, "-2I" description, Rightarrow] & \sfU_1 
\end{tikzcd}
    \]
    Here, the bigrading is relative to the underlined object of bigrading $(0, -2)$. 
\end{prop}

\begin{proof}
    First, we prove the odd case. From \Cref{prop:twist-tangle-retract,prop:twist-tangle-retract-op}, the complex $K_{2k + 1} = D([T^{\up\up}_2], [T^{\dn\up}_{-2k - 1}])$ strongly deformation retracts (over $R$) onto the complex $X_{2k + 1} = D(E^{\up\up}_2, E^{\dn\up}_{-2k - 1})$ given by 
    \[
\begin{tikzcd}
\underline{\sfU_2} \arrow[r, "a"] \arrow[d, "m"'] & \sfU_2 \arrow[r, "b"] \arrow[d, "m"'] & \sfU_2 \arrow[r, "a"] \arrow[d, "m"'] & \cdots \arrow[r, "a"]  & \sfU_2 \arrow[r, "b"] \arrow[d, "m"'] & \sfU_2 \arrow[r, "m"] \arrow[d, "m"'] & \sfU_1 \arrow[d, "\Delta_1"] \\
\sfU_1 \arrow[r, "-H"] \arrow[d, dotted] & \sfU_1 \arrow[d, dotted] \arrow[r, dotted] & \sfU_1 \arrow[r, "-H"] \arrow[d, dotted] & \cdots \arrow[r, "-H"] & \sfU_1 \arrow[d, dotted] \arrow[r, dotted] & \sfU_1 \arrow[d, dotted] \arrow[r, "-\Delta_2"] & \sfU'_2 \arrow[d, "b_1"] \\
\sfU_1 \arrow[r, "H"] & \sfU_1 \arrow[r, dotted] & \sfU_1 \arrow[r, "H"] & \cdots \arrow[r, "H"]  & \sfU_1 \arrow[r, dotted] & \sfU_1 \arrow[r, "\Delta_2"] & \sfU'_2. 
\end{tikzcd}
    \]
    Here, the vertical direction corresponds to the differential of $E^{\up\up}_2$, and the horizontal direction to that of $E^{\dn\up}_{-2k - 1}$. The diagram $\sfU'_2 = D(\sfE_1, \sfE_0)$ appearing at the right end denotes the nested $2$-component unlink diagram. Dotted arrows indicate the zero map, and the double arrows are suppressed for the moment. 

    One sees that the top row can be eliminated, leaving a single $\sfU_1$ at the left end. This transforms the vertical $\sfU_2 \xrightarrow{\ m\ } \sfU_1$ to $\sfU_1 \xrightarrow{\ 0\ } \sfU_1$. Furthermore, consider the bottom right box
    \[
\begin{tikzcd}
\sfU_1 \arrow[d, dotted] \arrow[r, "-\Delta_2"] & \sfU'_2 \arrow[d, "b_1"] \\
\sfU_1 \arrow[r, "\Delta_2"] & \sfU'_2.                 
\end{tikzcd}        
    \]
    By delooping the inner circle of $\sfU'_2$, it transforms into, 
    \[
\begin{tikzcd}[row sep=3em]
\sfU_1 \arrow[d, dotted] \arrow[rr, "-Y", bend left] \arrow[r, "-I"'] & \sfU^+_1 \arrow[d, "Y"'] \arrow[rd] & \sfU^-_1 \arrow[d, "X"] \arrow[ld, "-I" description, pos=.75] \arrow[l, "\oplus", phantom] \\
\sfU_1 \arrow[rr, "Y"', bend right] \arrow[r, "I"] & \sfU^+_1 & \sfU^+_- \arrow[l, "\oplus", phantom]
\end{tikzcd}
    \]
    and then, by applying Gaussian elimination to the top horizontal arrow $-I$ and to the diagonal arrow $-I$, reduces to 
    \[
\begin{tikzcd}
0 \arrow[d, dotted] \arrow[r, dotted] & 0 \arrow[d, dotted] \\
\sfU_1 \arrow[r, "H"] & U^+_1.
\end{tikzcd}
    \]
    Thus we obtain the underlying complex described in the statement. 
    
    Next, we describe $u = u(X_{2k + 1})$. Recall that
    \[
        u(E^{\up\up}_2)\colon \sfE_1 \xrightarrow{\ 2 
        s\ } \underline{\sfE_0}
    \]
    and $u(E^{\up\dn}_{-p}) = 0$. Put $\xi^1 = \xi_1(E^{\up\up}_2)$ and $\xi^2 = \xi_1(E^{\up\dn}_{-p})$. In the diagram of \Cref{fig:twist-knot}, by starting from the point where an arrow is drawn, the cross terms of $u$ are given by 
    \[
        \xi^1(\xi^2 - \xi^1 - \xi^2) + \xi^2(-\xi^1 - \xi^2) + (-\xi^1)(-\xi^2) = 2\xi^1\xi^2.
    \]
    Therefore, we have 
    \[
        u = u(E^{\up\up}_2) + 0 + 2\xi^1\xi^2.
    \]
    Consider the following repetitive block in $X_{2k + 1}$:
    \[
\begin{tikzcd}
\sfU_2 \arrow[r, "a"] \arrow[d] & \sfU_2 \arrow[r, "b"] \arrow[d] & \sfU_2 \arrow[d] \\
\sfU_1 \arrow[r] \arrow[d, dotted] & \sfU_1 \arrow[d, dotted] \arrow[r, dotted] & \sfU_1 \arrow[d, dotted] \\
\sfU_1 \arrow[r] \arrow[uu, "2s", Rightarrow, bend left] & \sfU_1 \arrow[r, dotted] \arrow[uu, Rightarrow, bend left] \arrow[lu, "2I" description, Rightarrow] & \sfU_1 \arrow[uu, Rightarrow, bend left] 
\end{tikzcd}
    \]
    The elimination of the top row transforms this block into
    \[
\begin{tikzcd}
\sfU^+_1 \arrow[r, dotted] \arrow[d, dotted] & 0 \arrow[r, dotted] \arrow[d, dotted] & 0 \arrow[d, dotted] \\
\sfU_1 \arrow[r] \arrow[d, dotted] & \sfU_1 \arrow[d, dotted] \arrow[r, dotted] & \sfU_1 \arrow[d, dotted] \\
\sfU_1 \arrow[r] \arrow[uu, "2I", Rightarrow, bend left] & \sfU_1 \arrow[r, dotted] \arrow[lu, "0" description, Rightarrow] & \sfU_1 \arrow[lu, "-2I" description, Rightarrow]
\end{tikzcd}
    \]
    Therefore, we obtain the description of $u$ as in the statement. 

    The proof for the even case proceeds similarly. The process of simplification is similar, and we obtain the strong deformation retract $X_{2k}$ of $K_{2k}$ described in the statement. In this case, we have $u(E^{\up\dn}_2) = u(E^{\dn\up}_{-2k}) = 0$. With $\xi^1 = \xi_1(E^{\up\dn}_2)$ and $\xi^2 = \xi_1(E^{\dn\up}_{-p})$, the cross terms of $u(X_{2k})$ are given by 
    \[
        \xi^1(\xi^2 - \xi^1 - \xi^2) + \xi^2(-\xi^1 - \xi^2) + (-\xi^1)(-\xi^2) = 2\xi^1\xi^2
    \]
    Therefore, we have 
    \[
        u(X_{2k}) = 2\xi^1\xi^2.
    \]
    The following repetitive block in $X_{2k}$:
    \[
\begin{tikzcd}
\sfU_2 \arrow[r, "b"] \arrow[d] & \sfU_2 \arrow[r, "a"] \arrow[d] & \sfU_2 \arrow[d] \\
\sfU_1 \arrow[r, dotted] \arrow[d, dotted] & \sfU_1 \arrow[d, dotted] \arrow[r, "-H"] & \sfU_1 \arrow[d, dotted] \\
\sfU_1 \arrow[r, dotted] & \sfU_1 \arrow[r, "H"] \arrow[lu, "2I" description, Rightarrow] & \sfU_1 
\end{tikzcd}
    \]
    transforms into
    \[
\begin{tikzcd}
\sfU^+_1 \arrow[r, dotted] \arrow[d, "H"'] & 0 \arrow[r, dotted] \arrow[d, dotted] & 0 \arrow[d, dotted] \\
\sfU_1 \arrow[r, dotted] \arrow[d, dotted] & \sfU_1 \arrow[d, dotted] \arrow[r, "-H"] & \sfU_1 \arrow[d, dotted] \\
\sfU_1 \arrow[r, dotted] & \sfU_1 \arrow[r, "H"] \arrow[lu, "2I" description, Rightarrow] & \sfU_1. 
\end{tikzcd}
    \]
\end{proof}

In the following corollary, we suppress the $q$-grading, since it can be easily read from \Cref{prop:twist-knot-reduced}. 

\begin{cor}
\label{cor:kh-twist-knot}
    For any $k \geq 0$, the following holds.
    \begin{enumerate}
        \item $\Kh(K_{2k}, \QQ)$ is supported on $\gr_\delta = \pm 1$, and $e$ acts as $0$. 
        \item $\Kh(K_{2k + 1}, \QQ)$ is supported on $\gr_\delta = -1, -3$, and has a single $e$-string $(\QQ \to \QQ)$ of dimension $2$ on $\delta = -1$; other $e$-strings are $1$-dimensional. 
        \item $\rKh(K_{2k}, \QQ)$ is supported on $\gr_\delta = 0$, and decomposes into a direct sum of $2k + 1$ copies of $\QQ$ and $k$ copies of $(\QQ \to \QQ)$. 
        \item $\rKh(K_{2k + 1}, \QQ)$ is supported on $\gr_\delta = -2$, and decomposes into a direct sum of $2k + 1$ copies of $\QQ$ and $k + 1$ copies of $(\QQ \to \QQ)$. 
    \end{enumerate}
\end{cor}

\begin{proof}
    Let $A$ denote the segment 
    \[
\begin{tikzcd}
\sfU_1 \arrow[r, "H"] & \sfU_1 & {} & \sfU_1 \arrow[r, "H"] \arrow[ll, "2I"', Rightarrow] & \sfU_1
\end{tikzcd}
    \]
    and $B$ the segment
    \[
\begin{tikzcd}
{\sfU_1} & {} & {\sfU_1} \arrow[r, "H"] \arrow[ll, "2I"', Rightarrow] & {\sfU_1}
\end{tikzcd}
    \]
    with appropriate bigrading shift. From \Cref{prop:twist-knot-reduced}, the complex $[K_{2k}]$ is homotopic to the direct sum of $\sfU_1$ and $k$ copies of $A$, and $[K_{2k+1}]$ is homotopic to the direct sum of $B$ and $k$ copies of $A$. With the TQFT $\mcF_0$ over $\QQ$, the segment $\mcF_0(A)$ is homotopic to 
    \[
\begin{tikzcd}
{\QQ\{X\}} \arrow[r, "0", dotted] & {\QQ\{1\}} & {} & {\QQ\{X\}} \arrow[r, "0", dotted] \arrow[ll, "0"', Rightarrow] & {\QQ\{1\}},
\end{tikzcd}
    \]
    and $\mcF_0(B)$ is homotopic to 
    \[
\begin{tikzcd}
{\begin{pmatrix} \QQ\{X\} \\ \QQ\{1\}\end{pmatrix}} & {} & {\QQ\{X\}} \arrow[r, "0", dotted] \arrow[ll, "{\begin{pmatrix} 2 \\ 0 \end{pmatrix}}"', Rightarrow] & {\QQ\{1\}}.
\end{tikzcd}
    \]
    Computation for the reduced case is similar. 
\end{proof}

\subsection{Comparison with HOMFLY--PT homology}

\begin{figure}[t]
    \centering
    \begin{subfigure}[b]{0.4\textwidth}
        \centering
        \input{tikzpictures/table-H-6_1}
        \caption{$\mcH(6_1^*)$}
        \label{fig:H-6_1}
    \end{subfigure}
    \begin{subfigure}[b]{0.4\textwidth}
        \centering
        \input{tikzpictures/table-Kh-6_1}
        \caption{$\rKh(6_1)$}
        \label{fig:Kh-6_1}
    \end{subfigure}
    \caption{$K_6 = 6_1^*$}
    \label{fig:6_1}
\end{figure}

Both $(2, q)$-torus knots and twist knots are \textit{$2$-bridge knots}, and it is proved in \cite[Corollary 1]{Rasmussen:2015} that $2$-bridge knots are KR-thin. Therefore, (i) their HOMFLY--PT homology is determined by the \textit{HOMFLY--PT polynomial} and the \textit{knot signature} (supported on $\Delta = \sigma(K)$), and (ii) their HOMFLY--PT homology is isomorphic to the reduced Khovanov homology (after collapsing the triple grading). Let us compare the results of \Cref{cor:kh-2q-torus,cor:kh-twist-knot} with the triply graded $\sl_2$-module structure of the HOMFLY--PT homology of these knots. 

The HOMFLY--PT polynomial of torus knots was found by Jones in \cite{Jones:1987}, and in particular the formula for $(2, 2k+1)$-torus knots is given in \cite[(72)]{DGR:2006} as 
\[
    P(T_{2, 2k+1})(q, a) = 
    a^{2k}\sum_{i=0}^{k} q^{4i-2k}
    \;-\;
    a^{2k+2}\sum_{i=1}^{k} q^{4i-2k-2}.
\]
One sees from the above formula that $\mcH(T_{2, 2k+1})$ has simple (irreducible) $\sl_2$-summands of dimensions $k + 1$ and $k$, lying on $a$-grading $2k$ and $2k + 2$ respectively. See \Cref{fig:7_1} in \Cref{sec:intro} for the case $k = 3$. 

The HOMFLY--PT polynomial of twist knots is computed by Ta{\c{s}}k{\"o}pr{\"u} and Sinan in \cite[Corollary 2.6]{TS:2021} as 
\[
    P(K_p)(q, a) = \begin{cases}
        \displaystyle
        -(q - q^{-1})^2 \left( \sum_{i = 0}^{k - 1} a^{2i} \right) + a^{2k} - a^{2k - 2} + a^{-2} & p = 2k,\\
        \displaystyle
        (q - q^{-1})^2 a^2 \left( \sum_{i = 0}^{k} a^{2i} \right) - a^{2k + 4} + a^{2k + 2} + a^{2} & p = 2k + 1.
    \end{cases}
\]
For the even $K_{2k}$, there are $k$ simple $2$-dimensional $\sl_2$-summands
\[
    q^{-2}a^{2i} \QQ \to q^{2}a^{2i} \QQ,
\]
and $2k + 1$ trivial $1$-dimensional $\sl_2$-summands lying on $\gr_q = 0$. Similarly, for the odd $K_{2k+1}$, there are $k+1$ simple $2$-dimensional $\sl_2$-summands
\[
    q^{-2}a^{2i + 2} \QQ \to q^{2}a^{2i + 2} \QQ,
\]
and the $2k + 1$ trivial $1$-dimensional lying on $\gr_q = 0$. See \Cref{fig:6_1} for the case $K_4 = 6_1^*$. Thus in all of these cases, we directly confirm that \Cref{mainthm:homfly} holds. 

%% file: tikzpictures/hollow-dot.tex
\tikzset{every picture/.style={line width=0.75pt}} 

\begin{tikzpicture}[x=0.75pt,y=0.75pt,yscale=-.75,xscale=.75]

\draw  [dash pattern={on 0.84pt off 2.51pt}] (10.93,11.56) -- (56.56,11.56) -- (56.56,55.49) -- (10.93,55.49) -- cycle ;
\draw  [dash pattern={on 0.84pt off 2.51pt}] (101.41,12.06) -- (148.08,12.06) -- (148.08,56.99) -- (101.41,56.99) -- cycle ;
\draw  [color={rgb, 255:red, 0; green, 0; blue, 0 }  ,draw opacity=1 ][fill={rgb, 255:red, 255; green, 255; blue, 255 }  ,fill opacity=1 ][line width=0.75]  (29.5,34.2) .. controls (29.5,32.43) and (30.93,31) .. (32.7,31) .. controls (34.47,31) and (35.91,32.43) .. (35.91,34.2) .. controls (35.91,35.97) and (34.47,37.41) .. (32.7,37.41) .. controls (30.93,37.41) and (29.5,35.97) .. (29.5,34.2) -- cycle ;
\draw  [fill={rgb, 255:red, 0; green, 0; blue, 0 }  ,fill opacity=1 ] (121.54,34.52) .. controls (121.54,32.75) and (122.97,31.32) .. (124.74,31.32) .. controls (126.51,31.32) and (127.95,32.75) .. (127.95,34.52) .. controls (127.95,36.29) and (126.51,37.73) .. (124.74,37.73) .. controls (122.97,37.73) and (121.54,36.29) .. (121.54,34.52) -- cycle ;
\draw  [dash pattern={on 0.84pt off 2.51pt}] (209.95,11.91) -- (254.54,11.91) -- (254.54,54.84) -- (209.95,54.84) -- cycle ;

\draw (68,25.32) node [anchor=north west][inner sep=0.75pt]    {$=$};
\draw (167,24.24) node [anchor=north west][inner sep=0.75pt]    {$-$};
\draw (194,24.17) node [anchor=north west][inner sep=0.75pt]    {$h$};

\end{tikzpicture}

%% file: tikzpictures/hollow-dot-NC.tex
\tikzset{every picture/.style={line width=0.75pt}} 

\begin{tikzpicture}[x=0.75pt,y=0.75pt,yscale=-1,xscale=1]

\draw  [color={rgb, 255:red, 0; green, 0; blue, 0 }  ,draw opacity=1 ] (13.43,26.19) .. controls (13.43,18.48) and (16,12.23) .. (19.18,12.23) .. controls (22.36,12.23) and (24.93,18.48) .. (24.93,26.19) .. controls (24.93,33.89) and (22.36,40.14) .. (19.18,40.14) .. controls (16,40.14) and (13.43,33.89) .. (13.43,26.19) -- cycle ;
\draw [color={rgb, 255:red, 0; green, 0; blue, 0 }  ,draw opacity=1 ]   (19.18,40.14) -- (78.44,40.14) ;
\draw [color={rgb, 255:red, 0; green, 0; blue, 0 }  ,draw opacity=1 ]   (19.18,12.23) -- (78.44,12.23) ;
\draw  [draw opacity=0][dash pattern={on 0.84pt off 2.51pt}] (76.15,39.93) .. controls (73.93,38.27) and (72.31,32.84) .. (72.31,26.39) .. controls (72.31,19.95) and (73.93,14.52) .. (76.15,12.86) -- (77.58,26.39) -- cycle ; \draw  [dash pattern={on 0.84pt off 2.51pt}] (76.15,39.93) .. controls (73.93,38.27) and (72.31,32.84) .. (72.31,26.39) .. controls (72.31,19.95) and (73.93,14.52) .. (76.15,12.86) ;  
\draw  [draw opacity=0] (77,12.41) .. controls (77.19,12.36) and (77.38,12.33) .. (77.58,12.33) .. controls (80.49,12.33) and (82.85,18.63) .. (82.85,26.39) .. controls (82.85,34.15) and (80.49,40.45) .. (77.58,40.45) .. controls (77.38,40.45) and (77.19,40.42) .. (77,40.37) -- (77.58,26.39) -- cycle ; \draw   (77,12.41) .. controls (77.19,12.36) and (77.38,12.33) .. (77.58,12.33) .. controls (80.49,12.33) and (82.85,18.63) .. (82.85,26.39) .. controls (82.85,34.15) and (80.49,40.45) .. (77.58,40.45) .. controls (77.38,40.45) and (77.19,40.42) .. (77,40.37) ;  

\draw  [color={rgb, 255:red, 0; green, 0; blue, 0 }  ,draw opacity=1 ] (124.2,27.41) .. controls (124.2,19.7) and (126.78,13.46) .. (129.95,13.46) .. controls (133.13,13.46) and (135.71,19.7) .. (135.71,27.41) .. controls (135.71,35.11) and (133.13,41.36) .. (129.95,41.36) .. controls (126.78,41.36) and (124.2,35.11) .. (124.2,27.41) -- cycle ;
\draw  [draw opacity=0] (129.37,13.36) .. controls (129.57,13.35) and (129.76,13.35) .. (129.95,13.35) .. controls (139.08,13.35) and (146.48,19.64) .. (146.48,27.41) .. controls (146.48,35.17) and (139.08,41.47) .. (129.95,41.47) .. controls (129.76,41.47) and (129.57,41.46) .. (129.37,41.46) -- (129.95,27.41) -- cycle ; \draw   (129.37,13.36) .. controls (129.57,13.35) and (129.76,13.35) .. (129.95,13.35) .. controls (139.08,13.35) and (146.48,19.64) .. (146.48,27.41) .. controls (146.48,35.17) and (139.08,41.47) .. (129.95,41.47) .. controls (129.76,41.47) and (129.57,41.46) .. (129.37,41.46) ;  
\draw  [draw opacity=0][dash pattern={on 0.84pt off 2.51pt}] (181.77,40.41) .. controls (179.55,38.76) and (177.92,33.32) .. (177.92,26.87) .. controls (177.92,20.43) and (179.55,15) .. (181.76,13.34) -- (183.19,26.87) -- cycle ; \draw  [dash pattern={on 0.84pt off 2.51pt}] (181.77,40.41) .. controls (179.55,38.76) and (177.92,33.32) .. (177.92,26.87) .. controls (177.92,20.43) and (179.55,15) .. (181.76,13.34) ;  
\draw  [draw opacity=0] (182.62,12.9) .. controls (182.81,12.84) and (183,12.82) .. (183.19,12.82) .. controls (186.1,12.82) and (188.46,19.11) .. (188.46,26.87) .. controls (188.46,34.64) and (186.1,40.93) .. (183.19,40.93) .. controls (183,40.93) and (182.81,40.91) .. (182.62,40.85) -- (183.19,26.87) -- cycle ; \draw   (182.62,12.9) .. controls (182.81,12.84) and (183,12.82) .. (183.19,12.82) .. controls (186.1,12.82) and (188.46,19.11) .. (188.46,26.87) .. controls (188.46,34.64) and (186.1,40.93) .. (183.19,40.93) .. controls (183,40.93) and (182.81,40.91) .. (182.62,40.85) ;  
\draw  [draw opacity=0] (183.77,12.82) .. controls (183.58,12.82) and (183.39,12.82) .. (183.19,12.82) .. controls (174.07,12.82) and (166.67,19.11) .. (166.67,26.87) .. controls (166.67,34.64) and (174.07,40.93) .. (183.19,40.93) .. controls (183.39,40.93) and (183.58,40.93) .. (183.77,40.92) -- (183.19,26.87) -- cycle ; \draw   (183.77,12.82) .. controls (183.58,12.82) and (183.39,12.82) .. (183.19,12.82) .. controls (174.07,12.82) and (166.67,19.11) .. (166.67,26.87) .. controls (166.67,34.64) and (174.07,40.93) .. (183.19,40.93) .. controls (183.39,40.93) and (183.58,40.93) .. (183.77,40.92) ;  

\draw  [color={rgb, 255:red, 0; green, 0; blue, 0 }  ,draw opacity=1 ] (222.73,27.41) .. controls (222.73,19.7) and (225.31,13.46) .. (228.49,13.46) .. controls (231.67,13.46) and (234.24,19.7) .. (234.24,27.41) .. controls (234.24,35.11) and (231.67,41.36) .. (228.49,41.36) .. controls (225.31,41.36) and (222.73,35.11) .. (222.73,27.41) -- cycle ;
\draw  [draw opacity=0] (227.91,13.36) .. controls (228.1,13.35) and (228.29,13.35) .. (228.49,13.35) .. controls (237.61,13.35) and (245.01,19.64) .. (245.01,27.41) .. controls (245.01,35.17) and (237.61,41.47) .. (228.49,41.47) .. controls (228.29,41.47) and (228.1,41.46) .. (227.91,41.46) -- (228.49,27.41) -- cycle ; \draw   (227.91,13.36) .. controls (228.1,13.35) and (228.29,13.35) .. (228.49,13.35) .. controls (237.61,13.35) and (245.01,19.64) .. (245.01,27.41) .. controls (245.01,35.17) and (237.61,41.47) .. (228.49,41.47) .. controls (228.29,41.47) and (228.1,41.46) .. (227.91,41.46) ;  
\draw  [draw opacity=0][dash pattern={on 0.84pt off 2.51pt}] (280.3,40.41) .. controls (278.09,38.76) and (276.46,33.32) .. (276.46,26.87) .. controls (276.46,20.43) and (278.08,15) .. (280.3,13.34) -- (281.73,26.87) -- cycle ; \draw  [dash pattern={on 0.84pt off 2.51pt}] (280.3,40.41) .. controls (278.09,38.76) and (276.46,33.32) .. (276.46,26.87) .. controls (276.46,20.43) and (278.08,15) .. (280.3,13.34) ;  
\draw  [draw opacity=0] (281.15,12.9) .. controls (281.34,12.84) and (281.53,12.82) .. (281.73,12.82) .. controls (284.64,12.82) and (287,19.11) .. (287,26.87) .. controls (287,34.64) and (284.64,40.93) .. (281.73,40.93) .. controls (281.53,40.93) and (281.34,40.91) .. (281.15,40.85) -- (281.73,26.87) -- cycle ; \draw   (281.15,12.9) .. controls (281.34,12.84) and (281.53,12.82) .. (281.73,12.82) .. controls (284.64,12.82) and (287,19.11) .. (287,26.87) .. controls (287,34.64) and (284.64,40.93) .. (281.73,40.93) .. controls (281.53,40.93) and (281.34,40.91) .. (281.15,40.85) ;  
\draw  [draw opacity=0] (282.31,12.82) .. controls (282.11,12.82) and (281.92,12.82) .. (281.73,12.82) .. controls (272.6,12.82) and (265.21,19.11) .. (265.21,26.87) .. controls (265.21,34.64) and (272.6,40.93) .. (281.73,40.93) .. controls (281.92,40.93) and (282.11,40.93) .. (282.31,40.92) -- (281.73,26.87) -- cycle ; \draw   (282.31,12.82) .. controls (282.11,12.82) and (281.92,12.82) .. (281.73,12.82) .. controls (272.6,12.82) and (265.21,19.11) .. (265.21,26.87) .. controls (265.21,34.64) and (272.6,40.93) .. (281.73,40.93) .. controls (281.92,40.93) and (282.11,40.93) .. (282.31,40.92) ;  

\draw  [fill={rgb, 255:red, 0; green, 0; blue, 0 }  ,fill opacity=1 ] (139.92,27.36) .. controls (139.92,26.5) and (140.61,25.8) .. (141.47,25.8) .. controls (142.32,25.8) and (143.02,26.5) .. (143.02,27.36) .. controls (143.02,28.21) and (142.32,28.91) .. (141.47,28.91) .. controls (140.61,28.91) and (139.92,28.21) .. (139.92,27.36) -- cycle ;
\draw  [color={rgb, 255:red, 0; green, 0; blue, 0 }  ,draw opacity=1 ][fill={rgb, 255:red, 255; green, 255; blue, 255 }  ,fill opacity=1 ] (267.09,26.62) .. controls (267.09,25.21) and (268.23,24.06) .. (269.64,24.06) .. controls (271.06,24.06) and (272.2,25.21) .. (272.2,26.62) .. controls (272.2,28.03) and (271.06,29.18) .. (269.64,29.18) .. controls (268.23,29.18) and (267.09,28.03) .. (267.09,26.62) -- cycle ;
\draw  [color={rgb, 255:red, 0; green, 0; blue, 0 }  ,draw opacity=1 ] (123.2,71.99) .. controls (123.2,64.29) and (125.78,58.04) .. (128.95,58.04) .. controls (132.13,58.04) and (134.71,64.29) .. (134.71,71.99) .. controls (134.71,79.7) and (132.13,85.94) .. (128.95,85.94) .. controls (125.78,85.94) and (123.2,79.7) .. (123.2,71.99) -- cycle ;
\draw  [draw opacity=0] (128.37,57.94) .. controls (128.57,57.94) and (128.76,57.93) .. (128.95,57.93) .. controls (138.08,57.93) and (145.48,64.23) .. (145.48,71.99) .. controls (145.48,79.76) and (138.08,86.05) .. (128.95,86.05) .. controls (128.76,86.05) and (128.57,86.05) .. (128.37,86.04) -- (128.95,71.99) -- cycle ; \draw   (128.37,57.94) .. controls (128.57,57.94) and (128.76,57.93) .. (128.95,57.93) .. controls (138.08,57.93) and (145.48,64.23) .. (145.48,71.99) .. controls (145.48,79.76) and (138.08,86.05) .. (128.95,86.05) .. controls (128.76,86.05) and (128.57,86.05) .. (128.37,86.04) ;  
\draw  [draw opacity=0][dash pattern={on 0.84pt off 2.51pt}] (180.77,85) .. controls (178.55,83.34) and (176.92,77.91) .. (176.92,71.46) .. controls (176.92,65.02) and (178.55,59.59) .. (180.76,57.92) -- (182.19,71.46) -- cycle ; \draw  [dash pattern={on 0.84pt off 2.51pt}] (180.77,85) .. controls (178.55,83.34) and (176.92,77.91) .. (176.92,71.46) .. controls (176.92,65.02) and (178.55,59.59) .. (180.76,57.92) ;  
\draw  [draw opacity=0] (181.62,57.48) .. controls (181.81,57.43) and (182,57.4) .. (182.19,57.4) .. controls (185.1,57.4) and (187.46,63.69) .. (187.46,71.46) .. controls (187.46,79.22) and (185.1,85.52) .. (182.19,85.52) .. controls (182,85.52) and (181.81,85.49) .. (181.62,85.43) -- (182.19,71.46) -- cycle ; \draw   (181.62,57.48) .. controls (181.81,57.43) and (182,57.4) .. (182.19,57.4) .. controls (185.1,57.4) and (187.46,63.69) .. (187.46,71.46) .. controls (187.46,79.22) and (185.1,85.52) .. (182.19,85.52) .. controls (182,85.52) and (181.81,85.49) .. (181.62,85.43) ;  
\draw  [draw opacity=0] (182.77,57.41) .. controls (182.58,57.4) and (182.39,57.4) .. (182.19,57.4) .. controls (173.07,57.4) and (165.67,63.69) .. (165.67,71.46) .. controls (165.67,79.22) and (173.07,85.52) .. (182.19,85.52) .. controls (182.39,85.52) and (182.58,85.51) .. (182.77,85.51) -- (182.19,71.46) -- cycle ; \draw   (182.77,57.41) .. controls (182.58,57.4) and (182.39,57.4) .. (182.19,57.4) .. controls (173.07,57.4) and (165.67,63.69) .. (165.67,71.46) .. controls (165.67,79.22) and (173.07,85.52) .. (182.19,85.52) .. controls (182.39,85.52) and (182.58,85.51) .. (182.77,85.51) ;  

\draw  [color={rgb, 255:red, 0; green, 0; blue, 0 }  ,draw opacity=1 ] (221.73,71.99) .. controls (221.73,64.29) and (224.31,58.04) .. (227.49,58.04) .. controls (230.67,58.04) and (233.24,64.29) .. (233.24,71.99) .. controls (233.24,79.7) and (230.67,85.94) .. (227.49,85.94) .. controls (224.31,85.94) and (221.73,79.7) .. (221.73,71.99) -- cycle ;
\draw  [draw opacity=0] (226.91,57.94) .. controls (227.1,57.94) and (227.29,57.93) .. (227.49,57.93) .. controls (236.61,57.93) and (244.01,64.23) .. (244.01,71.99) .. controls (244.01,79.76) and (236.61,86.05) .. (227.49,86.05) .. controls (227.29,86.05) and (227.1,86.05) .. (226.91,86.04) -- (227.49,71.99) -- cycle ; \draw   (226.91,57.94) .. controls (227.1,57.94) and (227.29,57.93) .. (227.49,57.93) .. controls (236.61,57.93) and (244.01,64.23) .. (244.01,71.99) .. controls (244.01,79.76) and (236.61,86.05) .. (227.49,86.05) .. controls (227.29,86.05) and (227.1,86.05) .. (226.91,86.04) ;  
\draw  [draw opacity=0][dash pattern={on 0.84pt off 2.51pt}] (279.3,85) .. controls (277.09,83.34) and (275.46,77.91) .. (275.46,71.46) .. controls (275.46,65.02) and (277.08,59.59) .. (279.3,57.92) -- (280.73,71.46) -- cycle ; \draw  [dash pattern={on 0.84pt off 2.51pt}] (279.3,85) .. controls (277.09,83.34) and (275.46,77.91) .. (275.46,71.46) .. controls (275.46,65.02) and (277.08,59.59) .. (279.3,57.92) ;  
\draw  [draw opacity=0] (280.15,57.48) .. controls (280.34,57.43) and (280.53,57.4) .. (280.73,57.4) .. controls (283.64,57.4) and (286,63.69) .. (286,71.46) .. controls (286,79.22) and (283.64,85.52) .. (280.73,85.52) .. controls (280.53,85.52) and (280.34,85.49) .. (280.15,85.43) -- (280.73,71.46) -- cycle ; \draw   (280.15,57.48) .. controls (280.34,57.43) and (280.53,57.4) .. (280.73,57.4) .. controls (283.64,57.4) and (286,63.69) .. (286,71.46) .. controls (286,79.22) and (283.64,85.52) .. (280.73,85.52) .. controls (280.53,85.52) and (280.34,85.49) .. (280.15,85.43) ;  
\draw  [draw opacity=0] (281.31,57.41) .. controls (281.11,57.4) and (280.92,57.4) .. (280.73,57.4) .. controls (271.6,57.4) and (264.21,63.69) .. (264.21,71.46) .. controls (264.21,79.22) and (271.6,85.52) .. (280.73,85.52) .. controls (280.92,85.52) and (281.11,85.51) .. (281.31,85.51) -- (280.73,71.46) -- cycle ; \draw   (281.31,57.41) .. controls (281.11,57.4) and (280.92,57.4) .. (280.73,57.4) .. controls (271.6,57.4) and (264.21,63.69) .. (264.21,71.46) .. controls (264.21,79.22) and (271.6,85.52) .. (280.73,85.52) .. controls (280.92,85.52) and (281.11,85.51) .. (281.31,85.51) ;  

\draw  [color={rgb, 255:red, 0; green, 0; blue, 0 }  ,draw opacity=1 ][fill={rgb, 255:red, 255; green, 255; blue, 255 }  ,fill opacity=1 ] (137.8,71.94) .. controls (137.8,70.47) and (138.99,69.27) .. (140.47,69.27) .. controls (141.94,69.27) and (143.13,70.47) .. (143.13,71.94) .. controls (143.13,73.41) and (141.94,74.61) .. (140.47,74.61) .. controls (138.99,74.61) and (137.8,73.41) .. (137.8,71.94) -- cycle ;
\draw  [fill={rgb, 255:red, 0; green, 0; blue, 0 }  ,fill opacity=1 ] (267.09,71.2) .. controls (267.09,70.35) and (267.79,69.65) .. (268.64,69.65) .. controls (269.5,69.65) and (270.2,70.35) .. (270.2,71.2) .. controls (270.2,72.06) and (269.5,72.75) .. (268.64,72.75) .. controls (267.79,72.75) and (267.09,72.06) .. (267.09,71.2) -- cycle ;

\draw (89.28,18.42) node [anchor=north west][inner sep=0.75pt]    {$=$};
\draw (199.67,19.22) node [anchor=north west][inner sep=0.75pt]    {$+$};
\draw (88.28,63) node [anchor=north west][inner sep=0.75pt]    {$=$};
\draw (198.67,63.8) node [anchor=north west][inner sep=0.75pt]    {$+$};

\end{tikzpicture}

%% file: tikzpictures/hollow-dot-relations.tex
\tikzset{every picture/.style={line width=0.75pt}} 

\begin{tikzpicture}[x=0.75pt,y=0.75pt,yscale=-.75,xscale=.75]

\draw  [dash pattern={on 0.84pt off 2.51pt}] (10.93,11.56) -- (56.56,11.56) -- (56.56,55.49) -- (10.93,55.49) -- cycle ;
\draw  [dash pattern={on 0.84pt off 2.51pt}] (128.41,11.06) -- (175.08,11.06) -- (175.08,55.99) -- (128.41,55.99) -- cycle ;
\draw  [dash pattern={on 0.84pt off 2.51pt}] (230.95,12.06) -- (275.54,12.06) -- (275.54,54.99) -- (230.95,54.99) -- cycle ;
\draw  [color={rgb, 255:red, 0; green, 0; blue, 0 }  ,draw opacity=1 ][fill={rgb, 255:red, 255; green, 255; blue, 255 }  ,fill opacity=1 ] (23.5,34.2) .. controls (23.5,32.43) and (24.93,31) .. (26.7,31) .. controls (28.47,31) and (29.91,32.43) .. (29.91,34.2) .. controls (29.91,35.97) and (28.47,37.41) .. (26.7,37.41) .. controls (24.93,37.41) and (23.5,35.97) .. (23.5,34.2) -- cycle ;
\draw  [color={rgb, 255:red, 0; green, 0; blue, 0 }  ,draw opacity=1 ][fill={rgb, 255:red, 255; green, 255; blue, 255 }  ,fill opacity=1 ] (36,34.2) .. controls (36,32.43) and (37.43,31) .. (39.2,31) .. controls (40.97,31) and (42.41,32.43) .. (42.41,34.2) .. controls (42.41,35.97) and (40.97,37.41) .. (39.2,37.41) .. controls (37.43,37.41) and (36,35.97) .. (36,34.2) -- cycle ;
\draw  [color={rgb, 255:red, 0; green, 0; blue, 0 }  ,draw opacity=1 ][fill={rgb, 255:red, 255; green, 255; blue, 255 }  ,fill opacity=1 ] (148.54,33.52) .. controls (148.54,31.75) and (149.97,30.32) .. (151.74,30.32) .. controls (153.51,30.32) and (154.95,31.75) .. (154.95,33.52) .. controls (154.95,35.29) and (153.51,36.73) .. (151.74,36.73) .. controls (149.97,36.73) and (148.54,35.29) .. (148.54,33.52) -- cycle ;
\draw  [dash pattern={on 0.84pt off 2.51pt}] (9.93,68.56) -- (55.56,68.56) -- (55.56,112.49) -- (9.93,112.49) -- cycle ;
\draw  [dash pattern={on 0.84pt off 2.51pt}] (129.95,70.06) -- (174.54,70.06) -- (174.54,112.99) -- (129.95,112.99) -- cycle ;
\draw  [color={rgb, 255:red, 0; green, 0; blue, 0 }  ,draw opacity=1 ][fill={rgb, 255:red, 0; green, 0; blue, 0 }  ,fill opacity=1 ] (22.5,91.2) .. controls (22.5,89.43) and (23.93,88) .. (25.7,88) .. controls (27.47,88) and (28.91,89.43) .. (28.91,91.2) .. controls (28.91,92.97) and (27.47,94.41) .. (25.7,94.41) .. controls (23.93,94.41) and (22.5,92.97) .. (22.5,91.2) -- cycle ;
\draw  [color={rgb, 255:red, 0; green, 0; blue, 0 }  ,draw opacity=1 ][fill={rgb, 255:red, 255; green, 255; blue, 255 }  ,fill opacity=1 ] (35,91.2) .. controls (35,89.43) and (36.43,88) .. (38.2,88) .. controls (39.97,88) and (41.41,89.43) .. (41.41,91.2) .. controls (41.41,92.97) and (39.97,94.41) .. (38.2,94.41) .. controls (36.43,94.41) and (35,92.97) .. (35,91.2) -- cycle ;
\draw  [draw opacity=0] (53.81,147.4) .. controls (53.07,150) and (45.16,152.04) .. (35.51,152.04) .. controls (25.37,152.04) and (17.15,149.78) .. (17.15,147) .. controls (17.15,146.83) and (17.18,146.66) .. (17.24,146.5) -- (35.51,147) -- cycle ; \draw  [color={rgb, 255:red, 128; green, 128; blue, 128 }  ,draw opacity=1 ] (53.81,147.4) .. controls (53.07,150) and (45.16,152.04) .. (35.51,152.04) .. controls (25.37,152.04) and (17.15,149.78) .. (17.15,147) .. controls (17.15,146.83) and (17.18,146.66) .. (17.24,146.5) ;  
\draw  [draw opacity=0][dash pattern={on 0.84pt off 2.51pt}] (17.6,145.75) .. controls (19.07,143.37) and (26.62,141.56) .. (35.71,141.56) .. controls (45.7,141.56) and (53.83,143.75) .. (54.06,146.48) -- (35.71,146.6) -- cycle ; \draw  [color={rgb, 255:red, 128; green, 128; blue, 128 }  ,draw opacity=1 ][dash pattern={on 0.84pt off 2.51pt}] (17.6,145.75) .. controls (19.07,143.37) and (26.62,141.56) .. (35.71,141.56) .. controls (45.7,141.56) and (53.83,143.75) .. (54.06,146.48) ;  
\draw   (17.01,147) .. controls (17.01,136.78) and (25.29,128.5) .. (35.51,128.5) .. controls (45.73,128.5) and (54.01,136.78) .. (54.01,147) .. controls (54.01,157.22) and (45.73,165.5) .. (35.51,165.5) .. controls (25.29,165.5) and (17.01,157.22) .. (17.01,147) -- cycle ;
\draw  [color={rgb, 255:red, 0; green, 0; blue, 0 }  ,draw opacity=1 ][fill={rgb, 255:red, 255; green, 255; blue, 255 }  ,fill opacity=1 ] (32.31,135.7) .. controls (32.31,133.93) and (33.74,132.5) .. (35.51,132.5) .. controls (37.28,132.5) and (38.71,133.93) .. (38.71,135.7) .. controls (38.71,137.47) and (37.28,138.91) .. (35.51,138.91) .. controls (33.74,138.91) and (32.31,137.47) .. (32.31,135.7) -- cycle ;

\draw (67,25.32) node [anchor=north west][inner sep=0.75pt]    {$=$};
\draw (189,25.32) node [anchor=north west][inner sep=0.75pt]    {$+$};
\draw (93,25.32) node [anchor=north west][inner sep=0.75pt]    {$\ -h$};
\draw (213,25.32) node [anchor=north west][inner sep=0.75pt]    {$t$};
\draw (67,82.32) node [anchor=north west][inner sep=0.75pt]    {$=$};
\draw (111,83.32) node [anchor=north west][inner sep=0.75pt]    {$t$};
\draw (67,137.4) node [anchor=north west][inner sep=0.75pt]    {$=\ 1$};

\end{tikzpicture}

%% file: tikzpictures/delooping.tex
\tikzset{every picture/.style={line width=0.75pt}} 

\begin{tikzpicture}[x=0.75pt,y=0.75pt,yscale=-.8,xscale=.8]

\draw  [dash pattern={on 0.84pt off 2.51pt}] (14,105) .. controls (14,90.64) and (25.64,79) .. (40,79) .. controls (54.36,79) and (66,90.64) .. (66,105) .. controls (66,119.36) and (54.36,131) .. (40,131) .. controls (25.64,131) and (14,119.36) .. (14,105) -- cycle ;
\draw    (86,85.49) -- (121.42,57.83) ;
\draw [shift={(123,56.6)}, rotate = 142.01] [color={rgb, 255:red, 0; green, 0; blue, 0 }  ][line width=0.75]    (10.93,-4.9) .. controls (6.95,-2.3) and (3.31,-0.67) .. (0,0) .. controls (3.31,0.67) and (6.95,2.3) .. (10.93,4.9)   ;
\draw    (223,146.4) -- (258.42,118.74) ;
\draw [shift={(260,117.51)}, rotate = 142.01] [color={rgb, 255:red, 0; green, 0; blue, 0 }  ][line width=0.75]    (10.93,-4.9) .. controls (6.95,-2.3) and (3.31,-0.67) .. (0,0) .. controls (3.31,0.67) and (6.95,2.3) .. (10.93,4.9)   ;
\draw    (86,117.51) -- (121.42,145.17) ;
\draw [shift={(123,146.4)}, rotate = 217.99] [color={rgb, 255:red, 0; green, 0; blue, 0 }  ][line width=0.75]    (10.93,-4.9) .. controls (6.95,-2.3) and (3.31,-0.67) .. (0,0) .. controls (3.31,0.67) and (6.95,2.3) .. (10.93,4.9)   ;
\draw    (223,56.6) -- (258.42,84.26) ;
\draw [shift={(260,85.49)}, rotate = 217.99] [color={rgb, 255:red, 0; green, 0; blue, 0 }  ][line width=0.75]    (10.93,-4.9) .. controls (6.95,-2.3) and (3.31,-0.67) .. (0,0) .. controls (3.31,0.67) and (6.95,2.3) .. (10.93,4.9)   ;
\draw  [color={rgb, 255:red, 0; green, 0; blue, 0 }  ,draw opacity=1 ] (82.6,152.34) .. controls (82.6,144.64) and (85.18,138.39) .. (88.35,138.39) .. controls (91.53,138.39) and (94.11,144.64) .. (94.11,152.34) .. controls (94.11,160.05) and (91.53,166.29) .. (88.35,166.29) .. controls (85.18,166.29) and (82.6,160.05) .. (82.6,152.34) -- cycle ;
\draw  [draw opacity=0] (87.77,138.29) .. controls (87.97,138.29) and (88.16,138.28) .. (88.35,138.28) .. controls (97.48,138.28) and (104.88,144.58) .. (104.88,152.34) .. controls (104.88,160.11) and (97.48,166.4) .. (88.35,166.4) .. controls (88.16,166.4) and (87.97,166.4) .. (87.77,166.39) -- (88.35,152.34) -- cycle ; \draw   (87.77,138.29) .. controls (87.97,138.29) and (88.16,138.28) .. (88.35,138.28) .. controls (97.48,138.28) and (104.88,144.58) .. (104.88,152.34) .. controls (104.88,160.11) and (97.48,166.4) .. (88.35,166.4) .. controls (88.16,166.4) and (87.97,166.4) .. (87.77,166.39) ;  
\draw  [color={rgb, 255:red, 0; green, 0; blue, 0 }  ,draw opacity=1 ][fill={rgb, 255:red, 255; green, 255; blue, 255 }  ,fill opacity=1 ] (97,152.29) .. controls (97,150.71) and (98.28,149.42) .. (99.87,149.42) .. controls (101.45,149.42) and (102.73,150.71) .. (102.73,152.29) .. controls (102.73,153.87) and (101.45,155.16) .. (99.87,155.16) .. controls (98.28,155.16) and (97,153.87) .. (97,152.29) -- cycle ;

\draw  [color={rgb, 255:red, 0; green, 0; blue, 0 }  ,draw opacity=1 ] (82.6,49.61) .. controls (82.6,41.9) and (85.18,35.66) .. (88.35,35.66) .. controls (91.53,35.66) and (94.11,41.9) .. (94.11,49.61) .. controls (94.11,57.31) and (91.53,63.56) .. (88.35,63.56) .. controls (85.18,63.56) and (82.6,57.31) .. (82.6,49.61) -- cycle ;
\draw  [draw opacity=0] (87.77,35.56) .. controls (87.97,35.55) and (88.16,35.55) .. (88.35,35.55) .. controls (97.48,35.55) and (104.88,41.84) .. (104.88,49.61) .. controls (104.88,57.37) and (97.48,63.67) .. (88.35,63.67) .. controls (88.16,63.67) and (87.97,63.66) .. (87.77,63.66) -- (88.35,49.61) -- cycle ; \draw   (87.77,35.56) .. controls (87.97,35.55) and (88.16,35.55) .. (88.35,35.55) .. controls (97.48,35.55) and (104.88,41.84) .. (104.88,49.61) .. controls (104.88,57.37) and (97.48,63.67) .. (88.35,63.67) .. controls (88.16,63.67) and (87.97,63.66) .. (87.77,63.66) ;  

\draw  [draw opacity=0][dash pattern={on 0.84pt off 2.51pt}] (256.7,165.88) .. controls (254.49,164.23) and (252.86,158.79) .. (252.86,152.34) .. controls (252.86,145.9) and (254.48,140.47) .. (256.7,138.81) -- (258.13,152.34) -- cycle ; \draw  [dash pattern={on 0.84pt off 2.51pt}] (256.7,165.88) .. controls (254.49,164.23) and (252.86,158.79) .. (252.86,152.34) .. controls (252.86,145.9) and (254.48,140.47) .. (256.7,138.81) ;  
\draw  [draw opacity=0] (257.55,138.37) .. controls (257.74,138.31) and (257.93,138.28) .. (258.13,138.28) .. controls (261.04,138.28) and (263.4,144.58) .. (263.4,152.34) .. controls (263.4,160.11) and (261.04,166.4) .. (258.13,166.4) .. controls (257.93,166.4) and (257.74,166.37) .. (257.55,166.32) -- (258.13,152.34) -- cycle ; \draw   (257.55,138.37) .. controls (257.74,138.31) and (257.93,138.28) .. (258.13,138.28) .. controls (261.04,138.28) and (263.4,144.58) .. (263.4,152.34) .. controls (263.4,160.11) and (261.04,166.4) .. (258.13,166.4) .. controls (257.93,166.4) and (257.74,166.37) .. (257.55,166.32) ;  
\draw  [draw opacity=0] (258.71,138.29) .. controls (258.51,138.29) and (258.32,138.28) .. (258.13,138.28) .. controls (249,138.28) and (241.61,144.58) .. (241.61,152.34) .. controls (241.61,160.11) and (249,166.4) .. (258.13,166.4) .. controls (258.32,166.4) and (258.51,166.4) .. (258.71,166.39) -- (258.13,152.34) -- cycle ; \draw   (258.71,138.29) .. controls (258.51,138.29) and (258.32,138.28) .. (258.13,138.28) .. controls (249,138.28) and (241.61,144.58) .. (241.61,152.34) .. controls (241.61,160.11) and (249,166.4) .. (258.13,166.4) .. controls (258.32,166.4) and (258.51,166.4) .. (258.71,166.39) ;  

\draw  [draw opacity=0][dash pattern={on 0.84pt off 2.51pt}] (256.7,60.61) .. controls (254.49,58.96) and (252.86,53.52) .. (252.86,47.07) .. controls (252.86,40.63) and (254.48,35.2) .. (256.7,33.54) -- (258.13,47.07) -- cycle ; \draw  [dash pattern={on 0.84pt off 2.51pt}] (256.7,60.61) .. controls (254.49,58.96) and (252.86,53.52) .. (252.86,47.07) .. controls (252.86,40.63) and (254.48,35.2) .. (256.7,33.54) ;  
\draw  [draw opacity=0] (257.55,33.1) .. controls (257.74,33.04) and (257.93,33.02) .. (258.13,33.02) .. controls (261.04,33.02) and (263.4,39.31) .. (263.4,47.07) .. controls (263.4,54.84) and (261.04,61.13) .. (258.13,61.13) .. controls (257.93,61.13) and (257.74,61.11) .. (257.55,61.05) -- (258.13,47.07) -- cycle ; \draw   (257.55,33.1) .. controls (257.74,33.04) and (257.93,33.02) .. (258.13,33.02) .. controls (261.04,33.02) and (263.4,39.31) .. (263.4,47.07) .. controls (263.4,54.84) and (261.04,61.13) .. (258.13,61.13) .. controls (257.93,61.13) and (257.74,61.11) .. (257.55,61.05) ;  
\draw  [draw opacity=0] (258.71,33.02) .. controls (258.51,33.02) and (258.32,33.02) .. (258.13,33.02) .. controls (249,33.02) and (241.61,39.31) .. (241.61,47.07) .. controls (241.61,54.84) and (249,61.13) .. (258.13,61.13) .. controls (258.32,61.13) and (258.51,61.13) .. (258.71,61.12) -- (258.13,47.07) -- cycle ; \draw   (258.71,33.02) .. controls (258.51,33.02) and (258.32,33.02) .. (258.13,33.02) .. controls (249,33.02) and (241.61,39.31) .. (241.61,47.07) .. controls (241.61,54.84) and (249,61.13) .. (258.13,61.13) .. controls (258.32,61.13) and (258.51,61.13) .. (258.71,61.12) ;  
\draw  [fill={rgb, 255:red, 0; green, 0; blue, 0 }  ,fill opacity=1 ] (244.49,46.82) .. controls (244.49,45.96) and (245.19,45.27) .. (246.04,45.27) .. controls (246.9,45.27) and (247.6,45.96) .. (247.6,46.82) .. controls (247.6,47.68) and (246.9,48.37) .. (246.04,48.37) .. controls (245.19,48.37) and (244.49,47.68) .. (244.49,46.82) -- cycle ;

\draw  [dash pattern={on 0.84pt off 2.51pt}] (279,103) .. controls (279,88.64) and (290.64,77) .. (305,77) .. controls (319.36,77) and (331,88.64) .. (331,103) .. controls (331,117.36) and (319.36,129) .. (305,129) .. controls (290.64,129) and (279,117.36) .. (279,103) -- cycle ;
\draw  [dash pattern={on 0.84pt off 2.51pt}] (158.5,157) .. controls (158.5,142.64) and (170.14,131) .. (184.5,131) .. controls (198.86,131) and (210.5,142.64) .. (210.5,157) .. controls (210.5,171.36) and (198.86,183) .. (184.5,183) .. controls (170.14,183) and (158.5,171.36) .. (158.5,157) -- cycle ;

\draw  [dash pattern={on 0.84pt off 2.51pt}] (158.5,43) .. controls (158.5,28.64) and (170.14,17) .. (184.5,17) .. controls (198.86,17) and (210.5,28.64) .. (210.5,43) .. controls (210.5,57.36) and (198.86,69) .. (184.5,69) .. controls (170.14,69) and (158.5,57.36) .. (158.5,43) -- cycle ;

\draw [color={rgb, 255:red, 255; green, 255; blue, 255 }  ,draw opacity=1 ][fill={rgb, 255:red, 255; green, 255; blue, 255 }  ,fill opacity=1 ]   (66,105) -- (279,103) ;
\draw  [line width=1.5]  (27,105) .. controls (27,97.82) and (32.82,92) .. (40,92) .. controls (47.18,92) and (53,97.82) .. (53,105) .. controls (53,112.18) and (47.18,118) .. (40,118) .. controls (32.82,118) and (27,112.18) .. (27,105) -- cycle ;
\draw  [line width=1.5]  (292,103) .. controls (292,95.82) and (297.82,90) .. (305,90) .. controls (312.18,90) and (318,95.82) .. (318,103) .. controls (318,110.18) and (312.18,116) .. (305,116) .. controls (297.82,116) and (292,110.18) .. (292,103) -- cycle ;

\draw (172.5,104) node  [font=\large]  {$\oplus $};
\draw (156.5,157) node [anchor=east] [inner sep=0.75pt]    {$q^{-1}$};
\draw (156.5,43) node [anchor=east] [inner sep=0.75pt]    {$q^{1}$};

\end{tikzpicture}

%% file: tikzpictures/T_q.tex
\tikzset{every picture/.style={line width=0.75pt}} 

\begin{tikzpicture}[x=0.75pt,y=0.75pt,yscale=-1,xscale=1]

\draw  [color={rgb, 255:red, 255; green, 255; blue, 255 }  ,draw opacity=1 ][fill={rgb, 255:red, 255; green, 255; blue, 255 }  ,fill opacity=1 ] (70,40) -- (100,40) -- (100,70) -- (70,70) -- cycle ;

\draw  [color={rgb, 255:red, 128; green, 128; blue, 128 }  ,draw opacity=1 ] (106.33,98.33) .. controls (111,98.33) and (113.33,96) .. (113.33,91.33) -- (113.33,65.08) .. controls (113.33,58.41) and (115.66,55.08) .. (120.33,55.08) .. controls (115.66,55.08) and (113.33,51.75) .. (113.33,45.08)(113.33,48.08) -- (113.33,18.83) .. controls (113.33,14.16) and (111,11.83) .. (106.33,11.83) ;
\draw    (70,70) -- (100,100) ;
\draw  [color={rgb, 255:red, 255; green, 255; blue, 255 }  ,draw opacity=1 ][fill={rgb, 255:red, 255; green, 255; blue, 255 }  ,fill opacity=1 ] (80,85) .. controls (80,82.24) and (82.24,80) .. (85,80) .. controls (87.76,80) and (90,82.24) .. (90,85) .. controls (90,87.76) and (87.76,90) .. (85,90) .. controls (82.24,90) and (80,87.76) .. (80,85) -- cycle ;
\draw    (100,70) -- (70,100) ;

\draw    (70,10) -- (100,40) ;
\draw  [color={rgb, 255:red, 255; green, 255; blue, 255 }  ,draw opacity=1 ][fill={rgb, 255:red, 255; green, 255; blue, 255 }  ,fill opacity=1 ] (80,25) .. controls (80,22.24) and (82.24,20) .. (85,20) .. controls (87.76,20) and (90,22.24) .. (90,25) .. controls (90,27.76) and (87.76,30) .. (85,30) .. controls (82.24,30) and (80,27.76) .. (80,25) -- cycle ;
\draw    (100,10) -- (70,40) ;

\draw (59,53.5) node [anchor=east] [inner sep=0.75pt]    {$T_{q} \ =$};
\draw (124.33,53.83) node [anchor=west] [inner sep=0.75pt]    {$q$};
\draw (85,55) node    {$\vdots $};

\end{tikzpicture}

%% file: tikzpictures/T_q-ori.tex
\tikzset{every picture/.style={line width=0.75pt}} 

\begin{tikzpicture}[x=0.75pt,y=0.75pt,yscale=-1,xscale=1]

\draw    (70,50) -- (100,80) ;
\draw  [color={rgb, 255:red, 255; green, 255; blue, 255 }  ,draw opacity=1 ][fill={rgb, 255:red, 255; green, 255; blue, 255 }  ,fill opacity=1 ] (80,65) .. controls (80,62.24) and (82.24,60) .. (85,60) .. controls (87.76,60) and (90,62.24) .. (90,65) .. controls (90,67.76) and (87.76,70) .. (85,70) .. controls (82.24,70) and (80,67.76) .. (80,65) -- cycle ;
\draw    (100,50) -- (70,80) ;

\draw    (71.41,21.41) -- (100,50) ;
\draw [shift={(70,20)}, rotate = 45] [color={rgb, 255:red, 0; green, 0; blue, 0 }  ][line width=0.75]    (6.56,-1.97) .. controls (4.17,-0.84) and (1.99,-0.18) .. (0,0) .. controls (1.99,0.18) and (4.17,0.84) .. (6.56,1.97)   ;
\draw  [color={rgb, 255:red, 255; green, 255; blue, 255 }  ,draw opacity=1 ][fill={rgb, 255:red, 255; green, 255; blue, 255 }  ,fill opacity=1 ] (80,35) .. controls (80,32.24) and (82.24,30) .. (85,30) .. controls (87.76,30) and (90,32.24) .. (90,35) .. controls (90,37.76) and (87.76,40) .. (85,40) .. controls (82.24,40) and (80,37.76) .. (80,35) -- cycle ;
\draw    (98.59,21.41) -- (70,50) ;
\draw [shift={(100,20)}, rotate = 135] [color={rgb, 255:red, 0; green, 0; blue, 0 }  ][line width=0.75]    (6.56,-1.97) .. controls (4.17,-0.84) and (1.99,-0.18) .. (0,0) .. controls (1.99,0.18) and (4.17,0.84) .. (6.56,1.97)   ;

\draw    (190,50) -- (218.59,78.59) ;
\draw [shift={(220,80)}, rotate = 225] [color={rgb, 255:red, 0; green, 0; blue, 0 }  ][line width=0.75]    (6.56,-1.97) .. controls (4.17,-0.84) and (1.99,-0.18) .. (0,0) .. controls (1.99,0.18) and (4.17,0.84) .. (6.56,1.97)   ;
\draw  [color={rgb, 255:red, 255; green, 255; blue, 255 }  ,draw opacity=1 ][fill={rgb, 255:red, 255; green, 255; blue, 255 }  ,fill opacity=1 ] (200,65) .. controls (200,62.24) and (202.24,60) .. (205,60) .. controls (207.76,60) and (210,62.24) .. (210,65) .. controls (210,67.76) and (207.76,70) .. (205,70) .. controls (202.24,70) and (200,67.76) .. (200,65) -- cycle ;
\draw    (220,50) -- (191.41,78.59) ;
\draw [shift={(190,80)}, rotate = 315] [color={rgb, 255:red, 0; green, 0; blue, 0 }  ][line width=0.75]    (6.56,-1.97) .. controls (4.17,-0.84) and (1.99,-0.18) .. (0,0) .. controls (1.99,0.18) and (4.17,0.84) .. (6.56,1.97)   ;
\draw    (190,20) -- (220,50) ;
\draw  [color={rgb, 255:red, 255; green, 255; blue, 255 }  ,draw opacity=1 ][fill={rgb, 255:red, 255; green, 255; blue, 255 }  ,fill opacity=1 ] (200,35) .. controls (200,32.24) and (202.24,30) .. (205,30) .. controls (207.76,30) and (210,32.24) .. (210,35) .. controls (210,37.76) and (207.76,40) .. (205,40) .. controls (202.24,40) and (200,37.76) .. (200,35) -- cycle ;
\draw    (220,20) -- (190,50) ;
\draw    (310,50) -- (338.59,78.59) ;
\draw [shift={(340,80)}, rotate = 225] [color={rgb, 255:red, 0; green, 0; blue, 0 }  ][line width=0.75]    (6.56,-1.97) .. controls (4.17,-0.84) and (1.99,-0.18) .. (0,0) .. controls (1.99,0.18) and (4.17,0.84) .. (6.56,1.97)   ;
\draw  [color={rgb, 255:red, 255; green, 255; blue, 255 }  ,draw opacity=1 ][fill={rgb, 255:red, 255; green, 255; blue, 255 }  ,fill opacity=1 ] (320,65) .. controls (320,62.24) and (322.24,60) .. (325,60) .. controls (327.76,60) and (330,62.24) .. (330,65) .. controls (330,67.76) and (327.76,70) .. (325,70) .. controls (322.24,70) and (320,67.76) .. (320,65) -- cycle ;
\draw    (340,50) -- (310,80) ;
\draw    (311.41,21.41) -- (340,50) ;
\draw [shift={(310,20)}, rotate = 45] [color={rgb, 255:red, 0; green, 0; blue, 0 }  ][line width=0.75]    (6.56,-1.97) .. controls (4.17,-0.84) and (1.99,-0.18) .. (0,0) .. controls (1.99,0.18) and (4.17,0.84) .. (6.56,1.97)   ;
\draw  [color={rgb, 255:red, 255; green, 255; blue, 255 }  ,draw opacity=1 ][fill={rgb, 255:red, 255; green, 255; blue, 255 }  ,fill opacity=1 ] (320,35) .. controls (320,32.24) and (322.24,30) .. (325,30) .. controls (327.76,30) and (330,32.24) .. (330,35) .. controls (330,37.76) and (327.76,40) .. (325,40) .. controls (322.24,40) and (320,37.76) .. (320,35) -- cycle ;
\draw    (340,20) -- (310,50) ;
\draw    (430,50) -- (460,80) ;
\draw  [color={rgb, 255:red, 255; green, 255; blue, 255 }  ,draw opacity=1 ][fill={rgb, 255:red, 255; green, 255; blue, 255 }  ,fill opacity=1 ] (440,65) .. controls (440,62.24) and (442.24,60) .. (445,60) .. controls (447.76,60) and (450,62.24) .. (450,65) .. controls (450,67.76) and (447.76,70) .. (445,70) .. controls (442.24,70) and (440,67.76) .. (440,65) -- cycle ;
\draw    (460,50) -- (431.41,78.59) ;
\draw [shift={(430,80)}, rotate = 315] [color={rgb, 255:red, 0; green, 0; blue, 0 }  ][line width=0.75]    (6.56,-1.97) .. controls (4.17,-0.84) and (1.99,-0.18) .. (0,0) .. controls (1.99,0.18) and (4.17,0.84) .. (6.56,1.97)   ;
\draw    (430,20) -- (460,50) ;
\draw  [color={rgb, 255:red, 255; green, 255; blue, 255 }  ,draw opacity=1 ][fill={rgb, 255:red, 255; green, 255; blue, 255 }  ,fill opacity=1 ] (440,35) .. controls (440,32.24) and (442.24,30) .. (445,30) .. controls (447.76,30) and (450,32.24) .. (450,35) .. controls (450,37.76) and (447.76,40) .. (445,40) .. controls (442.24,40) and (440,37.76) .. (440,35) -- cycle ;
\draw    (458.59,21.41) -- (430,50) ;
\draw [shift={(460,20)}, rotate = 135] [color={rgb, 255:red, 0; green, 0; blue, 0 }  ][line width=0.75]    (6.56,-1.97) .. controls (4.17,-0.84) and (1.99,-0.18) .. (0,0) .. controls (1.99,0.18) and (4.17,0.84) .. (6.56,1.97)   ;

\draw (68,50) node [anchor=east] [inner sep=0.75pt]    {$T_{2}^{\uparrow \uparrow } \ =$};
\draw (428,50) node [anchor=east] [inner sep=0.75pt]    {$T_{2}^{\downarrow \uparrow } \ =$};
\draw (188,50) node [anchor=east] [inner sep=0.75pt]    {$T_{2}^{\downarrow \downarrow } \ =$};
\draw (308,50) node [anchor=east] [inner sep=0.75pt]    {$T_{2}^{\uparrow \downarrow } \ =$};
\draw (70,83.4) node [anchor=north] [inner sep=0.75pt]  [font=\scriptsize]  {$P_{1}$};
\draw (100,83.4) node [anchor=north] [inner sep=0.75pt]  [font=\scriptsize]  {$P_{2}$};
\draw (220,16.6) node [anchor=south] [inner sep=0.75pt]  [font=\scriptsize]  {$P_{1}$};
\draw (190,16.6) node [anchor=south] [inner sep=0.75pt]  [font=\scriptsize]  {$P_{2}$};
\draw (310,83.4) node [anchor=north] [inner sep=0.75pt]  [font=\scriptsize]  {$P_{1}$};
\draw (340,16.6) node [anchor=south] [inner sep=0.75pt]  [font=\scriptsize]  {$P_{2}$};
\draw (460,83.4) node [anchor=north] [inner sep=0.75pt]  [font=\scriptsize]  {$P_{1}$};
\draw (430,16.6) node [anchor=south] [inner sep=0.75pt]  [font=\scriptsize]  {$P_{2}$};

\end{tikzpicture}

%% file: tikzpictures/E01.tex
\tikzset{every picture/.style={line width=0.75pt}} 

\begin{tikzpicture}[x=0.75pt,y=0.75pt,yscale=-.75,xscale=.75]

\draw  [dash pattern={on 4.5pt off 4.5pt}] (206,34.77) .. controls (206,20.71) and (217.4,9.31) .. (231.46,9.31) .. controls (245.52,9.31) and (256.91,20.71) .. (256.91,34.77) .. controls (256.91,48.83) and (245.52,60.22) .. (231.46,60.22) .. controls (217.4,60.22) and (206,48.83) .. (206,34.77) -- cycle ;
\draw  [draw opacity=0][line width=1.5]  (251.56,20.25) .. controls (246.19,23.83) and (239.17,26) .. (231.48,26) .. controls (223.83,26) and (216.84,23.85) .. (211.48,20.3) -- (231.48,2.73) -- cycle ; \draw  [line width=1.5]  (251.56,20.25) .. controls (246.19,23.83) and (239.17,26) .. (231.48,26) .. controls (223.83,26) and (216.84,23.85) .. (211.48,20.3) ;  
\draw  [draw opacity=0][line width=1.5]  (211.4,49.29) .. controls (216.77,45.71) and (223.79,43.54) .. (231.48,43.54) .. controls (239.13,43.54) and (246.12,45.69) .. (251.48,49.24) -- (231.48,66.81) -- cycle ; \draw  [line width=1.5]  (211.4,49.29) .. controls (216.77,45.71) and (223.79,43.54) .. (231.48,43.54) .. controls (239.13,43.54) and (246.12,45.69) .. (251.48,49.24) ;  

\draw  [dash pattern={on 4.5pt off 4.5pt}] (86.48,9.29) .. controls (100.54,9.29) and (111.94,20.69) .. (111.94,34.75) .. controls (111.94,48.81) and (100.54,60.21) .. (86.48,60.21) .. controls (72.42,60.21) and (61.02,48.81) .. (61.02,34.75) .. controls (61.02,20.69) and (72.42,9.29) .. (86.48,9.29) -- cycle ;
\draw  [draw opacity=0][line width=1.5]  (101,54.85) .. controls (97.42,49.48) and (95.25,42.46) .. (95.25,34.77) .. controls (95.25,27.12) and (97.4,20.13) .. (100.95,14.77) -- (118.52,34.77) -- cycle ; \draw  [line width=1.5]  (101,54.85) .. controls (97.42,49.48) and (95.25,42.46) .. (95.25,34.77) .. controls (95.25,27.12) and (97.4,20.13) .. (100.95,14.77) ;  
\draw  [draw opacity=0][line width=1.5]  (71.96,14.69) .. controls (75.54,20.06) and (77.71,27.08) .. (77.71,34.77) .. controls (77.71,42.42) and (75.56,49.41) .. (72.01,54.77) -- (54.43,34.77) -- cycle ; \draw  [line width=1.5]  (71.96,14.69) .. controls (75.54,20.06) and (77.71,27.08) .. (77.71,34.77) .. controls (77.71,42.42) and (75.56,49.41) .. (72.01,54.77) ;

\draw (7,22.4) node [anchor=north west][inner sep=0.75pt]    {$\mathsf{E}_{0} \ =\ $};
\draw (149,23.4) node [anchor=north west][inner sep=0.75pt]    {$\mathsf{E}_{1} \ =\ $};
\draw (121,24) node [anchor=north west][inner sep=0.75pt]   [align=left] {,};

\end{tikzpicture}

%% file: tikzpictures/E1_uxlx.tex
\tikzset{every picture/.style={line width=0.75pt}} 

\begin{tikzpicture}[x=0.75pt,y=0.75pt,yscale=-.8,xscale=.8]

\draw  [dash pattern={on 4.5pt off 4.5pt}] (69,32.27) .. controls (69,18.21) and (80.4,6.81) .. (94.46,6.81) .. controls (108.52,6.81) and (119.91,18.21) .. (119.91,32.27) .. controls (119.91,46.33) and (108.52,57.72) .. (94.46,57.72) .. controls (80.4,57.72) and (69,46.33) .. (69,32.27) -- cycle ;
\draw  [draw opacity=0][line width=1.5]  (114.56,17.75) .. controls (109.19,21.33) and (102.17,23.5) .. (94.48,23.5) .. controls (86.83,23.5) and (79.84,21.35) .. (74.48,17.8) -- (94.48,0.23) -- cycle ; \draw  [line width=1.5]  (114.56,17.75) .. controls (109.19,21.33) and (102.17,23.5) .. (94.48,23.5) .. controls (86.83,23.5) and (79.84,21.35) .. (74.48,17.8) ;  
\draw  [draw opacity=0][line width=1.5]  (74.4,46.79) .. controls (79.77,43.21) and (86.79,41.04) .. (94.48,41.04) .. controls (102.13,41.04) and (109.12,43.19) .. (114.48,46.74) -- (94.48,64.31) -- cycle ; \draw  [line width=1.5]  (74.4,46.79) .. controls (79.77,43.21) and (86.79,41.04) .. (94.48,41.04) .. controls (102.13,41.04) and (109.12,43.19) .. (114.48,46.74) ;  

\draw  [fill={rgb, 255:red, 0; green, 0; blue, 0 }  ,fill opacity=1 ] (91,23.27) .. controls (91,21.5) and (92.43,20.07) .. (94.2,20.07) .. controls (95.97,20.07) and (97.41,21.5) .. (97.41,23.27) .. controls (97.41,25.04) and (95.97,26.47) .. (94.2,26.47) .. controls (92.43,26.47) and (91,25.04) .. (91,23.27) -- cycle ;
\draw  [dash pattern={on 4.5pt off 4.5pt}] (209,31.27) .. controls (209,17.21) and (220.4,5.81) .. (234.46,5.81) .. controls (248.52,5.81) and (259.91,17.21) .. (259.91,31.27) .. controls (259.91,45.33) and (248.52,56.72) .. (234.46,56.72) .. controls (220.4,56.72) and (209,45.33) .. (209,31.27) -- cycle ;
\draw  [draw opacity=0][line width=1.5]  (254.56,16.75) .. controls (249.19,20.33) and (242.17,22.5) .. (234.48,22.5) .. controls (226.83,22.5) and (219.84,20.35) .. (214.48,16.8) -- (234.48,-0.77) -- cycle ; \draw  [line width=1.5]  (254.56,16.75) .. controls (249.19,20.33) and (242.17,22.5) .. (234.48,22.5) .. controls (226.83,22.5) and (219.84,20.35) .. (214.48,16.8) ;  
\draw  [draw opacity=0][line width=1.5]  (214.4,45.79) .. controls (219.77,42.21) and (226.79,40.04) .. (234.48,40.04) .. controls (242.13,40.04) and (249.12,42.19) .. (254.48,45.74) -- (234.48,63.31) -- cycle ; \draw  [line width=1.5]  (214.4,45.79) .. controls (219.77,42.21) and (226.79,40.04) .. (234.48,40.04) .. controls (242.13,40.04) and (249.12,42.19) .. (254.48,45.74) ;  

\draw  [fill={rgb, 255:red, 0; green, 0; blue, 0 }  ,fill opacity=1 ] (231.25,40.47) .. controls (231.25,38.7) and (232.69,37.27) .. (234.46,37.27) .. controls (236.23,37.27) and (237.66,38.7) .. (237.66,40.47) .. controls (237.66,42.24) and (236.23,43.68) .. (234.46,43.68) .. controls (232.69,43.68) and (231.25,42.24) .. (231.25,40.47) -- cycle ;

\draw (17,24.67) node [anchor=north west][inner sep=0.75pt]    {$u_{X} \ =\ $};
\draw (127,25.77) node [anchor=north west][inner sep=0.75pt]   [align=left] {,};
\draw (157,23.67) node [anchor=north west][inner sep=0.75pt]    {$l_{X} \ =\ $};

\end{tikzpicture}

%% file: tikzpictures/endo-ab.tex
\tikzset{every picture/.style={line width=0.75pt}} 

\begin{tikzpicture}[x=0.75pt,y=0.75pt,yscale=-.9,xscale=.9]

\draw  [color={rgb, 255:red, 0; green, 0; blue, 0 }  ,draw opacity=1 ][line width=1.5]  (177.91,28.77) .. controls (177.91,23.46) and (182.41,19.16) .. (187.96,19.16) .. controls (193.5,19.16) and (198,23.46) .. (198,28.77) .. controls (198,34.07) and (193.5,38.37) .. (187.96,38.37) .. controls (182.41,38.37) and (177.91,34.07) .. (177.91,28.77) -- cycle ;
\draw  [dash pattern={on 4.5pt off 4.5pt}] (162.5,33.77) .. controls (162.5,19.71) and (173.9,8.31) .. (187.96,8.31) .. controls (202.02,8.31) and (213.41,19.71) .. (213.41,33.77) .. controls (213.41,47.83) and (202.02,59.22) .. (187.96,59.22) .. controls (173.9,59.22) and (162.5,47.83) .. (162.5,33.77) -- cycle ;
\draw  [draw opacity=0][line width=1.5]  (168.89,50.71) .. controls (174.12,48.05) and (180.76,46.45) .. (187.98,46.45) .. controls (194.91,46.45) and (201.3,47.92) .. (206.43,50.39) -- (187.98,65.81) -- cycle ; \draw  [line width=1.5]  (168.89,50.71) .. controls (174.12,48.05) and (180.76,46.45) .. (187.98,46.45) .. controls (194.91,46.45) and (201.3,47.92) .. (206.43,50.39) ;  
\draw  [dash pattern={on 4.5pt off 4.5pt}] (37,33.77) .. controls (37,19.71) and (48.4,8.31) .. (62.46,8.31) .. controls (76.52,8.31) and (87.91,19.71) .. (87.91,33.77) .. controls (87.91,47.83) and (76.52,59.22) .. (62.46,59.22) .. controls (48.4,59.22) and (37,47.83) .. (37,33.77) -- cycle ;
\draw  [draw opacity=0][line width=1.5]  (42.4,48.29) .. controls (47.77,44.71) and (54.79,42.54) .. (62.48,42.54) .. controls (70.13,42.54) and (77.12,44.69) .. (82.48,48.24) -- (62.48,65.81) -- cycle ; \draw  [line width=1.5]  (42.4,48.29) .. controls (47.77,44.71) and (54.79,42.54) .. (62.48,42.54) .. controls (70.13,42.54) and (77.12,44.69) .. (82.48,48.24) ;  
\draw    (100,33) -- (147,33) ;
\draw [shift={(149,33)}, rotate = 180] [color={rgb, 255:red, 0; green, 0; blue, 0 }  ][line width=0.75]    (10.93,-4.9) .. controls (6.95,-2.3) and (3.31,-0.67) .. (0,0) .. controls (3.31,0.67) and (6.95,2.3) .. (10.93,4.9)   ;
\draw    (151,44) -- (100,44) ;
\draw [shift={(98,44)}, rotate = 360] [color={rgb, 255:red, 0; green, 0; blue, 0 }  ][line width=0.75]    (10.93,-4.9) .. controls (6.95,-2.3) and (3.31,-0.67) .. (0,0) .. controls (3.31,0.67) and (6.95,2.3) .. (10.93,4.9)   ;
\draw    (223.77,25.32) .. controls (232.08,17.65) and (256.35,7.43) .. (256.99,33.62) .. controls (257.6,58.63) and (237.61,49.89) .. (225.44,44.01) ;
\draw [shift={(223.77,43.2)}, rotate = 25.94] [color={rgb, 255:red, 0; green, 0; blue, 0 }  ][line width=0.75]    (6.56,-1.97) .. controls (4.17,-0.84) and (1.99,-0.18) .. (0,0) .. controls (1.99,0.18) and (4.17,0.84) .. (6.56,1.97)   ;

\draw (179,69.65) node [anchor=north west][inner sep=0.75pt]    {$T'$};
\draw (53,69.65) node [anchor=north west][inner sep=0.75pt]    {$T$};
\draw (117,11.4) node [anchor=north west][inner sep=0.75pt]    {$\Delta$};
\draw (117,47.4) node [anchor=north west][inner sep=0.75pt]    {$m$};
\draw (263,23.4) node [anchor=north west][inner sep=0.75pt]    {$a,\ b$};

\end{tikzpicture}

%% file: tikzpictures/twist-knot.tex
\tikzset{every picture/.style={line width=0.75pt}} 

\begin{tikzpicture}[x=0.75pt,y=0.75pt,yscale=-1,xscale=1]

\draw    (130,140) -- (100,110) ;
\draw    (190,140) -- (160,110) ;
\draw    (134.33,50) -- (155.67,71.33) ;
\draw  [color={rgb, 255:red, 255; green, 255; blue, 255 }  ,draw opacity=1 ][fill={rgb, 255:red, 255; green, 255; blue, 255 }  ,fill opacity=1 ] (141.44,60.67) .. controls (141.44,58.7) and (143.04,57.11) .. (145,57.11) .. controls (146.96,57.11) and (148.56,58.7) .. (148.56,60.67) .. controls (148.56,62.63) and (146.96,64.22) .. (145,64.22) .. controls (143.04,64.22) and (141.44,62.63) .. (141.44,60.67) -- cycle ;
\draw    (155.67,50) -- (134.33,71.33) ;

\draw    (134.33,28.67) -- (155.67,50) ;
\draw  [color={rgb, 255:red, 255; green, 255; blue, 255 }  ,draw opacity=1 ][fill={rgb, 255:red, 255; green, 255; blue, 255 }  ,fill opacity=1 ] (141.44,39.33) .. controls (141.44,37.37) and (143.04,35.78) .. (145,35.78) .. controls (146.96,35.78) and (148.56,37.37) .. (148.56,39.33) .. controls (148.56,41.3) and (146.96,42.89) .. (145,42.89) .. controls (143.04,42.89) and (141.44,41.3) .. (141.44,39.33) -- cycle ;
\draw    (155.67,28.67) -- (134.33,50) ;
\draw  [color={rgb, 255:red, 155; green, 155; blue, 155 }  ,draw opacity=1 ][dash pattern={on 4.5pt off 4.5pt}] (134.33,28.67) -- (155.67,28.67) -- (155.67,71.33) -- (134.33,71.33) -- cycle ;

\draw  [color={rgb, 255:red, 255; green, 255; blue, 255 }  ,draw opacity=1 ][fill={rgb, 255:red, 255; green, 255; blue, 255 }  ,fill opacity=1 ] (160,110) -- (160,140) -- (130,140) -- (130,110) -- cycle ;

\draw  [color={rgb, 255:red, 255; green, 255; blue, 255 }  ,draw opacity=1 ][fill={rgb, 255:red, 255; green, 255; blue, 255 }  ,fill opacity=1 ] (115,120) .. controls (117.76,120) and (120,122.24) .. (120,125) .. controls (120,127.76) and (117.76,130) .. (115,130) .. controls (112.24,130) and (110,127.76) .. (110,125) .. controls (110,122.24) and (112.24,120) .. (115,120) -- cycle ;
\draw  [color={rgb, 255:red, 255; green, 255; blue, 255 }  ,draw opacity=1 ][fill={rgb, 255:red, 255; green, 255; blue, 255 }  ,fill opacity=1 ] (175,120) .. controls (177.76,120) and (180,122.24) .. (180,125) .. controls (180,127.76) and (177.76,130) .. (175,130) .. controls (172.24,130) and (170,127.76) .. (170,125) .. controls (170,122.24) and (172.24,120) .. (175,120) -- cycle ;
\draw    (130,110) -- (100,140) ;
\draw    (190,110) -- (160,140) ;
\draw  [color={rgb, 255:red, 155; green, 155; blue, 155 }  ,draw opacity=1 ][dash pattern={on 4.5pt off 4.5pt}] (190,110) -- (190,140) -- (100,140) -- (100,110) -- cycle ;
\draw [color={rgb, 255:red, 128; green, 128; blue, 128 }  ,draw opacity=1 ]   (100,140) .. controls (26.34,132) and (80.32,8.33) .. (134.33,28.67) ;
\draw [color={rgb, 255:red, 128; green, 128; blue, 128 }  ,draw opacity=1 ]   (134.33,71.33) .. controls (118.85,86.67) and (79.98,91.67) .. (100,110) ;
\draw [color={rgb, 255:red, 128; green, 128; blue, 128 }  ,draw opacity=1 ]   (190,140) .. controls (263.66,132) and (209.68,8.33) .. (155.67,28.67) ;
\draw [color={rgb, 255:red, 128; green, 128; blue, 128 }  ,draw opacity=1 ]   (155.67,71.33) .. controls (171.15,86.67) and (210.02,91.67) .. (190,110) ;
\draw [color={rgb, 255:red, 128; green, 128; blue, 128 }  ,draw opacity=1 ]   (114.33,84) -- (116.54,82.89) ;
\draw [shift={(118.33,82)}, rotate = 153.43] [color={rgb, 255:red, 128; green, 128; blue, 128 }  ,draw opacity=1 ][line width=0.75]    (10.93,-3.29) .. controls (6.95,-1.4) and (3.31,-0.3) .. (0,0) .. controls (3.31,0.3) and (6.95,1.4) .. (10.93,3.29)   ;

\draw (53.67,78.17) node [anchor=east] [inner sep=0.75pt]    {$K_{p} \ =$};
\draw (145,125) node  [rotate=-90]  {$\vdots $};
\draw (146.51,74.99) node [anchor=north] [inner sep=0.75pt]  [font=\footnotesize,color={rgb, 255:red, 128; green, 128; blue, 128 }  ,opacity=1 ]  {$T_{2}$};
\draw (90.33,113.67) node [anchor=west] [inner sep=0.75pt]  [font=\footnotesize,color={rgb, 255:red, 128; green, 128; blue, 128 }  ,opacity=1 ,rotate=-90]  {$T_{-p}$};

\end{tikzpicture}

%% file: tikzpictures/table-H-6_1.tex
\tikzset{every picture/.style={line width=0.75pt}} 

\begin{tikzpicture}[x=0.75pt,y=0.75pt,yscale=-1,xscale=1]

\draw  [draw opacity=0] (90,20) -- (150.33,20) -- (150.33,100.23) -- (90,100.23) -- cycle ; \draw  [color={rgb, 255:red, 155; green, 155; blue, 155 }  ,draw opacity=1 ] (90,20) -- (90,100.23)(110,20) -- (110,100.23)(130,20) -- (130,100.23)(150,20) -- (150,100.23) ; \draw  [color={rgb, 255:red, 155; green, 155; blue, 155 }  ,draw opacity=1 ] (90,20) -- (150.33,20)(90,40) -- (150.33,40)(90,60) -- (150.33,60)(90,80) -- (150.33,80)(90,100) -- (150.33,100) ; \draw  [color={rgb, 255:red, 155; green, 155; blue, 155 }  ,draw opacity=1 ]  ;
\draw  [color={rgb, 255:red, 245; green, 166; blue, 35 }  ,draw opacity=1 ][fill={rgb, 255:red, 0; green, 0; blue, 0 }  ,fill opacity=1 ] (99,46.8) .. controls (99,45.58) and (99.98,44.6) .. (101.2,44.6) .. controls (102.42,44.6) and (103.4,45.58) .. (103.4,46.8) .. controls (103.4,48.02) and (102.42,49) .. (101.2,49) .. controls (99.98,49) and (99,48.02) .. (99,46.8) -- cycle ;
\draw  [color={rgb, 255:red, 245; green, 166; blue, 35 }  ,draw opacity=1 ][fill={rgb, 255:red, 0; green, 0; blue, 0 }  ,fill opacity=1 ] (138.8,46.8) .. controls (138.8,45.58) and (139.78,44.6) .. (141,44.6) .. controls (142.22,44.6) and (143.2,45.58) .. (143.2,46.8) .. controls (143.2,48.02) and (142.22,49) .. (141,49) .. controls (139.78,49) and (138.8,48.02) .. (138.8,46.8) -- cycle ;
\draw [color={rgb, 255:red, 245; green, 166; blue, 35 }  ,draw opacity=1 ]   (103.4,46.8) .. controls (116.52,40.14) and (122.38,39.91) .. (137.16,46.1) ;
\draw [shift={(138.8,46.8)}, rotate = 203.33] [color={rgb, 255:red, 245; green, 166; blue, 35 }  ,draw opacity=1 ][line width=0.75]    (6.56,-1.97) .. controls (4.17,-0.84) and (1.99,-0.18) .. (0,0) .. controls (1.99,0.18) and (4.17,0.84) .. (6.56,1.97)   ;

\draw  [color={rgb, 255:red, 208; green, 2; blue, 27 }  ,draw opacity=1 ][fill={rgb, 255:red, 0; green, 0; blue, 0 }  ,fill opacity=1 ] (118.2,30.7) .. controls (118.2,29.48) and (119.18,28.5) .. (120.4,28.5) .. controls (121.62,28.5) and (122.6,29.48) .. (122.6,30.7) .. controls (122.6,31.92) and (121.62,32.9) .. (120.4,32.9) .. controls (119.18,32.9) and (118.2,31.92) .. (118.2,30.7) -- cycle ;
\draw  [color={rgb, 255:red, 126; green, 211; blue, 33 }  ,draw opacity=1 ][fill={rgb, 255:red, 0; green, 0; blue, 0 }  ,fill opacity=1 ] (99,69.7) .. controls (99,68.48) and (99.98,67.5) .. (101.2,67.5) .. controls (102.42,67.5) and (103.4,68.48) .. (103.4,69.7) .. controls (103.4,70.92) and (102.42,71.9) .. (101.2,71.9) .. controls (99.98,71.9) and (99,70.92) .. (99,69.7) -- cycle ;
\draw  [color={rgb, 255:red, 126; green, 211; blue, 33 }  ,draw opacity=1 ][fill={rgb, 255:red, 0; green, 0; blue, 0 }  ,fill opacity=1 ] (138.8,69.7) .. controls (138.8,68.48) and (139.78,67.5) .. (141,67.5) .. controls (142.22,67.5) and (143.2,68.48) .. (143.2,69.7) .. controls (143.2,70.92) and (142.22,71.9) .. (141,71.9) .. controls (139.78,71.9) and (138.8,70.92) .. (138.8,69.7) -- cycle ;
\draw [color={rgb, 255:red, 126; green, 211; blue, 33 }  ,draw opacity=1 ]   (103.4,69.7) .. controls (116.52,63.04) and (122.38,62.81) .. (137.16,69) ;
\draw [shift={(138.8,69.7)}, rotate = 203.33] [color={rgb, 255:red, 126; green, 211; blue, 33 }  ,draw opacity=1 ][line width=0.75]    (6.56,-1.97) .. controls (4.17,-0.84) and (1.99,-0.18) .. (0,0) .. controls (1.99,0.18) and (4.17,0.84) .. (6.56,1.97)   ;

\draw  [color={rgb, 255:red, 245; green, 166; blue, 35 }  ,draw opacity=1 ][fill={rgb, 255:red, 0; green, 0; blue, 0 }  ,fill opacity=1 ] (118.2,50.7) .. controls (118.2,49.48) and (119.18,48.5) .. (120.4,48.5) .. controls (121.62,48.5) and (122.6,49.48) .. (122.6,50.7) .. controls (122.6,51.92) and (121.62,52.9) .. (120.4,52.9) .. controls (119.18,52.9) and (118.2,51.92) .. (118.2,50.7) -- cycle ;
\draw  [color={rgb, 255:red, 126; green, 211; blue, 33 }  ,draw opacity=1 ][fill={rgb, 255:red, 0; green, 0; blue, 0 }  ,fill opacity=1 ] (113.2,71.7) .. controls (113.2,70.48) and (114.18,69.5) .. (115.4,69.5) .. controls (116.62,69.5) and (117.6,70.48) .. (117.6,71.7) .. controls (117.6,72.92) and (116.62,73.9) .. (115.4,73.9) .. controls (114.18,73.9) and (113.2,72.92) .. (113.2,71.7) -- cycle ;
\draw  [color={rgb, 255:red, 126; green, 211; blue, 33 }  ,draw opacity=1 ][fill={rgb, 255:red, 0; green, 0; blue, 0 }  ,fill opacity=1 ] (123.2,71.7) .. controls (123.2,70.48) and (124.18,69.5) .. (125.4,69.5) .. controls (126.62,69.5) and (127.6,70.48) .. (127.6,71.7) .. controls (127.6,72.92) and (126.62,73.9) .. (125.4,73.9) .. controls (124.18,73.9) and (123.2,72.92) .. (123.2,71.7) -- cycle ;
\draw  [color={rgb, 255:red, 74; green, 144; blue, 226 }  ,draw opacity=1 ][fill={rgb, 255:red, 0; green, 0; blue, 0 }  ,fill opacity=1 ] (118.2,89.7) .. controls (118.2,88.48) and (119.18,87.5) .. (120.4,87.5) .. controls (121.62,87.5) and (122.6,88.48) .. (122.6,89.7) .. controls (122.6,90.92) and (121.62,91.9) .. (120.4,91.9) .. controls (119.18,91.9) and (118.2,90.92) .. (118.2,89.7) -- cycle ;

\draw (82,31.6) node  [font=\scriptsize]  {$4$};
\draw (82,51.6) node  [font=\scriptsize]  {$2$};
\draw (116.5,130.5) node    {$\Delta \ =\ 0$};
\draw (82.5,11.5) node  [font=\footnotesize]  {$a$};
\draw (163.5,108.5) node  [font=\footnotesize]  {$q$};
\draw (120.5,110.8) node  [font=\scriptsize]  {$0$};
\draw (98,110.8) node  [font=\scriptsize]  {$-2$};
\draw (140.5,110.8) node  [font=\scriptsize]  {$2$};
\draw (82,72.5) node  [font=\scriptsize]  {$0$};
\draw (79.5,92.5) node  [font=\scriptsize]  {$-2$};

\end{tikzpicture}

%% file: tikzpictures/table-Kh-6_1.tex
\tikzset{every picture/.style={line width=0.75pt}} 

\begin{tikzpicture}[x=0.75pt,y=0.75pt,yscale=-1,xscale=1]

\draw  [draw opacity=0] (50,30) -- (190.2,30) -- (190.2,50.1) -- (50,50.1) -- cycle ; \draw  [color={rgb, 255:red, 155; green, 155; blue, 155 }  ,draw opacity=1 ] (50,30) -- (50,50.1)(70,30) -- (70,50.1)(90,30) -- (90,50.1)(110,30) -- (110,50.1)(130,30) -- (130,50.1)(150,30) -- (150,50.1)(170,30) -- (170,50.1)(190,30) -- (190,50.1) ; \draw  [color={rgb, 255:red, 155; green, 155; blue, 155 }  ,draw opacity=1 ] (50,30) -- (190.2,30)(50,50) -- (190.2,50) ; \draw  [color={rgb, 255:red, 155; green, 155; blue, 155 }  ,draw opacity=1 ]  ;
\draw  [color={rgb, 255:red, 245; green, 166; blue, 35 }  ,draw opacity=1 ][fill={rgb, 255:red, 0; green, 0; blue, 0 }  ,fill opacity=1 ] (117.6,44.33) .. controls (117.6,45.54) and (118.58,46.53) .. (119.8,46.53) .. controls (121.02,46.53) and (122,45.54) .. (122,44.33) .. controls (122,43.11) and (121.02,42.13) .. (119.8,42.13) .. controls (118.58,42.13) and (117.6,43.11) .. (117.6,44.33) -- cycle ;
\draw  [color={rgb, 255:red, 245; green, 166; blue, 35 }  ,draw opacity=1 ][fill={rgb, 255:red, 0; green, 0; blue, 0 }  ,fill opacity=1 ] (157.4,44.33) .. controls (157.4,45.54) and (158.38,46.53) .. (159.6,46.53) .. controls (160.82,46.53) and (161.8,45.54) .. (161.8,44.33) .. controls (161.8,43.11) and (160.82,42.13) .. (159.6,42.13) .. controls (158.38,42.13) and (157.4,43.11) .. (157.4,44.33) -- cycle ;
\draw [color={rgb, 255:red, 245; green, 166; blue, 35 }  ,draw opacity=1 ]   (122,44.33) .. controls (135.12,50.98) and (140.98,51.22) .. (155.76,45.02) ;
\draw [shift={(157.4,44.33)}, rotate = 156.67] [color={rgb, 255:red, 245; green, 166; blue, 35 }  ,draw opacity=1 ][line width=0.75]    (6.56,-1.97) .. controls (4.17,-0.84) and (1.99,-0.18) .. (0,0) .. controls (1.99,0.18) and (4.17,0.84) .. (6.56,1.97)   ;

\draw  [color={rgb, 255:red, 208; green, 2; blue, 27 }  ,draw opacity=1 ][fill={rgb, 255:red, 0; green, 0; blue, 0 }  ,fill opacity=1 ] (178,40) .. controls (178,38.78) and (178.98,37.8) .. (180.2,37.8) .. controls (181.42,37.8) and (182.4,38.78) .. (182.4,40) .. controls (182.4,41.22) and (181.42,42.2) .. (180.2,42.2) .. controls (178.98,42.2) and (178,41.22) .. (178,40) -- cycle ;
\draw  [color={rgb, 255:red, 126; green, 211; blue, 33 }  ,draw opacity=1 ][fill={rgb, 255:red, 0; green, 0; blue, 0 }  ,fill opacity=1 ] (78,37.8) .. controls (78,36.58) and (78.98,35.6) .. (80.2,35.6) .. controls (81.42,35.6) and (82.4,36.58) .. (82.4,37.8) .. controls (82.4,39.02) and (81.42,40) .. (80.2,40) .. controls (78.98,40) and (78,39.02) .. (78,37.8) -- cycle ;
\draw  [color={rgb, 255:red, 126; green, 211; blue, 33 }  ,draw opacity=1 ][fill={rgb, 255:red, 0; green, 0; blue, 0 }  ,fill opacity=1 ] (117.8,37.8) .. controls (117.8,36.58) and (118.78,35.6) .. (120,35.6) .. controls (121.22,35.6) and (122.2,36.58) .. (122.2,37.8) .. controls (122.2,39.02) and (121.22,40) .. (120,40) .. controls (118.78,40) and (117.8,39.02) .. (117.8,37.8) -- cycle ;
\draw [color={rgb, 255:red, 126; green, 211; blue, 33 }  ,draw opacity=1 ]   (82.4,37.8) .. controls (95.52,31.14) and (101.38,30.91) .. (116.16,37.1) ;
\draw [shift={(117.8,37.8)}, rotate = 203.33] [color={rgb, 255:red, 126; green, 211; blue, 33 }  ,draw opacity=1 ][line width=0.75]    (6.56,-1.97) .. controls (4.17,-0.84) and (1.99,-0.18) .. (0,0) .. controls (1.99,0.18) and (4.17,0.84) .. (6.56,1.97)   ;

\draw  [color={rgb, 255:red, 245; green, 166; blue, 35 }  ,draw opacity=1 ][fill={rgb, 255:red, 0; green, 0; blue, 0 }  ,fill opacity=1 ] (138.4,40) .. controls (138.4,38.78) and (139.38,37.8) .. (140.6,37.8) .. controls (141.82,37.8) and (142.8,38.78) .. (142.8,40) .. controls (142.8,41.22) and (141.82,42.2) .. (140.6,42.2) .. controls (139.38,42.2) and (138.4,41.22) .. (138.4,40) -- cycle ;
\draw  [color={rgb, 255:red, 126; green, 211; blue, 33 }  ,draw opacity=1 ][fill={rgb, 255:red, 0; green, 0; blue, 0 }  ,fill opacity=1 ] (92.6,40) .. controls (92.6,38.78) and (93.58,37.8) .. (94.8,37.8) .. controls (96.02,37.8) and (97,38.78) .. (97,40) .. controls (97,41.22) and (96.02,42.2) .. (94.8,42.2) .. controls (93.58,42.2) and (92.6,41.22) .. (92.6,40) -- cycle ;
\draw  [color={rgb, 255:red, 126; green, 211; blue, 33 }  ,draw opacity=1 ][fill={rgb, 255:red, 0; green, 0; blue, 0 }  ,fill opacity=1 ] (102.6,40) .. controls (102.6,38.78) and (103.58,37.8) .. (104.8,37.8) .. controls (106.02,37.8) and (107,38.78) .. (107,40) .. controls (107,41.22) and (106.02,42.2) .. (104.8,42.2) .. controls (103.58,42.2) and (102.6,41.22) .. (102.6,40) -- cycle ;
\draw  [color={rgb, 255:red, 74; green, 144; blue, 226 }  ,draw opacity=1 ][fill={rgb, 255:red, 0; green, 0; blue, 0 }  ,fill opacity=1 ] (57.6,40) .. controls (57.6,38.78) and (58.58,37.8) .. (59.8,37.8) .. controls (61.02,37.8) and (62,38.78) .. (62,40) .. controls (62,41.22) and (61.02,42.2) .. (59.8,42.2) .. controls (58.58,42.2) and (57.6,41.22) .. (57.6,40) -- cycle ;

\draw (119.5,60) node  [font=\scriptsize]  {$2$};
\draw (99.5,60) node  [font=\scriptsize]  {$0$};
\draw (159.5,60) node  [font=\scriptsize]  {$6$};
\draw (139.5,60) node  [font=\scriptsize]  {$4$};
\draw (179.5,60) node  [font=\scriptsize]  {$8$};
\draw (78.5,60) node  [font=\scriptsize]  {$-2$};
\draw (58.5,60) node  [font=\scriptsize]  {$-4$};
\draw (43,18) node  [font=\small]  {$\delta $};
\draw (200.8,58.2) node  [font=\small]  {$q_{2}$};
\draw (43,39) node  [font=\scriptsize]  {$0$};

\end{tikzpicture}

%% file: 6.tex
\section{Algorithmic computations}
\label{sec:direct-computation}

Finally, we describe an algorithm that computes the bigraded $k[e]$-module structure of $\Kh(L; k)$ over any computable field $k$, and present computational results performed by the program \cite{Sano:YUI} developed by the author. 

\subsection{Algorithm and implementation}

Recall that, for a finitely generated graded $k[e]$-module $M$ with $e$ acting nilpotently on $M$, the structure theorem states that $M$ decomposes as a finite direct sum
\begin{equation}
\label{eqn:M-decomp}
    M \isom \bigoplus_i q^{j_i} (k[e] / (e^{l_i}))
\end{equation}
with $j_i \in \ZZ$ and $0 < l_1 \leq l_2 \leq \cdots$. 

First, the (non-graded) $k[e]$-module structure of $M$ can be computed as follows. Fix a $k$-basis $\{m_i\}$ of $M$, and let $\varphi\colon k^n \to M$ be the map sending each unit vector $e_i$ to $m_i$. Let $A$ be the matrix representing the action of $e$ on $M$. By extending both $\varphi$ and $A$ over $k[e]$, there is a short exact sequence 
\[
\begin{tikzcd}
    0 \arrow[r] & {k[e]^n} \arrow[r, "eI - A"] & {k[e]^n} \arrow[r, "\varphi", two heads] & M \arrow[r] & 0
\end{tikzcd}
\]
giving the isomorphism $M \isom k[e]^n / (eI - A)$. The computation of the \textit{Smith normal form} (SNF) of $eI - A$ 
\[
    S = \operatorname{diag}(1, \ldots, 1, e^{l_1}, e^{l_2}, \ldots)
\]
by elementary row and column operations gives invertible matrices $P, Q \in GL(n; k[e])$ such that the following diagram commutes:
\[
\begin{tikzcd}[row sep = 3em, column sep = 3em]
{k[e]^n} \arrow[d, "Q"'] \arrow[r, "eI - A"] & {k[e]^n} \arrow[d, "P"] \\
{k[e]^n} \arrow[r, "S"] & {k[e]^n}.    
\end{tikzcd}
\]
The cokernel of the bottom row gives the decomposition \eqref{eqn:M-decomp} (without the grading). Furthermore, let $v_i$ denote the preimage of the $i$-th unit vector $e_i$ by $P$, 
\[
    v_i := P^{-1}e_i \in k[n]^n.
\]
Then $z_i := \varphi(v_i) \in M$ gives a generator of the summand $k[e] / (e^{l_i})$. 

Next, assume that $M$ is graded, $e$ is homogeneous, and the $k$-basis of $M$ is chosen to be homogeneous. Then we may perform the SNF computation within the category of graded $k[e]$-modules (by restricting to grading preserving elementary row and column operations). In this case each generator $z_i$ is homogeneous, and is expressed as a linear combination $z_i = \sum_j c_j m_j$ of the homogeneous basis with monomial coefficients. Its grading can be read off from a term $c_j m_j$ having a non-zero constant coefficient $c_j \in k^\times$, in which case $\deg(z_i) = \deg(m_j)$.

Now let $K$ be a knot diagram. For $M = \Kh(K; k)$ as a bigraded $k$-module, it obviously satisfies the above conditions. To perform the actual computation, first compute the homology $\Kh(K; k)$ together with explicit homogeneous generators and transformation matrices. Then compute the action of $e$ on each of the representatives, which gives the matrix representing the action of $e$ on $\Kh(K; k)$. Perform graded SNF computation to obtain the $e$-string decomposition of $\Kh(K; k)$. 


Here, we introduce a more compact form of presenting the $e$-string decomposition of $\Kh(L; k)$. Let 
\[
    e(l) := 1 + q^{4} e + \cdots + q^{4(l - 1)} e^{l - 1}
\]
and identify it with the graded $k[e]$-module $k[e] / (e^l)$ of dimension $l$. The total bigraded $k[e]$-module structure of $\Kh(L; k)$, given in \eqref{eqn:kh-e-string-decomp}, will be presented as a polynomial in $(\delta, q, e)$, as 
\[
    \Kh(L; k) = \sum_{i, j, l} \delta^i q^j e(l).
\]

The above explained algorithm has been implemented in the program \cite{Sano:YUI}. It is capable of computing the $e$-string decomposition of both unreduced and reduced Khovanov homology for knots, over $k = \QQ$ and other computable fields. It can optionally show $\sle$ in the form of matrices, which even works over $k = \ZZ$. (An $e$-string decomposition over $\ZZ$ does not make sense, since $\ZZ[e]$ is not a PID.) 

Here, we demonstrate the computation of an $e$-string decomposition for the left-handed trefoil $3_1$. With the command,

\begin{verbatim}
> ykh sl2 3_1 -t Q
\end{verbatim}

\noindent
the program outputs

\begin{verbatim}
 j\i  -3  -2  -1  0 
  9   Q   .   .   . 
  7   .   .   .   . 
  5   .   Q   .   . 
  3   .   .   .   Q 
  1   .   .   .   Q 

  d^1 q^1 e(2) + d^3 q^3 e(1) + d^3 q^9 e(1)
\end{verbatim}

\noindent
It shows the bigraded module $\Kh(3_1; \QQ)$ with the usual bigrading $(\gr_t, \gr_q)$, together with an $e$-string decomposition
\[
    \Kh(3_1; \QQ) = \delta^1 q^1 e(2) + \delta^3(q^3 e(1) + q^9 e(1)).
\]
The result agrees with \Cref{cor:kh-2q-torus} for $k = 1$, after mirroring. 

Below is another computation for the reduced Khovanov homology of $7_1$, this time computed over $\ZZ$. Instead of showing the $e$-string decomposition, it shows the explicit description of $\sle$ in the form of matrices. 

\begin{verbatim}
> ykh sl2 7_1 -t Z
\end{verbatim}

\begin{verbatim}
 j\i  -7  -6  -5  -4  -3  -2  -1  0 
  20  Z   .   .   .   .   .   .   . 
  18  .   Z   .   .   .   .   .   . 
  16  .   .   Z   .   .   .   .   . 
  14  .   .   .   Z   .   .   .   . 
  12  .   .   .   .   Z   .   .   . 
  10  .   .   .   .   .   Z   .   . 
  8   .   .   .   .   .   .   .   . 
  6   .   .   .   .   .   .   .   Z 

delta: 6

  Z(0, 6) -> Z(-2, 10); rank: 1
  [-6]
  
  Z(-2, 10) -> Z(-4, 14); rank: 1
  [4]

  Z(-3, 12) -> Z(-5, 16); rank: 1
  [-4]

  Z(-4, 14) -> Z(-6, 18); rank: 1
  [-2]

  Z(-5, 16) -> Z(-7, 20); rank: 1
  [2]
\end{verbatim}

\noindent
Again, we see that the result agrees with \Cref{cor:kh-2q-torus} for $k = 3$. 

The current implementation is not very efficient, as it naively computes $\Kh(K)$ (without Bar-Natan's fast computation algorithm \cite{BarNatan:2007}), and then computes the matrix for $\sle$. Obviously, the reduction technique explained in \Cref{sec:diagrammatic-computation} can be implemented in cooperation with Bar-Natan's algorithm, and would drastically improve the performance. The improvement is left as future work. 

\subsection{Computation results}

We computed the $\QQ[e]$-module structure of reduced Khovanov homology for all prime knots with up to $11$ crossings, using the program \cite{Sano:YUI} together with the dataset provided by the \textit{Knot Atlas}~\cite{KnotAtlas}. Below, we summarize the results, which also combines computational results for the (reduced) HOMFLY--PT homology computed in \cite{Nakagane-Sano:2025}.

The following proposition contains the result for the \textit{Conway knot} $K_C = 11n_{34}$ and the \textit{Kinoshita--Terasaka knot} $K_{KT} = 11n_{42}$ featured in \Cref{mainthm:conway-and-kt}. 

\begin{prop}
\label{prop:stronger-than-homfly}
    Of the 44 groups of prime knots given in \cite[Proposition 4.4]{Nakagane-Sano:2025}, each of which contains two or more distinct knots that have identical HOMFLY--PT homology and identical reduced Khovanov homology, the following four pairs have distinct $e$-operator on reduced Khovanov homology. 
    \[
        (11n_{34}, 11n_{42}), \quad
        (11n_{39}, 11n_{45}), \quad
        (11n_{73}, 11n_{74}), \quad
        (11n_{151}, 11n_{152}).
    \]
    Indeed, the $\QQ[e]$-module structures are given by:
    {\allowdisplaybreaks
    \begin{align*}
        \rKh&(11n_{34}) \\ 
            &= \delta^{-2} (q^{-6}e(2) + q^{-4}e(2) + q^{-4}e(3) + q^{-2}e(3) + e(1) + q^2e(1) + q^2e(2) + q^4e(2)) \\
            &+ \delta^{0}  (q^{-12}e(2) + q^{-10}e(2) + q^{-10}e(3) + q^{-8}e(3) + q^{-6}e(1) + q^{-4}e(1) + q^{-4}e(2) + q^{-2}e(2) + e(1) ) \\ 
        \rKh&(11n_{42})  \\ 
            &= \delta^{-2} (q^{-6}e(2) + q^{-4}e(1) + q^{-4}e(3) + q^{-2}e(2) + 2e(1) + 2q^2e(2) + q^4e(1) + q^8e(1)) \\ 
            &+ \delta^{0}  (q^{-12}e(1) + 2q^{-10}e(2) + q^{-8}e(1) + q^{-8}e(3) + q^{-6}e(2) + 2q^{-4}e(1) + q^{-2}e(2) + 2e(1)), \\
        \rKh&(11n_{39}) \\
            &= \delta^{0} (q^{-4}e(1) + 3q^{-2}e(2) + 3e(1) + 2e(3) + 2q^2e(2) + 4q^4e(1) + 3q^6e(2) + 2q^8e(1) + q^{12}e(1)) \\
            &+ \delta^{2} (q^{-8}e(1) + q^{-6}e(2) + q^{-4}e(1) + q^{-2}e(2) + e(1) + q^4e(1)), \\
        \rKh&(11n_{45}) \\
            &= \delta^{0} (q^{-4}e(1) + 3q^{-2}e(2) + 2e(1) + e(2) + 2e(3) + q^2e(2) + q^2e(3) + 3q^4e(1) + q^6e(1) + 2q^6e(2) \\ 
            &+ q^8e(1) + q^8e(2)) + \delta^{2} (q^{-8}e(2) + q^{-6}e(3) + q^{-2}e(1) + e(2)) \\ 
        \rKh&(11n_{73}) \\
            &= \delta^{0} (q^{-2}e(2) + e(1) + e(2) + e(3) + q^2e(3) + q^4e(1) + q^6e(1) + q^6e(2) + q^8e(2)) \\
            &+ \delta^{2} (q^{-8}e(2) + q^{-6}e(3) + q^{-2}e(1) + e(2)), \\
        \rKh&(11n_{74}) \\
            &= \delta^{0} (q^{-2}e(2) + 2e(1) + e(3) + q^2e(2) + 2q^4e(1) + 2q^6e(2) + q^8e(1) + q^{12}e(1)) \\
            &+ \delta^{2} (q^{-8}e(1) + q^{-6}e(2) + q^{-4}e(1) + q^{-2}e(2) + e(1) + q^4e(1)),\\
        \rKh&(11n_{151}) \\
            &= \delta^{2} (3q^2e(2) + 2q^4e(1) + 2q^4e(3) + 2q^6e(2) + 4q^8e(1) + 3q^{10}e(2) + 2q^{12}e(1) + q^{16}e(1)  \\
            &+ \delta^{4} (q^{-4}e(1) + q^{-2}e(2) + e(1) + q^2e(2) + q^4e(1) + q^8e(1)),\\
        \rKh&(11n_{152}) \\
            &= \delta^{2} (3q^2e(2) + q^4e(1) + q^4e(2) + 2q^4e(3) + q^6e(2) + q^6e(3) + 3q^8e(1) + q^{10}e(1) + 2q^{10}e(2) \\
            &+ q^{12}e(1) + q^{12}e(2)) + \delta^{4} (q^{-4}e(2) + q^{-2}e(3) + q^2e(1) + q^4e(2)).
    \end{align*}
    }
\end{prop}

\begin{prop}
\label{prop:stronger-than-kh}
    Of the 90 groups of prime knots with up to $11$ crossings, each of which contains two or more knots with identical reduced Khovanov homology, the following 57 groups are ruled out as having distinct $\sle$-operator. 
    \[
\begin{array}{l}
\{9_{8}, m(11n_{60})\},\ 
\{9_{12}, m(11n_{15})\},\ 
\{9_{14}, m(11n_{53})\},\ 
\{9_{20}, m(11n_{90})\},\ 
\{9_{21}, m(11n_{129})\},\\ 
\{9_{22}, m(11n_{3})\},\ 
\{9_{25}, m(11n_{25})\},\ 
\{9_{27}, 11n_{83}\},\ 
\{9_{30}, m(11n_{114})\},\ 
\{9_{36}, 11n_{16}\},\\ 
\{9_{39}, 11n_{11} = 11n_{112}\},\ 
\{9_{41}, m(11n_{4}) = m(11n_{21})\},\ 
\{10_{22}, 10_{35}\},\ 
\{10_{41}, m(10_{94})\},\ 
\{10_{43}, 10_{91}\},\\ 
\{10_{59}, 10_{106}\},\ 
\{10_{60}, m(10_{86})\},\ 
\{10_{71}, 10_{104})\},\ 
\{10_{73}, m(10_{83})\},\
\{10_{81}, m(10_{109})\},\\
\{10_{137}, m(10_{155}) = 11n_{37}\},\
\{10_{138}, 11n_{117}\},\ 
\{11a_{3}, 11a_{51}, 11a_{331}\},\
\{11a_{5}, 11a_{112}\},\ 
\{11a_{7}, 11a_{325}\},\\ 
\{11a_{9}, 11a_{140}\},\ 
\{11a_{10}, 11a_{42}\},\ 
\{11a_{16}, 11a_{280}\},\ 
\{11a_{28}, m(11a_{87}), 11a_{96}\},\\ 
\{11a_{30} = 11a_{272}, m(11a_{189})\},\ 
\{11a_{34}, 11a_{89}\},\
\{11a_{35} = 11a_{316}, m(11a_{36})\},\ 
\{11a_{45}, 11a_{118}\},\\ 
\{11a_{60}, 11a_{220}\},\ 
\{11a_{66}, 11a_{121}\},\ 
\{11a_{67}, 11a_{317}\},\ 
\{11a_{68}, 11a_{111}\},\ 
\{11a_{80}, m(11a_{270})\},\\ 
\{11a_{84}, 11a_{88}\},\ 
\{11a_{86}, m(11a_{205})\},\ 
\{11a_{110}, 11a_{257}\},\ 
\{11a_{120}, 11a_{295}\},\ 
\{11a_{125}, m(11a_{297})\},\\ 
\{11a_{131}, 11a_{218}\},\ 
\{11a_{146}, 11a_{294}\},\ 
\{11a_{147}, m(11a_{322})\},\ 
\{11a_{159}, 11a_{347}\},\ 
\{11a_{228}, 11a_{251} = 11a_{253}\},\\
\{11n_{7}, 11n_{36} = 11n_{44}\},\ 
\{11n_{34}, 11n_{42}\},\ 
\{11n_{39}, 11n_{45}\},\ 
\{11n_{66}, 11n_{150}\},\ 
\{11n_{73}, 11n_{74}\},\\ 
\{11n_{103}, m(11n_{175})\},\ 
\{11n_{124}, 11n_{166}\},\\ 
\{11n_{148}, 11n_{168}\},\ 
\{11n_{151}, 11n_{152}\}.
\end{array}
    \]
    Here, notations such as $\{K_1, K_2 = K_3\}$ indicate that the $\sle$-operators of $K_1$ and $K_2$ are distinct, but those of $K_2$ and $K_3$ are identical. 
\end{prop}

It is also interesting to ask whether $\sle$ can detect chiral knots $K$ (i.e. $K \neq K^*$) that cannot be detected by the HOMFLY--PT homology and reduced Khovanov homology. However, no such knots were found within our computation. 

\begin{prop}
\label{prop:cannot-detect-chiral}
    The following four are the only chiral prime knots with up to $11$ crossings whose HOMFLY--PT homology is identical to that of their mirror (\cite[Proposition 4.3]{Nakagane-Sano:2025}). These knots also have $\sle$-operator on reduced Khovanov homology identical to that of their mirror.
    \[
        10_{48},\  
        10_{71},\  
        10_{91},\ 
        10_{104}.
    \]
\end{prop}

We also performed computations for the \textit{unreduced} Khovanov homology. As we have seen in \Cref{cor:kh-twist-knot}, the strength of $\sle$ is not equal between the two (in fact, for twist knots, the $\sle$-operator on the unreduced homology is almost trivial). Within our computations, there were no differences of strength between the reduced and the unreduced homology. To be precise, the above \Cref{prop:stronger-than-homfly,prop:stronger-than-kh,prop:cannot-detect-chiral} hold verbatim if we replace \textit{reduced} with \textit{unreduced}. 

Recently, large scale computations of Khovanov homology and HOMFLY--PT homology have been performed in \cite{CM:2024,KLTVZ:2025}. It is interesting to see how much $\sle$ improves their results. 

\subsection
[Analysis of Proposition \ref{prop:stronger-than-homfly}]
{Analysis of \Cref{prop:stronger-than-homfly}}
\label{subsec:conway-and-kt}

Finally, we give a deeper analysis of \Cref{prop:stronger-than-homfly}, especially for the \textit{Conway knot} $K_C = 11n_{34}$ and the \textit{Kinoshita--Terasaka knot} $K_{KT} = 11n_{42}$. Hereafter, we let $K$ denote either $K_C$ or $K_{KT}$. First, from the computational result of \cite{Nakagane-Sano:2025}, the triply-graded $\QQ$-module structure of $\mcH(K)$ is given by 
\begin{center}
\begin{tabular}{r|lllllll}
$a \setminus q$ 
      & $-6$       & $-4$       & $-2$       & $0$        & $2$        & $4$        & $6$ \\
\hline
$0$   & .          & $1_{(-4)}$ & $1_{(-2)}$ & $3_{(0)}$  & $1_{(2)}$  & $1_{(4)}$  & . \\
$-2$  & $1_{(-10)}$& $1_{(-8)}$ & $3_{(-6)}$ & $2_{(-4)}$ & $3_{(-2)}$ & $1_{(0)}$  & $1_{(2)}$ \\
$-4$  & .          & $1_{(-12)}$& $1_{(-10)}$& $2_{(-8)}$ & $1_{(-6)}$ & $1_{(-4)}$ & .
\end{tabular}
$$ \Delta = 0$$

\begin{tabular}{r|lllllll}
$a \setminus q$ 
      & $-6$       & $-4$       & $-2$       & $0$        & $2$        & $4$       & $6$ \\
\hline
$2$   & .          & $1_{(0)}$  & $1_{(2)}$  & $2_{(4)}$  & $1_{(6)}$  & $1_{(8)}$ & . \\
$0$   & $1_{(-6)}$ & $1_{(-4)}$ & $3_{(-2)}$ & $2_{(0)}$  & $3_{(2)}$  & $1_{(4)}$ & $1_{(6)}$ \\
$-2$  & .          & $1_{(-8)}$ & $1_{(-6)}$ & $2_{(-4)}$ & $1_{(-2)}$ & $1_{(0)}$ & .
\end{tabular}
$$ \Delta = -2$$
\end{center}
Here, the two non-trivial $\Delta$-slices $\Delta = -2, 0$ are shown as bigraded $\QQ$-modules with bigrading $(q, a)$. Each number in the cell at position $(q, a) = (i, j)$ denotes the dimension of $\mcH^{i, j,\ \Delta - (i + j)}(K)$, and the subscripted number indicates the $\sl_2$-modified $q$-grading $q + 2a$. As reviewed in \Cref{subsec:homfly-overview}, the first differential $d^1$ of the spectral sequence from $\mcH(K)$ to $\overline{H}_2(K) \isom \rKh(K^*)$ has $(\deg_\Delta, \deg_q, \deg_a) = (2, 4, -2)$, and the higher differentials are all trivial. Therefore, we have 
\[
    E^2(K) = H(\mcH(K), d^1) \isom \rKh(K^*).
\]
The bigraded $\QQ$-module structure of $\rKh(K^*)$ is given by 
\begin{center}
\begin{tabular}{r|lllllllllll}
$\delta \setminus q$ 
     & $-12$ & $-10$ & $-8$ & $-6$ & $-4$ & $-2$ & $0$ & $2$ & $4$ & $6$ & $8$ \\
\hline
$0$  & $1$   & $2$   & $2$  & $3$  & $3$  & $2$  & $3$ & $1$ & .   & .   &  . \\
$-2$ & .     & .     & .    & $1$  & $2$  & $2$  & $3$ & $3$ & $2$ & $2$ & $1$ 
\end{tabular}\end{center}
From $E^2 \isom \rKh(K^*)$, its triply graded structure is partially determined as 
\begin{center}
\begin{tabular}{r|lllllll}
$a \setminus q$ 
      & $-6$       & $-4$       & $-2$       & $0$        & $2$        & $4$       & $6$ \\
\hline
$2$   & .          & $1^?_{(0)}$  & $1^?_{(2)}$  & $1_{(4)}$  & $1_{(6)}$  & $1_{(8)}$ & . \\
$0$   & $1^?_{(-6)}$ & $1^?_{(-4)}$ & $1_{(-2)}$ & $2^?_{(0)}$  & $3^?_{(2)}$  & $1_{(4)}$ & $1_{(6)}$ \\
$-2$  & .          & . & $1^?_{(-6)}$ & $2^?_{(-4)}$ & $1_{(-2)}$ & $1_{(0)}$ & .
\end{tabular}
$$ \Delta = -2$$

\begin{tabular}{r|lllllll}
$a \setminus q$ 
      & $-6$       & $-4$       & $-2$       & $0$        & $2$        & $4$        & $6$ \\
\hline
$0$   & .          & $1_{(-4)}$ & $1_{(-2)}$ & $3^?_{(0)}$  & $1^?_{(2)}$  & . & . \\
$-2$  & $1_{(-10)}$& $1_{(-8)}$ & $3^?_{(-6)}$ & $2^?_{(-4)}$ & $1_{(-2)}$ & $1^?_{(0)}$  & $1^?_{(2)}$ \\
$-4$  & .          & $1_{(-12)}$& $1_{(-10)}$& $1_{(-8)}$ & $1^?_{(-6)}$ & $1^?_{(-4)}$ & .
\end{tabular}
$$ \Delta = 0$$
\end{center}
Here, notation such as $n^?_{(q)}$ implies that the dimension is undetermined and is bounded above by $n$. 

Now, the bigraded $\QQ[e]$-module structure of $\rKh(K_C^*)$ is given by 
\begin{center}
\begin{tabular}{r|lllllllllll}
$\delta \setminus q$ 
     & $-12$ & $-10$       & $-8$  & $-6$  & $-4$        & $-2$   & $0$   & $2$         & $4$   & $6$ & $8$ \\
\hline
$0$  & $e(2)$ & $e(2) + e(3)$ & $e(3)$ & $e(1)$ & $e(1) + e(2)$ & $e(2)$  & $e(1)$ & .           & .     & .   & . \\
$-2$ & .     & .           & .     & $e(2)$ & $e(2) + e(3)$ & $e(3)$  & $e(1)$ & $e(1) + e(2)$ & $e(2)$ & .   & . 
\end{tabular}\end{center}
In general, $e$-strings with longer length give stronger constraints. For example, consider the $e$-string $\delta^{-2}q^{-4}e(3)$. Since the $\sle$-operator is compatible with the spectral sequence (\Cref{prop:sl2-action-spectral-seq}), there is only one place where it can sit in $E^2(K_C)$, which is 
\[
    \Delta^{-2}q^{-4}a^0e(3)\colon 1_{(-4)} \to 1_{(0)} \to 1_{(4)}.
\]
In this way, we may fully determine the triply graded $\QQ[e]$-module structure of $E^2(K_C)$ as
\begin{center}
\begin{tabular}{r|lllllll}
$a \setminus q$ 
      & $-6$ & $-4$  & $-2$  & $0$   & $2$   & $4$ & $6$ \\
\hline
$2$   & .    & .     & $e(2)$ & $e(2)$ & .     & .   & . \\
$0$   & .    & $e(3)$ & $e(3)$ & $e(1)$ & $e(1)$ & .   & . \\
$-2$  & .    & .     & $e(2)$ & $e(2)$ & .     & .   & .
\end{tabular}
$$ \Delta = -2$$

\begin{tabular}{r|lllllll}
$a \setminus q$ 
      & $-6$  & $-4$  & $-2$  & $0$   & $2$ & $4$ & $6$ \\
\hline
$0$   & .     & $e(2)$ & $e(2)$ & $e(1)$ & .   & .   & . \\
$-2$  & $e(3)$ & $e(3)$ & $e(1)$ & $e(1)$ & .   & .   & . \\
$-4$  & .     & $e(2)$ & $e(2)$ & .     & .   & .   & .
\end{tabular}
$$ \Delta = 0$$
\end{center}

\noindent
Next, the bigraded $\QQ[e]$-module structure of $\rKh(K_{KT}^*)$ is given by 
\begin{center}
\begin{tabular}{r|lllllllllll}
$\delta \setminus q$ 
     & $-12$ & $-10$  & $-8$        & $-6$  & $-4$        & $-2$   & $0$    & $2$    & $4$   & $6$ & $8$  \\
\hline
$0$  & $e(1)$ & $2e(2)$ & $e(1) + e(3)$ & $e(2)$ & $2e(1)$      & $e(2)$  & $2e(1)$ & .      & .     & .   &  . \\
$-2$ & .     & .      & .           & $e(2)$ & $e(1) + e(3)$ & $e(2)$  & $2e(1)$ & $2e(2)$ & $e(1)$ & .   & $e(1)$ 
\end{tabular}\end{center}
Proceeding as before, we determine the triply graded $\QQ[e]$-module structure of $E^2(K_{KT})$ as
\begin{center}
\begin{tabular}{r|lllllll}
$a \setminus q$ 
      & $-6$ & $-4$  & $-2$  & $0$   & $2$   & $4$     & $6$ \\
\hline
$2$   & .    & .     & $e(2)$ & $e(1)$ & .     & $e(1)$   & . \\
$0$   & .    & $e(3)$ & $e(2)$ & $e(1)$ & $e(2)$ & .       & . \\
$-2$  & .    & .     & $e(2)$ & $e(1)$ & .     & $e(1)$   & .
\end{tabular}
$$ \Delta = -2$$

\begin{tabular}{r|lllllll}
$a \setminus q$ 
      & $-6$  & $-4$  & $-2$  & $0$    & $2$ & $4$ & $6$ \\
\hline
$0$   & .     & $e(1)$ & $e(2)$ & $2e(1)$ & .   & .   & . \\
$-2$  & $e(2)$ & $e(3)$ & $e(2)$ & $e(1)$  & .   & .   & . \\
$-4$  & .     & $e(1)$ & $e(2)$ & $e(1)$  & .   & .   & .
\end{tabular}
$$ \Delta = 0$$
\end{center}
These results are visualized as \Cref{fig:Kh-conway,fig:Kh-KT} in \Cref{sec:intro}. 
